# SMARANDACHE NEAR-RINGS

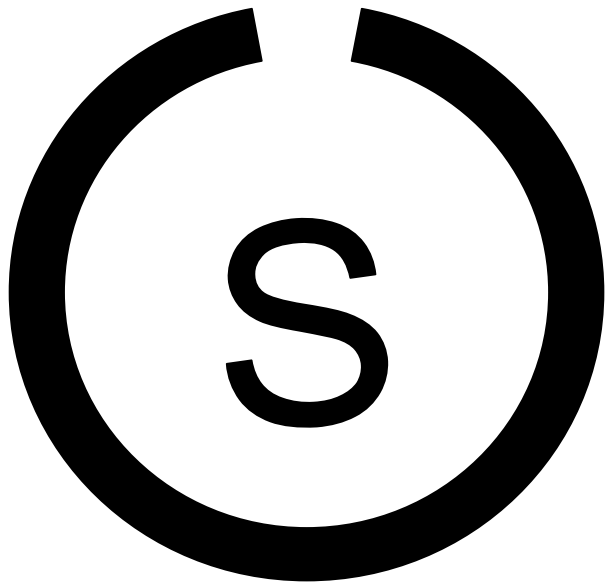

**w. b. vasantha kandasamy**
Department of Mathematics
Indian Institute of Technology, Madras
Chennai – 600 036
e-mail: vasantha@iitm.ac.in
web: **http://mat.iitm.ac.in/~wbv**

2002

*The picture on the cover is a representation of the Smarandache near-ring (S-near-ring). The international near-ring community represents the near-ring by a structure which is nearly a ring, but not exactly a ring, i.e. a near-ring. For the S-near-ring, of course, we have a "S" within the near-ring structure!*

# CONTENTS

# PREFACE

Near-rings are one of the generalized structures of rings. The study and research on near-rings is very systematic and continuous. Near-ring newsletters containing complete and updated bibliography on the subject are published periodically by a team of mathematicians (Editors: Yuen Fong, Alan Oswald, Gunter Pilz and K. C. Smith) with financial assistance from the National Cheng Kung University, Taiwan. These newsletters give an overall picture of the research carried out and the recent advancements and new concepts in the field. Conferences devoted solely to near-rings are held once every two years. There are about half a dozen books on near-rings apart from the conference proceedings. Above all there is a online searchable database and bibliography on near-rings. As a result the author feels it is very essential to have a book on Smarandache near-rings where the Smarandache analogues of the near-ring concepts are developed. The reader is expected to have a good background both in algebra and in near-rings; for, several results are to be proved by the reader as an exercise.

This book is organized into ten chapters: chapter one recalls some of the basic notions on groups, semigroups, groupoids semirings and lattices. The basic notions on near-rings are dealt in chapter two. In chapter three several definitions from available researchers are restated with the main motivation of constructing their Smarandache analogues. The main concern of this book is the study of Smarandache analogue properties of near-rings and Smarandache near-rings; so it does not promise to cover all concepts or the proofs of all results. Chapter four introduces the concept of S-near-rings and some of its basic properties. In chapter five all associative Smarandache near-rings built using group near-rings/semi near-rings and semigroup near-rings and semigroup seminear-rings are discussed. This alone has helped the author to build several new and innovative classes of near-rings and seminear-rings, both commutative and non-commutative, finite and infinite; using the vast classes of groups and semigroups. It has become important to mention here that the study of near-rings or seminear-rings using groups or semigroup is void except for a less than half a dozen papers. So when Smarandache notions are introduced the study in this direction would certainly be perused by innovative researchers/students.

Chapter six introduces and studies the concepts about seminear-rings and Smarandache seminear-rings. It is worthwhile to mention here that even the study of seminear-rings seems to be very meagre. The two applications of seminear-rings in case of group automatons and Balanced Incomplete Block Designs (BIBD) are given in chapter seven and the corresponding applications of the Smarandache structures are introduced. One of the major contribution is that by defining Smarandache planar near-rings we see for a given S-planar near-ring we can have several BIBD's which is impossible in case of planar near-rings. So, the introduction of several BIBD's from a single S-planar near-ring will lead to several error-correcting codes. Codes using BIBD may be having common properties that will certainly be of immense use to a coding theorist.

The study of non-associative structures in algebraic structures has become a separate entity; for, in the case of groups, their corresponding non-associative structure i.e. loops is dealt with separately. Similarly there is vast amount of research on the non-associative structures of semigroups i.e. groupoids and that of rings i.e. non-associative rings. However it is unfortunate that we do not have a parallel notions or study of non-associative near-rings. The only known concept is loop near-rings where the additive group structures of a near-ring is replaced by a loop. Though this study was started in 1978, further development and research is very little. Further, this definition is not in similar lines with rings. So in chapter eight we have defined non-associative near-rings and given methods of building non-associative near-rings and seminear-rings using loops and groupoids which we call as groupoid near-rings, near loop rings, groupoid seminear-rings and loop seminear-ring. For all these concepts a Smarandache analogue is defined and several Smarandache properties are introduced and studied. The ninth chapter deals with fuzzy concepts in near-rings and gives 5 new unconventional class of fuzzy near-rings; we also define their Smarandache analogues. The final chapter gives 145 suggested problems that will be of interest to researchers. It is worthwhile to mention that in the course of this book we have introduced over 260 Smarandache notions pertaining to near-ring theory.

My first thanks are due to Dr. Minh Perez, whose constant encouragement and intellectual support has made me to write this book.

I am indebted to my devoted husband Dr. Kandasamy and my very understanding children Meena and Kama without whose active support this book would have been impossible.

I humbly dedicate this book to the heroic social revolutionary Chhatrapathi Shahuji Maharaj, King of Kolhapur state, who, a century ago, began the policy of affirmative action in India, thereby giving greater opportunities for the socially dispriveleged and traditionally discriminated oppressed castes. His life is an inspiration to all of us, and his deeds are benchmarks in the history of the struggle against the caste system.

**Chapter One**

# PREREQUISITES

In this chapter we just recall the definition of groups, groupoids, loops, semigroups, semirings and lattices. We only give the basic definitions and some of its properties essential for the study of this book. We at the outset expect the reader to have a good knowledge in algebra and in near-rings. This chapter has five sections. In the first section we recall the definition of groups and some of its basic properties. In section two we define groupoids and loops and introduce a new class of loops and groupoids built using the modulo integers $Z_n$. Section three is devoted for giving the notions about semigroups and its substructures. In section four we recall the definition of semirings and its properties. In the final section we define lattices and give some properties about them.

## 1.1 Groups with examples

In this section we just recapitulate the definition of groups and some of its properties as the concept of groups is made use of in studying several properties of near-rings.

**DEFINITION 1.1.1**: *Let G be a non-empty set. G is said to form a group if in G there is defined a binary operation, called the product and denoted by '.' such that*

1. *a, b ∈ G implies a . b ∈ G.*
2. *a, b, c ∈ G implies that a . (b . c) = (a . b) . c.*
3. *There exists an element e ∈ G such that a . e = e . a = a for all a ∈ G.*
4. *For every a ∈ G there exists an element $a^{-1}$ ∈ G such that a . $a^{-1}$ = $a^{-1}$ . a = e.*

*A group G is said to be abelian (or commutative) if for every a, b∈ G, a . b = b . a. A group which is not abelian is called naturally enough non-abelian or non-commutative.*

*The number of elements in a group G is called the order of G and it is denoted by o(G) or |G|. When the number of elements in G is finite we say G is a finite group, otherwise the group is said to be an infinite group.*

**DEFINITION 1.1.2**: *A non-empty subset H of a group G is said to be a subgroup of G if under the product in G, H itself forms a group.*

**Example 1.1.1**: Let G = Q \ {0} be the set of rationals barring 0. G is an abelian group under multiplication and is of infinite order.

**Example 1.1.2**: Let S = {1, −1}, S is a group under multiplication. S is a finite group of order 2.

***Example 1.1.3***: Let X = { 1, 2 ,…, n } be set with n elements. Denote by $S_n$ the set of all one to one mappings of the set X to itself. Define a binary operation on $S_n$ as the composition of maps. $S_n$ is a group of finite order. $|S_n| = n!$ and $S_n$ is a non-commutative group.

Throughout this paper $S_n$ will be called the symmetric group of degree n or a permutation group on n elements. $A_n$ is a subgroup of $S_n$ which will be called as the alternating group and its order is n!/2.

**DEFINITION 1.1.3**: *A subgroup N of a group G is said to be a normal subgroup of G if for every g ∈ G and n ∈ N, $gng^{-1}$ ∈ N or equivalently if by $gNg^{-1}$ we mean the set of all $gng^{-1}$, n ∈ N then N is a normal subgroup of G if and only if $gNg^{-1}$ ⊂ N for every g ∈ G.*

*N is a normal subgroup of G if and only if $gNg^{-1}$ = N for every g ∈ G.*

**DEFINITION 1.1.4**: *Let G be a group, we say G is cyclic if G is generated by a single element from G.*

***Example 1.1.4***: Let G = $Z_5$ \ {0} be the set of integers modulo 5. G is a cyclic group. For 3 ∈ G generates G. G = {3, $3^2$, $3^3$, $3^4$} = {3, 4, 2, 1}, $3^2$ ≡ 4 (mod 5), $3^3$ ≡ 2 (mod 5) and $3^4$ ≡ 1 (mod 5). So G is a cyclic group of order 4.

It is obvious that all cyclic groups are abelian.

**DEFINITION 1.1.5**: *Let G and H be two groups. A map $\phi$ from G to H is a group homomorphism if for all a, b ∈ G we have $\phi(ab) = \phi(a)\,\phi(b)$.*

*The concept of isomorphism, epimorphism and automorphism can be defined in a similar way for groups.*

## 1.2 Definition of Groupoids and Loops with examples

In this section we introduce the notion of half-groupoids, groupoids and loops. Further we include here the loops and groupoids built, using the modulo integers $Z_n$. We illustrate these with examples.

**DEFINITION 1.2.1**: *A half-groupoid G is a system consisting of a non-empty set G and a binary operation '.' on G such that, for a, b ∈ G, a . b may be in G or may not be in G if a . b ∈ G we say the product under '.' is defined in G, otherwise the product under '.' is undefined in G.*

Trivially all groupoids are half-groupoids.

***Example 1.2.1***: Let G = {0, 1, 4, 3, 2}. Clearly 3 . 4 , 2 . 4, 3 . 2 are some elements not defined in G where '.' is the usual multiplication.

**DEFINITION 1.2.2**: *Let A be a non-empty set, A is said to be a groupoid if on A is defined a binary operation '∗' such that for all a, b ∈ A; a ∗ b ∈ A and in general a ∗ (b ∗ c) ≠ (a ∗ b) ∗ c for a, b, c ∈ A.*

*A groupoid A can be finite or infinite according as the number of elements in A are finite or infinite respectively. If a groupoid A has finite number of elements we say A is of finite order and denote it by o(A) or |A|.*

**Example 1.2.2**: The following table A = {a, b, c, d, e} gives a groupoid under the operation '∗'.

| ∗ | a | b | c | d | e |
|---|---|---|---|---|---|
| a | a | c | b | d | a |
| b | c | d | e | a | b |
| c | b | b | c | b | b |
| d | d | e | e | d | a |
| e | e | a | e | a | d |

Clearly (A, ∗) is a groupoid and the operation '∗' is non-associative and A is of finite order.

If a groupoid (A, ∗) contains an element i, called the identity such that a ∗ i = i ∗ a = a for all a ∈ A, then we say A is a groupoid with identity. If in particular a ∗ b = b ∗ a for all a, b ∈ A we say the groupoid is commutative.

Now we proceed on to define the new classes of groupoid using $Z_n$.

**DEFINITION 1.2.3**: *Let $Z_n$ = {0, 1, 2, ..., n – 1}, n a positive integer. Define an operation '∗' on $Z_n$ by a ∗ b = ta + ub (mod n) where t, u ∈ $Z_n$ \ {0} with (t, u) = 1 for all a, b ∈ $Z_n$. Clearly {$Z_n$, '∗', (t, u)} is a groupoid.*

Now for any fixed n and for varying t, u ∈ $Z_n$ \{0} with (t, n) = 1 we get a class of groupoids. This class of groupoids is called as the new class of groupoids built using $Z_n$. For instance for the set $Z_5$ alone we can have 10 groupoids. Thus the number of groupoids in the class of groupoids built using $Z_5$ is 10.

**Example 1.2.3**: Let ($Z_5$, ∗) be a groupoid given by the following table:

| ∗ | 0 | 1 | 2 | 3 | 4 |
|---|---|---|---|---|---|
| 0 | 0 | 4 | 3 | 2 | 1 |
| 1 | 3 | 2 | 1 | 0 | 4 |
| 2 | 1 | 0 | 4 | 3 | 2 |
| 3 | 4 | 3 | 2 | 1 | 0 |
| 4 | 2 | 1 | 0 | 4 | 3 |

Here t = 3 and u = 4, a ∗ b = 3a + 4b (mod 5). Clearly this is a non-commutative groupoid with no identity element having just five elements. We can also define groupoids using $Z_n$ where (t, u) = d, d > 1, t, u ∈ $Z_n$ \ {0}.

Now the groupoids in which (t, u) = d or t, u ∈ $Z_n$ \ {0} from a generalized new class of groupoids which contain the above-mentioned class. They have more number of elements in them. Just for comparison using $Z_5$ we have, if we assume (t, u) = d we have in this class 12 groupoids and if just t, u ∈ $Z_n$\{0} with no condition on t and u we have the number of groupoids in this new class of groupoids using $Z_5$ to be 16.

Now we proceed on to define Smarandache groupoid (S-groupoid).

**DEFINITION 1.2.4**: *A Smarandache groupoid (S-groupoid) G is a groupoid which has a proper subset S, S ⊂ G such that S under the operations of G is a semigroup.*

**Example 1.2.4**: Let (G, ∗) be a groupoid given by the following table:

| ∗ | 0 | 1 | 2 | 3 | 4 | 5 |
|---|---|---|---|---|---|---|
| 0 | 0 | 3 | 0 | 3 | 0 | 3 |
| 1 | 1 | 4 | 1 | 4 | 1 | 4 |
| 2 | 2 | 5 | 2 | 5 | 2 | 5 |
| 3 | 3 | 0 | 3 | 0 | 3 | 0 |
| 4 | 4 | 1 | 4 | 1 | 4 | 1 |
| 5 | 5 | 2 | 5 | 2 | 5 | 2 |

Clearly (G, ∗) is a S-groupoid for it has A = {2, 5} to be a semigroup.

**DEFINITION 1.2.5**: *Let (G, ∗) be any group. A proper subset P of G is said to be a subgroupoid of G if (P, ∗) is a groupoid.*

*Similarly (P, ∗) will be a Smarandache subgroupoid (S-subgroupoid) of G if (P, ∗) is a S-groupoid under the operation of ∗ ; i.e. P has a proper subset which is a semigroup under the operation '∗'.*

Now we proceed on to define free groupoid and Smarandache free groupoid as the notion of these concepts will be used in the construction of Smarandache semi-automaton and Smarandache automaton.

**DEFINITION 1.2.6**: *Let S be a non-empty set. Generate a free groupoid using S and denote it by ⟨S⟩. Clearly the free semigroup generated by the set S is properly contained in ⟨S⟩; as in ⟨S⟩, we may or may not have the associative law to be true.*

<u>Note</u>: If S = {a, b, c, …} we see even (ab)c ≠ a(bc) in general for all a, b, c ∈ ⟨S⟩. Thus unlike a free semigroup where the operations is assumed to be associative, in case of free groupoids we do not even assume the associativity while placing them in juxtaposition.

**THEOREM [98]**: *Every free groupoid is a S-free groupoid.*

*Proof*: Obvious from the fact that the set S which generates the free groupoid will certainly contain the free semigroup generated by S. Hence the claim.

Now we proceed on to define loops and give the major important identities.

**DEFINITION 1.2.7**: *(L, ∗) is said to be a loop where, L is a non-empty set and '∗' a binary operation called product defined on L satisfying the following conditions:*

1. *For a, b ∈ L we have a ∗ b ∈ L.*
2. *There exists an element e ∈ L such that a ∗ e = e ∗ a = a for all a ∈ L. e is called the identity element of L.*
3. *For every ordered pair (a, b) ∈ L × L there exists a unique pair (x, y) ∈ L × L such that ax = b and ya = b.*

*We say L is a commutative loop if a ∗ b = b ∗ a for all a, b ∈ L, otherwise L is non-commutative.*

*The number of elements in a loop is the order of the loop denoted by o(L) or |L|; if |L| < ∞ it is finite, otherwise infinite.*

***Example 1.2.5:*** (L, ∗) is a loop given by the following table:

| ∗ | e | $a_1$ | $a_2$ | $a_3$ | $a_4$ | $a_5$ |
|---|---|---|---|---|---|---|
| e | e | $a_1$ | $a_2$ | $a_3$ | $a_4$ | $a_5$ |
| $a_1$ | $a_1$ | e | $a_3$ | $a_5$ | $a_2$ | $a_4$ |
| $a_2$ | $a_2$ | $a_5$ | e | $a_4$ | $a_1$ | $a_3$ |
| $a_3$ | $a_3$ | $a_4$ | $a_1$ | e | $a_5$ | $a_2$ |
| $a_4$ | $a_4$ | $a_3$ | $a_5$ | $a_2$ | e | $a_1$ |
| $a_5$ | $a_5$ | $a_2$ | $a_4$ | $a_1$ | $a_3$ | e |

It is easily verified this loop is non-commutative and is of finite order.

**DEFINITION 1.2.8**: *Let (L, ∗) be a loop. A proper subset P of L is said to be a subloop of L if (P, ∗) is a loop.*

**DEFINITION 1.2.9**: *A loop L is said to be a Moufang loop if it satisfies any one of the following identities:*

1. *(xy) (zx) = (x(yz)) x.*
2. *((xy)z) y = x(y (zy)).*
3. *x(y(xz)) = ((xy)x)z for all x, y, z ∈ L.*

**DEFINITION 1.2.10**: *Let L be a loop. L is called a Bruck loop if (x(yx)) z = (x(y(xz)) and (xy)⁻¹ = x⁻¹y⁻¹ for all x, y, z ∈ L.*

**DEFINITION 1.2.11**: *A loop (L, ∗) is called a Bol loop if ((xy)z)y = x((yz)y) for all x, y, z ∈ L.*

**DEFINITION 1.2.12**: *A loop L is said to be right alternative if (xy) y = x (yy) for x, y ∈ L and L is left alternative if (xx) y = x(xy) for all x, y ∈ L. L is said to be an alternative loop if it is both right and left alternative loop.*

**DEFINITION 1.2.13**: *A loop (L, .) is called a weak inverse property loop (WIP-loop) if (xy) z = e imply x(yz) = e for all x, y, z ∈ L.*

**DEFINITION 1.2.14**: *A loop L is said to be semialternative if (x, y, z) = (y, z, x) for all x, y, z ∈ L where (x, y, z) denotes the associator of elements x, y, z ∈ L.*

We mainly introduce these concepts as we are interested in defining non-associative near-ring in chapter six and the concept of loops and groupoids will play a major role in that chapter.

Now we just proceed on to define Smarandache loops and conclude this section.

**DEFINITION 1.2.15**: *The Smarandache loop (S-loop) is defined to be a loop L such that a proper subset A of L is a subgroup with respect to the same induced operation that is $\phi \neq A \subset L$.*

***Example 1.2.6***: The loop (L, $*$) given by the following table is a S-loop.

| $*$ | e | $a_1$ | $a_2$ | $a_3$ | $a_4$ | $a_5$ | $a_6$ | $a_7$ |
|---|---|---|---|---|---|---|---|---|
| e | e | $a_1$ | $a_2$ | $a_3$ | $a_4$ | $a_5$ | $a_6$ | $a_7$ |
| $a_1$ | $a_1$ | e | $a_4$ | $a_7$ | $a_3$ | $a_6$ | $a_2$ | $a_5$ |
| $a_2$ | $a_2$ | $a_6$ | e | $a_5$ | $a_1$ | $a_4$ | $a_7$ | $a_3$ |
| $a_3$ | $a_3$ | $a_4$ | $a_7$ | e | $a_6$ | $a_2$ | $a_5$ | $a_1$ |
| $a_4$ | $a_4$ | $a_2$ | $a_5$ | $a_1$ | e | $a_7$ | $a_3$ | $a_6$ |
| $a_5$ | $a_5$ | $a_7$ | $a_3$ | $a_6$ | $a_2$ | e | $a_1$ | $a_4$ |
| $a_6$ | $a_6$ | $a_5$ | $a_1$ | $a_4$ | $a_7$ | $a_3$ | e | $a_2$ |
| $a_7$ | $a_7$ | $a_3$ | $a_6$ | $a_2$ | $a_5$ | $a_1$ | $a_4$ | e |

The subset $H_i$ = {e, $a_i$} for all $a_i \in$ L is a subgroup of L, hence L is S-loop.

Now we just recall the new class of loops constructed using $Z_n$.

**DEFINITION 1.2.16**: *Let $L_n(m)$ = {e, 1, 2, 3 ,…, n} be the set where n > 3, n is odd and m is a positive integer such that (m, n) = 1 and (m − 1, n) = 1 with m < n. Define on $L_n(m)$ a binary operation '.' as follows:*

*i)*    *e . i = i . e = i for all i ∈ $L_n(m)$.*

*ii)*   *i . i = $i^2$ = e for all i ∈ $L_n(m)$.*

*iii)*  *i . j = t where t = (mj − (m − 1)i) (mod n) for all i, j ∈ $L_n(m)$, i ≠ j, i ≠ e and j ≠ e.*

Clearly $L_n(m)$ is a loop. Let $L_n$ denote the new class of loops, $L_n(m)$ for a fixed n and varying m satisfying the conditions; m < n, (m, n) = 1 and (m − 1, n) = 1, that is $L_n$ = {$L_n(m)$ / n > 3, n odd, m < n, (m, n) = 1, (m − 1, n) = 1}.

We see each loop in this new class is of even order. For in $L_5$ we have only 3 loops given by $L_5(2)$, $L_5(3)$ and $L_5(4)$ none of them are groups. Just for illustrative purpose we give the loop $L_5(3)$. The loop $L_5(3)$ in $L_5$ is given by the following table:

| . | e | 1 | 2 | 3 | 4 | 5 |
|---|---|---|---|---|---|---|
| e | e | 1 | 2 | 3 | 4 | 5 |
| 1 | 1 | e | 4 | 2 | 5 | 3 |
| 2 | 2 | 4 | e | 5 | 3 | 1 |
| 3 | 3 | 2 | 5 | e | 1 | 4 |
| 4 | 4 | 5 | 3 | 1 | e | 2 |
| 5 | 5 | 3 | 1 | 4 | 2 | e |

This loop is commutative and is of order 6.

## 1.3 Semigroups

In this section we recall the concept of semigroups, Smarandache semigroups and some of its basic properties as this concept would be used to construct semigroup near-rings; a class of near-rings having a variety of properties.

**DEFINITION 1.3.1**: *Let S be a non-empty set. S is said to be a semigroup if on S is defined a binary operation '.' such that for all a, b ∈ S, a . b ∈ S and (a . b) . c = a . (b . c) for all a, b, c ∈ S.*

*A semigroup S is said to be commutative if a . b = b . a for all a, b ∈ S. A semigroup in which ab ≠ ba for at least a pair is said to be a non-commutative semigroup.*

*A semigroup which has an identity element e ∈ S is called a monoid, if e is such that a . e = e . a = a for all a ∈ S. A semigroup S with finite number of elements is called a finite semigroup and its order is finite and it is denoted by o(S) = |S|. If |S| is infinite we say S is a semigroup of infinite order.*

**DEFINITION 1.3.2**: *Let (S, .) be a semigroup. P a non-empty proper subset of S is said to be a subsemigroup if (P, .) is a semigroup.*

**DEFINITION 1.3.3**: *Let (S, .) be a semigroup. A non-empty subset I of S is said to be a right ideal (left ideal) of S if s ∈ S and a ∈ I then a s ∈ I (s a ∈ I). An ideal is thus simultaneously a left and a right ideal of S.*

***Example 1.3.1***: Let $Z_9$ = {0, 1, 2 ,…, 8} be the integers modulo 9. $Z_9$ is a semigroup under multiplication modulo 9. P = {0, 3, 6} is an ideal of $Z_9$.

**DEFINITION 1.3.4**: *Let S be a semigroup, with 0, an element 0 ≠ a ∈ S is said to be a zero divisor in S if there exists b ≠ 0 in S such that a.b = 0.*

**DEFINITION 1.3.5**: *Let S be a semigroup an element x ∈ S is called an idempotent if $x^2$ = x.*

**DEFINITION 1.3.6**: *Let S be a semigroup with unit 1. We say an element x ∈ S is said to be invertible if there exist y ∈ S such that xy = 1.*

Now we proceed on to define the concept of Smarandache semigroup.

**DEFINITION 1.3.7**: *A Smarandache semigroup (S-semigroup) is defined to be a semigroup A such that a proper subset of A is a group with respect to the operation of the semigroup A.*

**Example 1.3.2**: Let $Z_{12}$ = {0, 1, 2 ,…, 11} be the semigroup under multiplication modulo 12. $Z_{12}$ is a S-semigroup for {4, 8} = P is a group under the operations of $Z_{12}$.

**Example 1.3.3**: Let $Z_4$ = {0, 1, 2, 3} is a semigroup under multiplication modulo 4. $Z_4$ is a S-semigroup for P = {1, 3} is a group.

**DEFINITION 1.3.8**: *Let S be a S-semigroup. A proper subset P of S is said to be a Smarandache subsemigroup (S-subsemigroup) if P is itself a S-semigroup under the operations of S.*

**Example 1.3.4**: Let $Z_{24}$ = {0, 1, 2 ,…, 23} be the S-semigroup under multiplication modulo 24. P = {0, 2, 4 ,…, 22} is a S-subsemigroup of $Z_{24}$ as A = {8, 16} is a subgroup in P. Hence the claim.

## 1.4 Semirings

In this section we give the definition of semirings. As books on semirings is very rare we felt it would be appropriate to give this concept as it would be used in defining Smarandache seminear-rings and its Smarandache analogue.

**DEFINITION 1.4.1**: *Let (S, +, .) be a non-empty set on which is defined two binary operation '+' and '.' satisfying the following conditions:*

1. *(S, +) is an additive abelian group with the identity zero.*
2. *(S, .) is a semigroup.*
3. *a . (b + c) = a . b + a . c and (a + c) . b = a . b + c . b for all a, b, c ∈ S. We call (S, +, .) a semiring.*

*The semiring S is commutative if a . b = b . a for all a, b ∈ S. The semiring in which ab ≠ ba for some a, b ∈ S, we say S is a non-commutative semiring. The semiring has unit if e ∈ S is such that a . e = e . a = a for all a ∈ S.*

*The semiring S is said to be a strict semiring if a + b = 0 implies a = b = 0 for all a, b ∈ S.*

*The semiring has zero divisors if a . b = 0 without a or b being equal to zero i.e. a, b ∈ S \ {0}. A commutative semiring with unit which is strict and does not contain zero divisors is called a semifield.*

***Example 1.4.1***: S = $Z^+ \cup \{0\}$ the set of positive integers with 0 is a semifield under usual addition '+' and usual multiplication '.'.

***Example 1.4.2***: Let S = $Z^+ \cup \{0\}$; take P = S × S; clearly P is a semiring which is not a semifield under componentwise addition and multiplication.

**DEFINITION 1.4.2**: *Let (S, +, .) be a semiring. A proper subset P of S is said to be a subsemiring if (P, +, .) is itself a semiring.*

**DEFINITION 1.4.3**: *Let S be a semiring, S is said to be a Smarandache semiring (S-semiring) if S has a proper subset B which is a semifield with respect to the operations of S.*

**DEFINITION 1.4.4**: *Let S be a semiring. A non-empty proper subset A of S is said to be a Smarandache subsemiring (S-subsemiring) if A is a S-semiring i.e. A has a proper subset P such that P is a semifield under the operations of S.*

***Example 1.4.3***: Let $Z^o[x]$ be the polynomial semiring. $Z^o = Z^+ \cup \{0\}$. $Z^o$ [x] is a S-semiring and $Z^o$ is a S-subsemiring.

***Example 1.4.4***: $Q^o = Q^+ \cup \{0\}$ is a S-semiring under usual addition '+' and multiplication '.'.

**DEFINITION 1.4.5**: *Let S be a S-semiring. We say S is Smarandache commutative semiring (S-commutative semiring) if S has a S-subsemiring A where A is a commutative semiring.*

For more about semirings and S-semirings please refer [100, 101].

## 1.5 Lattices and its properties

This section is devoted to the introduction of lattices and some of its basic properties as the concept of lattices is very much used in near-rings.

**DEFINITION 1.5.1**: *Let A and B be two non-empty sets. A relation R from A to B is a subset of A × B. Relations from A to A are called relation on A, for short. If (a, b)∈R then we write aRb and say that 'a is in relation to b'. Also if a is not in relation to b we write a ℟ b.*

*A relation R on a nonempty set A may have some of the following properties:*

*R is reflexive if for all a in A we have aRa.*

*R is symmetric if for all a and b in A aRb implies bRa. R is transitive if for all a, b, c in A aRb and bRc imply aRc.*

*A relation R on A is an equivalence relation, if R is reflexive, symmetric and transitive. In case [a]: {b ∈ A / aRb} is called the equivalence class of a, for any a ∈ A.*

**DEFINITION 1.5.2**: *A relation R on a set A is called partial order (relation) if R is reflexive, anti-symmetric and transitive.*

*In this case (A, R) is called a partial ordered set or a poset.*

**DEFINITION 1.5.3**: *A partial order relation '≤' on A is called a total order (or linear order) if for each a, b ∈ A either a ≤ b or b ≤ a; (A, ≤) is then called a chain or a totally ordered set.*

**DEFINITION 1.5.4**: *Let (A, ≤) be a poset and B ⊆ A.*

  a) *a ∈ A is called an upper bound of B ⇔ ∀ b ∈ B, b ≤ a.*
  b) *a ∈ A is called a lower bound of B ⇔ ∀ b ∈ B a ≤ b.*
  c) *The greatest amongst the lower bounds whenever it exists, is called the infimum of B and is denoted by inf B.*
  d) *The least upper bound of B whenever it exists is called the supremum of B and is denoted by sup B.*

**DEFINITION 1.5.5**: *A poset (L, ≤) is called lattice ordered if for every pair of elements x, y in L the sup (x, y) and inf (x, y) exists.*

**_Remark_**:

i) Every ordered set is lattice ordered.
ii) If a lattice ordered set (L, ≤) the following statement are equivalent for all x and y in L

  a)  x ≤ y.
  b)  sup (x, y) = y.
  c)  inf (x, y) = x .

**DEFINITION 1.5.6**: *An algebraic lattice (L, ∩, ∪) is a non-empty set L with two binary operations ∩ (meet) and ∪ (join) (also called intersection or product and union or sum respectively) which satisfy the following conditions for all x, y, z ∈ L.*

  $L_1$     $x \cap y = y \cap x;$ *and* $x \cup y = y \cup x$
  $L_2$     $x \cap (y \cap z) = (x \cap y) \cap z$ *and* $x \cup (y \cup z) = (x \cup y) \cup z$
  $L_3$     $x \cap (x \cup y) = x$ *and* $x \cup (x \cap y) = x$

*Two application of $L_3$ namely* $x \cap x = x \cap (x \cap x)$ *lead to idempotent law. viz.* $x \cap x = x$ *and* $x \cup x = x.$

**DEFINITION 1.5.7**: *A nonempty subset S of a lattice L is called the sublattice of L if S is a lattice with respect to the restriction of ∪ and ∩ of L on to S.*

**DEFINITION 1.5.8**: *A lattice L is called modular if for all  x, y, z ∈ L x ≤ z implies x ∪ (y ∩ z) = (x ∪ y)  ∩ z this equation is known as the modular equation.*

***Example 1.5.1***: The lattice given by the following figure:

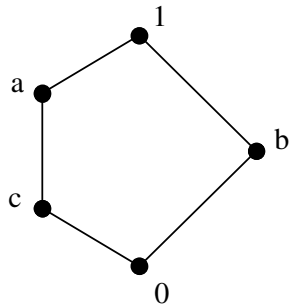

**Figure 1.5.1**

is called the pentagon lattice which is the smallest non-modular lattice.

***Example 1.5.2***: The lattice given by the following figure is a modular lattice.

This lattice will be called as the diamond lattice in this book.

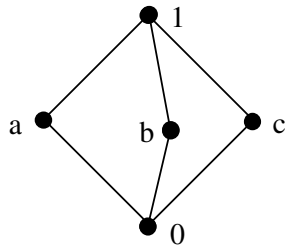

**Figure 1.5.2**

**DEFINITION 1.5.9**: *A lattice L is distributive if either of the following conditions hold for all x, y, z in L*

$$x \cup (y \cap z) = (x \cup y) \cap (x \cup z) \ or$$
$$x \cap (y \cup z) = (x \cap y) \cup (x \cap z).$$

*These equations are known as distributive equations.*

**Remark:**

1. The pentagon lattice is non-distributive.
2. All distributive lattices are modular.
3. The diamond lattice is non-distributive but modular. It is the smallest non-distributive lattice which is modular.

**DEFINITION 1.5.10**: *A lattice in which every pair of elements is comparable i.e. if a, b ∈ L we have a ≤ b or b ≤ a then we call such lattices as chain lattices. All chain lattices are distributive.*

**DEFINITION 1.5.11**: *Let L be a lattice with 0 and 1. Let x ∈ L we say x has a complement y in L if x ∩ y = 0 and x ∪ y = 1. A lattice L with 0 and 1 is called complemented if for each x ∈ L there is at least one element y in L such that x ∩ y = 0 and x ∪ y = 1; y is called a complement of x.*

**DEFINITION 1.5.12**: *A complemented distributive lattice is a Boolean algebra.*

***Example 1.5.3***: The following lattice is a Boolean algebra.

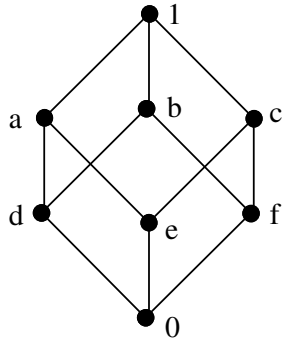

**Figure 1.5.3**

<u>**Remark**</u>: We have defined in the section 1.4 the concept of semirings and we see all Boolean algebras are semirings. In fact all distributive lattices are semirings.

**Chapter Two**

# NEAR-RINGS AND ITS PROPERTIES

This chapter has five sections; in the first section we just recall the definition of near-rings and give some examples. Section two is devoted to the introduction of N-groups, homomorphism and ideal like substructures. Direct products, and subdirect products in near-rings are introduced in section three. Section four is utilized in giving properties about ideals in near-rings. The fifth section is devoted to the recollection of modularity in near-rings. Near-polynomial rings are introduced in section six and the final section is devoted to introduce near-matrix rings. We have only recalled the definitions of these concepts and as we expect the reader to be well versed in the theory of near-rings the results and properties about near-rings should be known to the reader. Secondly as our main purpose is the study of Smarandache near-ring and Smarandache property in near-rings we do not even give problems for the reader in these sections.

## 2.1 Definition of Near-ring and some of its basic properties

In this section we just recall the definition of near-ring and illustrate it with examples. Several of its basic properties like homomorphism, direct products, N-groups etc are given in this section mainly to make it self contained. Certainly this does not mean that all properties about near-rings are given but we demand the reader to be well versed in near-ring theory. Most of the results given here are from the book of G.Pilz [61] who can be called as the father of near-rings.

**DEFINITION 2.1.1**: *A near-ring is a set N together with two binary operations '+' and '.' such that*

   a. *(N, +) is a group (not necessarily abelian).*
   b. *(N, .) is a semigroup.*
   c. *For all $n_1, n_2, n_3, \in N$; $(n_1 + n_2) . n_3 = n_1 . n_3 + n_2 . n_3$ (right distributive law).*

*This near-ring will be termed as right near-ring. If $n_1 (n_2 + n_3) = n_1 . n_2 + n_1 . n_3$ instead of condition (c) the set N satisfies, then we call N a left near-ring.*

In this text by a near-ring we will mean only a right near-ring. This convention will be strictly followed unless otherwise stated explicitly. Further we will write ab for a . b just for simplicity of notation.

***Example 2.1.1***: Z be the set of positive and negative integers with 0. (Z, +) is a group. Define '.' on Z by a . b = a for all a, b ∈ Z.   Clearly (Z, +, .) is a near-ring.

***Example 2.1.2***: Let $Z_{12}$ = {0, 1, 2 ,..., 11}. ($Z_{12}$, +) is a group under '+' modulo 12. Define '.' on $Z_{12}$ by a . b = a for all a ∈ $Z_{12}$. Clearly ($Z_{12}$, +, .) is a near-ring.

***Example 2.1.3***: Let $M_{2\times 2} = \{(a_{ij}) \,/\, Z$ ; Z is treated as a near-ring$\}$. $M_{2\times 2}$ under the operation of '+' and matrix multiplication '.' is defined by the following:

$$\begin{pmatrix} a & b \\ c & d \end{pmatrix} \cdot \begin{pmatrix} e & f \\ g & h \end{pmatrix} = \begin{pmatrix} a+b & a+b \\ c+d & c+d \end{pmatrix},$$

because we use the multiplication in Z i.e. a.b = a. So

$$\begin{pmatrix} a & b \\ c & d \end{pmatrix} \begin{pmatrix} e & f \\ g & h \end{pmatrix} = \begin{bmatrix} ae+bg & af+bh \\ ce+dg & cf+dh \end{bmatrix} = \begin{pmatrix} a+b & a+b \\ c+d & c+d \end{pmatrix}.$$

It is easily verified $M_{2\times 2}$ is a near-ring.

**DEFINITION 2.1.2**: *(N, +, .) is a near-ring. |N| or o(N) denotes the order of the near-ring N i.e. the number of elements in N. If |N| < ∞ we say the near-ring is finite if |N| = ∞ we call N to be an infinite near-ring. A non-empty set R is said to be a ring if on R is defined two binary operations '+' and '.' such that*

i)      *(R, +) is an abelian group.*
ii)     *(R, .) is a semigroup.*
iii)    *a . (b + c) = a . b + a . c and (a + b) . c = a . c + b . c for all a, b, c ∈ R.*

*If a . b = b . a for all a, b ∈ R we say R is a commutative ring. If 1 ∈ R is such that 1 . a = a . 1 = a for all a ∈ R; then R is a ring with unit. If R \ {0} is an abelian group under '.' we say R is a field. Thus every ring is a field.*

The following concepts are defined in near-rings as in rings : Left (right) identities, left (right) zero divisors, left (right) invertible elements, left (right) cancellable elements, idempotents and nilpotent elements.

***Remark***: Every ring is a near-ring.

*Notation*: Let N be a near-ring. d ∈ N is called distributive if for all n, $n_1$ ∈ N we have d(n + $n_1$) = dn + d$n_1$. $N_d$ = {d ∈ N / d is distributive}

**DEFINITION 2.1.3**: *Let N be a near-ring. If (N, +) is abelian, we call N an abelian near-ring. If (N, .) is commutative we call N itself a commutative near-ring. If N = $N_d$, N is said to be distributive. If all non-zero elements of N are left (right) cancelable, we say that N fulfills the left (right) cancellation law. N is an integral domain if N has no non-zero divisors of zero. If $N^* = N \setminus \{0\}$ is a group N is called a near-field. A near-ring, which is not, a ring will be referred to as non-ring. Similarly a non-field is a nf which is not a field. A near ring with the property that $N_d$ generates (N, +) is called a distributively generated near-ring.*

***Example 2.1.4***: $(Z_2, +, .)$ be the set. $Z_2 = \{0, 1\}$ under '+' usual addition '.' defined as 1.1 = 1 = 1.0 and 0.0 = 0.1 = 0 is a near-field.

Throughout this book η will denote the class of all near-rings.

## 2.2 N-groups, Homomorphism and ideal-like subsets

In this section we recall the definition of N-groups, N-homomorphism and ideals in a near-ring and illustrate with examples. For more about these concepts the reader is requested to refer [61].

**DEFINITION 2.2.1**: *Let $(P, +)$ be a group with 0 and let N be a near-ring. Let $\mu$ ; $N \times P \rightarrow P$; $(P, \mu)$ is called an N-group if for all $p \in P$ and for all $n, n_1, \in N$ we have $(n + n_1) p = np + n_1p$ and $(nn_1) p = n(n_1p)$. $N^P$ stands for N-groups.*

**DEFINITION 2.2.2**: *A subgroup M of a near-ring N with $M.M \subseteq M$ is called a subnear-ring of N. A subgroup S of $N^P$ with $NS \subset S$ is a N-subgroup of P.*

**Example 2.2.1**: Let $(Z, +, .)$ be the near-ring. $(2Z, +, .)$ is a subnear-ring.

**Example 2.2.2**: Let $(R, +, .)$ be a near-ring. R denote the set of reals. R is a N-group. Q the set of rationals is a N-subgroup of R.

**DEFINITION 2.2.3**: *Let N and $N_1$ be two near-rings. P and $P_1$ be two N-groups.*

  a. *$h : N \rightarrow N_1$ is called a near-ring homomorphism if for all $m, n \in N$*
     *$h(m + n) = h(m) + h(n)$ and $h(mn) = h(m) h(n)$.*
  b. *$h : P \rightarrow P_1$ is called N-homomorphism if for all $p, q \in P$ and for all $n \in N$*
     *$h(p + q) = h(p) + h(q)$ and $h(np) = nh(p)$.*

**DEFINITION 2.2.4**: *Let N be a near-ring and P a N-group. A normal subgroup I of $(N, +)$ is called an ideal of N if*

  i.   *$IN \subseteq I$.*
  ii.  *$\forall n, n_1 \in N$ and for all $i \in I$ , $n(n_1 + i) - nn_1 \in I$.*

*Normal subgroup T of $(N, +)$ with (i) is called right ideals of N while normal subgroup L of $(N, +)$ with (ii) are called left ideals. A normal subgroup S of $N^P$ is called ideal of $N^P$ if for all $s \in P$; $s_1 \in S$ and for all $n \in N$, $n(s + s_1) - ns \in S$. Factor near-ring $N / I$ and factor N-groups $P / S$ are defined as in case of rings.*

**DEFINITION 2.2.5**: *A subnear-ring M of N is called invariant if $MN \subseteq M$ and $NM \subseteq M$.*

It is pertinent to note that subnear-rings and ideals coincide in case of rings but does not in case of near-rings. A near-ring is simple if it has no ideals. $N^P$ is called N-simple if it has no N-subgroups except 0 and P.

**DEFINITION 2.2.6**: *A minimal ideal N is an ideal which is minimal in the set of all non-zero ideals. Similarly, one defines minimal right ideals, left ideals, N-subgroups. Dually, one gets the concept of maximal ideals etc.*

**DEFINITION 2.2.7**: *Let X, Y be subsets of $N^P$.*

$$(X : Y) = \{n \in N / nY \subseteq X\}.$$
$$(0 : X) \text{ is called the annihilator of } X.$$

*We denote it by $(X : Y)_N$ as N is the related near-ring. $N^P$ is faithful if $(0:P) = \{0\}$.*

## 2.3 Products, Direct sums and subdirect products in near-rings

In this section we recall the concepts of products, direct sums and subdirect products in near-rings and illustrate them with examples.

**DEFINITION 2.3.1**: *Let $\{N_i\}$ be a family of near-rings ($i \in I$, I an indexing set). $N_1 \times \ldots \times N_r = \underset{i \in I}{X} N_i$ with the component wise defined operations '+' and '.', N is called the direct product of the near-rings $N_i$.*

**Example 2.3.1**: Let $Z_2 = \{0, 1\}$ be the near-field. Z be the near-ring. $N = Z \times Z_2$ is a direct product of near-rings. Here $N = \{(z, p) / z \in Z \text{ and } p \in Z_2\}$ is a near-ring under component wise addition '+' and '.'.

**DEFINITION 2.3.2**: *The subnear-ring of $\underset{i \in I}{X} N_i$ consisting of those elements where all components except a finite number of them equal to zero is called the external direct sum $\underset{i \in I}{\bigoplus} N_i$ of the $N_i$'s. More generally, every subnear-ring N of $\underset{i \in I}{X} N_i$ where all projection maps $\pi_I$ ($i \in I$) are surjective (in other words for all $i \in I$ and for all $i_n \in N_i$, $n_i$ is the $i^{th}$ component of some element of N) is called the subdirect product of the $N_i$'s.*

**DEFINITION 2.3.3**: *Let N be a near-ring and S a subsemigroup of (N, +). A near-ring $N_s$ is called a near-ring of left (right) quotients of N with respect to S if*

   a. *$N_s$ has identity.*
   b. *N is embeddable in $N_s$, by a homomorphism, h.*
   c. *For all $s \in S$; $h(s)$ is invertible in $(N_s, .)$.*
   d. *For all $q \in N_s$ there exists $s \in S$ and there exist $n \in N$ such that $q = h(n) h(s)^{-1}$ ($q = h(s)^{-1} h(n)$).*

**DEFINITION 2.3.4**: *The near-ring N is said to fulfil the left (right) ore condition (ore (1)) with respect to a given subsemigroup S of (N, .) if for $(s, n) \in S \times N$ there exists $ns_1 = sn_1$ ($s_1n = n_1s$).*

**DEFINITION 2.3.5**: *If $S = \{s \in N / s \text{ is cancellable}\}$ then $N_s$ (if it exists) is called the left (right) quotient near-ring of N.*

*Notation:* Let V denote the set of all variety of near-rings (example - all near-rings, all abelian near-rings, all near-rings with unity, or all distributive near-rings). Let X be any non-empty subset.

**DEFINITION 2.3.6**: *A near-ring $F_X \in V$ is called a free near-ring in V over X if there exists $f : X \rightarrow F$ (where X is any non-empty set) for all $N \in V$ and for all $g : X \rightarrow N$ there exists a homomorphism $h \in Hom (F_X, N)$: h o f = g.*

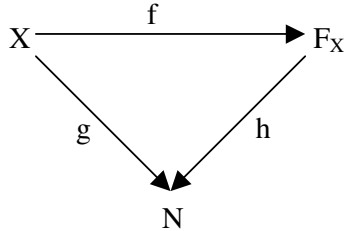

**Figure 2.3.1**

**DEFINITION 2.3.7**: *Let V = {set of all near-rings}, we simply speak about free near-ring on X. A near-ring is called free if it is free over some set X.*

On similar lines one can define free N-groups. For more about these properties refer [61].

Finally it is interesting to note that if G is a group. M(G) be the set of all self maps of G and define operations on M(G) as additions of maps '+' and compositions of maps 'o' clearly (M(G), +, o) is a near-ring.

## 2.4 Ideals in Near-rings

Here we just recall the definitions related to ideals in near-rings. All proofs and theorems related to them are left as exercise for the reader to study and refer.

**DEFINITION 2.4.1**: *Let N be a near-ring. $\{I_k\}_{k \in K}$ be the collection of all ideals of the near-ring N. Then the following are equivalent :*

a. *The set of all finite sums of the elements of the $I_k$'s.*
b. *The set of all finite sums of elements of different $I_k$'s.*
c. *The sum of the normal subgroups $(I_k, +)$.*
d. *The subgroup of (N, +) generated by $\bigcup_{k \in K} I_k$ .*
e. *The normal subgroups of (N, +) generated by $\bigcup_{k \in K} I_k$ .*
f. *The ideal of N generated by $\bigcup_{k \in K} I_k$ .*

*The set (a) – (f) above is called the sum of ideals $I_k$ ($k \in K$) and denoted by $\sum_{k \in K} I_k$ (for k = {1, 2, …} also by $I_1 + I_2$ …).*

**DEFINITION 2.4.2**: *Let $\{I_k\}_{k \in K}$ be the ideals of the near-ring N. Their sum $\sum_{k \in K} I_k$ is called an (internal) direct sum if each element of $\sum_{k \in K} I_k$ has a unique representation as a finite sum of elements of different $I_k$'s.*

*In this case we write for the sum $\sum_{k \in I} * I_k$ (or $I_1 + I_2 + I_3 \ldots$).*

**DEFINITION 2.4.3**: *Let I be an ideal of the near-ring N. I is called a direct summand of N if there exists an ideal J in N such that $N = I \dot{+} J$.*

**DEFINITION 2.4.4**: *A finite sequence $N = N_0 \supset N_1 \supset N_2 \supset \ldots \supset N_n = \{0\}$ of subnear-rings $N_i$ of N is called a normal sequence of N if and only if for all $i \in \{1, 2, \ldots, n\}$ $N_i$ is an ideal of $N_{i-1}$.*

In the special case that all $N_i$ is an ideal of N call the normal sequence an invariant sequence. n is called the length of the sequence and $N_{i-1} / N_i$ are called the factors of the finite sequence. Furthermore properties can be had from [61]. As in case of rings in case of near-rings also D.C.C and A.C.C can be defined.

**DEFINITION 2.4.5**: *Let N be a near-ring ($N^p$ a N-group). N is called decomposable if it is the direct sum of nontrivial ideals or if it has nontrivial direct summand, otherwise indecomposable. N is completely reducible if it is the direct sum of simple ideals.*

***Example 2.4.1***: The ring Z is indecomposable but is not simple. There exists nice relationship between completely reducible and ACCI and DCCI.

<u>*Notation*</u>: If S and T are ideals of N then ST = {st / s $\in$ S and t $\in$ T}. For n $\in$ N we can define $S^n$ in a similar way when S is an ideal of N.

**DEFINITION 2.4.6**: *Let P be an ideal of the near-ring N. P is called a prime ideal if for all ideals I and J of N ; $IJ \subseteq P$ implies $I \subset P$ or $J \subset P$.*

**DEFINITION 2.4.7**: *N is called a prime near-ring, if {0} is a prime ideal.*

**DEFINITION 2.4.8**: *Let $M \subset N$, N a near-ring. M is called a m – system if for all a, b $\in$ M there exists $a_1 \in (a)$ and $b_1 \in (b)$: $a_1 b_1 \in M$.*

**DEFINITION 2.4.9**: *Let S be an ideal in N, S is semiprime if and only if for all ideals I of N $I^2 \subseteq S$ implies $I \subset S$.*

**DEFINITION 2.4.10**: *N is called a semiprime near-ring if {0} is a semiprime ideal.*

**DEFINITION 2.4.11**: *Let I be an ideal of N; the prime radical of I denoted by P(I) = $\bigcap_{I \subseteq P} P$.*

**DEFINITION 2.4.12**: *n ∈ N is nilpotent if $n^k = 0$ for some k > 0, k an integer. S ⊂ N is nilpotent if $S^k = \{0\}$. S ⊂ N is called nil if all s ∈ S is nilpotent.*

## 2.5 Modularity in near-rings

In this section we recall the concept of modularity, quasi regularity, idempotents and centralizer in near-ring. We give only the basic definitions.

**DEFINITION 2.5.1**: *Let L be a left ideal in a near-ring N. L is called modular if and only if there exists e ∈ N and for all n ∈ N, n − ne ∈ L. In this case we also say that L is modular by e and that e is a right identity modulo L since for all n ∈ N; ne = n (mod L).*

***Results:***

1. Each modular left ideal L ≠ N is contained in a maximal one.
2. If N = $L_1$ + $L_2$ where $L_1$, $L_2$ are modular left ideals, then $L_1 ∩ L_2$ is again modular.

*Notation*: For z ∈ N, denote the left ideal generated by the set {n- nz / n ∈ N} by $L_z$. Note that for $L_Z$ = N, z = 0.

**DEFINITION 2.5.2**: *z ∈ N is called quasi regular if z ∈ $L_z$, S ⊂ N is called quasi regular if and only if for all s ∈ S, s is quasi regular.*

***Result***: Each nil ideal I of a near-ring N is quasi regular.

**DEFINITION 2.5.3**: *Let E be the set of idempotents in N. A set E of idempotents is called orthogonal if for all e, f ∈ E, e ≠ f implies e.f = 0.*

**DEFINITION 2.5.4**: *An idempotent e ∈ N is called central if it is in the center of (N, .) i.e. if for all n ∈ N; en = ne.*

**DEFINITION 2.5.5**: *A ring R is called biregular if each principal ideal is generated by an idempotent.*

**DEFINITION 2.5.6**: *A near-ring is said to be biregular if there exists some set E of central idempotents with*

a. *For all e ∈ N, Ne is an ideal of N.*
b. *For all n ∈ N there exist e ∈ E; Ne = (n) (principal ideal generated by n).*
c. *For all e, f ∈ E, e + f = f + e.*
d. *For all e, f ∈ E, ef ∈ E and e + f − ef ∈ E.*

[63] showed that each commutative biregular near-ring is isomorphic to a subdirect product of fields and hence is a biregular ring.

**DEFINITION 2.5.7**: *Let I be an ideal of N. The intersection of all prime ideals P which contain I is called the prime radical of I.*

We have just recalled some of the definitions, which we are interested in developing to Smarandache notions. We do not claim to have exhausted all the Smarandache notions about near-ring as in certain places we are not in a position to define Smarandache analogue. We have developed only those concepts which has interested us. The innovative reader can develop Smarandache notions for other concepts about near-rings which we have not dealt in this book.

## 2.6 Near Polynomial rings

In this section we define a new notion called near polynomial rings which are near-rings only when the near-ring N is of the special form, otherwise the near polynomial ring will fail to be a near-ring.

Throughout this section (N, +, .) will be a near-ring satisfying the following postulates (where N is a non-empty set).

1. (N, +) is a group with 0 as the additive identity.
2. For all m, n ∈ N we define m.n = m. Hence (N, .) is a semigroup.
3. 1 ∈ N is such that 1.x = 1 for all x ∈ N
4. (a + b)c = ac + bc is the only distributive law satisfied by N for all a, b, c ∈ N.

Clearly N is a near-ring. We use only this near-ring to built near polynomial rings.

**DEFINITION 2.6.1**: *Let (N, +, .) be a near-ring and x an indeterminate. N(x) will be called the near polynomial ring where*

$$N(x) = \left\{ \sum_{i=1}^{\alpha} n_i x^i \ / \ n_i \in N \right\}.$$

*In N(x) two polynomials p(x), q(x) will be equal; where*

$p(x) = p_0 + p_1 x + \ldots + p_n x^n$
$q(x) = q_0 + q_1 x + \ldots + q_m x^m$ *if and only if $p_i = q_i$ and m = n, for each i.*

*The addition of polynomials is given by*

*$p(x) + q(x) = (p_0 + q_0) + (p_1 + q_1) x + \ldots + p_n x^n (q_m x^m)$ if n > m (if m > n). Thus (N(x), +) is a group with $0 = 0 + 0x + \ldots + 0x^n + \ldots$ which acts as the 0 polynomial. For p(x), q(x) ∈ N(x) define '.' as p(x).q(x) = p(x). (N(x), .) is a semigroup. Clearly (p(x) + q(x)) r(x) = p(x) r(x) + q(x) r(x). Thus (N(x), +, .) is a near-ring called the near polynomial ring.*

**DEFINITION 2.6.2**: *Let N(x) be a near polynomial ring. For $0 \neq f(x) = a_0 + a_1 x + \ldots a_m x^m$ where $a_m \neq 0$; in N we define degree of f(x) to be m, thus deg f(x) is the index of the highest non-zero coefficient of f(x).*

$$deg\ [f(x)\ g(x)] = deg\ f(x).$$

*This is the marked difference between the polynomial ring and the near polynomial ring.*

***Example 2.6.1***: Let $Z_3$ [x] be the near polynomial ring of degree less than or equal to 2. That is $Z_3[x] = \{0, 1, 2, x, 2x, x + 2, x + 1, 2x + 1\ 2x + 2, x^2, 2x^2, x^2 + 2, x^2 + 1, x^2 + x, x + 2x^2, x + x^2 + 1, x + 2x^2 + 1, 2x + x^2 + 1\ 2 + x + x^2, 2x + 2x^2 + 2, 2x + 2 + x^2, 2x^2 + 2x + 1, x + 2x^2 + 2, 2 + 2x^2, 2x^2 + 2x, 2x^2 + 1\ 2x + x^2\}$. $Z_3[x]$ is a near polynomial ring with 27 elements in it.

***Example 2.6.2***: Let $Z_2 = \{0, 1\}$ be the near-ring. $Z_2[x]$ is a near polynomial ring. Thus we get a class of near polynomial rings using $Z_n$, Q or R but only restriction is a.b = a for all a, b $\in$ $Z_n$ (or Q or R).

Thus this gives us a new class of near-rings which will certainly have more interesting properties.

## 2.7 Near Matrix rings

In this section we introduce a new class of near-rings called near matrix ring. We will show near matrix rings are also near-rings. To build near matrix rings we take only near-rings N, defined in the earlier section i.e. in section 2.6 in which we have a.b = a for all a, b $\in$ N and (a + b).c = a.c + b.c is the distributive law satisfied by N. We define near matrix rings as follows:

**DEFINITION 2.7.1**: *Let N be a near-ring. Consider the set of all n × n matrices with entries from N is denoted by $M_{n \times n}$. Let M $\in$ $M_{n \times n}$ then*

$$M = (\ a_{ij}\ ) = \begin{bmatrix} a_{11} & a_{12} & \dots & a_{1n} \\ a_{21} & a_{22} & \dots & a_{2n} \\ \vdots & \vdots & & \vdots \\ a_{n1} & a_{n2} & \dots & a_{nn} \end{bmatrix}$$

*we define '+' on $M_{n \times n}$ as in case of rings. The product of two n × n matrices A and B in $M_{n \times n}$ is defined as follows A = $(a_{ij})$ and B = $(b_{ij})$*

$$A.B = \begin{bmatrix} a_{11} & \dots & a_{1n} \\ a_{21} & \dots & a_{2n} \\ \vdots & & \vdots \\ a_{n1} & \dots & a_{nn} \end{bmatrix} . \begin{bmatrix} b_{11} & \dots & b_{1n} \\ b_{21} & \dots & b_{2n} \\ \vdots & & \vdots \\ b_{n1} & \dots & b_{nn} \end{bmatrix} = \begin{bmatrix} p_{11} & \dots & p_{1n} \\ p_{12} & \dots & p_{2n} \\ \vdots & & \vdots \\ p_{n1} & \dots & p_{nn} \end{bmatrix}$$

*where*

$$p_{11} = a_{11} + a_{12} + \dots + a_{1n}$$
$$\vdots$$
$$p_{1n} = a_{11} + a_{12} + \dots + a_{in}$$
$$\vdots$$
$$p_{n1} = a_{n1} + a_{n2} + \dots + a_{nn}$$
$$\vdots$$
$$p_{nn} = a_{n1} + a_{n2} + \dots + a_{nn}$$

It is easily verified $(M_{n \times n}, +, .)$ is a near-ring called the near matrix ring.

*Now let*

$$I_{n \times n} = \begin{bmatrix} 1 & 0 & \dots & 0 \\ 0 & 1 & \dots & 0 \\ \vdots & \vdots & & \vdots \\ 0 & 0 & \dots & 1 \end{bmatrix}$$

*$I_{n \times n} A \neq I_{n \times n}$ and $A\, I_{n \times n} \neq A$ for $A \in M_{n \times n}$ also $AB \neq A$ for $A, B \in M_{n \times n}$.*

***Example 2.7.1***: Let $M_{2 \times 2} = \{(a_{ij}) \,/\, a_{ij} \in Z\}$. $Z$ is a near-ring with usual addition but a . b = a for all a, b $\in$ Z. This near-ring $M_{2 \times 2}$ has zero divisors for take $A = \begin{pmatrix} 1 & -1 \\ -5 & 5 \end{pmatrix}$.

Clearly AB = (0) for all B $\in$ $M_{2 \times 2}$. This new class of near-rings will prove to be useful in our study of Smarandache-near-rings. This near-ring has $I_{n \times n}$ such that $I_{n \times n} \times A \neq I_{n \times n}$ for all A $\in$ $M_{n \times n}$. In fact

$$I_{n \times n} \times A = \begin{pmatrix} 1 & 1 & \dots & 1 \\ 1 & 1 & \dots & 1 \\ \vdots & \vdots & & \vdots \\ 1 & 1 & \dots & 1 \end{pmatrix}.$$

for all A $\in$ $M_{n \times n}$. These near-rings enjoy several interesting properties and is a class of near-rings different from the usual near-rings.

**Chapter Three**

# SPECIAL CLASSES OF NEAR-RINGS AND THEIR GENERALIZATIONS

This chapter has 5 sections. As the name of the chapter suggests it deals with the introduction of several properties of near-rings, seminear-rings. We do not promise to prove any of the properties about these near-rings or seminear-rings. We in fact recall the definition of group near-rings and give the classical problems to be solved about these structures.

In fact till date no one has introduced or studied the concept of semigroup near-ring. In this chapter we define and introduce those concepts. Section one recalls the definition of several properties about near-rings. Section two is solely devoted to the introduction and study of the concepts like group near-rings, semigroup near-rings, group seminear-rings and semigroup seminear-rings. Section three is completely used to define a new concept called near loop-rings. Loop near-rings have been already studied in [64] but the structure of loop near-rings is constructed by replacing a group by a loop in the definition of near-ring, so this study is interesting as it gives a class of non-associative near-rings and its generalizations. The fourth section is a still generalization of section three for it introduces a new class of non-associative near-rings using groupoids over near-rings which the author chooses to call as groupoid near-rings. The study of this structure gives a non-abstract class of non-associative near-rings. The final section gives many special properties in near-rings and concepts related to near-rings.

## 3.1 IFP-near-rings

In this section we just recall the definition of IFP-near-rings as given by [61] and give some of the basic properties enjoyed by these near-rings.

**DEFINITION 3.1.1**: *Let N be a near-ring. N is said to fulfil the insertion of factors properly (IFP) provided that for all a, b, n ∈ N we have ab = 0 implies a n b = 0. N has strong IFP property if every homomorphic image of N has the IFP.*

**THEOREM 3.1.1**: *N has the strong IFP if and only if for all I in N and for all a, b, n ∈ N. ab ∈ I implies a n b ∈ I.*

*Proof*: The proof is straightforward, hence omitted.

**THEOREM 3.1.2**: *The following assertions are equivalent*

    a) *N has the IFP property.*
    b) *for all n ∈ N (0:n) is an ideal of N.*
    c) *for all S ⊆ N; (0:S) is an ideal of N.*

*Proof*: Left for the reader to prove.

**DEFINITION 3.1.2**: *Let p be a prime. A near-ring N is called a p-near-ring provided that for all $x \in N$, $x^p = x$ and $px = 0$.*

**DEFINITION 3.1.3**: *A near-ring N is Boolean if and only if for all $x \in N$, $x^2 = x$.*

**DEFINITION 3.1.4**: *Let (G, +) be a group. Define '*' on G by $a * b = 0$ for all a, b $\in$ G. Then (G, +, *) is a IFP near-ring.*

**DEFINITION 3.1.5**: *Let (G, +) be a group. Define '*' on G by $a * b = a$ if $b \neq 0$, $a * b = 0$ if b = 0. Then (G, +, *) is a IFP near-ring.*

Now we get loop near-rings which are IFP as follows:

**DEFINITION 3.1.6**: *Let (L, +) be a loop. Define * on L by $a * b = 0$ for all a, b $\in$ L; then (L, +, *) is a IFP-loop near-ring.*

**DEFINITION 3.1.7:** *Let (L, +) be a loop. Define * on L by $a * b = 0$ if $b \neq 0$ and $a * b = 0$ if b = 0. Now (L, +, *) is a IFP-loop near-ring.*

**DEFINITION [80]**: *A right ideal I of a near-ring N is called right quasi reflexive if whenever A and B are ideals of N with $A.B \subset I$ then b (b' + a) − bb' $\in$ I for all a $\in$ A and for all b, b' $\in$ B.*

**DEFINITION [80]**: *A near-ring N is strongly subcommutative if every right ideal of it is right quasi reflexive; any near-ring in this class of rings is therefore right quasi reflexive.*

***Example 3.1.1:*** The class of semiprime and prime near-rings are strongly subcommutative.

**THEOREM [80]:** *Let N be a strongly commutative near-ring and e be an idempotent in N. If eN is a right ideal of N, then eN is a ideal of N.*

*Proof*: Since N is strongly subcommutative near-ring, eN is right quasi reflexive ideal of N. Now by definition we see eN is an ideal of N as $eN.N \subset eN$ we must have $Ne.N \subset eN$.

**THEOREM [80]:** *Let N be a strongly subcommutative near-ring and e a non-zero idempotent in N. If eN is a minimal right ideal of N then eNe is a division subnear-ring of N.*

*Proof*: Straightforward.

**<u>*Result [80]:*</u>** If eN is a minimal right ideal with e a central idempotent in N then

1) eN is a division near-ring.

2) If the right annihilator of e in N is zero then N is a division near-ring.

**DEFINITION [80]:** *Let N be a near-ring. S a subnormal subgroup of (N, +). S is called a quasi ideal of N if $SN \subseteq NS \subset S$ where by NS we mean elements of the form {n(n' + s) – nn' / for all s ∈ S and for all n, n' ∈ N} = NS.*

**THEOREM [80]:** *Let N be any strongly subcommutative near-ring and e a non-zero idempotent in N. If eNe is a minimal quasi ideal of N then eN is a minimal quasi reflexive ideal of N.*

*Proof*: Left for the reader as an exercise.

**DEFINITION [8]:** *A left near-ring N is said to satisfy left self-distributive identity if abc = abac for a, b, c ∈ N.*

**DEFINITION [8]:** *A near-ring N is said to be*

> *left permutable if abc = bac,*
> *right permutable if abc = acb*
> *medial if abcd = acbd*
> *right self-distributive abc = acbc*

*for all a, b, c, d ∈ N.*

**DEFINITION [19]:** *Let (S, +, .) be a ring and M a right S-module. Let R = S × M and define (α, s) ⊙ (β, t) = (αβ, sβ + t). Then (R, +, ⊙) is a left near-ring, the abstract affine near-ring inducted by S and M. If S is a right self-distributive near-ring and $MS^2 = 0$ then R is a right self-distributive near-ring. If MS = 0 and S is a medial, left permutable or left-self distributive then (R, ⊙) is a medial left-permutable or left-self distributive respectively.*

***Example 3.1.2:*** Let (G, +) be a group and let h be any idempotent endomorphism on (G, +), define a ⊗ b = h(b) for each a, b ∈ G. Then (G, +, ⊗) is both a right self-distributive and a left self-distributive near-ring.

***Example 3.1.3:*** Let F be a free near-ring and J the ideal of F generated as an ideal by {abc – abac / a, b, c ∈ F} then F/J is a left-self distributive near-ring.

**THEOREM 3.1.3**: *Let R be a left-distributive near-ring such that a, b,c ∈ R and X is a non-empty subset of R. $abc = aba^nc = ab^nc$ for $n \geq 1$. Hence $a^3 = a^n$ is an idempotent for $n \geq 3$ and if abc = 0 then bac = 0.*

*Proof*: By the very definition, the result follows.

Several results can be got in this direction [8].

**DEFINITION 3.1.8**: *A near-ring R is defined to be equiprime if for all $0 \neq a \in R$ and x, y ∈ R, arx = ary for all r ∈ R implies x = y. If P is an ideal of R, P is called an equiprime ideal of R if R/P is an equiprime near-ring. It can be proved if P is an ideal* *of r if and only if for all $a \in R \setminus P$, $x, y \in P$, $arx - ary \in P$ for all $r \in R$ implies $x - y \in P$.*

**DEFINITION [10]:** *An infra near-ring (INR) is a triple $(N, +, .)$ where*

> *i)    $(N, +)$ is a group.*
> *ii)   $(N, .)$ is a semigroup.*
> *iii)  $(x + y)z = xz - 0z + yz$ for all $x, y, z \in N$.*

**DEFINITION [10]:** *Let $I$ be a left-ideal of $N$. Suppose that $I$ satisfies the following conditions:*

> *i)    $a, x, y \in N$, $anx - any \in I$ for all $n \in N$ implies $x - y \in I$.*
> *ii)   $I$ is left invariant.*
> *iii)  $0N \subset I$.*

*Then $I$ is called an equiprime left-ideal of $N$. We call an idempotent $e$ of a near-ring to be central if it is in the center of $(N, .)$ i.e. for all $n \in N$, $ne = ex$.*

**DEFINITION 3.1.9**: *Let $N$ be a near-ring. $N$ is said to be partially ordered by $\leq$ if*

> *a)  $\leq$ makes $(N, +)$ into a partially (fully ordered group).*
> *b)  $\forall n, n^1 \in N$, $n \geq 0$ and $n' > 0$ implies $nn' \geq 0$.*

**DEFINITION [87]:** *Let $N$ be a near-ring. Let $I_1, I_2, \ldots, I_m$ be the collection of all right ideals of $N$. We say $N$ is a n-ideal near-ring (n an integer less than or equal to m) if for every n right ideals $I_1, I_2, \ldots, I_n$ of $N$ and for every n elements $X_1, X_2, \ldots, X_n$ in $N \setminus (I_1 \cup I_2 \cup \ldots \cup I_n)$ we have $\langle X_1 \cup I_1 \cup I_2 \ldots \cup I_n \rangle = \langle X_2 \cup I_1 \cup I_2 \cup \ldots \cup I_n \rangle = \ldots = \langle X_n \cup I_1 \cup \ldots \cup I_n \rangle$ where $\langle \rangle$ denotes the right ideal generated by $X_i \cup I_1 \cup \ldots \cup I_n$; $1 \leq i \leq n$.*

***Example 3.1.4***: Let $Z_8 = \{0, 1, 2, \ldots, 7\}$ be the group under addition and multiplication as

$$0.0 = 0.1 = 0.2 = \ldots = 0.7 = 0$$
$$1.0 = 1.1 = 1.2 = \ldots = 1.7 = 1$$
$$\vdots$$
$$7.0 = 7.1 = 7.2 = \ldots = 7.7 = 7.$$

$I_1 = \{2, 4, 6, 0\}$ and $I_2 = \{0, 4\}$ are the only ideals of $Z_8$. Clearly $Z_8$ is a 2-ideal near-ring.

***Example 3.1.5***: Let $Z_{12} = \{0, 1, 2, \ldots, 11\}$. Define usual '+' on $Z_{12}$ and '.' by a.b = a for all a, b $\in Z_{12}$. Clearly $(Z_{12}, +, .)$ is a near-ring. $I_1 = \{0, 6\}$, $I_2 = \{0, 4, 8\}$, $I_3 = \{0, 3, 6, 9\}$, $I_4 = \{0, 2, 4, 6, 8, 10\}$ be the set of all ideals of $Z_{12}$. $Z_{12}$ is not a 2-ideal near-ring for we take $I_1$ and $I_2$ and $3, 2 \in Z_{12} \setminus (I_1 \cup I_2)$, $\langle 3 \cup I_1 \rangle = I_3$ and $\langle 2 \cup I_2 \rangle = I_4$. But it is easily verified $Z_{12}$ is a 3-ideal near-ring and 4-ideal near-ring.

***Example 3.1.6***: Let $Z_7 = \{0, 1, 2, 3, \ldots, 6\}$. Clearly $(Z_7, +, .)$ is a near-ring which is not an n-ideal near-ring for any positive integer n.

**THEOREM 3.1.4**: *Let N be a near-ring. Let $I_1, I_2, \ldots, I_m$ be the collection of all right ideals of N. If $I_i \cap I_j = \{0\}$ for $i \neq j$ for all $1 \leq i, j \leq n$. Then N is a n-ideal near-ring $n \leq m$.*

*Proof*: Let $I_1, \ldots, I_m$ be the set of right ideals of N. If $L_i \cap L_j = \{0\}$ if $i \neq j$ for all $1 \leq i, j \leq n$. Let $X_1, \ldots, X_m \in N \setminus \{I_1 \cup I_2 \cup \ldots \cup I_m\}$. Then $\langle X_1 \cup I_1 \cup I_2 \cup \ldots \cup I_m \rangle = \langle X_2 \cup I_1 \cup \ldots \cup I_m \rangle = \ldots = \langle X_m \cup I_1 \cup \ldots \cup I_m \rangle = N$.

**THEOREM [87]**: *Let N be a near-ring. Let $I_1, I_2, \ldots, I_m$ be the collection of all right ideals of N. N is not an n-ideal near-ring $(n < m)$, if $I_1 \cup I_2 \cup \ldots \cup I_n \subseteq I_k$ where $I_k$ is a proper right ideal of N.*

*Proof*: Let $I_1, \ldots, I_m$ be any arbitrary collection of right ideals of N. Take $X \in I_k \setminus (I_1 \cup \ldots \cup I_n)$ and $Y \in N \setminus (I_1 \cup \ldots \cup I_n \cup I_k)$ then $\langle X \cup I_1 \cup I_2 \cup \ldots \cup I_n \rangle = I_k$ and $\langle Y \cup I_1 \cup \ldots \cup I_n \rangle \neq I_k$, hence N is not a n-ideal near-ring.

The concept of $\Gamma$-near-ring, a generalization of ring was introduced by Nobusawa [56] and generalized by [3]. A generalization of both the concepts of near-ring and the $\Gamma$-ring, namely $\Gamma$-near-ring was introduced by [67] and later studied by [9].

**DEFINITION [67]**: *Let $(M, +)$ be a group (not necessarily abelian) and $\Gamma$ be a non-empty set. Then M is said to be a $\Gamma$-near-ring if there exists a mapping $M \times \Gamma \times M \rightarrow$ (the images of $(a, \alpha, b)$ is denoted by $a\alpha b$) satisfying the following conditions :*

      *a.  $(a + b)\alpha c = a\alpha c + b\alpha c$   and*
      *b.  $(a\alpha b)\beta c = a\alpha(b\beta c)$*

*for all a, b, c $\in$ M and $\alpha, \beta \in \Gamma$. M is said to be a zero symmetric $\Gamma$ near-ring if $a\alpha 0 = 0$ for all $a \in M$ and $\alpha \in \Gamma$ where 0 denotes the additive identity in M.*

***Example [67]***: Let $(G, +)$ be a group and X be a non-empty set. Let $M = \{f / f: X \rightarrow G\}$. Then M is a group under pointwise addition. If G is non-abelian, then $(M, +)$ is non-abelian. To see this let $a, b \in \Gamma$ such that $a + b \neq b + a$. Now define $f_a$, $f_b$ from x to G by $f_a(x) = b$ for all $x \in X$, $f_a, f_b \in M$ and $f_a + f_b \neq f_b + f_a$. Thus if G is non-abelian then M is also non-abelian.

Let $\Gamma$ be a set of all mappings of G into X. If $f_1, f_2 \in M$ and $g \in \Gamma$, then, obviously $f_1 g f_2 \in M$ for all $f_1, f_2, f_3 \in M$ and $g_1, g_2 \in M$ it is clear that

      i)      $(f_1 g_1 f_2)g_2 f_3 = f_1 g_1 (f_2 g_2 f_3)$ and
      ii)     $(f_1 + f_2)g_1 f_3 = f_1 g_1 f_3 + f_2 g_1 f_3$.

But $f_1 g_1 (f_2 + f_3)$ need not be equal to $f_1 g_1 f_2 + f_1 g_1 f_3$. To see this fix $0 \neq z \in M$ and $u \in X$. Define $G_\mu : G \rightarrow X$ by $g_\mu(x) = \mu$ for all $x \in G$ and $f_z : X \rightarrow G$ by $f_z(x) = z$ for all $x \in X$. Now for any two elements $f_2, f_3 \in M$, consider $f_z g_\mu (f_2 + f_3)$ and $f_z g_\mu f_2 + f_z g_\mu f_3$.

For all x ∈ X $[f_z g_\mu (f_2 + f_3)](x) = f_z[g_\mu(f_2(x) + f_3(x))] = f_z(\mu) = z$ and $[f_z g_\mu f_2 + f_z g_\mu f_3](x)$ $= f_z g_\mu f_2(x) + f_z g_\mu f_3(x) = f_z(\mu) + f_z(\mu) = z + z$. Since $z \neq 0$ we have $z \neq z + z$ and hence $f_z g_\mu (f_2 + f_3) \neq f_z g_\mu f_2 + f_z g_\mu f_3$. Now we have the following:

If $(G, +)$ is non-abelian and X is a non-empty set then M = {f / f: X → G} is a non-abelian group under pointwise addition and there exists a mapping M × Γ × M → M where Γ = {g / g: G → X} satisfying the following conditions:

        i)        $(f_1 g_1 f_2) g_2 f_3 = f_1 g_1 (f_2 g_2 f_3)$ and
        ii)      $(f_1 + f_2) g_1 f_3 = f_1 g_1 f_3 + f_2 g_1 f_3$.

for all $f_1$, $f_2$, $f_3 \in$ M and for all $g_1$, $g_2 \in \Gamma$. Therefore M is a Γ-near-ring.

**DEFINITION [67]:** *Let M be a Γ-near-ring. Then a normal subgroup I of (M, +) is called*

        *i)*       *a left ideal if $a\,\alpha\,(b + i) - a\,\alpha\,b \in I$ for all a, b ∈ M, $\alpha \in \Gamma$ and i ∈ I.*
        *ii)*      *a right ideal if $i\,\alpha\,a \in I$ for all a ∈ M, i ∈ I and*
        *iii)*     *an ideal if it is both a left and a right ideal.*

**DEFINITION [67]:** *An ideal A of M is said to be prime if B and C are ideals of M such that $B\,\Gamma C \subseteq A$ implies $B \subseteq A$ or $C \subseteq A$.*

**DEFINITION [67]:** *Let $M_1$ and $M_2$ be Γ-near-rings. A group homomorphism f of $(M_1, +)$ into $(M_2, +)$ is called a Γ-homomorphism if $f(x\,\alpha\,y) = f(x)\,\alpha\,f(y)$ for all x, y ∈ M and $\alpha \in \Gamma$. We say f is a Γ-isomorphism if f is one to one and onto.*

For an ideal I of a Γ-near-ring the quotient Γ-near-ring M/I is defined as usual:

**THEOREM [67]:** *Let I be an ideal of M and f, the canonical group epimorphism from M onto M/I. Then f is a Γ-homomorphism of M onto M/I with kernel I. Conversely if f is a Γ-epimorhism of $M_1$ onto $M_2$ and I is the kernel of f then $M_1/I$ is isomorphic to $M_2$.*

*Proof*: Left for the reader to prove.

There is yet another theorem the proof of which is left as an exercise to the reader.

**THEOREM [67]:** *Let f be a Γ-homomorphism of $M_1$ onto $M_2$ with ker f = I and $T^*$ a non-empty subset of $M_2$. Then $T^*$ is an ideal of $M_2$ if and only if $f'(T^*) = T$ is an ideal of $M_1$ containing I. In this case we have $M_1/T$, $M_2/T^*$ and $(M_1/I)(T/I)$ are Γ-isomorphic.*

*Example [67]:* Let G be a non-trivial group and X be a non-empty set. If M is the set of all mappings from X into G and Γ be the set of all mappings from G into X, then M is a Γ-near-ring. Let y be a non-zero fixed element of G. Define $\phi$ : X → G by $\phi(x) =$ y for every x ∈ X. Then $0 \neq \phi \in$ M, where 0 is the additive identity in M and $\phi g 0 = \phi$ $\neq 0$ for any g ∈ Γ. Therefore M is a Γ-near-ring which is not zero-symmetric.

**PROBLEMS:**

1.  Give an example of a IFP-near-ring.
2.  Give an example of a finite loop near-ring.
3.  Let $Z_p$ be the near-field. L a loop say $L = L_7(3)$. Find zero-divisors of the near loop ring $Z_p L_7(3)$. What is the major difference between a loop near-ring and a near loop ring.
4.  Can $Z_{30}$ have quasi-ideals?
5.  Show $(Z_{11}, +, .)$ is not a n-ideal near-ring.
6.  Will $(Z_{15}, +, .)$ be a n-ideal near-ring for some suitable integer n?
7.  Prove $Z_p$ when p is a prime is never a n-ideal near-ring.
8.  Prove if R is a Boolean near-ring, the following are equivalent
    i)      R is LSD.
    ii)     R is left-permutable.
    iii)    R is medial.
    iv)     R is self-distributive.
9.  Prove if the near-ring R is self-distributive and if R is simple, then $R^2 = 0$ and $R^+$ is a simple group.

## 3.2 Group near-rings and its generalizations

In this section we introduce the concept of group near-rings analogous to group rings where rings are replaced by near-rings. Such study was started only in the year 1992. When group near-rings find so much of place why group near-rings never grew into a theory is really incomprehensible. As there are several books on group rings there is not even a single book on group near-rings neither is the study become very popular. One of the reasons for it may be the following [61] concentrates only on the study of ideals, distributivity and idempotents. Apart from this one fails to study the isomorphism problem or semisimplicity problem or other famous problems about group near-rings.

Just for the sake of completeness we recall the definition of group near-rings and develop the subject analogous to group rings.

**DEFINITION 3.2.1**: *Let G be any group under multiplication. N a near-ring with 1 which is assumed to be commutative or a near-field. The group near-ring NG consists of all finite formal sums of the form $\Sigma \alpha_i g_i$ where i runs over a finite index and $\alpha_i \in N$ and $g_i \in G$; satisfying the following conditions.*

> 1. $\Sigma \alpha_i g_i = \Sigma \beta_i g_i \Leftrightarrow \alpha_i = \beta_i$ where $\alpha_i, \beta_i \in N$.
> 2. $\Sigma \alpha_i g_i + \Sigma \beta_i g_i = \Sigma (\alpha_i + \beta_i) g_i$.
> 3. $\alpha_i g_i = g_i \alpha_i$.
> 4. $(\Sigma \alpha_i g_i)(\Sigma \beta_j h_j) = \Sigma \lambda_k k$ and $g_i h_j = k$ and $\lambda_k = \Sigma \alpha_i \beta_j$.

*whenever the elements $\alpha_i$'s and $\beta_j$'s are distributive, if they are not distributive we assume ( $\Sigma \alpha_i g_i = (\alpha_1 g_1 + ... + \alpha_n g_n)$ and $\Sigma \beta_j h_j = (\beta_1 h_1 + ... + \beta_m h_m)) \Sigma \alpha_i g_i \Sigma \beta_j h_j = \alpha_1 g_1 (\beta_1 h_1 + ... + \beta_m h_m) + \alpha_2 g_2 (\beta_1 h_1 + ... + \beta_m h_m) + ... + \alpha_n g_n (\beta_1 h_1$* $+ \dots + \beta_m h_m) = \alpha_1 (\beta_1 g_1 h_1 + \dots + \beta_m g_1 h_m) + \alpha_2 (\beta_2 g_2 h_2 + \dots + \beta_m g_2 h_m) + \dots + \alpha_n (\beta_1 g_n h_1 + \dots + \beta_m g_n h_m).$

*This will be defined in this text in two ways:*

1. $\Sigma \alpha_i g_i \Sigma \beta_j h_j = \Sigma \alpha_l g_i$ or
2. $\alpha_1 (g_1 h_1 + \dots + g_1 h_m) + \dots + \alpha_n (g_h h_1 + \dots + g_n h_m).$

*Here 'or' in the mutually exclusive sense. We assume $1 \cdot g = g$ for $g \in G$.*

**Example 3.2.1**: Let $G = \langle g / g^n = 1 \rangle$ and $Z_2 = \{0, 1\}$ be the near-field. The group near-ring; $Z_2 G = \left\{ 1, g, g^2, \dots, g^{n-1}, \sum_{i=1}^{n} \alpha_1 g^i \text{ where } \alpha_i \in \{0,1\} \right\}.$

**Example 3.2.2**: Let $G = \langle g / g^2 = 1 \rangle$ and $Z_2 = \{0, 1\}$. The group near-ring $Z_2 G = \{0, 1, g, 1 + g\}$.

Study of zero divisor, units, idempotents, nilpotent elements happen to be an interesting feature but to the best of my knowledge there are only 5 papers on group near-rings where the study is made only relative to ideals, modules, direct decomposition and pseudo distributivity. The literature was got from the international near-ring community on searchable database [104]. Thus the author concludes if the study is taken online with the development and research as for the group ring, the group near-ring certainly will have more interesting results.

**DEFINITION 3.2.2**: *Let G be a group and N a near-ring. NG the group near-ring*

1. *has zero divisors if $\alpha\beta = 0$ if $\alpha, \beta \in NG \setminus \{0\}$.*
2. *$\alpha^2 = \alpha$ for $\alpha \in NG \setminus \{0, 1\}$ will be called idempotents of NG.*
3. *$\alpha\beta = 1$ where $\alpha, \beta \in NG \setminus \{1\}$ will be called units of NG.*
4. *$\alpha^n = 0$ where $\alpha \in NG \setminus \{0\}$ will be called as nilpotent elements of NG.*

We know if G is a finite group the group ring KG is artinian. We propose the analogous problems; when is NG artinian given G is a finite group and N a near-ring? Such study is classically important as near-rings are the generalizations of rings it is all the more relevant. Analogous to the classical theorem on group rings which says if K is a field of characteristic 0, then KG is semiprime for any group G; what is the analogous question and solution in case of group near-ring? Now the study of the question: If $Z_m$ be a near-ring and G be a group with no elements of order m, m a prime. Can we say the group near-ring $Z_m G$ has no nonzero nil ideals?

In case of group rings we have $Z_p G$ has no non-zero nil ideals if G has no elements of order p.

**Example 3.2.3**: Let $Z_3 = \{0, 1, 2\}$ be a near-ring under usual '+' modulo 3 and '.' defined by $a \cdot b = a$ for all $a, b \in Z_3$. $G = \langle g / g^2 = 1 \rangle$ be any group. The group near-ring $Z_3 G = \{0, 1, 2, g, 2g, 1 + g, 1 + 2g, 2 + g\}$. Does $Z_3 G$ have units, zero divisors or idempotents?

Is the group near-ring $Z_3G$ a p-near-ring? These two questions are left for the reader as exercise.

The concept of modular group near-ring is briefly described here following the definition of [39]. Let $Z_pG$ be the group near-ring of a group G of prime order p over the near-ring $Z_p$. We define the mod p-envelop of G denoted by $G^* = 1 + U$ where $U = \{\Sigma\ \alpha_i\ g_i\ /\ \Sigma\ \alpha_i = 0\}$; $G^*$ may be a group or a semigroup.

***Example 3.2.4***: Let $G = \langle g\ /\ g^2 = 1\rangle$ be a group of order two, $Z_2 = [0, 1]$ be the near-ring with 2 elements and $Z_2G$ be the group near-ring. $G^* = 1 + U = 1 + \{0, 1 + g\} = \{1, g\}$, so $G = G^*$.

***Example 3.2.5***: Let $G = \langle g\ /\ g^2 = 1\rangle$ and $Z_3 = \{0, 1, 2\}$ be the near-ring. For the group near-ring $Z_3G$ the mod p-envelop of G is $G^* = \{1, g, 2 + 2g, 1 + g\}$ is a semigroup under multiplication.

***Example 3.2.6***: Let $G = \langle g\ /\ g^2 = 1\rangle$ be a group and $Z_p = \{0, 1, 2, 3, \ldots, p - 1\}$ be a near-ring. Find $G^*$ for the group near-ring $Z_pG$.

What is the order of $G^*$? Is $G^*$ a semigroup?

It is still more surprising to see that the notion of semigroup near-rings have not even been introduced, we on the contrary see that semigroup rings have become a very developed concept. Now semigroup near-rings are constructed in an analogous way as group near-rings. If we replace the group by a semigroup in the definition of group near-rings we get semigroup near-rings.

All this study is significant as the class of semigroup near-rings is a more generalized concept than that of the group near-rings or to be more precise the class of group near-rings is strictly and completely contained in the class of semigroup near-rings.

We proceed on to study some famous and important concepts about semigroup near-rings and at the same time we bring to the notice of the reader that a varied and a vast research can be carried out in this direction. Further we suggest several problems for the reader to solve in chapter ten.

The study and characterization of zero divisors, units, idempotents, nilpotents in semigroup near-rings happens to be a routine one left for the reader to develop by putting conditions on both the semigroups and near-rings. We first give the condition for a semigroup near-ring to be a p-near-ring, secondly for it to be a Marot near-ring and finally for it to be the mod p-envelop of a semigroup S when $Z_p$ is the integral near-ring of order p.

Throughout this section we assume $Z_p$ to be a near-ring of the form $Z_p$ where $Z_p$ is a group under '+' modulo p and multiplication '.' defined by a . b = a for all a, b $\in Z_p$.

**THEOREM 3.2.1**: *Let S be a semigroup with identity 1. If the semigroup near-ring NS is a p-near-ring then S is a commutative semigroup with $x^p = x$ for every $x \in S$ and N is a p-near ring without identity.*

*Proof*: Given NS the semigroup near-ring is a p-near-ring. As $1 \in S$ we have $N.1 \subseteq NS$ so N is also a p-near-ring. Let $n \neq 0$ be in N and $s \in S$, $ns \in NS$; $(ns)^p = n^p s^p = ns^p = ns$ (given) forcing $s^p = s$. Since N is commutative we see S is also commutative.

**THEOREM 3.2.2**: *Let N be a finite p-near-field and S a commutative semigroup such that $s^p = s$ for every s in S and s . $s_1 = 0$ if $s \neq s_1$, for all s, $s_1 \in S$ with multiplicative identity 1, then NS is a p-near-ring.*

*Proof*: A finite p-near-field is isomorphic to the field $Z_p$ so we have $N = Z_p$; hence $Z_p$ S is a p-near-ring. For if $\alpha = \Sigma \, \alpha_i \, s_i \in Z_p$ S we have $\alpha^p = \alpha$ (using the fact $\alpha_I^p = \alpha_I$ for all $\alpha_I \in Z_p$ and $s^p = s$ for all s in S and S is a commutative semigroup with 1.

**THEOREM 3.2.3**: *Let N be a finite near-field. The semigroup near-ring NS is a p-near-ring if and only if N is a p-near-ring and S is a commutative semigroup with identity in which every $s \in S$; $s^p = s$ and $ss_1 = 0$ if $s \neq s_1$ for all s, $s_1 \in S$.*

*Proof*: It is left for the reader to prove using the results given in the earlier theorems.

**THEOREM 3.2.4**: *Let S be a commutative idempotent semigroup in which $xy = 0$ if $x \neq y$ and $Z_p$ be the near-field of integers modulo p, p a prime. Then the semigroup near-ring $Z_p S$ is a p-near-ring.*

*Proof*: Obvious from the fact $x^p = x$ for every $x \in S$. So $Z_p S$ is a p-near-ring.

Finally it is left for the reader to prove.

**THEOREM 3.2.5**: *Let N be a Boolean near-ring and S is an idempotent semigroup with identity in which $x.y = 0$ if $x \neq y$. Then the semigroup near-ring NS is a Boolean near-ring.*

Recall by a regular element of a near-ring N, we mean a nonzero divisor in N. An ideal containing regular elements is called regular. A commutative near-ring with identity is called a Marot near-ring if each regular ideal of N is generated by a regular element of N.

Clearly every near-field is a Marot near-ring.

**THEOREM 3.2.6**: *Let S be an ordered semigroup without divisors of zero and with no elements of finite order and N be a near-field. The semigroup near-ring NS is a Marot near-ring.*

*Proof*: The proof follows from the fact that the semigroup near-ring NS has no zero divisors; using the fact S is ordered and has no elements of finite order and N is a near-field.

**THEOREM 3.2.7**: *Suppose NS is a Marot semigroup near-ring with nontrivial divisors of zero then*

1. *N has divisors of zero.*
2. *S has divisors of zero.*
3. *S has elements of finite order.*
4. *N has divisors of zero and S has divisors of zero and elements of finite order.*

*Proof*: Given the semigroup near-ring is a Marot near-ring, so both N and S are commutative as $N \subseteq NS$ and $S \subset NS$, N and S can have divisors of zero or S can have elements of finite order which obviously contribute to divisors of zero.

**THEOREM 3.2.8**: *Let N be a commutative Marot near-ring and $N \setminus \{0\} = N_d$. S be a commutative semigroup with identity, having no divisors of zero and $0 \notin S$. Then NS is a Marot semigroup near-ring.*

*Proof*: NS is a commutative near-ring as N and S are given to be commutative. To prove NS is a Marot near-ring we need to only prove that every regular ideal is generated by a regular element of NS; to prove this we have to prove.

  i. If a regular ideal I is generated by a regular element, I contains no proper divisors of zero.

  ii. If a regular ideal I is generated by a divisor of zero then I cannot contain regular elements.

Proof of i. Let I be generated by the regular element $\alpha$. Let $\beta \neq 0 \in I$ such that $\beta\lambda = 0$ ($\lambda \neq 0$). Now

$$\beta = \sum_i \alpha s_i, \quad \beta\lambda = 0 \quad = \sum \alpha s_i \lambda = \alpha \sum \lambda s_i = 0$$

a contradiction to our assumption $\alpha$ is regular $\therefore \beta \notin I$ with $\beta\lambda = 0$.

Proof of ii. Let I be generated by a divisor of zero say $\beta$ ($\beta \neq 0$, $\lambda \neq 0$) with $\beta\lambda = 0$. Let $\alpha \in I$ be a regular element

$$\alpha \sum_i \beta\delta_i; \; \alpha\lambda = \sum \beta\lambda\delta_i = 0.$$

A contradiction $\therefore \alpha \notin I$. Hence the theorem.

The following result is obvious and is left as an exercise for the reader.

**THEOREM 3.2.9**: *If NS is a Marot semigroup near-ring then N is a commutative near-ring and S is a commutative semigroup.*

**THEOREM 3.2.10**: *Let N be a near-field. If NS is a Marot semigroup near-ring with divisors of zero then S is a commutative semigroup with divisors of zero or S has elements of finite order.*

*Proof*: Straightforward and hence left for the reader to prove.

The study of mod p-envelop of group of order p and K a field of characteristic p have been studied by [39]. Such study for group near-rings have not been carried out by any researcher. Here we just mention about this analysis and leave the rest for the reader to prove. We study this presently for the semigroup near-rings $Z_pS$; where $Z_p$ is the integral near-ring and the mod p envelope of the semigroup S is denoted by S*.

**DEFINITION 3.2.3**: *Let $Z_pS$ denote the integral semigroup near-ring of the semigroup S over the near-ring $Z_p$. Let $U = \{\Sigma \alpha_i s_i / \Sigma \alpha_i = 0\}$ then $1 + U = S*$ is called the mod p envelope of the semigroup S over the integral near-ring $Z_p$.*

*Example 3.2.7*: Let $S = \{s / s^2 = s\}$ be a semigroup and $Z_2 = \{0, 1\}$ be the integral near-ring with usual '+' and in which '.' is defined by $a.b = a \; \forall \; a, b \in Z_2$. $Z_2S = \{0, 1, 0.s, 1.s, 1 + 0.s, \ldots\}$; $U = \{0, 0s, 1 + 1.s, 1 + 0.s + 1s\}$; $S* = 1 + U = \{1, 1 + 0s, 1.s, 0.s + 1.s\}$. S* can be easily verified to be a semigroup under multiplication and $|S*| = 4$. It is left for the reader to verify the following example.

*Example 3.2.8*: Let $Z_2 = \{0, 1\}$ and $S = \{a, b, ab / a^2 = a, b^2 = b$ and $ab = ba\}$ be the semigroup. The semigroup near-ring $Z_2S$ have its mod p envelope to be a semigroup S*, where $|S*| = 2^6$.

We only propose several problems for the reader. Now the study of the following classical problems would certainly be of interest to any researcher. Now suppose given two semigroup near-rings NS and $N_1S_1$, when are they isomorphic. This problem for the case of group rings remains open. The study of such structures will give many concrete examples of near-rings with varied and similar properties.

<u>**PROBLEMS:**</u>

1. Find for the group near-ring $Z_2S_3$, the zero divisors, units and idempotents? Does $Z_2S_3$ have right ideals?
2. For the group near-ring $Z_{20}S_3$ find
   1. ideals.
   2. subnear-rings.
   3. right ideals which are not ideals.
3. Find the order of the semigroup near-ring $Z_2S(3)$. Does $Z_2S(3)$ have ideals?
4. Is $Z_2S_5 \cong Z_5S_3$? Justify your answer! Obtain a nontrivial homomorphism between these two group near-rings.
5. Does the group near-ring ZG, (where $(Z, +, '.')$ is a near-ring in which $a.b = a$ for all $a, b \in Z$. G a finite group) have
   1. Zero divisors?
   2. idempotents?
   3. semi-idempotents?
   4. units?
   Find them when $G = S_4$.

### 3.3 Loop near-rings and its generalization

In this section we introduce the concept of non-associative near-rings and define two new classes of near-rings using groupoids and loops. They form a nice class of near-rings satisfying certain special properties as well as satisfying some famous identities like Moufang. Bol, alternative etc. But the study of loop near rings was existent even from the year 1978 [64].

Here we bring in the notions of bipotent loop near-ring, regular loop near-ring and several of the properties enjoyed by them.

Study is also done about groupoid near-rings and half groupoid near-rings. Two types of loop near-rings are studied: we use loop near-rings analogous to group rings, which we denote by near loop rings. Another type of loop near-rings are got when we replace the group in the near-ring by a loop and the semigroup remains so.

**DEFINITION [75]**: *The system $N = (N, +, ., 0)$ is called a right loop-half groupoid near-ring provided.*

1. *$(N, +, 0)$ is a loop.*
2. *$(N, .)$ is a half groupoid.*
3. *$(n_1.n_2).n_3 = n_1 (n_2.n_3)$ for all $n_1, n_2, n_3 \in N$ for which $n_1.n_2$, $n_2.n_3$, $(n_1.n_2).n_3$ and $n_1.(n_2.n_3) \in N$.*
4. *$(n_1 + n_2) n_3 = n_1.n_3 + n_2.n_3$ for all $n_1, n_2, n_3 \in N$ for which $(n_1 + n_2)n_3$, $n_2.n_3 \in N$; we just denote $n_1.n_3$ by $n_1n_3$. If instead of (4) N satisfies $n_1(n_2 + n_3) = n_1n_2 + n_1n_3$ for every $n_1, n_2, n_3 \in N$ then we say N is a left loop-half groupoid near-ring.*

*For any a belonging to the additive loop L we shall denote the unique left and right inverse of a by $a_l$ and $a_r$ respectively.*

**DEFINITION [75]**: *A right loop near-ring N is a system $(N, +, .)$ of double composition '+' and '.' such that*

1. *$(N, +)$ is a loop.*
2. *$(N, .)$ is a semigroup.*
3. *The multiplication '.' is right distributive over addition + i.e. for all $n_1, n_2, n_3 \in N$; $(n_1 + n_2).n_3 = n_1.n_3 + n_2.n_3$; here after by convention by a loop near-ring we mean only a right loop near-ring.*

**_Result [75]:_** A loop-half groupoid near-ring N is a loop near-ring if and only if $(N, .)$ is a groupoid.

**_Example 3.3.1_**: Every near ring is a loop near-ring.

**_Example 3.3.2_**: Let $(G, +, \overline{0})$ be an additive loop where $\overline{0}$ is the identity element of G. Let $\Delta$ be a proper subset of G. Define $a.b = a$ for all $a \in G$ and $b \in \Delta$. Then $(G, \text{'+'}, \text{'.'}, \overline{0})$ is a loop half groupoid near-ring.

**DEFINITION [63]**: *Let* $(G, +, \overline{0})$ *be a loop and* $\Delta$ *a subset of G. A set S of endomorphism of G is called a* $\Delta$ *centralize of G provided.*

1. *The zero endomorphism* $\overline{0} \in S$.
2. *$S / \hat{0}$ (component of $\overline{0}$ in S) is a group of automorphism of G.*
3. *$\phi(\Delta) \subseteq \Delta$ for all $\phi \in S$.*
4. *$\phi, \psi \in S$ and $\phi(\omega) = \psi(\omega)$ for some $\overline{0} \neq \omega \in \Delta$ imply $\phi = \psi$.*

**DEFINITION [63]**: *Let* $(G, +, \overline{0})$ *be a loop. Let* $\Delta$ *be a subset of G and S a* $\Delta$ *-centralizer of G. A mapping T of G into itself is called a* $\Delta$ *-transformation of G over S provided $T[\phi(\omega)] = \phi[T(\omega)]$ for all $\omega \in \Delta$ and $\phi \in S$. It is interesting to note that if $\overline{0} \in \Delta$ and T is a $\Delta$ -transformation of G over S then T fixes $\overline{0}$ i.e.*

$$T(\overline{0}) = (\overline{0}).$$

*We shall denote the set of all $\Delta$-transformations of G over S by N (S, $\Delta$). Further we see for any endomorphism $\phi$ of a loop G, $[\phi(g)]_r = \phi(g_r)$ for all $g \in G$.*

**_Result [63]:_** Let $(G, +, 0)$ be a loop. $\Delta$ a subset of G containing $\overline{0}$ and S a $\Delta$ - centralizer of G. Then N (S, $\Delta$) is a loop-half groupoid near-ring under addition and composition of mappings but N (S, $\Delta$) is not a loop near-ring.

**DEFINITION [63]**: *A non-empty subset M of a loop near-ring (N, '+', '.', 0) is said to be a subloop near-ring of N if and only if (M, '+', '.', 0) is a loop near-ring.*

**DEFINITION [63]**: *A loop near-ring; N is said to be zero symmetric if and only if $n0 = 0$ for every $n \in N$ where '0' is the additive identity. Zero symmetric loop near-ring will be denoted by $N_0$.*

**DEFINITION [63]**: *A loop near-ring is said to be a near-ring if the additive loop of the loop near-ring is an additive group.*

**DEFINITION [63]**: *An element a of a loop near-ring N is said to be a left (right) zero divisor in N, if there exists an element $b \neq 0$ in N such that $ab = 0$ $(ba = 0)$. Element which is both a left and a right zero divisor in N is called a two sided zero divisor or simply a zero divisor in N.*

**DEFINITION [63]**: *A map $\theta$ from a loop near-ring N into a loop near-ring $N_1$ is called a homomorphism if*

$$\theta(n_1 + n_2) = \theta(n_1) + \theta(n_2),$$
$$\theta(n_1 n_2) = \theta(n_1) \theta(n_2) \text{ for every } n_1, n_2 \in N.$$

**DEFINITION [37]**: *Let N be a loop near-ring. An additive subloop A of N is a N-subloop (right N-subloop) if $NA \subset A$ $(AN \subset A)$ where $NA = \{na / n \in N, a \in A\}$.*

**DEFINITION [63]**: *A non-empty subset I of N is called left ideal in N if*

1.  *(I, +) is a normal subloop of (N, +).*
2.  *n (n₁ + i) + nᵣn₁ ∈ I for each i ∈ I and n, n₁ ∈ N, where nᵣ denotes the unique right inverse of n.*

**DEFINITION [37]:** *A non-empty subset I of N is called an ideal in N if*

1.  *I is a left ideal.*
2.  *IN ⊆ I.*

**DEFINITION [38]:** *A loop near-ring N is said to be left bipotent if Na = Na² for every a in N.*

***Example [38]:*** Let $N_1 = \{0, 1, 2, 3\}$, $N_2 = \{0, 1, 2, 3, 4\}$ and $N_3 = \{0, 1, 2, 3, 4, 5, 6\}$ with additions defined as modulo 4, 5 and 7 respectively and multiplication by the following tables:

$(N_1, .)$

| . | 0 | 1 | 2 | 3 |
|---|---|---|---|---|
| 0 | 0 | 0 | 0 | 0 |
| 1 | 0 | 3 | 0 | 1 |
| 2 | 0 | 2 | 0 | 2 |
| 3 | 0 | 1 | 0 | 3 |

$(N_2, .)$

| . | 0 | 1 | 2 | 3 | 4 |
|---|---|---|---|---|---|
| 0 | 0 | 0 | 0 | 0 | 0 |
| 1 | 0 | 0 | 4 | 1 | 0 |
| 2 | 0 | 0 | 3 | 2 | 0 |
| 3 | 0 | 0 | 2 | 3 | 0 |
| 4 | 0 | 0 | 1 | 4 | 0 |

and $(N_3, .)$ is given by

| . | 0 | 1 | 2 | 3 | 4 | 5 | 6 |
|---|---|---|---|---|---|---|---|
| 0 | 0 | 0 | 0 | 0 | 0 | 0 | 0 |
| 1 | 0 | 1 | 2 | 4 | 4 | 2 | 1 |
| 2 | 0 | 2 | 4 | 1 | 1 | 4 | 2 |
| 3 | 0 | 3 | 6 | 5 | 5 | 6 | 3 |
| 4 | 0 | 4 | 1 | 2 | 2 | 1 | 4 |
| 5 | 0 | 5 | 3 | 6 | 6 | 3 | 5 |
| 6 | 0 | 6 | 5 | 3 | 3 | 5 | 6 |

Clearly $N_1$, $N_2$, $N_3$ are left bipotent near-rings, hence left bipotent loop near-rings further we see in general a left bipotent loop near-ring need not be a left bipotent near-ring by the following example [38].

***Example [38]:*** Consider the set N = {e, a, b, c, d} define addition '+' as follows:

| + | e | a | b | c | d |
|---|---|---|---|---|---|
| e | e | a | b | c | d |
| a | a | b | e | d | c |
| b | b | c | d | a | e |
| c | c | d | a | e | b |
| d | d | e | c | b | a |

Let the multiplication be defined by x.y = x for every x, y ∈ N. Since a + (b + d) ≠ (a + b) + d; (N, +) is a loop and not a group. One can easily see that $Na^2 = N$ for every a ∈ N. Thus N is a left bipotent loop near-ring which is not a left bipotent near-ring.

**DEFINITION [37, 38]:** *A loop near-ring N is said to be a s-loop near-ring if a ∈ Na for every a in N.*

**Example [37, 38]:** The following loop near ring (N, +, .) given by the following table where N = {0, a, b, c} is s-loop near-ring.

| + | 0 | a | b | c |
|---|---|---|---|---|
| 0 | 0 | a | b | c |
| a | a | 0 | c | b |
| b | b | c | 0 | a |
| c | c | b | a | 0 |

| . | 0 | a | b | c |
|---|---|---|---|---|
| 0 | 0 | 0 | 0 | 0 |
| a | 0 | a | b | c |
| b | 0 | b | 0 | 0 |
| c | 0 | c | b | c |

**DEFINITION [38]:** *A loop near-ring N is said to be regular if for each a in N there exists x in N such that a = axa.*

**THEOREM [38]:** *Let N be a s-loop near-ring, then N is regular if and only if for each a (≠0) in N there is an idempotent 'e' such that Na = Ne.*

*Proof*: Proof is left for the reader.

**DEFINITION [38]:** *A loop near-ring N is said to be strictly duo (duo) if every N-subloop (left ideal) is also a right N-subloop (right ideal).*

**DEFINITION [38]:** *For any subset A of a loop near-ring N we define $\sqrt{A} = \{ x \in N / x^n \in A$ for some n $\}$.*

**DEFINITION [38]:** *A loop near-ring N is called irreducible (simple) if it contains only the trivial N-subloops (ideals) (0) and N itself.*

**DEFINITION [38]:** *A loop near-ring N is called a loop near-field if it contains an identity and each non-zero element has a multiplicative inverse.*

**DEFINITION [38]:** *An ideal (P ≠ N) is called strictly prime (prime) if for any two N-subloops (ideals) A and B of N such that AB ⊂ P then A ⊂ P or B ⊂ P.*

**Definition [38]:** *A left ideal B of a loop near-ring N is called strictly essential if B ∩ K ≠ {0} for any zero N-subloop K of N. An element x in N is said to be singular if there exists a nonzero strictly essential left ideal A in N such that Ax = {0}.*

**DEFINITION [38]:** *An element x of a loop near-ring N is said to be central if xy = yx for all y ∈ N.*

**<u>*Result*</u>**: If $N^1$ is a homomorphic image of a regular loop near-ring N whose idempotents are central then the idempotents of $N^1$ are also central.

**DEFINITION [38]:** *A nonzero loop near-ring N is said to be subdirectly irreducible if the intersection of all the non-zero ideals of N is non-zero.*

Now the study of loop rings was existent soon after the definition of group rings. Loop rings are nothing but loops over rings analogous to group rings which are groups over rings. But the study of "Loop near-rings" i.e. loops over near-rings was non-existent till 1990. The author was the first one to study this concept [78]. As loop near-ring was not a loop over near-ring but in the definition of near-rings the additive group was replaced by a additive loop, so the author cannot call the study of loops over near-rings as loop near-rings so has decided to call these structure as near-loop rings which are loops over near-rings. Before we proceed on to define near loop rings we just give the definition of non-associative near-ring.

**DEFINITION [78]:** *Let (N, +, .) be a nonempty set. N is said to be a non-associative right near-ring if the following conditions are true.*

> *1. (N, +) is a group with identity 0.*
> *2. (N, .) is a non-associative semigroup (i.e. a groupoid).*
> *3. (a + b)c = ac + bc for all a, b, c ∈ N.*

*If (N, +, .) has an element 1 such that a.1 = 1.a = a for all a ∈ N, we call N a non-associative right near-ring with unit.*

*If (N, .) is a loop; barring the zero element from N, we call N a non-associative division right near-ring. If N \ {0} is a commutative loop under the operation '.' we call N a non-associative near-field.*

*We say (N, +, .) is a non-associative left-near ring if instead of 3 we have a(b + c) = a . b + a . c for all a, b, c ∈ N. We study in this text only non-associative right near-ring which we simply term as non-associative near-rings.*

**DEFINITION [78]:** *Let (L, ., 1) be a loop and N a near-ring with a multiplicative identity. The loop over near-ring, that is the near loop ring denoted by NL with identity is a non-associative near-ring consisting of all finite formal sums of the form α = Σ α (m)m; m ∈ L and α (m) ∈ N.*

*(where supp $\alpha$ = {m / $\alpha$(m) ≠ 0}, the support of $\alpha$ is finite) satisfying the following operational rules.*

1. $\Sigma\,\alpha(m)\,m = \Sigma\,\mu(m)\,m \Leftrightarrow \alpha(m) = \mu(m)$ *for all* $m \in L$.
2. $\Sigma\,\alpha(m)\,m + \Sigma\,\mu(m)\,m = \Sigma\,(\alpha(m) + \mu(m))\,m$; $m \in L$; $\alpha(m), \mu(m) \in N$.
3. $(\Sigma\,\alpha(m)\,m)\,(\Sigma\,\mu(m)\,n) = \Sigma\,\gamma(k)\,k$; $m, n, k \in L$ *where* $\gamma(k) = \Sigma\,\alpha(m)\,\mu(n)$; $mn = k \in L$ *(whenever they are distributive, otherwise it is assumed to be in NL).*
4. $n\,(\Sigma\,\alpha(m)\,m) = \Sigma\,n\,\alpha(m)\,m$ *for all* $n \in N$ *and* $m \in L$.
5. $\alpha(m)\,m = m\,\alpha(m)$ *for all* $\alpha(m) \in N$ *and* $m \in L$.

Dropping the zero components of the formal sum we may write $\alpha = \Sigma\,\alpha_i\,m_i$; $i = 1, 2$ ,…, n. Thus m → m . 1 is an embedding of N in NL. After the identification of N with N . 1 we shall assume that N is contained NL. It is easily verified NL is a near-ring and it is a non-associative near-ring. Throughout this section N is right near-ring with identity 1.

***Result [78]***: Let L be a Moufang loop. N a near-ring with identity. Then the near loop ring NL is a N-group.

*Proof*: Left for the reader to prove.

**THEOREM [78]:** *Let L be a Moufang loop and N a near-ring; NL the loop near-ring. The set of all ideals, (N-subgroups) of an N-group of NL forms a inductive Moore systems of NL.*

**THEOREM [78]:** *Let N be a commutative near-ring without divisors of zero and L a Moufang loop with no elements of finite order. Then the near loop ring NL contains nontrivial associative near domains.*

*Proof*: Left for the reader to prove.

We now proceed to define Marot loop near-rings.

**DEFINITION [83]:** *By a regular element of a near-ring N we mean a non-zero divisor in N. An ideal containing regular elements is called a regular ideal. A commutative near-ring with identity is called a Marot near-ring if each regular ideal of N is generated by a regular element of N.*

***Example 3.3.3***: Every near-field is a Marot near-ring.

**DEFINITION [83]:** *We call a commutative near loop ring NL with identity to be a Marot near-ring if each regular ideal of NL is generated by a regular element of NL; by a regular element we mean a nonzero divisor in NL.*

**THEOREM [83]:** *Let N be a near-field and L a commutative ordered loop without elements of finite order then the near loop ring NL is a Marot near loop ring.*

*Proof*: To prove NL is a Marot near loop ring it is sufficient if we prove NL has no zero divisors i.e. all elements in NL will be regular elements. To this end we make use of the following facts:

1. N is a near-field so mn = 0 if and only if m = 0 or n = 0.
2. L is ordered and has no elements of finite order. We prove this by assuming NL has zero divisors and arrive at a contradiction. Let

$$\alpha = \sum_{i=1}^{n} n_i g_i \ \text{ and } \ \beta = \sum_{j=1}^{m} m_j h_j \ \ (n_i, m_j \in N \setminus \{0\})$$

belong to NL; with $g_1 < g_2 < \ldots < g_n$ and $h_1 < h_2 < \ldots < h_m$ suppose $\alpha\beta = \Sigma n_i g_i \Sigma m_j h_j = \Sigma n_i m_j g_i h_j = 0$; now $n_i m_j \neq 0$ for all $1 \leq i \leq n$ and $1 \leq j \leq m$ with $g_1 h_1$ to be the smallest element and $g_n h_m$ to be the largest element. So $\alpha\beta = 0$ is impossible as $n_1 m_1 \neq 0$ and $n_n m_m \neq 0$. Thus $\alpha\beta = 0$ is impossible, hence NL has no zero divisors so all elements in NL are regular, thus NL is a Marot near loop ring.

The reader is requested to develop more properties for Marot near-rings.

**THEOREM [82]**: *Let L be a loop and N a near-ring, NL the near loop ring (loop L over the near-ring N). If $M \subset N$ is an invariant subnear-ring then the near loop ring ML an invariant subloop near-ring.*

*Proof*: Follows from the fact $M \subset N$, $MN \subset M$ and $NM \subset M$. Hence ML. NL $\subseteq$ ML and NL. ML $\subseteq$ ML.

**THEOREM [78]**: *Let N be a near-ring and L a disassociative loop then the near loop ring NL has proper near-ring.*

*Proof*: Follows from the fact L is diassociative; hence has subgroups and the group near-rings are near-rings.

Now using the definition of [63] and [64] we take the near-ring to be the non-associative near-ring with respect to '+' given by them. Now using this loop near-ring as a near-ring we take loop L under multiplication and form the loop-loop near-ring.

**DEFINITION [91]**: *Let L be a loop (under multiplication) N a loop near-ring which is taken as a loop near-field. The loop – loop near-ring NL of the loop L over the loop near-ring N with identity is not even an associative loop but a doubly non-associative near-ring as (NL, '+', '.') is non-associative with respect to both '+' and '.'. NL consists of all finite formal sums of the form $\alpha = \sum_{m \in L} \alpha(m)m;\ \alpha(m) \in N$ such that support $\alpha = \{m \mid \alpha(m) \neq 0\}$ is a finite set, satisfying the following operational rules:*

*1. $\sum_{m \in L} \alpha(m)m = \sum_{m \in L} \beta(m)m \Leftrightarrow \alpha(m) = \beta(m)$ for all $m \in L$ (we assume in the finite formal sums the coefficient of each element in L occurs only once in the summation). (For otherwise $(ag + bg) + cg \neq ag + (bg + cg)$ if $(a + b) + c \neq a + (b + c)$).*

2. $\sum_{m \in L} \alpha(m)m + \sum_{m \in L} \beta(m)m = \sum_{m \in L} (\alpha(m) + \beta(m))m$.

3. $(\Sigma \alpha(m)m) (\Sigma \beta(n)n) = \Sigma \gamma(k) k$; $\gamma(k) = \Sigma \alpha(m) \beta(n)$; $mn = k$; only when $\alpha(m)$'s are distributive in N if $\alpha(m)$'s are not distributive in N we assume the product is in NL.

4. $n \Sigma \alpha(m) m = \Sigma n \alpha(m)m$ if and only if $n \in N_d$. Dropping the zero components of the formal sum we may write $\alpha = \Sigma \alpha_i m_i$. Thus $n \to n.1_m$, $1_m$ is the identity of the loop L. So $N.1 = N$ so $N.1 = N \subseteq NL$. The element $1.1m$ acts as the identity element of NL.

This new concept has lead to define a new algebraic structure called Lield. Please refer [91] for more properties.

**DEFINITION [91]**: *A nonempty set (L, 0, '+', '.') is called a Lield if*

1. *L is an abelian loop under '+' with 0 acting as the identity.*
2. *L is commutative loop with respect to '.' ; 1 acts as the multiplicative identity.*
3. *The distributive laws $a (b + c) = ab + ac$ , $(a + b) c = ac + bc$ for all a, b, c ∈ L holds good.*

NL defined above will be a Lield only if L is a commutative loop and NL has no zero divisors and NL = $N_d$.

Near domains are defined in a different way by [42].

**DEFINITION [42]:** *The algebraic structure (D, '+', '.' 0, 1) is called a near domain if it satisfies the following axioms:*

1. *(D, +, 0) is a loop.*
2. *$a + b = 0 \Rightarrow b + a = 0$ for all a, b ∈ D.*
3. *$\langle D^*, '.' 1 \rangle$ is a group where $D^* = D \setminus \{0\}$.*
4. *$0.a = a.0 = 0$ for all a ∈ D.*
5. *$a (b + c) = ab + bc$ for all a, b, c ∈ D.*
6. *For every pair a, b ∈ D there exist $d_{a,b} \in D^*$ such that for every x ∈ D, $a + (b + x) = (a + b) + d_{a,b} x$.*

For more about near domains refer [42].

Using this we define loop near domains which are analogous to groups over rings.

**DEFINITION [78]:** *Let L be a finite loop under '+' and D be a near domain, the loop near domain DL contains elements generated by $d_i m_i$ where $d_i \in D$ and $m_i \in L$ where we admit only finite formal sums satisfying the following conditions:*

1. *(DL, 0, '+') is a loop under '+'.*
2. *$d_i m_i \cdot d_j m_j = d_i d_j m_i m_j$ where we assume $m_i m_j \in DL$ as an element and $d_i d_j = d_k d_i$, $d_j$, $d_k \in D$.*

3. *We assume $(m_i \cdot m_j) \cdot m_k = m_i (m_j \cdot m_k)$ where just the elements $m_i$ $m_j$ $m_k$ are juxtaposed; $m_i$ , $m_j$ , $m_k \in L$.*
4. *We assume only finite sequence of elements of the form $m_1$ , $m_2$ ,..., $m_k \in DL$.*
5. *$m_i m_j = m_j m_i$ for all $m_i$, $m_j \in L$.*
6. *$m_i$ $m_j$ ... $m_k = m_1{}^I$ $m_2{}^J$ ... $m_k{}^l$ if and only if $m_i = m_i{}^l$ for $i = 1, 2, 3, ..., k$.*
7. *Since L is a loop under '+' and D is also a loop under '+'; DL inherits the operation of '+' and '0' serves as the additive identity of the loop DL.*
8. *Since $1 \in D$ we have only $1.L \subset DL$ so L is a subset of DL and $1.m_i = m_i.1 = m_i$ for all $i \in L$.*

Here the study, under what conditions will DL contain a near domain is an interesting one.

**DEFINITION 3.3.1**: *We call an additive loop L to be conditionally associative if*

1. *L is closed under some product operation '.'.*
2. *For every pair a, b $\in$ L there exists $d_{a,b} \in L \setminus \{0\}$ such that for every $x \in L$; a + (b + x) = (a + b) + $d_{a.b}.x$.*

***Example 3.3.4***: Near domains are trivially conditionally associative loops. The following theorem is left for the reader as an exercise.

**THEOREM 3.3.1**: *Let L be a conditionally associative loop and D a near domain. The loop near domain DL is a conditionally associative loop.*

When we take for the near-ring the Boolean near-ring we can define the loop Boolean near-ring analogous to the loop near-ring. We just by default of notion call loops over Boolean near-rings as loop Boolean near-rings.

**THEOREM 3.3.2**: Let N be a Boolean near-ring and L a loop. The loop Boolean near-ring is not a Boolean near-ring.

*Proof*: Straightforward, hence left for the reader as an exercise.

**THEOREM 3.3.3:** *Let N be a Boolean near-ring such that $x.y = 0$ for all x, y $\in N$ , x $\neq$ y and L a finite loop. The loop Boolean near-ring NL has nontrivial idempotents.*

*Proof*: Let $|L| = n$ and $L = \{m_1 = 1, m_2 ,..., m_n\}$. Suppose $\alpha = \sum_{m_i \in L} \alpha_i m_i$ with $\alpha_i \in N$ with $\alpha_i$'s distinct in N then, $\alpha^2 = \alpha$ is a nontrivial idempotent of NL.

**THEOREM 3.3.4**: *Let L be a finite loop and $Z_2 = (0, 1)$ be the Boolean near-ring. The loop Boolean near-ring $Z_2L$ has nontrivial subnear-rings which are Boolean near-rings.*

*Proof*: $|L| = n$ with $L = \{1, m_2 ,..., m_n \}$, n odd. Take $S = \{0, 1 + m_2 + ... + m_n\}$. S is a Boolean subnear-ring of NL.

**Example 3.3.5**: Let $Z_p = \{0, 1, \ldots, p-1\}$, be a near ring and p be a prime. $Z_p$ is a group under usual addition, '+' modulo p.

$$
\begin{array}{ccccccc}
0.0 & = & 0.1 & = & \cdots & = & 0.p-1 & = & 0 \\
1.0 & = & 1.1 & = & \cdots & = & 1.p-1 & = & 1 \\
\vdots & & & \ddots & & & & & \vdots \\
p-1.0 & = & p-1.1 & = & \cdots & = & p-1.p-1 & = & p-1
\end{array}
$$

is a near-ring. L any loop. $Z_pL$ is a near loop ring which is a non-associative ring.

**Example 3.3.6**: Let $Z_2 = (0, 1)$ be a near-ring and L be a loop given by the following table:

| . | 1 | a | b | c | d | e |
|---|---|---|---|---|---|---|
| 1 | 1 | a | b | c | d | e |
| a | a | 1 | e | b | c | d |
| b | b | c | d | a | e | 1 |
| c | c | d | 1 | e | a | b |
| d | d | e | c | 1 | b | a |
| e | e | b | a | d | 1 | c |

$Z_2L$ is a near loop ring; $(1 + a) \in Z_2L$ is a zero divisor which can be easily verified.

$(1 + a) (1 + a + b + c + d + e) = 0$ but $(1 + a + b + c + d + e) (1 + a) \neq 0$; as $Z_2$ is not left distributive. $(1 + a + b + c + d + e)^2 = 0$.

**Example 3.3.7**: Let $Z_2 = (0, 1)$ be the near-ring and L be a loop given by the following table:

| . | 1 | a | b | c | d |
|---|---|---|---|---|---|
| 1 | 1 | a | b | c | d |
| a | a | 1 | c | d | b |
| b | b | d | a | 1 | c |
| c | c | b | d | a | 1 |
| d | d | c | 1 | b | a |

$Z_2L$ the near loop ring is a non-associative ring. $(1 + a + b + c + d)^2 = 1 + a + b + c + d$ is an idempotent of $Z_2L$.

In view of the two examples we give the following theorem:

**THEOREM 3.3.5:** *Let $Z_2 = \{0, 1\}$ be the near-ring. L be a loop of finite order. The near loop ring has atleast nontrivial idempotent if |L| is odd and has zero divisors in |L| is even.*

*Proof*: Let $|L| = n + 1$ with $L = \{1, m_1, \ldots, m_n\}$. Choose $\alpha = \sum_{m_i \in L} m_i$ if $|L|$ is odd then we get $\alpha^2 = \alpha$ and if $|L|$ is even we get $\alpha^2 = 0$; hence the claim.

Certainly the near loop ring $Z_2L$ may have more idempotents and zero divisors than the ones mentioned in the theorems. Now we proceed on to define a new concept called modular near loop rings; modular group rings of a finite p-group was studied by [39].

**DEFINITION [39]:** *Let $Z_pG$ be the group ring of a finite group G over the field $Z_p$ or the ring of integer $Z_n$. The mod p envelop of G is defined to be set $G^* = 1 + U$ where*

$$U = \left\{ \sum_i \alpha_i g_i \ / \ \Sigma\alpha_i = 0 \ \text{ where } \ \alpha_i \in Z_p \right\}.$$

*Example 3.3.8*: Let $Z_2 = \{0, 1\}$ be the prime field of characteristic 2 and $G = \{g \ / \ g^3 = 1\}$ be the group. The group ring $Z_2G$ has $G^* = 1 + U$ where $U = \{0, 1 + g, 1 + g^2, g + g^2\}$.

Thus $G^* = 1 + U = \{1, g, g^2, 1 + g + g^2\}$ is the mod p envelop of G.
We proceed on to define this concept in case of near loop rings.

**DEFINITION 3.3.2**: *Let L be a loop of finite order and $Z_p$ be the near-ring. $Z_pL$ be the near loop ring of the loop L over the near-ring $Z_p$. The mod p envelop of the loop, $L^* = 1 + U$ where $U = \{\Sigma\alpha_i m_i \ / \ \Sigma\alpha_i = 0, \ \alpha_i \in Z_p\}.$*

*Example 3.3.9*: Let L be a loop given by the following table with 5 elements.

|   | 1 | a | b | c | d |
|---|---|---|---|---|---|
| 1 | 1 | a | b | c | d |
| a | a | 1 | c | d | b |
| b | b | d | a | 1 | c |
| c | c | b | d | a | 1 |
| d | d | c | 1 | b | a |

$Z_2 = \{0, 1\}$ be the near-ring. $Z_2L$ be the near loop ring.

$L^* = 1 + U = 1 + \{0, 1 + a, 1 + b, 1 + c, 1 + d, a + b, a + c, a + d, b + c, b + d, c + d, 1 + a + b + c, 1 + a + b + d, 1 + a + d + c, 1 + b + c + d, a + b + c + d\} = \{1, a, b, c, d, 1 + a + b, 1 + a + c, 1 + a + d, 1 + b + c, 1 + b + d, 1 + c + d, a + b + c, a + b + d, a + c + d, b + c + d, 1 + a + b + c + d\}$. $L^*$ is not a loop; $L^*$ is a groupoid. It is easily verified $L^*$ is closed under the product.

In view of this we have the following result; the proof of which is left for the reader as an exercise.

**THEOREM 3.3.6:** *Let L be a loop of odd order n and $Z_2 = \{0, 1\}$ be the near-ring. $Z_2L$ be the near loop ring. Then $L^*$ is a groupoid with 1 and $|L^*| = 2^{n-1}$.*

All near-rings N in the following definitions are assumed to be non-associative for otherwise all identities will be trivially true.

**DEFINITION 3.3.3**: *Let N be a non-associative near-ring we call N a moufang near-ring if it satisfies any one of the following identities.*

     i.        *(xy) (zx) = (x (yz)) x.*
     ii.       *((xy) z) y = x (y (zy)).*
     iii.     *x (y (xz)) = ((xy) x) z for all x, y, z in N.*

**DEFINITION 3.3.4**: *We call a non-associative near-field N to be a Bruck near-ring if*

$$(x (yx)) z = x (y (xz))$$
*and* $(xy)^{-1} = x^{-1} y^{-1}$ *for all x, y, z* $\in$ *N.*

**DEFINITION 3.3.5**: *Let N be a non-associative near-ring. N is said to be a Bol near-ring if ((xy) z) y = x ((yz) y) for all x, y, z in N.*

**DEFINITION 3.3.6**: *We say the non-associative near-ring N to be right alternative near-ring if (xy) y = x (yy) for all x, y in N and left alternative if x (xy) = (xx) y. N is alternative if N is simultaneously right and left alternative.*

**DEFINITION 3.3.7**: *A non-associative near-ring N is said to have weak inverse property (WIP) if (xy) z = 1 implies x (yz) = 1 for all x, y, z in N.*

**THEOREM 3.3.7**: *Let $L_n(2)$ be the loop n > 3 and n odd and N any near-ring. The loop near-ring $NL_2(2)$ is a right alternative near-ring.*

*Proof*: Left for the reader to prove.

**THEOREM 3.3.8**: *Let $L_n(n-1)$ be the power associative loop and N any near-ring. $NL_n(n-1)$ is a left alternative near-ring.*

*Proof*: Left as an exercise to the reader to prove.

**Theorem 3.3.9**: *Let $L_n(m)$ be any loop in $L_n$. $NL_n(m)$ for no near-ring N is a alternative near-ring.*

*Proof*: Straightforward by definition.

**THEOREM 3.3.10**: $L_n(m) \in L_n$; *$NL_n(m)$ for no $L_n(m)$ is a Moufang near-ring.*

*Proof*: Left for the reader to prove.

**THEOREM 3.3.11**: *Let N be a near-ring. $L_n(m) \in L_n$ be any loop. The near loop ring $NL_n(m)$ is never a Bol near-ring.*

*Proof*: Follows from the fact $L_n(m)$ are never Bol loops.

Now we proceed on to define the concept of quasi regular elements in near loop rings.

**THEOREM 3.3.12**: *Let L be a loop and N any near-field. NL the near loop ring. Every element of the form $1 \pm m_j$, $m_j \in L$ is quasi regular in the near loop ring NL.*

*Proof*: Let $x = 1 \pm m_j$ choose $y = \pm m_i$ such that $m_i\, m_j = 1$. Then $y\,(1-x) = \pm m_i\,(1 - 1 \mp , m_j) = 1$ (Choose $-m_j$ if $x = 1 + m_i$ and $m_j$ if $x = 1 - m_i$), hence x is quasi regular in NL.

**THEOREM 3.3.13**: *Let N be a near-field and L = {1, $m_2$ ,…, $m_n$} be a loop. NL the near loop ring. Elements in L can never be quasi regular.*

*Proof*: Suppose $m_i \in L$ be quasi regular then it implies $\alpha = \Sigma a_i m_i \in NL$ is such that $\alpha\,(1 - m_i) = 1$ but $\alpha\,(1 - m_i) = \alpha - \alpha = 0$. Hence the claim.

Several interesting results can be obtained in this direction. It is left for the reader to prove or disprove if L = {$m_1 = 1$, $m_2$, …, $m_n$} be loop of finite order and N a near-field. $\alpha = \Sigma a_i\, m_i \in NL$ is quasi regular if and only if $\Sigma a_i$ is quasi regular in N. (Hint: $m_i$ runs over all elements of L and $a_i \neq 0$ for i = 1, …, n.).

**PROBLEMS**:

1.  Let N be the near-ring ($Z_8$, +, .) i.e. ($Z_8$, +) an additive group modulo 8 a.b = a for all a $\in$ N. Let L be the loop given by the following table:

    | .     | e     | $a_1$ | $a_2$ | $a_3$ | $a_4$ | $a_5$ |
    |-------|-------|-------|-------|-------|-------|-------|
    | e     | e     | $a_1$ | $a_2$ | $a_3$ | $a_4$ | $a_5$ |
    | $a_1$ | $a_1$ | e     | $a_3$ | $a_5$ | $a_2$ | $a_4$ |
    | $a_2$ | $a_2$ | $a_5$ | e     | $a_4$ | $a_1$ | $a_3$ |
    | $a_3$ | $a_3$ | $a_4$ | $a_1$ | e     | $a_5$ | $a_2$ |
    | $a_4$ | $a_4$ | $a_3$ | $a_5$ | $a_2$ | e     | $a_1$ |
    | $a_5$ | $a_5$ | $a_2$ | $a_4$ | $a_1$ | $a_3$ | e     |

    a.  Define the near loop ring $Z_8 L$.
    b.   How many elements does $Z_8 L$ have? Find zero divisors in $Z_8 L$.
    c.  Does $Z_8 L$ have subnear loop rings?
    d.  Can $Z_8 L$ have subnear-rings which are associative near-rings?

2.  Replace $Z_8$ in problem (1) by the loop near-ring given by the following table for '+' and '.' respectively

    | + | 0 | a | b | c |
    |---|---|---|---|---|
    | 0 | 0 | a | b | c |
    | a | a | 0 | c | b |
    | b | b | c | 0 | a |
    | c | c | b | a | 0 |

| . | 0 | a | b | c |
|---|---|---|---|---|
| 0 | 0 | 0 | 0 | 0 |
| a | 0 | a | b | c |
| b | 0 | b | 0 | 0 |
| c | 0 | c | b | c |

Let (N, '+', '.') be the loop near-ring for the same loop L given in problem 1. Find the loop near-ring. Find the difference between the near loop ring $Z_8L$ and NL.

3. Find whether NL given in Problem 2 is a s-loop near-ring.
4. Prove $Z_2L_5$ (2) is an alternative near loop ring.
5. Find the mod p envelope of $Z_9 L_7(3)$.
6. Find zero divisors, units and idempotents of the near loop ring $Z_3L_9(5)$.
7. Find quasi regular elements in $Z_{12}L_5(3)$.
8. Find the smallest Lield?
9. Is NL given in problem 2 left bipotent? Justify your claim.

## 3.4 Groupoid Near-rings and its Properties

In this section we introduce a new type of non-associative near-rings using a groupoid and a near-ring and obtain several nice properties about these near-rings. The generalization of this is nothing but groupoids over seminear-rings these groupoid near-rings would be a larger class of non-associative near-rings which will contain the class of seminear-rings. In fact this class will be a class of non-associative seminear-rings. Throughout this section we assume N to be a right seminear-ring with 1 and G a groupoid; G may have one of may not have one.

**DEFINITION 3.4.1**: *Let N be a seminear-ring with 1 and G any groupoid. The groupoid seminear-ring NG of the groupoid G over the seminear-ring N consists of all finite formal sums of the form $\Sigma\ \alpha_i\ g_i$; $\alpha_i \in N$; $g_i \in G$ (i – running over a finite index) satisfying the following properties:*

1. *We assume NG consists of all elements of the groupoid G i.e. $G \subseteq N$ by the assumption $1.g = g = g.1$ where $1 \in N$ and $g \in N$.*
2. *We have*
   *$\Sigma\ \alpha_{li}\ g_i = \Sigma\ \beta_i\ g_i \Leftrightarrow \alpha_i = \beta_i$ for all i.*
3. *$\Sigma\ \alpha_{li}\ g_i + \Sigma\ \beta_i\ g_i = \Sigma\ (\alpha_i + \beta_i)\ g_i$.*
4. *$g_i\ \Sigma\ \alpha_i\ h_j = \Sigma\ \alpha_i\ g_i\ h_j$ (Clearly it can be checked NG is a semigroup under '+' (NG, .) is a groupoid).*
5. *$(\Sigma\alpha_i\ g_i)\ (\Sigma\ \beta_j\ h_j) = \Sigma\ \lambda_k\ s_k$ where $\lambda_k = \Sigma\alpha_i\ \beta_j$ and $g_i\ h_j = s_k$ if the elements $\alpha_i$ distributes over $\beta_j$ otherwise we see $(\alpha_i\ g_i + \ldots + \alpha_n\ g_n)\ \Sigma\ \beta_j\ h_j = \alpha_1\ g_1\ \Sigma\ \beta_j\ h_j + \ldots + \alpha_n\ g_n\ \Sigma\ \beta_j\ h_j$ (under assumption all our near-rings and seminear-rings are right distributive).*
   *$= \alpha_i\ \Sigma\ \beta_j\ g_i\ h_j + \ldots + \alpha_n\ \Sigma\ \beta_j\ g_n\ h_j$*
   *$= \alpha\ \Sigma\ \beta_j\ s_{ij} + \ldots + \alpha_n\ \Sigma\ \beta_j\ s_{nj}$, $\quad g_i\ h_j = s_{ij} \in G$.*

*6.*   $g_i \, \alpha = \alpha \, g_i$ *for all* $\alpha \in N$ *and* $g_i \in G$.

*Thus (NG, +, .) is a non-associative seminear-ring. This is the only known method to us by which we get non-associative seminear-rings.*

Since the seminear-rings happens to be non-associative we can define analogously the identities enjoyed by non-associative structures like Moufang, Bol, right (left) alternative, WIP.

All other properties which are true in case of seminear-rings will continue to be true even in case of non-associative seminear-rings as in places where we have used product of 3 elements, except in the definition of IFP and regular elements.

***Example 3.4.1***: Let $Z_4 = \{0, 1, 2, 3\}$ define '×' and '.' on $Z_4$ as '×' is the usual multiplication modulo 4, '.' is a.b = a for all a $\in Z_4$. Clearly $(Z_4, \times, .)$ is a right seminear-ring which is associative. G be a groupoid given by the following table:

| .     | $a_0$ | $a_1$ | $a_2$ | $a_3$ | $a_4$ | $a_5$ |
|-------|-------|-------|-------|-------|-------|-------|
| $a_0$ | $a_0$ | $a_3$ | $a_0$ | $a_3$ | $a_0$ | $a_3$ |
| $a_1$ | $a_2$ | $a_5$ | $a_2$ | $a_5$ | $a_2$ | $a_5$ |
| $a_2$ | $a_4$ | $a_1$ | $a_4$ | $a_1$ | $a_4$ | $a_1$ |
| $a_3$ | $a_0$ | $a_3$ | $a_0$ | $a_3$ | $a_0$ | $a_3$ |
| $a_4$ | $a_2$ | $a_5$ | $a_2$ | $a_5$ | $a_2$ | $a_5$ |
| $a_5$ | $a_4$ | $a_1$ | $a_4$ | $a_1$ | $a_4$ | $a_1$ |

$Z_4G$ is a groupoid seminear-ring which is a non-associative seminear-ring. Just we will call a non-associative seminear-ring N to be a Bol seminear-ring if ((xy)z)y = x(yz)y) for all x, y, z $\in$ N.

Let

$$\alpha \quad = \quad 2a_1 + 3a_3 + a_0$$
$$\beta \quad = \quad 3a_5 + 2a_2 + 2a_4$$

$$
\begin{aligned}
\alpha \, \beta \quad &= \quad (2a_1 + 3a_3 + a_0)\,(3a_5 + 2a_2 + 2a_4) \\
&= \quad 2a_1\,[\,3a_5 + 2a_2 + 2a_4\,] + 3a_3\,[3a_5 + 2a_2 + 2a_4\,] + a_0\,[3a_5 + 2a_2 + 2a_4\,] \\
&= \quad 2\,[3a_1a_5 + 2a_1a_2 + 2a_1a_4] + 3\,[3a_3a_5 + 2a_3a_2 + 2a_3a_4] + 1\,[a_0a_5 + 2a_0a_2 + \\
&\qquad 2a_0a_4] \\
&= \quad 2\,[3a_5] + 3\,[3a_3] + 1\,[3a_3] \\
&= \quad 2a_5 + 3a_3 + a_3 \\
&= \quad 2a_5 + 3a_3.
\end{aligned}
$$

This example clearly illustrates how the product operation is performed in $Z_4G$.

The above groupoid G is a Bol groupoid. It can be easily verified by the reader.

***Example 3.4.2***: Let $Z_6 = \{0, 1, 2, 3, 4, 5\}$ with '×' and '.' the usual multiplication modulo 6 and '.' is a.b = a for all a, b in $Z_6$ be the seminear-ring. G be a groupoid given by the following table:

|     | $x_1$ | $x_2$ | $x_3$ |
| --- | --- | --- | --- |
| $x_1$ | $x_1$ | $x_3$ | $x_2$ |
| $x_2$ | $x_2$ | $x_1$ | $x_3$ |
| $x_3$ | $x_3$ | $x_2$ | $x_1$ |

$Z_6G$ is a groupoid seminear-ring if $\alpha = 3x_1 + 2x_2$ we see $\alpha^2 = 3x_3 + 2x_2$ . $\alpha = x_1 + x_2 + x_3$ is such that $\alpha^2 = (x_1 + x_2 + x_3)^2 = x_1 + x_2 + x_3 = \alpha$. Now as we have constructed seminear-rings which are non-associative we can obtain zero divisor, idempotents, units from these groupoid seminear-rings.

Thus finding these special elements in groupoid seminear-rings happens to be a very interesting problem.

***Example 3.4.3***: Let $Z_4 = \{0, 1, 2, 3\}$ be a seminear-ring and G be a groupoid given by the following table:

| . | $a_0$ | $a_1$ | $a_2$ | $a_3$ |
| --- | --- | --- | --- | --- |
| $a_0$ | $a_0$ | $a_2$ | $a_0$ | $a_2$ |
| $a_1$ | $a_1$ | $a_3$ | $a_1$ | $a_3$ |
| $a_2$ | $a_2$ | $a_0$ | $a_2$ | $a_0$ |
| $a_3$ | $a_3$ | $a_1$ | $a_3$ | $a_1$ |

$Z_4G$ be the groupoid seminear-ring. Let $\alpha = 2 (a_0 + a_1 + a_2 + a_3)$ and $\beta = (a_0 + a_1 + a_2 + a_3)$; it is easily verified that $\alpha\beta = 0$. $I = Z_4H$ where $H = \{a_0, a_2\}$ is not a left ideal of the seminear-ring. Recall in a seminear-ring N, I is a left ideal in N if

1. $(I, +)$ is a subsemigroup of $(N, +)$.
2. $n (n^1 + i) + nn^1 \in I$ (for each $i \in I$ and $n, n_1 \in N$). But this seminear-ring $Z_4G$ has several subsemigroups.

Just like groupoids can be used to built automaton we see using group automaton we can construct syntactic near-rings. Chapter seven is completely reserved for dealing with the probable applications of Smarandache near-rings and Smarandache seminear-rings both associative and non-associative.

The study of the groupoid seminear-rings when we take the groupoid from the new class of groupoids and seminear-ring as $Z_n$, n non-prime or $Z^+ \cup \{0\}$ happens to be a very concrete study. We see there are a lot of openings to find under what conditions these groupoid seminear-rings happen to be Bol, Moufang, right (left) alternative or WIP.

Such study is at a very dormant stage as this is the first place we have made the definition non-associative seminear-rings to be satisfying several of these identities. Before we propose problems we make some more new definitions.

**DEFINITION 3.4.2**: *Let N be a non-associative seminear-ring. We say N is a weak Bol seminear-ring if N has a proper subseminear-ring M which is non-associative and all elements of M satisfies the Bol identity.*

<u>Note</u>: N need not be a Bol seminear-ring still N can be a weak Bol seminear-ring i.e. all elements of N need not satisfy Bol identity.

**DEFINITION 3.4.3**: *Let N be a non-associative seminear-ring. If every subseminear-ring M of N which is also a non-associative subseminear-ring satisfies the Bol identity or is a Bol subseminear-ring then we say N is a strong Bol seminear-ring.*

**THEOREM 3.4.1**: *Let N be a strong Bol seminear-ring then N is a weak Bol seminear-ring.*

*Proof*: Obvious by the very definition. The reader is requested to give an example of a weak Bol seminear-ring which is not a strong Bol seminear-ring.

***Remark***: In a similar way as we have defined weak and strong Bol seminear-rings we can define weak Moufang seminear-rings, strong Moufang seminear-ring, weak alternative seminear-rings and strong alternative seminear-rings and so on. The reader is requested to develop these ideas and characterize these groupoid seminear-rings. We just define P-seminear-rings where N is a non-associative seminear-rings.

**DEFINITION 3.4.4**: *Let N be a non-associative seminear-ring; we say N is a P-seminear-ring if (x.y). x = x.(y.x) for all x, y ∈ N.*

One can define P-weak seminear-ring and P-strong seminear-ring and obtain interesting results about them. Thus using groupoids and seminear-rings we can construct a class of non-associative seminear-rings.

Apart from these class of non-associative seminear-rings the author is not in a position to get any non-associative seminear-rings.

<u>**PROBLEMS:**</u>

1.    Let $(Z_8, \text{'×'}, .)$ be a seminear-ring. G be a groupoid given by the following table:

| . | $a_0$ | $a_1$ | $a_2$ | $a_3$ |
|---|---|---|---|---|
| $a_0$ | $a_0$ | $a_3$ | $a_2$ | $a_1$ |
| $a_1$ | $a_2$ | $a_1$ | $a_0$ | $a_3$ |
| $a_2$ | $a_0$ | $a_3$ | $a_2$ | $a_1$ |
| $a_3$ | $a_2$ | $a_1$ | $a_0$ | $a_3$ |

Find for the groupoid seminear-ring $Z_8 G$

    a.  Subseminear-rings.
    b.  left ideals.
    c.  ideals.
    d.  Does this seminear-ring satisfy any of the well known identies.

2. Let {$Z_3$, +, '.'} be a near-ring and G be the grouped given in problem 1.
   Is $Z_3$G a non-associative seminear-ring which is not a near-ring? Justify your answer.
3. What is the order of the smallest non-associative seminear-ring?
4. Is the order of the smallest associative seminear-ring 4? Justify your answer.
5. Give an example of a strong moufang seminear-ring.
6. Does there exist a Bol seminear-ring?
7. Give an example of a weak alternative seminear-ring.
8. Does there exists a P seminear-ring?
9. From the groupoids built using $Z_n$ and using an appropriate $Z_m$; is it possible to find an example of strong Bol seminear-ring, strong WIP seminear-ring?
10. Find all ideals in the groupoid seminear-ring $Z^o$G where the groupoid G is given in problem 1.

## 3.5 Special properties of Near-rings

This section has very unique feature, for here we just recall several definitions of different authors and have listed them. The main purpose of this book is to study and introduce as many as Smarandache concepts in near-rings; so we avoid dealing with near-rings; whereever it is a dire necessity we give the definitions. This section gives the definitions from several authors.

**DEFINITION [22]**: *Let N be a zero symmetric left near-ring. A mapping $F : N \rightarrow N$ is said to be commuting on the subset S of N if [x, F(x)] = 0 for all x ∈ S, where [x, y] = xy − yx.*

*A mapping $D : N \rightarrow N$ is called a multiplicative derivation if D(xy) = x D(y) + D(x)y for all x, y ∈ N; if D is also additive, D is called a derivation.*

Several interesting properties can be derived in case of S-near-ring and all the more on non-associative near-rings.

**DEFINITION [5]**: *Near-rings here are left near-rings. A near-ring N is weakly divisible if for all x, y ∈ N there exists z ∈ N, such that xz = y or yz = x.*

For more about these near-rings please refer [5].

**DEFINITION [14]**: *The right near-ring N is called strongly prime if for each a ∈ N \ {0} there exists a finite subset F such that aFx ≠ {0} for all x ∈ N \ {0} [respectively such that afx = afy for all f ∈ F implies x = y].*

**DEFINITION [51]**: *Let S be a non-empty subset of a right near-ring N and define the left and right polar subsets of N, respectively by L(S) = {x / xNs = (0) ∀ s ∈ S} and R(S) = {y / sNy = (0) ∀ s ∈ S} If $P_L$(N) is the set of left polar subsets of N and $P_R$(N) is the set of nonempty right polar subsets of N, then $P_L$(N) and $P_R$(N) are complete bounded lattices.*

**DEFINITION [4]**: *A generalized p-near-field (gp near-field) is a near-ring with 1 in which for each x there is an $n = n(x)$ with $x^{n(x)} = 1$. N is called an gp-near-ring if it has nonzero nilpotent elements and its one sided proper non-nil ideals are gp-near-fields.*

**DEFINITION [43]**: *If R is a ring and M is an R-module then there is exactly one way to extend the module multiplication 'o' in $(N, +) = (R, +) \oplus (M, +)$, to get an abstract affine near-ring namely $(r, m) o (s, n) = (rs, rn + m)$.*

**DEFINITION [69]**: *Let N be a near-ring, S a non-empty subset of an N-group V. An N-subgroup H of V centralizes S if $(h + u) a = ha + ua$ for all $h \in H$, $u \in S$ and $a \in N$.*

**DEFINITION [69]**: *An involution in a near-ring N (with identity) is a proper unit $\lambda$ that satisfies $\lambda^2 = 1$. If N is generated by some single involution then N is an involution near-ring.*

**DEFINITION [18]**: *A group semiautomaton is an ordered triple $S = (G, I, \delta)$ where G is a group, I is a set of inputs and $\delta : G \times I \to G$ is a transition function. The syntactic near-ring of S denoted by N (S) is a subnear-ring of M(G) generated by the maps $f_x$ where $af_x = \delta(a, x)$ together with the identity map of G.*

**DEFINITION [102]**: *A $\theta$-symmetric near-ring N is called reliable if whenever $\theta : I \to N$ is an epimorphism and I is an ideal of a near-ring A then for $(x, y \in I, a \in A)$; $x - y \in$ Ker $\theta$ implies $ax - ay \in$ Ker $\theta$.*

**DEFINITION [30]**: *Let N be a near-ring one can define an equivalence relation on N by $a \equiv b$ if and only if $na = nb$ for all $n \in N$.*

*N is called planar if there are more than three equivalence classes and every equation of the form $xa = xb + c$ has a unique solution x when $a \not\equiv b$.*

A planar near-ring can be used to construct balanced incomplete block designs (BIBD) of high efficiency. By taking either rows or columns of the incidence matrix of such a BIBD one can obtain error correcting codes with several nice features.

**DEFINITION [49]**: *An I-E group is a group which satisfies $I (G) = E(G)$ where $I(G)$ is the near-ring generated by all the inner automorphisms of G and $E(G)$ is the near-ring generated by all the endomorphisms of G. The authors show that a semidirect product $G = H.K$ of I-E groups where $(|H|, |K|) = 1$ and H is normal in G in an I-E-group if and only if projection of G onto K is in $I(G)$.*

**DEFINITION [50]**: *A group G written additively is an I-E group if $I(G) = E(G)$. Here $I(G)$ and $E(G)$ are near-rings additively generated by the inner automorphisms and the endomorphisms of G respectively.*

**DEFINITION [73]**: *Let $G = S \oplus Z$ be the direct sum of a simple non-abelian group S of order n, and a cyclic group Z of order q, where S contains a cyclic group of order q. Let $E(G)$ be the near-ring generated by all the endomorphism of G. Let $N = \{a \in E(G) / Ga \subseteq S, Sa = \{0\}\}$ Then N is an ideal of order $n^{n(q-1)}$.*

**DEFINITION [20]**: *S is a unitary non-associative zero symmetric right near-ring. L denotes the set of all elements in S without left inverses. U = S \ L and (A, B, C) = {(ab) c − a (bc) / a ∈ A, b ∈ B, c ∈ C} where (A, B, C) ∈ {L, S, U}.*

**DEFINITION [27]**: *Let (G, +) be a group written additively not necessarily abelian, with identity element denoted by 0. Let T ⊂ G \ {0} and define . $\overline{\Gamma}$ by a. $\overline{\Gamma}$ b = a if b ∈ T; a $\dot{\Gamma}$ b = 0 if b ∉ T. Then {G, +, $\overline{\Gamma}$ } is a right zero symmetric near-ring.*

Several interesting results have been carried out in this direction.

**DEFINITION [15]**: *Let N be a near-ring and let $N_0$ be a subset of N such that $N_0$ generates N and $NN_0 \subseteq N_0$; suppose that $A_0$ is a subset of $N_0$ then, $NA_0 \subset A_0$ and $A_0 N_0 \subset A_0$. Then $A_0$ is called the σ-subset of N and the subnear-ring generated by $A_0$ is called a σ -subnear-ring N. A class B of near-rings is called σ - hereditary of N ∈ B. A σ - subnear-ring E of N is a called essential in N if 0 ≠ A is an ideal of N implies E ∩ A ≠ 0.*

**DEFINITION [34]**: *A subset S of a near-ring N is integral if S has no nonzero divisors of zero. If H is an integral subset of N with the further property that $N^2 \subset H$ then N is said to be an H- integral near-ring.*

**DEFINITION [13]**: *Let N be a right near-ring and G a N-group. G is called equiprime if NG ≠ 0, N0 = 0 and if for each a ∈ N and g, $g^l$ ∈ G with Ga ≠ 0; ang = ang$^l$ for all n ∈ N implies g = $g^l$. The equiprime ideals of N are precisely the annihilators of equiprime N-groups.*

**DEFINITION [27]**: *A semiautomaton S = (Q, X, δ) with state set Q, input monoid X and state transition function δ : Q × X → Q is called a group semiautomaton if Q is an additive group. S is called additive if there is some $x_0$ ∈ X with δ(q, x) = δ(q, $x_0$) + δ (0, x) and δ(q − $q^l$, $x_0$) = δ(q, $x_0$) - δ($q^l$, $x_0$) for all q, $q^l$ ∈ Q and x ∈ X. Then there is some homomorphism*

$$\psi : Q \rightarrow Q \text{ and some map}$$
$$\alpha : X \rightarrow Q \text{ with } \phi(x_0) = 0$$

*and δ(q, x) = ψ(q) + α(x). Let $δ_x$ for a fixed x ∈ X be the map Q → Q, q → δ(q, x) then {$δ_x$ / x ∈ X} generates a subnear-ring N(P) of the near-ring (m(Q), +, .) of all mappings on Q. The near-ring N(P) is called the syntactic near-ring of P.*

**DEFINITION [13]**: *Let N be a right near-ring and A an ideal or a left ideal of N. We define three properties as follows:*

a.   *A is equiprime if for any a, x, y ∈ N such that anx − any ∈ A for all n ∈ N, we have a ∈ A or x − y ∈ A.*

b.   *A is strongly semiprime if for each a ∈ N \ A there exists a finite subset F of N such that if x, y ∈ N and afx − afy ∈ A for all f ∈ F then x − y ∈ A.*

c.   *A is completely equiprime if a ∈ N \ A and ax − ay ∈ A imply x − y ∈ A.*

**DEFINITION [55]**: *A right near-ring N is a left infra near-ring if (N, +, .) is a triple where N is a non-empty set + and . are binary operations satisfying*

    *1.  (N, +) is a group.*
    *2.  (N, .) is a semigroup.*
    *3.  For all x, y, z in N, x.(y + z) = x.y − x.0 + x.z.*

**DEFINITION [41]**: *A near-ring M which satisfies the property, A an ideal of B, B an ideal of N and B / A ≅ M implies A is an ideal of N for all near-rings N and subnear-ring A and B of N are called F near-ring.*

**DEFINITION [20]**: *The authors consider a left zero symmetrical near-ring R and subsemigroup S of R under multiplication. Let D be the normal subgroup of R under addition generated by the set of differences { −(xs + ys) + (x + y)s / x, y ∈ R; s ∈ S}. Call D the defect of distributivity. D is a minimal defect if D is either {0} or minimal as a normal subgroup of R under addition.*

**DEFINITION [23]**: *The author defines weak-weakly regular near-ring. A near-ring N is left weak-weakly regular if every x ∈ N is of the form x = ux for some u in the principal ideal generated by x and is left weakly regular if every a ∈ N is in the product of {Na} * {Na} = {finite sums $\Sigma x_k\, y_k\, / x_k, y_k \in Na$}.*

**DEFINITION [60]**: *A composition near-ring is a quadruple {C, +, 0, .) where (C, +,.) and (C, +, o) are near-rings such that (a.b) o c = (a o b) .c for all a, b, c ∈ C.*

**DEFINITION 3.5.1**: *A nonzero ideal H of G is said to be uniform if for each pair of ideals $K_1$ and $K_2$ of G such that $K_1 \cap K_2 = (0)$; $K_1 \subset H$, $K_2 \subset H$ implies $K_1 = (0)$ or $K_2 = (0)$.*

**DEFINITION 3.5.2**: *An ideal H of G is said to have finite Goldie dimension (written as FGD) if H does not contain an infinite number of non-zero ideals of G whose sum is direct.*

**THEOREM [67]**: *Let G be a N-group with FGD then every nonzero ideal of G contains a uniform ideal.*

The proof is left for the reader. Now we proceed on to introduce a new notion called normal near-rings.

**DEFINITION 3.5.3**: *Let N be a near-ring; an element a ∈ N is called normal element of N if aN = Na. If aN = Na for every a ∈ N then N is called a normal near-ring. Let n(N) denote the set of all normal elements of N. N is a normal near-ring if and only if n(N) = N.*

**THEOREM 3.5.1**: *Let N be a near-ring we can have n(N) = ϕ.*

*Proof*: Take $Z_m$ = {0, 1, 2, …, m-1} the set of integers modulo m. {$Z_m$, +, .} is a near-ring under addition '+'modulo m and '.' defined by a.b = a for all a, b ∈ $Z_m$, n(N) = ϕ.

**DEFINITION 3.5.4**: *Let N be a near-ring. A subnear-ring S of N is called a normal subnear-ring with respect to a non-empty subset T of N if ts = st for every t ∈ T. In particular we have nS = Sn for every n ∈ N we say S is a normal subnear-ring of N.*

**DEFINITION 3.5.5**: *Let N be a near-ring, a subsemigroup S of N is called a normal subsemigroup of S if nS = Sn for all n ∈ N.*

**PROPOSITION 3.5.1:** *Let N be a near-ring and S be a normal subsemigroup of (N, .) then N satisfies the left ore condition with respect to S of (N, .)*

*Proof*: N is said to fulfill the left ore condition with respect to a given subsemigroup S of (N, .) if for all (s, n) ∈ S × N there exists $(s_1, n_1)$ ∈ S ×N; $ns_1 = sn_1$. Now if S is a normal subsemigroup then we have nS = Sn that is $ns_1 = s_2n$ for n ∈ N; hence N satisfies the left ore condition with respect to normal subsemigroup S of (N, .).

**DEFINITION 3.5.6**: *Let N be a near-ring a subnear-ring M of N is invariant with respect to an element a ∈ N \ M if a M ⊂ M and Ma ⊂ M and strictly invariant with respect to a if aM = M and Ma = M.*

**THEOREM 3.5.2**: *A subnear-ring M of N is invariant with respect to every element a of N; then M is invariant subnear-ring of N.*

*Proof*: Straightforward; left for the reader to prove.

**THEOREM 3.5.3**: *Let N be a near-ring. A subnear-ring M is strictly invariant with respect to every element T in a subset in N. Then M is a normal subnear-ring with respect to T.*

*Proof*: Left for the reader to prove.

**DEFINITION 3.5.7**: *A commutative near-ring N with identity is called a Marot near-ring if each regular ideal of N is generated by regular elements, where by regular ideal we mean an ideal with regular elements and a non-zero divisor of a near-ring N.*

***Example 3.5.1***: All near-fields are Marot near-rings.

**THEOREM 3.5.4**: *Let G be a torsion free abelian group and N a near-field. The group near-ring is a Marot near-ring.*

*Proof*: Straightforward, hence left for the reader to prove. A near-ring N is Boolean if and only if for all x ∈ N, $x^2 = x$. A near-ring N is a p-near-ring if and only if for all x ∈ N; $x^p = x$ and px = 0.

**THEOREM 3.5.5**: *Let N be a near-ring which is Boolean and is such that x.y = 0 for every distinct pair of elements x, y ∈ N and S a semigroup. The semigroup near-ring NS is a Boolean near-ring if and only S is a semigroup such that in S for any $a_i$, $a_j$ ∈ S, $a_i ≠ a_j$ we have $a_i a_j = 0$ and $a_i a_j = a_i$ if $a_i = a_j$.*

*Proof*: Left as an exercise for the reader to prove. It is important to note that if we replace the semigroup by a group and even if N is a Boolean near-ring with $x^2 = x$ and $xy = 0$ if $x \neq y$ for all $x, y \in$ N; NG the group near-ring is never a Boolean near-ring as in G we can never have $g^2 = g$.

**THEOREM 3.5.6**: *Let G be a group of finite order say n, n-odd $Z_2 = \{0, 1\}$ be the near-ring. Then the group near-ring $Z_2G$ has nontrivial subBoolean near-rings.*

*Proof*: Clearly S = $\{0, (1 + g_2 + \ldots + g_n)\}$ is a nontrivial Boolean near-ring.

**DEFINITION [66]**: *Let K and I be ideals of a N-group G, K is said to be the complement of I if the following two conditions are true.*

    1. *$K \cap I = (0)$ and*
    2. *$K_1$ is an ideal of G such that $K \subset K_1$, $K_1 \neq K$ imply $K_1 \cap I \neq 0$.*

**THEOREM [66]**: *Let G be a N-group. Fix an ideal I of G, by Zorn's lemma, the set of all ideals H of N satisfying $H \cap I = (0)$ contains a maximal element say K. Again by Zorn's lemma the set of all ideals X of G satisfying $I \subset X$ and $X \cap K = (0)$ contains a maximal element say $K^0$ then we have*

    a. *K is complement of I,*
    b. *$K^0$ is a complement of K;*
    c. *$K + I$ and $K + K^0$ are essential ideals and*
    d. *I is essential in $K^0$.*

*Proof:* The proof can be had from [66].

**DEFINITIONS [66]**:

    a. *A subset S of the N-group G is said to be small in G if $S + K = G$ and K is an ideal of G imply K = G.*
    b. *G is said to be hollow if every proper ideal of G is small in G and*
    c. *G is said to have finite spanning dimension (FSD in short) if for any decreasing sequence of N-subgroups.*

*$X_0 \supset X_1 \supset X_2 \supset \ldots$ of G such that $X_i$ an ideal of $X_{i-1}$ ; there exists an integer k such that $X_j$ is small in G for all $j > k$.*

**DEFINITION [65]**: *Suppose H and K are two N-subgroups of G. Then K is said to be supplement for H if $H + K = G$ and $H + K_1 \neq G$ for any proper ideal $K_1$ of K.*

**THEOREM [65]**: *Let H be a non-small N-subgroup of G. If every proper ideal of H is small in G then H is hollow.*

*Proof*: Left for the reader to refer [65]

**THEOREM [65]**: *If G has FSD, then every non-small N-subgroup H of G contains a hollow N-subgroup which is non-small in G.*

Several more results can be had from [65]. Now we proceed on to define a new class of near-rings called near semilattices [85].

**DEFINITION [85]**: *Let L be a lattice with 0 and 1. The new product is defined on L $\times_p$ L as*

1. *$a_1, a, b_1, b \in L$ we have $(a, b) \cup (a_1, b_1) = (a_1 \cup a_1, b \cup b_1)$.*
2. *$(a, b) \cap (a_1, b_1) = (a \cap b_1, b \cap a_1)$ ('$\cup$' and '$\cap$' can also be inter changed).*
   *Clearly (0, 0) and (1, 1) are the least and the greatest elements of L $\times_p$ L. Obviously L $\times_p$ L is not a lattice; we call this structures (L $\times_p$ L, $\cup$, $\cap$) a near-semilattice.*

**DEFINITION [85]**: *Let L be a lattice, if we define a new product on L $\times_q$ L as*

1. *$L \times_n L = \{(a, b) / a, b \in L\}$.*
2. *For $(a, b), (c, d) \in L \times_n L$ $(a, b) \cup (c, d) = (a \cup d, b \cup c)$.*
3. *$(a, b) \cap (c, d) = (a \cap d, b \cap c)$.*

*then we call L $\times_q$ L a quasi near semilattice. For this also (0, 0) acts as the least element and (1, 1) acts as the greatest element. If lattice L is distributive we get a near semilattice to be a seminear-ring.*

**Example 3.5.2**: Let L = {0, 1} be the two element lattice L $\times_p$ L = {(0,0) (1, 0) (0, 1), (1, 1)}. Clearly L $\times_p$ L is a semilattice with respect to $\cap$ (or '$\cup$') but not even a semilattice under '$\cup$' (or $\cap$) L $\times_p$ L = {(0, 0), (1, 0), (0, 1), (1, 1)} is the quasi near-semilattice.

**THEOREM [85]**: *Let L $\times_p$ L, L $\times_q$ L be the near semilattice and quasi near semilattice where L has 0 and 1. Then L $\times_p$ L and L $\times_q$ L have nontrivial idempotents.*

*Proof*: Obvious from the fact (x, x) under the operation $\cup$, $\cap$ are idempotents. It is very important to note no other elements other than elements of the form (x, x) are the only idempotents.

**Example [85]:** Let L = {0, 1}; L $\times_p$ L, be the near semilattice {L $\times_p$ L, $\cup$, $\cap$} has the following diagram:

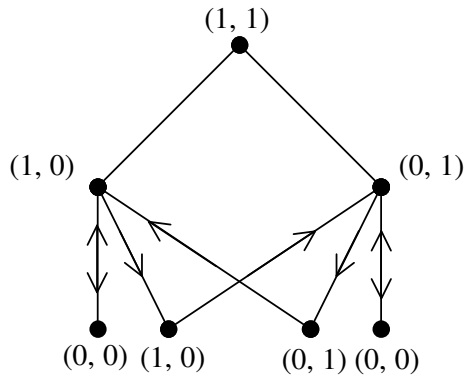

**Figure 3.5.1**

**DEFINITION [85]**: *Let L and $L^1$ be any two lattices L $\times_p$ L and $L^1$ $\times_p$ $L^1$ be near semilattices. A map $\phi$ from L $\times_p$ L to $L^1$ $\times_p$ $L^1$ is called a near semilattice homomorphism if*

$$\phi((a, b) \cup (c, d)) = \phi(a^1\ b^1) \cup \phi(c^1\ d^1)\ and$$
$$\phi((a, b) \cap (c, d)) = \phi(a^1\ b^1) \cap \phi(c^1\ d^1)$$

*for all (a, b), (c, d) $\in$ L $\times_p$ L and $(a^1, b^1)$, $(c^1, d^1) \in L^1 \times_p L^1$. Clearly $\phi(0, 0) = (0, 0)$ and $\phi(1, 1) = (1, 1)$.*

**DEFINITION [85]**: *Let L be a lattice L $\times_p$ L be a near semilattice of L. S a proper subset of L $\times_p$ L is called a subnear-semilattice if $(S, \cup, \cap)$ itself a near semilattice.*

*Also if $S^1$ is a sublattice of L then $S = S^1 \times_p S^1$ is also a subnear semilattice. This study becomes interesting when the lattices under study are distributive so that L $\times_p$ L becomes a seminear-ring.*

Finally we introduce yet another structure about seminear-rings, which we call as seminear pseudo rings (SNP-ring).

**DEFINITION 3.5.8**: *Let P be a non-empty set on which is defined two binary operations $\oplus$ and $\odot$ satisfying the following conditions:*

a. *For all p, q $\in$ P, p $\oplus$ q $\in$ P.*
b. *p $\oplus$ p = e for all p $\in$ P, e is called the one sided identity in P with respect to the operation $\oplus$.*
c. *p $\oplus$ e = p for all p $\in$ P.*
d. *p $\oplus$ q $\neq$ q $\oplus$ p for p, q $\in$ P.*
e. *(p $\oplus$ q) $\oplus$ r $\neq$ p $\oplus$ (q $\oplus$ r), in general for p, q, r $\in$ P.*
f. *p $\odot$ q $\in$ P for all p, q $\in$ P.*
g. *p $\odot$ p = $e^1$, $e^1 \in$ P is the one sided identity with respect to $\odot$ for all p $\in$ P.*
h. *p $\odot$ q $\neq$ q $\odot$ p for p, q $\in$ P.*
i. *p $\odot$ (q $\odot$ r) $\neq$ (p $\odot$ q) $\odot$ r for p, q, r $\in$ P.*
j. *(p $\oplus$ q) $\odot$ r = p $\odot$ r $\oplus$ q $\odot$ r for all p, q, r $\in$ P.*

*Then (P, $\oplus$, $\odot$) is called a seminear pseudo ring.*

**Example 3.5.3**: Q be the field of rationals; define on Q '$\oplus$' as the usual '$-$' and $\odot$ as p $\odot$ q = p $\odot$ $q^{-1}$. It is easily verified (Q, $\oplus$, $\odot$) is a SNP-ring.

**Example 3.5.4**: Let D be a division ring of characteristic 0. Define p $\oplus$ q = p $-$ q , p $\odot$ q = $\dfrac{p}{q}$ , (q $\neq$ 0), p $\odot$ 0 = 0 . Clearly (D, $\oplus$, $\odot$) is a SNP-ring.

**DEFINITION 3.5.9**: *Let (R, $\oplus$, $\odot$) be a SNP-ring. A proper subset S $\subset$ R is said to be a SNP-subring of R if (S, $\oplus$, $\odot$) is a SNP-ring.*

**DEFINITION 3.5.10**: *Let (R, ⊕, ⊙) be a SNP-ring. A nonempty subset I of R is called the SNP-ideal if*

1. *For all p, q ∈ I, p ⊕ q ∈ I.*
2. *0 ∈ I.*
3. *For all p ∈ I and r ∈ R we have p ⊙ r or r ⊙ p ∈ I.*

**DEFINITION 3.5.11**: *Let (R, ⊕, ⊙) be a SNP-ring if (a ⊕ b) ⊙ c = a ⊙ c ⊕ b ⊙ c for all a, b, c ∈ R; then R is a quasi SNP-ring.*

**THEOREM 3.5.7**: *Every SNP-ring is a quasi SNP-ring and not conversely.*

*Proof*: Straightforward by the very definition.

**DEFINITION 3.5.12**: *Let (R, ⊕, ⊙) and (Q, ⊕, ⊙) be two SNP-rings. A SNP-ring homomorphism $\phi$ form (R, ⊕, ⊙) to (Q, $⊕^1$, $⊙^1$) satisfies the following operational rules*

1. *$\phi (r ⊕ s) = \phi (r) ⊕^1 \phi (s)$.*
2. *$\phi (0) = 0$.*
3. *$\phi (r ⊙ s) = \phi (r) ⊙^1 \phi (s)$.*
4. *$\phi (e) = e^1$ where e and $e^1$ are the identities of R and Q respectively.*

*The kernel of the SNP-ring homomorphism and the concept of isomorphism are defined as for in case of rings.*

<u>**PROBLEMS:**</u>

1. Give an example of a weakly divisible near-ring.
2. Does the near-ring $Z_{14}$ have a normal element?
3. Does $Z_{24}$ have a normal subnear-ring?
4. Find an example of a near-ring which is normal.
5. Is $Z_{32}$ a Marot near-ring?
6. Give an example of a near-ring which is Boolean.
7. Does there exists a finite SNP-ring? Substantiate your claim.

**Chapter Four**

# SMARANDACHE NEAR-RINGS

In this chapter we introduce the concept of Smarandache near-rings (S-near-ring) and introduce several of the properties in near-rings to the S-near-rings. The study of S-near-rings is very recent [99]. The near-rings introduced in this chapter will be called as S-near-rings of level I. S-near-rings of level II will also be introduced in the next chapter. Throughout this chapter by a S-near-ring we mean only S-near-rings of level I. All near rings and S-near rings in this book will only be right near-rings. This chapter has five sections. In section one we define S-near-rings and give examples; Smarandache N-groups are introduced in section two. In section three Smarandache direct product and Smarandache free near-rings defined and described. The concept of Smarandache ideals in near-rings is introduced in section four and in the final section we give the notion of S-idempotents and Smarandache modularity.

## 4.1 Definition of S-near-ring with examples

In this section we define S-near-rings, Smarandache subnear-rings and illustrate them with examples. We do not claim by any way that we have exhausted all Smarandache analogous for near-rings. The only aim of this book is to introduce S-near-rings and define some properties and leave the rest of the work for any innovative researcher/ student to pursue the research on Smarandache notions about near-rings.

**DEFINITION 4.1.1**: *N is said to be a Smarandache near-ring (S-near-ring) if (N, +, .) is a near-ring and has a proper subset A such that (A, +, .) is a near-field.*

**Example 4.1.1**: Let $Z_2 = \{0, 1\}$ be the near-field. Take any group G such that $Z_2G$ is the group near-ring of the group G over the near-field $Z_2$. $Z_2G$ is a S-near-ring as $Z_2 \subset Z_2G$ and $Z_2$ is a near-field.

It is important to note that in case of S-near-rings we can get several or a new class of S-near-ring by using group near-rings and semigroup near-rings using basically the near-field $Z_2$ or any $Z_p$.

**Example 4.1.2**: Let $Z_2 = \{0, 1\}$ be a near-field and S(n) any symmetric semigroup. The semigroup near-ring $Z_2S(n)$ is a S-near-ring. This S-near-ring can be of any order finite or infinite depending on the order of the group or the semigroup which is used. Similarly they can be commutative or non-commutative which depends basically on the groups and semigroups.

**DEFINITION 4.1.2**: *Let (N, +, .) be a S-near-ring. A non-empty proper subset T of N is said to be a Smarandache subnear-ring (S-subnear ring) if (T, +, .) is a S-near-ring; i.e. T has a proper subset which is a nearfield.*

***Example 4.1.3***: Let $Z_2$ be the near-field. $S_n$ be the symmetric group of degree n. $Z_2 S_n$ is the group near-ring. Clearly $Z_2 H$ is a S-subnear ring where

$$H = \left\{ 1 = \begin{pmatrix} 1 & 2 & 3 \\ 1 & 2 & 3 \end{pmatrix} , \ p_4 = \begin{pmatrix} 1 & 2 & 3 \\ 2 & 3 & 1 \end{pmatrix} , \ p_5 = \begin{pmatrix} 1 & 2 & 3 \\ 3 & 1 & 2 \end{pmatrix} \right\}$$

is a subgroup of $S_n$ and $Z_2 \subset Z_2 H$. Hence the claim.

***Example 4.1.4***: Let $M_{3 \times 3} = \{(a_{ij}) \ / \ a_{ij} \in Z_2 = \{0, 1\}$ the near-field$\}$ be the collection of all $3 \times 3$ matrices. $M_{3 \times 3}$ under the usual matrix addition, and multiplication defined as from $Z_2$ makes $M_{3 \times 3}$ a near-ring which is also a S-near-ring. This has also S-subnear-rings which are nontrivial.

***Example 4.1.5***: $Z_2$ be the near-field. $Z_2[x]$ is a S-near-field. By taking $p(x).q(x) = p(x)$ for all $p(x), q(x) \in Z_2[x]$ we see $(Z_2[x], +, '.')$ is a S-near-ring. $Z_2[x]$ has nontrivial S-subnear-rings.

<u>**PROBLEMS:**</u>

1.  Find a S-near-ring of order 32.
2.  Does their exist a S-near-ring of order 17? Justify your answer.
3.  Is the group near-ring $Z_2 G$ where $Z_2$ is a near-field and $G = \langle g \ / \ g^8 = 1 \rangle$, a S-near-ring?
4.  Does there exist at least two nontrivial S-subnear-rings of the near-ring given in problem 3?
5.  Let $M_{4 \times 4} = \{(a_{ij}) \ / \ a_{ij} \in Z_2 - Z_2$ a near-field$\}$ prove $M_{4 \times 4}$ is a S-near-ring. Find a S-subnear-ring in $M_{4 \times 4}$. What is the order of $M_{4 \times 4}$?
6.  Prove $S = Z_2 \times Z_2$ is a S-near-ring where $Z_2$ is the near-field. Can $Z_2 \times Z_2$ have S-subnear-rings?
7.  Let $Z_n = \{0, 1, 2, \dots, n\text{-}1\}$, $(Z_n, +)$ is a group and define product '.' as a.b = a for all a, b $\in Z_n$. Is $(Z_n, +, .)$ a S-near-ring?

## 4.2 Smarandache N-groups

In this section we introduce a concept analogous to N-groups in a near-ring the Smarandache N-groups (S-N-groups). We define and illustrate them with examples and proceed onto define S-ideals of a near-ring and Smarandache homomorphism (S-homomorphism) of near-rings.

**DEFINITION 4.2.1**: *Let (P, +) be a S-semigroup with zero 0 and let N be a S-near-ring. Let $\mu : N \times Y \rightarrow Y$ where Y is a proper subset of P which is a group under the operations of P. (P, $\mu$) is called the Smarandache N-group (S-N-group) if for all $y \in Y$ and for all n, $n_1 \in N$ we have $(n + n_1) \ y = ny + n_1 y$ and $(nn_1) \ y = n(n, \ y)$. $S \ (N^p)$ stands for S-N-groups.*

**DEFINITION 4.2.2**: *A S-semigroup M of a near-ring N is called a Smarandache-quasi subnear-ring (S-quasi subnear-ring) of N if X ⊂ M where X is a subgroup of M which is such that X.X ⊆ X.*

**THEOREM 4.2.1**: *If N has a S-quasi subnear-ring then N has a subnear-ring.*

*Proof*: Obvious by the very definitions of these concepts.

It is noteworthy to state that the existence of a subnear-ring in a near-ring need not in general be a S-quasi subnear-ring.

**DEFINITION 4.2.3**: *A S-subsemigroup Y of S($N^P$) with NY ⊂ Y is said to be a Smarandache N-subgroup (S-N-subgroup) of P.*

**DEFINITION 4.2.4**: *Let N and $N_1$ be two S-near-rings. P and $P_1$ be S-N-subgroups.*

a. *h : N → $N_1$ is called a Smarandache near-ring homomorphism (S-near-ring homomorphism) if for all m, n ∈ M (M is a proper subset of N which is a near-field) we have*

$$h\,(m + n) = h(m) + h(n),$$
$$h\,(mn) = h(m)\,h(n)\ where\ h\,(m)\ and$$

*h (n) ∈ $M_1$ ($M_1$ is a proper subset of $N_1$ which is a near-field).It is to be noted that h need not even be defined on whole of N.*

b. *h : P → $P_1$ is called the Smarandache N-subgroup-homomorphism (S-N-subgroup-homomorphism) if for all p, q in S (S the proper subset of P which is a S-N-subgroup of the S-semigroup P) and for all m ∈ M ⊂ N (M a subfield of N) h (p + q) = h (p) + h (q)  and h (mp) = mh (p); h(p), h(q) and mh(p) ∈ $S_1$ ($S_1$ a proper subset of $P_1$ which is a S-$N_1$ subgroup of S-semigroup $P_1$).*

*Here also we do not demand h to be defined on whole of P.*

**DEFINITION 4.2.5**: *Let N be a S-near-ring. A normal subgroup I of (N, +) is called a Smarandache ideal (S-ideal) of N related to X if*

    a.   *IX ⊆ I.*
    b.   *∀ x, y ∈ X and for all i ∈ I; x (y + i) − xy ∈ I,*

*where X is the near-field contained in N.*

A subgroup may or may not be a S-ideal related to all near-fields. Thus while defining S-ideal it is important to mention the related near-field.

Normal subgroup T of (N, +) with (a) is called S-right ideals of N related to X while normal subgroups L of (N, +) with (b) are called S-left ideals of N related to X.

**DEFINITION 4.2.6**: *A proper subset S of P is called a Smarandache ideal of $S(N^p)$ (S-ideal of $S(N)^P$) related to Y if*

1. *S is a S-normal subgroup of the S-semigroup P.*
2. *For all $s_1 \in S$ and $s \in Y$ (Y is the subgroup of P) and for all $m \in M$ (M the near-field of N)*
   *$n ( s + s_1 ) - ns \in S$.*

*A S-near ring is Smarandache simple (S-simple) if it has no S-ideals. $S(N^p)$ is called Smarandache N-simple (S-N-simple) if it has no S-normal subgroups except 0 and P.*

**DEFINITION 4.2.7**: *A S-subnear-ring M of a near-ring N is called Smarandache invariant (S-invariant) related to the near-field X in N if $M X \subset M$ and $X M \subset M$ where X is a S-near-field of N. Thus in case of Smarandache invariance it is only a relative concept as a S-subnear ring M may not be invariant related to every near-fields in the near-ring N.*

**DEFINITION 4.2.8**: *Let X and Y be S-semigroups of $S(N^P)$. $(X{:}Y) = \{n \in M \mid nY \subset X\}$ where M is a near-field contained in N. $(0, x)$ is called the Smarandache annihilator (S-annihilator) of X.*

<u>**PROBLEMS**</u>:

1. Is the near-ring $Z_n$ under usual addition and multiplication defined by $m.p = m$ for all $m, p \in Z_n$.
   - a. A-S-near-ring?
   - b. Does $Z_n$ have S-quasi subnear-ring?
   - c. Find the S-N-subgroup for $Z_n$.
2. Find the order of the smallest S-near-ring (It is well known $Z_2$ is the smallest near-ring).
3. Prove $N = Z_2 \times Z_9 \times Z_5$ is a S-near-ring.
4. Construct a S-homomorphism between N and $N_1$ where $N = Z_2 \times Z_7$ and $N_2 = Z \times Z_2$.
5. Find the S-ideal for N given in problem 3.

## 4.3 Smarandache direct product and Smarandache free near-rings

In this section we define Smarandache direct product (S-direct product) of near-rings and introduce the concept of Smarandache free near rings (S-free near rings) and several related concepts are defined and it is left for the reader to develop more results and theorems based on these definitions.

**DEFINITION 4.3.1**: *Let $\{N_i\}$ be a family of near-rings which has at least one S-near-ring. Then this direct product $N_1 \times \ldots \times N_r = \underset{i \in I}{X} N_i$ with component wise defined operations '+' and '.' is called Smarandache direct product (S-direct product) of near-rings.*

**THEOREM 4.3.1**: *The S-direct product of near ring is a S-near-ring.*

*Proof*: Follows from the fact in the product at least one of the $N_i$ is a S-near-ring.

***Example 4.3.1***: Let $N = Z_2 \times Z \times Z_3 \times Z_7$ be the S-direct product of near-rings. N is a S-near-ring.

**DEFINITION 4.3.2**: *The S-subnear-ring of the S-direct product $\underset{i \in I}{X} N_i$ consisting of those elements where all components except a finite number; equal to zero, is called the Smarandache (external) direct sum (S-direct sum) $\underset{i \in I}{\oplus} N_i$ of the $N_i$'s.*

More generally every S-subnear-ring N of $\underset{i \in I}{X} N_i$ where all projection maps $\pi_I$ ($i \in I$) are surjective (in other words for all $i \in I$ and for all $n_i \in N_i$; $n_i$ is the $i^{th}$ component of some element of N) is called the Smarandache subdirect product (S-subdirect product) of the $N_i$'s.

Note we would be defining Smarandache mixed direct product of near-rings which will be different from these; as the sole motivation of defining Smarandache mixed direct product is to define Smarandache near-rings of several levels.

**THEOREM 4.3.2**: *If $N = \oplus N_i$ is the S-external direct sum of near-rings then N is a S-near-ring.*

*Proof*: Follows from the very definition of S-external direct sum.

**DEFINITION 4.3.3**: *Let N be a near-ring and S a S-subsemigroup of (N, +). The near-ring $N_s$ is called the Smarandache-near-ring of left (right) quotients (S-near-ring of left (right) quotients) of N with respect to S if*

    *a.  $N_s$ has identity.*
    *b.  N is embeddable in $N_s$, by a homomorphism h.*
    *c.  For all $s \in S$; h(s) is invertible in ($N_s$, .).*
    *d.  For all $q \in N_s$, there exists $s \in S$ and there exist $n \in N$ such that $q = h(n)$ $h(s)^{-1}$ ($q = h(s)^{-1} h(n)$).*

The problem whether $N_s$ is a S-near-ring is left as an open problem.

**DEFINITION 4.3.4**: *The near-ring N is said to fulfil the Smarandache left (right) ore conditions (S-left(right) ore condition) (ore (1)) with respect to a given S-subsemigroup P of (N, .) if for (s, n) $\in$ S $\times$ N there exists $ns_1 = sn_1$ ($s_1 n = n_1 s$).*

**THEOREM 4.3.3**: *If a near-ring N satisfies S-left (right) ore condition with respect to S-subsemigroup P then the near-ring N fulfils left (right) ore condition (ore (1)) with respect to a given subsemigroup P of (N, .).*

*Proof*: Obvious by the very definitions. The reader is advised to find an example of near-ring N which fulfills left (right) ore condition (ore (1)) with respect to subsemigroup P of (N, .) but does not satisfy S-left (right) ore condition.

**DEFINITION 4.3.5**: *If $S = \{s \in N \,/\, s$ is cancellable\}, the $N_S$ if it exists and if $N_S$ is a S-near-ring then $N_s$ is called the Smarandache left (right) quotient (S-left(right) quotient) near-ring of N.*

**THEOREM 4.3.4**: *If $N_S$ is a S-quotient near-ring then $N_S$ is a quotient near-ring.*

*Proof*: Obvious by the very definition.

It is once again an interesting problem to obtain a necessary and sufficient condition for a left (right) quotient near-ring N to be a S-quotient near-ring.
Let $S(V)$ denote the collection of all S-near-rings and X be any nonempty subset which is S-semigroup under '+' or '.'.

We define Smarandache free near-ring as follows.

**DEFINITION 4.3.6**: *A S-near-ring $F_x \in S(V)$ is called a Smarandache free near-ring (S-free near ring) in V over X if there exist $f : X \rightarrow F$ for all $N \in V$ and for all $g : X \rightarrow N$ there exists a S-near-ring homomorphism $h \in S\,(Hom\,F_x,\,N)$ [Here $S\,(Hom\,(F_x,\,N))$ denotes the collection of S-homomorphism form $F_x$ to N] such that $h \circ f = g$.*

**THEOREM 4.3.5**: *Let $F_x$ be a S-free near-ring then $F_x$ is a free near-ring.*

*Proof*: Obvious by the very definition.

<u>**PROBLEMS:**</u>

1. Give an example of a near-ring which is not a S-near-ring.
2. Is $N = Z_5 \times Z_7 \times Z_3$ a S-direct product? Justify your answer.
3. Prove $N = Z_2 \times Z$ is a S-direct product.
4. Is $Z_2 = \{0, 1\}$ a S-near-ring?
5. Let $N = Z_2 \times Z_{12}$? Choose $S = Z_2 \times \{0, 2, 4, 6, 8, 10\}$ which is a S-subsemigroup of $(N, +)$. Find $N_S$.
6. Is $N_s$ given in problem 5 a S-near-ring?
7. Give an example of a S-free near-ring.

## 4.4 Smarandache ideals in Near-rings

In this section we introduce the notion of Smarandache internal direct sum, Smarandache normal sequence, Smarandache prime ideal, Smarandache semiprime ideal and the concept of Smarandache nilpotent elements and Smarandache nilpotent ideals in a near-ring. We illustrate this by examples and give some interesting related results, we propose several problems as we except the reader to get involved in this Smarandache concepts about near-rings.

**DEFINITION 4.4.1**: *N be a near-ring {$S_k$}$_{k \in K}$ be the collection of all S-ideals of N. $\sum_{k \in K} S_k = S_1 + \ldots + S_n + \ldots$ is called the Smarandache internal direct sum (S-internal direct sum) if each element of $\Sigma S_k$ has a unique representation as a finite sum of elements for different $S_k$'s.*

**THEOREM 4.4.1**: *If the near-ring has S-internal direct sum then it has internal direct sum.*

*Proof*: Follows by the very definition of the two concepts.

**DEFINITION 4.4.2**: *Let I be an S-ideal of the near-ring N. I is called a Smarandache direct summand (S-direct summand) of N if there exists an ideal J such that $N = I \dotplus J$. I is called a Smarandache strong direct summand (S-strong direct summand) of N if there exists a S-ideal J such that $N = I \dotplus J$.*

**THEOREM 4.4.2**: *If the S-ideal I of N is a S-strong direct summand then it is a S-direct summand.*

*Proof*: Left for the reader as an exercise.

**THEOREM 4.4.3**: *If the S-ideal I of N is a S-strong direct summand or S-direct summand then obviously I is a direct summand.*

*Proof*: The reader is requested to prove using the definitions.

**DEFINITION 4.4.3**: *A finite sequence $N = N_0 \supset N_1 \supset N_2 \supset \ldots \supset N_n = \{0\}$ of S-subnear-rings $N_i$ of N is called a Smarandache normal sequence (S-normal sequence) of N if and only if for all $i \in \{1, 2, \ldots, n\}$, $N_i$ is an S-ideal of $N_{i-1}$. In the special case that all $N_i$ is an S-ideal of N then we call the normal sequence an Smarandache invariant sequence; (S-invariant sequence) and n is called the Smarandache length (S-length) of the sequence and $N_{i-1} / N_i$ are called the Smarandache factors (S-factors) of the sequence.*

**DEFINITION 4.4.4**: *Let N be a near-ring. N is called Smarandache decomposable (S-decomposable) if it is the direct sum of nontrivial S-ideals or if it has nontrivial direct summand. Otherwise it is S-indecomposable. Thus even if it is indecomposable yet not as S-ideals then it is still S-indecomposable.*

**DEFINITION 4.4.5**: *Let P be a S-ideal of the near-ring N. P is called Smarandache prime ideal (S-prime ideal) if for all S-ideals I and J of N, $IJ \subseteq P$ implies $I \subset P$ or $J \subset P$.*

**THEOREM 4.4.4**: *Let N be a near-ring if P is a S-prime ideal then P is a prime ideal.*

*Proof*: Easily proved.

If P is not in general a S-prime ideal even if P is a prime ideal. The study of finding these would be an interesting problem.

**DEFINITION 4.4.6**: *Let P be an S-ideal in N, P is Smarandache semiprime (S-semiprime) if and only if for all S-ideals I of N $I^2 \subset P$ implies $I \subset P$.*

**DEFINITION 4.4.7**: *Let I be a S-ideal of N; the Smarandache prime radical (S-prime radical) of I denoted by S (P (I)) = $\bigcap_{I \subseteq P} P$. Thus if S (P (I)) is a S- prime radical then it has prime radical.*

**DEFINITION 4.4.8**: *Let N be a near-ring. A nilpotent element $0 \neq x \in N$ is said to be a Smarandache nilpotent element (S-nilpotent element) if $x^n = 0$ and there exists $y \in R \setminus \{0, x\}$ such that $x^r y = 0$ or $yx^r = 0$, r, s, > 0 and $y^m \neq 0$ for an integer m > 1.*

**THEOREM 4.4.5**: *Let N be a near-ring $0 \neq x \in N$ is S-nilpotent then x is nilpotent.*

*Proof*: Obvious by the very definition.

**DEFINITION 4.4.9**: *Let S be a S-ideal we say $S \subset N$ is nilpotent if $S^k = \{0\}$. $S \subset N$ is called Smarandache nil (S-nil) if for all $s \in S$, s is S-nilpotent.*

**PROBLEMS:**

1. Give an example of a near-ring which has S-nilpotents and S-nil ideals.
2. Let N = $Z_2 \times Z_7 \times Z$.
    a. Find an S-ideal in N.
    b. Will N be a S-strong direct summand?
3. Given N = $M_{3 \times 3}$ = {$(a_{ij})$ / $a_{ij} \in Z_2$} where $Z_2$ is a near-ring. N is near-ring:
    a. Can N have a S-normal sequence?
    b. Is N, S-decomposable?
    c. Can N have S-prime ideals?
4. Let N = $M_{5 \times 5}$= {$(a_{ij})$ / $a_{ij} \in Z_2 \times Z_2$} be a near-ring.
    a. Does N have any S-semiprime ideal?
    b. Find a semiprime ideal in N which is not S-semiprime.
5. Find for the near-ring N given in problem 4.
    a. S-ideals.
    b. Ideals which are not S-ideals.
    c. Nilpotent elements which are not S-nilpotents.
    d. Prime ideals which are not S-prime ideals.
    e. Is N decomposable or S-decomposable?

## 4.5 Smarandache Modularity in near-rings

In this section we introduce the concept of Smarandache idempotents (S-idempotents) as this is the basic notion needed for the definition of Smarandache modular ideals and other related concepts. Using the concept of S-idempotents we develop Smarandache biregular and Smarandache prime radical and obtain some interesting results about them.

**DEFINITION 4.5.1**: *Let N be a near-ring. An element $0 \neq x \in N$ is a Smarandache idempotent (S-idempotent) of R if*

      *a.  $x^2 = x$.*
      *b.  There exist $a \in N \setminus \{x, 1, 0\}$ (if N has 1). With*

            *1.  $a^2 = x$ and*
            *2.  xa = a (ax = a) or ax = x (xa = x)*

*or in (b,2) is in the mutually exclusive sense.*

*The element $a \in N$ is called the Smarandache co-idempotent (S-co-idempotent).*

**DEFINITION 4.5.2**: *Let P be a S-left ideal of a near ring N. P is called Smarandache modular (S-modular) if and only if there exists a S-idempotent $e \in N$ and for all $n \in N$; $n - ne \in P$.*

In this case we also say that P is S-modular by e and e is right identity modulo P since for all $n \in N$, ne = n (mod p).

**THEOREM 4.5.1**: *Let N be a near-ring. If $x \in N$ is a S-idempotent then x is an idempotent.*

*Proof*: Obvious by the very definition of these concepts.

**THEOREM 4.5.2**: *Let N be a near-ring. If P is S-modular than P is modular.*

*Proof*: Follows from the definition and theorem 4.5.1.

**DEFINITION 4.5.3**: *For $z \in N$ denote the S-left ideal generated by set $\{n - nz \mid n \in N\}$ by $S(L_z)$. Note that for L = N, z = 0.*

**DEFINITION 4.5.4**: *$z \in N$ is called Smarandache quasi regular (S-quasi regular) if $z \in L_z$. An S-ideal $P \subset N$ is called S-quasi regular if and only if for all $s \in S$, s is S-quasi regular.*

**DEFINITION 4.5.5**: *Let E be the set of all S-idempotents on N. A set E of S-idempotents is called orthogonal if e.f = 0 $e \neq f$. $e, f \in E$.*

**DEFINITION 4.5.6**: *Let N be a ring. An S-ideal is said to be Smarandache principal ideal (S-principal ideal) if it is generated by a single element.*

**DEFINITION 4.5.7**: *A ring R is said to be Smarandache biregular (S-biregular) if each S-principal ideal is generated by an idempotent.*

**DEFINITION 4.5.8**: *Let N be a near-ring. N is said to by Smarandache biregular (S-biregular) if there exists some set E of central S-idempotents with*

      *a.  For all $e \in N$; Ne is an S-ideal of N.*

b. *For all $n \in N$ there exists $e \in E$; $Ne = (n)$ (principal ideal generated by n).*
c. *For all $e, f \in E$   $e + f = f + e$.*
d. *For all $e, f \in E$, $ef \in E$ and $e + f - ef \in E$.*

**DEFINITION 4.5.9**: *Let I be an S-ideal. The intersection of all S-prime ideal P, such that $I \subseteq P$ is called the Smarandache prime radical (S-prime radical) of I, i.e. $S(I) = \bigcap_{I \subseteq P} P$.*

**THEOREM 4.5.3**: *Let I be a S-ideal of N. S(I) be the S-prime radical of I then S(I) is a prime radical of I.*

*Proof*: Obvious by the very definition of these two concepts.

The study of when will a prime radical be a S-prime radical is an innovative one.

**PROBLEMS:**

1.   Give an example of a near-ring N in which no idempotent is an S-idempotent.
2.   Find a near-ring N so that every idempotent of N is a S-idempotent.
3.   Let $N = Z_2 \times Z_2 \times Z_2$
    a. Does N contain a S-ideal P which is S-modular?
    b. Does N contain an ideal P which is modular?
4.   Can $M_{3 \times 3} = \{(a_{ij}) \ / \ a_i \in Z_2\}$ have S-ideals which are S-modular?
5.   a) Find an ideal I in $M_{3 \times 3}$ given in problem 4 so that S(I) is a S-prime radical of I.
    b) Can $M_{3 \times 3}$ be S-biregular?

**Chapter Five**

# SPECIAL PROPERTIES OF CLASSES OF SMARANDACHE NEAR-RINGS AND ITS GENERALIZATIONS

In this chapter we introduce several new classes of Smarandache near-rings and their generalizations; all new concepts in near-rings and seminear-rings are introduced to S-near-rings and S-seminear-rings. This chapter has five sections. In section one we introduce five types of Smarandache mixed direct products which will pave way for several types of S-near-rings and S-seminear-rings. Section two is used to define special classes of S-near-rings like S-IFP near-ring etc. The study of Smarandache group near-rings is introduced in section three, section four studies S-seminear-ring and their generalizations. The final section is devoted to the study of special properties in S-near-rings.

## 5.1 Smarandache mixed direct product of near-rings and seminear rings

This section is completely devoted to the introduction of Smarandache mixed direct product of near-rings. Only by defining these Smarandache mixed direct product we are able to define several levels of S-near rings and S-seminear-rings and S-quasi seminear-rings. We illustrate them with examples and denote how each of these S-mixed products yield different levels of S-near-rings and S-seminear-rings.

**DEFINITION 5.1.1**: *Let $N = N_1 \times N_2 \times \ldots \times N_n$ be the Smarandache mixed direct product I (S-mixed direct product I) of seminear-rings and near-rings where $N_i$'s are from the classes of near-rings and seminear-rings which we call as S-mixed direct product I.*

The use of S-mixed direct product I is two-fold.

**THEOREM 5.1.1**: *Let $N = N_1 \times N_2 \times \ldots \times N_n$ be S-mixed direct product I, then N is a S-seminear-ring.*

*Proof*: Obvious by the very definition.

**THEOREM 5.1.2**: *Let $N = N_1 \times \ldots \times N_r$ be S-mixed direct product I, then N has S-subnear-rings and S-ideals.*

*Proof*: Direct by using definitions.

**DEFINITION 5.1.2**: *Let $N = N_1 \times N_2 \times \dots \times N_r$ where some $N_i$'s are near-rings and some of them are near-fields. Then this Smarandache mixed direct product of near-rings is called S-mixed direct product II.*

It is easily verified that S-mixed direct product II gives S-near-rings. Further these S-mixed direct product II when r > 2 can give S-subnear-rings. It is also useful in getting S-N groups. It paves way for S-quasi subnear-rings.

It is also useful to define notions like S-invariant subnear-rings related to near-field as we can take in the S-mixed direct product II more than one near-field. It is pertinent to mention here that S-mixed direct product near-rings I and II are entirely different from S-direct product of near-rings. Thus we can say a S-mixed direct product of near-rings is a S-direct product if and only if one of the $N_i$ is a S-near-ring.

Thus if we assume none of the $N_i$ is a S-near-ring then certainly the S-mixed direct product I or II is a S-near ring but the product is never a S-direct product of near-rings. But they are S-near-rings or S-seminear-ring but they are different in their own ways. Now we call these S-near-rings and S-seminear-rings as level I. Now we proceed on to define Smarandache near-rings and Smarandache seminear-rings of level II.

**DEFINITION 5.1.3**: *Let $N = N_1 \times \dots \times N_r$ be the Smarandache mixed direct product III (S-mixed direct product III) of near-rings and rings i.e. some of the $N_i$'s are rings and certainly some $N_j$'s are near-rings. N is not a S-near-ring unless one of the $N_i$'s is a near-field.*

We now define Smarandache near-ring of level II.

**DEFINITION 5.1.4**: *Let N be a near-ring. N is said to be a Smarandache near-ring of level II (S-near-ring II) if N has proper subset P, $P \subset N$ where P is a ring.*

Thus we see S-mixed direct product II is always a S-near-ring II. We see only by this we can get S-near-ring II for it is not easy to obtain subset in a near-rings which are rings. We see S-near-ring I and S-near-ring II are in general not related; they happen to be two distinct classes.

**DEFINITION 5.1.5**: *Let $N = N_1 \times \dots \times N_t$ where some of the $N_i$'s are the seminear-ring and some of them are semirings then N is called the Smarandache mixed direct product IV (S-mixed direct product IV) of seminear-rings.*

The S-mixed direct product IV paves way for the following definition.

**DEFINITION 5.1.6**: *Let N be a seminear-ring. N is said to be a Smarandache seminear-ring of level II (S-seminear-ring of level II) if N contains a proper subset P which is a semiring. Clearly the S-mixed direct product IV gives S-seminear-ring II.*

Here also we wish to state that by no means in general S-seminear-ring I and S-seminear-ring II are related in fact they are disjoint notions but not totally disjoint which is seen by the following example.

***Example 5.1.1***: Let $N = Z^o \times Z_{15}$ where $Z^o = Z^+ \cup \{0\}$ is a semiring under usual '+' and '.'. $(Z_{15}, \times, '.')$ is a seminear-ring. Clearly N is a S-seminear-ring II. This is also a S-seminear-ring I as $Z_{15}$ has a proper subset P where P is a near-ring.

Thus we cannot say in general that a S-seminear-ring I is never a S-seminear-ring II.

To obtain a Smarandache pseudo seminear-ring we are not in any position to construct S-mixed direct products. It is interesting to note that in case of non-associative near-rings (NA-near-rings) and non-associative seminear-rings (NA-seminear-rings) every S-mixed direct product of various types helps us to build S-NA-seminear-rings and S-NA-near-rings.

**DEFINITION 5.1.7**: *Let $N = N_1 \times N_2 \times \ldots \times N_n$ be the Smarandache mixed direct product V (S-mixed direct product V) of non-associative seminear-rings, rings and non-associative near-rings. Clearly N is a S-quasi near-ring.*

*Thus using this S-mixed direct product V we get a class of S-quasi near-rings. If we include semirings in the S-mixed direct product we can get S-quasi seminear-rings.*

Thus whenever we need these S-mixed direct products we would be using them, but what is more important is that we are not able to construct by any natural means NA-near-ring or NA-seminear-rings. The only way is that we get NA-near-rings and NA-seminear-rings by building near loop-rings and groupoid near-rings, loop seminear-rings and groupoid seminear-rings.

## PROBLEMS:

1. Let $N = Q \times Z_{15} \times Z_2$ (where $Z_2$ is the near-field, Q the field of rationals and $Z_{15}$ a seminear-ring) be the S-mixed direct product. Find whether N is a S-near-ring I or a S-seminear-ring I or S-near-ring II or S-seminear-ring II?

2. Let $N = Z_3 \times Z^o \times Z_{18}$ be the S-mixed direct product, where $Z_3$ is the prime field of characteristic 3, $Z^o$ the semiring and $Z_{18}$ a seminear-ring. Is N a S-mixed direct product I, or II, or III, or IV, or V or none? Does N have S-seminear-rings and S-near-rings? Justify your answers.

3. Let $Z_{16}G$ be the groupoid seminear-ring of the groupoid G over the seminear-ring $Z_{16}$. The groupoid G is given by the following table:

| * | a | b | c | d |
|---|---|---|---|---|
| a | a | c | a | c |
| b | a | c | a | c |
| c | a | c | a | c |
| d | a | c | a | c |

Is $Z_{16}G$ a S-near-ring? Is $Z_{16}G$ a S-seminear-ring ? Is $Z_{16}G$ at least a S-quasi seminear-ring?

4.  Let $N = Z_{10} \times Z_8 \times Z_2$ be the S-mixed direct product of near-ring $Z_{10}$, $Z_8$ seminear-rings and $Z_2$ the field of characteristic 2. Is N a S-near-ring? Is N a S-seminear-ring?

5.  Let N be a S-mixed direct product of $N = Z^o \times Z_{12} \times Z_p$ where $Z^o$ is a semiring, $Z_p$ a field and $Z_{12}$ a seminear-ring. Is N a S-quasi near-ring? Is N a S-quasi seminear-ring? Substantiate your claim.

6.  Give an example of a S-quasi seminear-ring which is not a S-quasi near-ring.

7.  Give an example of a S-seminear-ring which is not a S-near-ring.

## 5.2 Special Classes of Smarandache near-rings

In this section we introduce some basic well-known classes of near-rings which are S-near-rings. We do not claim to give all the results pertaining to them. We provide only the main definitions and leave the rest for the reader to develop. Concepts like S-IFP, S-n-ideal near-rings, S-LSD near-rings and S-equiprime near-rings are introduced.

**DEFINITION 5.2.1**: *Let N be a S-near-ring. N is said to fulfil the Smarandache insertion of factors property (S-IFP for short) if for all a, b $\in$ N we have a.b = 0 implies anb = 0 for all n $\in$ P, P $\subset$ N, where P is a near-field.*

We say N has S-strong IFP property if every homomorphic image of N has the IFP property or we redefine this as:

**DEFINITION 5.2.2**: *Let N be a S-near-ring and I a S-ideal of N. N is said to fulfil Smarandache strong IFP property (S-strong IFP property) if and only if for all a, b $\in$ N, ab $\in$ I implies anb $\in$ I where n $\in$ P, P $\subset$ N and P is a near-field.*

**DEFINITION 5.2.3**: *Let p be a prime. A S-near-ring N is called a Smarandache p-near-ring (S-p-near-ring) provided for all x $\in$ P, $x^p = x$ and px = 0, where P $\subset$ N and P is a semifield. Clearly if N is itself a p-near-ring and N is a S-near-ring then N will trivially be a S-p-near-ring. Conversely if N is a S-p-near-ring N in general need not be a S-p-near-ring.*

**DEFINITION 5.2.4**: *Let N be a near-ring. A S-right ideal I of N is called Smarandache right quasi reflexive (S-right quasi reflexive) if whenever A and B are S-ideals of N with AB $\subset$ I, then b(b' + a) − bb' $\in$ I for all a $\in$ A and for all b, b' $\in$ B.*

**DEFINITION 5.2.5**: *Let N be a near-ring. N is said to be Smarandache strongly subcommutative (S-strongly subcommutative) if every S-right ideal of it is S-right quasi reflexive.*

**DEFINITION 5.2.6**: *Let N be a near-ring. S a S-subnormal subgroup of (N, +). S is called a Smarandache quasi-ideal (S-quasi ideal) of N if SN $\subset$ N and NS $\subset$ S where by NS we mean elements of the form {n(n' + s) − nn'/ for all s $\in$ S and for all n, n' $\in$ N} = NS.*

**DEFINITION 5.2.7**: *A left-near-ring N is said to be Smarandache left-self-distributive (S-left self distributive) if the identity abc = abac is satisfied for a, b, c ∈ A, A ⊂ N and A is a S-subnear-ring of N.*

**DEFINITION 5.2.8**: *Let N be a near-ring. N is said to be Smarandache left permutable (S-left permutable) if abc = bac for all a, b, c ∈ A, A a subnear-ring of N.*

*The near-ring N is said to be Smarandache right permutable (S-right permutable) if abc = acb for all a, b, c ∈ A; A ⊂ N, A a subnear-ring of N. The near-ring N is said to be Smarandache medial (S-medial) if abcd = acbd for all a, b, c, d in A where A is a S-subnear-ring of N.*

*The near-ring N is said to be Smarandache right self-distributive (S-right self distributive) if abc = acbc for all a, b, c in A, A a S-subnear-ring of N.*

**DEFINITION 5.2.9**: *Let (R, +, .) be a S-ring and M a S-right R-module. Let W = R × M and define (α, s) ⊙ (β, t) = (αβ, sβ + t). Then (W, +, ⊙) is a S-left near-ring, the abstract Smarandache affine near-ring (S-affine near-ring) inducted by R and M. If R is S-right self-distributive and $MR^2 = 0$ then W is a S-right self-distributive near-ring. If MR = 0 and R is S-medial, S-left permutable or S-left self-distributive then (W, ⊙) is S-meidal, S-left permutable or S-left self-distributive respectively.*

**DEFINITION 5.2.10**: *Let R be a S-near-ring. R is said to be Smarandache equiprime (S-equiprime) if for all 0 ≠ a ∈ P ⊂ R. P the near-field in R and for x, y ∈ R, arx = ary for all r ∈ R implies x = y. If B is a S-ideal of R, B is called a Smarandache equiprime ideal (S-equiprime ideal) if R/B is an S-equiprime near-ring.*

**DEFINITION 5.2.11**: *Let (N, +, .) be a triple. N is said to be a Smarandache infra-near-ring (S-INR) where*

    *1. (S, +) is a S-semigroup.*
    *2. (N, .) is a semigroup.*
    *3. (x + y)z = xz − 0z + yz for all x, y, z ∈ N.*

**DEFINITION 5.2.12**: *Let I be a S-left ideal of N. Suppose I satisfies the following conditions:*

    *1. a, x, y ∈ N, anx − any ∈ I for all n ∈ N implies x − y ∈ I.*
    *2. I is left invariant.*
    *3. 0N ⊂ I.*

*Then I is called a Smarandache equiprime left-ideal (S-equiprime left ideal) of N.*

**DEFINITION 5.2.13**: *Let N be a S-near-ring. N is said to be Smarandache partially ordered (S-partially ordered) by ≤ if*

    *1. ≤ makes (N, +) into a partially ordered S-semigroup.*
    *2. for all n, n' ∈ N, n ≥ 0 and n' > 0 implies nn' ≥ 0.*

**DEFINITION 5.2.14**: *Let N be a near-ring. Let $I_1, I_2, \ldots, I_m$ be the collection of all S-right ideals of N. We say N is a Smarandache n-ideal near-ring (S-n-ideal near-ring) (n an integer less than or equal to m) if for every n, S-right ideals $I_1, \ldots, I_n$ and for every n elements $X_1, X_2, \ldots, X_n$ in $N \setminus (I_1 \cup I_2 \cup \ldots \cup I_n)$ we have $\langle X_1 \cup I_1 \cup \ldots \cup I_n \rangle = \langle X_2 \cup I_1 \cup \ldots \cup I_n \rangle = \ldots = \langle X_n \cup I_1 \cup \ldots \cup I_n \rangle$ where $\langle\ \rangle$ denotes the S-right ideal generated by $X_i \cup I_1 \cup \ldots \cup I_n$; $1 \leq i \leq n$.*

**DEFINITION 5.2.15**: *Let (M, +) be a S-semigroup (not necessarily abelian) and $\Gamma$ a non-empty set. Then M is said to be a Smarandache $\Gamma$-near ring (S-$\Gamma$-near-ring) if there exists a mapping $M \times \Gamma \times M \to M$ the images (a, $\alpha$, b) denoted by $a\alpha b$ satisfying the following conditions:*

1. *$(a + b)\alpha c = a\alpha c + b\alpha c$.*
2. *$(a\alpha b)\beta c = a\alpha(b\beta c)$.*
   *for all a, b, c $\in$ M and $\alpha$, $\beta \in \Gamma$.*

*M is said to be Smarandache zero symmetric $\Gamma$-near-ring (S- zero symmetric $\Gamma$-near-ring) if $a\alpha 0 = 0$ for all $a \in M$ and $\alpha \in \Gamma$ where 0 denotes the additive identity in M.*

**DEFINITION 5.2.16**: *Let M be a S-$\Gamma$-near-ring, then a S-normal subsemigroup I of (M, +) is called*

1. *a S-left ideal if $a\alpha(b + i) - a\alpha b \in I$ for all a, b $\in$ M, $\alpha \in \Gamma$ and $i \in I$.*
2. *a S-right ideal if $i\ \alpha\ a \in I$ for all $a \in M$, $\alpha \in \Gamma$ and $i \in I$ and*
3. *a S-ideal if it is both S-left and S-right ideal.*

**DEFINITION 5.2.17**: *An S-ideal A of M is said to be Smarandache prime (S-prime) if B and C are S-ideals of M such that $B\Gamma C \subset A$ implies $B \subseteq A$ or $C \subseteq A$.*

**DEFINITION 5.2.18**: *Let $M_1$ and $M_2$ be any two S-$\Gamma$-near-rings. A S-S-semigroup homomorphism f of $(M_1, +)$ into $(M_2, +)$ is called a Smarandache $\Gamma$-homomorphism (S-$\Gamma$-homomorphism) if $\Gamma(x\alpha y) = f(x)\ \alpha\ f(y)$ for all x, y $\in$ M and $\alpha \in \Gamma$. We say f is a S-$\Gamma$-isomorphism if f is one to one and onto. For an S-ideal I of a S-$\Gamma$ near-ring, the S-quotient $\Gamma$near-ring M/I is defined in a similar way as quotients are defined for any near-ring.*

Now several interesting Smarandache results can be obtained on S-$\Gamma$-near-rings. The reader is entrusted with that work. Now we proceed on to define the concept of Smarandache left-self-distributive near-ring.

**DEFINITION 5.2.19**: *Let N be a near-ring. We say N is a Smarandache left self-distributive near-ring (S-LSD near-ring) if N has a proper S-subnear-ring which is a LSD, we do not need the totality of the near-ring to be a LSD near-ring.*

**THEOREM 5.2.1**: *If N is a LSD near-ring having a S-subnear-ring then N is a S-LSD near-ring.*

*Proof*: Straightforward.

Give an example of a LSD near-ring which is not a S-LSD near-ring. Also give an example of S-LSD near-ring which is not a LSD near-ring.

<u>**PROBLEMS:**</u>

1. Let $(G, +)$ be a S-semigroup, let h be an S-idempotent endomorphism on $(G, +)$ define $a \otimes b = h(b)$ for each $a, b \in A$, A a subgroup of G. Is $(G, +, \otimes)$ a S-right self-distributive near-ring? Is $((G, +, \otimes)$ a S-left self-distributive near-ring? Substantiate your claim.

2. Let R be a S-left distributive near-ring such that $a, b, c \in R$ and X a non-empty subset of R. Is $abc = aba^nc = ab^nc$ for $n \geq 1$? If this is true can we say $a^3 = a$ is an idempotent for $n \geq 3$ and if $abc = 0$ then $bac = 0$?

3. Is the near-ring $Z_{15} = \{0, 1, 2, \ldots, 14\}$ an S-n-ideal near-ring?

4. Test whether $N = Z_4 \times Z_{21}$ is a S-n-ideal near-ring where $Z_4$ and $Z_{21}$ are near-rings.

5. Let $(G, +)$ be a group, X a non-empty set. Let $M = \{f/f : X \to F\}$. Then M is a group under pointwise addition. Prove by suitable definition M is S-$\Gamma$-near-ring, where $\Gamma$ is a set of all mappings of G into X.

6. Let I be an S-ideal of M and f the canonical S-semigroup, S-epimorphism from M on to M/I.

   i)   Prove or disprove f is a Smarandache $\Gamma$-homomorphism of M on to M/I.
   ii)  What is its kernel?

## 5.3 Smarandache Group near-rings and their generalization

When the study of group rings is very vast with over a dozen texts on them we see group near-rings is very rare only less than half a dozen research papers and all the more the introduction of even semigroup near-rings is totally absent. But as we don't have several classes or concrete examples of near-rings except those near-rings got from $Z_n$ or Z and some as group morphisms; the concept of group near-rings and the introduction of semigroup near-rings will certainly serve to provide us with a vast stretch of examples as both the algebraic structures groups and semigroups are very non-abstract with several concrete examples.

Thus in this section we introduce the concept of Smarandache group near-rings and Smarandache semigroup near-rings and generalize them and built Smarandache group seminear-rings and Smarandache semigroup seminear-rings. We obtain several interesting and innovative illustrations and results and give four distinct classes of S-near-rings and S-seminear-rings. Several of the open problems about group rings are also reformulated in the case of group near-rings and semigroup near-rings.

**DEFINITION 5.3.1**: *Let G be a group and N a S-near-ring. Then we define the group near- ring NS to be a Smarandache group near-ring (S-group near-ring).*

It is important to note that all group near-rings are not S-group near-rings. By an example we illustrate this.

*Example 5.3.1*: Let $Z_2 = (0, 1)$ be the near-field and $G = \langle g \, / \, g^3 = 1 \rangle$. The group near-ring $Z_2G$ is not a S-group near-ring.

Thus we make it very clear that if N is not a S-near-ring but the group near-ring NG happens to be a S-near-ring still we do not call NS the S-group near-ring, we can say group near-ring is a S-near-ring. Thus we can redefine S-group near-ring NS as:

**DEFINITION 5.3.2**: *Let N be a near-ring and G any group. The group near-ring NS is a S-group near-ring if and only if N is a S-near-ring.*

Thus in this section we will be mainly interested in the study of S-group near-rings than the S-near-rings.

**THEOREM 5.3.1**: *Let NG be a S-group near-ring. Then NG has a S-subnear-ring.*

*Proof*: Since it is given N has a proper subset P which is a near-field we can take PG and PG $\subset$ NG and P $\subset$ PG so PG is a S-subnear-ring. Thus if NG is a S-group near-ring we can always say NG contains a S-subnear-ring.

Now the natural question would be if the group near-ring NG has a S-subnear-ring; can we say NG is a S-group near-ring? The answer to this question in general is true. For we may have a S-subnear-ring of NG got from a proper subset of NG and not of N so NG with the properties of G may yield a S-subnear-ring yet NG may fail to be a S-group near-ring.

We just recall if G is a group, $\Delta = \Delta (G) = \{x \in G \, / \, [G : C_G (x)] < \infty\}$ since the conjugates of x are in one to one correspondence with the right cosets of $C_G(x)$, it follows that x has only finitely many conjugates if and only if x $\in$ $\Delta$. We observe $\Delta$ is a normal subgroup of G. For clearly $1 \in \Delta$ and since $C_G (x^{-1}) = C_G(x)$ we see that x $\in$ $\Delta$ implies $x^{-1} \in \Delta$. $\Delta$ is called the finite conjugate subgroup of G. The importance of $\Delta$ is two fold. Firstly reduce the problem studied from NG to N$\Delta$ and second we are able to handle the much simpler group $\Delta$.

Clearly if NG is a S-group near-ring we see N$\Delta$ is also a S-group near-ring. But N$\Delta \subseteq$ NG. Let $\theta$ denote the projection.

$\theta : NG \rightarrow N\Delta$ given by $\alpha = \sum_{g_i \in G} n_i g_i \rightarrow \theta(\alpha) = \sum_{g_i \in \Delta} n_i g_i$, then $\theta$ is clearly a map and not a near-ring homomorphism.

a.  Suppose A and B be S-ideals of NG; then can we say θ(A) is a S-ideal of NΔ?

b.  If A ≠ 0 can we say θ (A) ≠ 0 and conversely?

c.  Can we prove if AB = 0 then it implies that θ (A) θ(B) = 0?

The above 3 questions are left for the reader to settle.

The next natural question would be if N is any near-field and G a torsion free abelian group can NG have zero divisors or suppose if H is a torsion free abelian subgroup of G, G any group can we say if α ∈ NH ⊂ NG with α ≠ 0 then α is not a zero divisor of NG? The study of these problems is related to the zero divisor problem for group rings, here we have only changed rings to near-rings.

Now we go for the introduction of the several of the properties enjoyed by group rings to group near-rings.

Now let M ⊂ N be a near-field in the near-ring N and G any group. Let MG be the subgroup near-ring of the S-group near-ring NG. Prove or disprove MG is prime if and only if Δ(G) is torsion free abelian. Yet another important question is can we say for any finite group G the S-group near-ring NG is artinian? This study is also dormant, a host of research can be done in this direction also. Another important study is if NG is a S-group near-ring, G a finite group; will NG satisfy N-artinian condition? To this end we introduce the S-artinian condition for S-near-rings N.

**DEFINITION 5.3.3**: *Let N be a near-ring. If the S-ideals of the near-ring N fulfils D.C.C conditions then we say N satisfies Smarandache DCC (S-DCC) conditions.*

*If the S-ideal of a near-ring N satisfies ACC condition then we say N satisfies Smarandache ACC (S-ACC) condition. If SDCC is satisfied by S-left ideals we distinguish it by S-DCC L. Similarly in case of S-right ideals we denote it by S-DCC R. If we have collection of S-N-subsemigroups of a near-ring N and if it satisfies DCC or ACC we denote it by S DCC N or S ACC N.*

*Thus if N is a near-field and G a finite group will the group near-ring NG satisfy*

       *1.  S-DCC R or*
       *2.  S-DCC N or*
       *3.  S-DCC L or*
       *4.  S-ACC R or*
       *5.  S-ACC L or*
       *6.  S-ACC N.*

This study is important and may yield a lot of results. For example if we take the near-ring to be the near-field $Z_2$ and G a solvable group $Z_2G$; the group near-ring will satisfy several of the conditions mentioned. Finally as the study in this direction is totally absent by researchers it may give any researcher a lot of scope for research in this direction.

We recall if N is a near-ring and P a S-prime ideal of N then N / P is a S-prime near-ring, we will say the near-ring N is S-semiprime if and only if the intersection of all S-prime ideals is 0.

Now let G be any group and N a near-field. When will the group near-ring NG be S-semiprime. This question too is not very simple it needs condition on N as well as G to be discussed. We have in case of group rings, KG is semiprime if K is a field of characteristic zero.

Another natural question about group near-ring and S-group near-ring is: Will the group near-ring or S-group near-ring contain nil ideals. The problem in case of prime fields of characteristic p and when G has no elements of order p is dealt by [58]. They have proved KG the group ring has no nonzero nil ideals. This question in case of group near-rings remains as an open problem.

The study of near-rings satisfying polynomial identity is a study in another direction as such study has never been contemplated even for near-rings we propose some problems in chapter X and accept our inability to discuss about it in this book. We now proceed on to define the concept of Smarandache near module. The near module is nothing but the N-group so the Smarandache near module is defined as the Smarandache N-subsemigroup. We end this study of group near-rings with

1. When is the S-group near-ring NG, S-regular?
2. When is the group near-ring regular?

It is well known [58], the group ring KG is regular if and only if G is a locally finite group and has no elements of order p in case K has characteristic p > 0.

Now we have already mentioned that the study of semigroup near-rings has not even started till date. Hence we in chapter 3 have mentioned some properties and have introduced the notion of semigroup near-rings. In this section we define Smarandache semigroup near-rings and indicate some more research in this direction.

**DEFINITION 5.3.4**: *Let N be a any near-ring and S a semigroup. The semigroup near-ring NS is said to be a Smarandache semigroup near-ring (S-semigroup near-ring) if S is a S-semigroup.*

Now even if N is a S-near-ring and S is not a S-semigroup we call NS a semigroup near-ring which is a S-near-ring and it is not a S-semigroup near-ring. It is pertinent to mention here that the S-semigroup near-ring may not in general be a S-near-ring. Even if NS is a S-near-ring; NS need not be a S-semigroup near-ring.

**THEOREM 5.3.2**: *Let NS be a S-semigroup near-ring; then NS has a S-group near-ring if and only if N is a S-near-ring.*

*Proof*: Given NS is a S-semigroup near-ring so S is a S-semigroup i.e. S has proper subset which is a subgroup say H. Now NH is a S-group near-ring if and only if N is a S-near-ring.

**THEOREM 5.3.3**: *Let N be a S-near-ring and S a S-semigroup. The S-semigroup near-ring NS has nontrivial S-subnear-rings.*

*Proof*: Obvious from the fact S is a S-semigroup so; G ⊂ S is a group; NG is a S-subnear-ring of NS, hence the claim.

Obtain any other interesting result about them.

*Example 5.3.2*: Let $Z_2 = \{0, 1\}$ be a near-ring and S(n) be the S-semigroup. $Z_2S(n)$ is a S-semigroup near-ring. Find S-ideals in $Z_2S(n)$. Is $Z_2S(n)$ a S-near-ring? Find zero divisors if any in $Z_2S(n)$. Does $Z_2S(n)$ have S-zero divisor?

<u>**PROBLEMS :**</u>

1. Let NG be a S-group near-ring Let A be an S-ideal of NG. Prove θ (A) is an S-ideal of NΔ, where θ : NG → NΔ is a projection.
2. In problem 1 can you prove. "A ≠ 0 if and only if θ (A) ≠ 0?
3. In problem 1 if AB = 0 prove θ (A) θ (B) = 0.
4. Study the same problems 1 to 3 when NG is not a S-group near-ring but only a S-near ring.
5. Prove or disprove the following: Let N be a S-near-ring i.e. M ⊂ N; M a near-field; G any group such that Δ(G) is torsion free abelian. NG is a S-group near-ring.
   a. Does it imply MG is prime if and only if Δ (G) is torsion free abelian?
   b. Does it imply MG is prime if and only if G has no non-identity finite normal subgroup?
6. Let G be a finite group; N a S-near-ring;
   a. Can the S-group near-ring NG satisfy Artinian conditions?
   b. Can the S-group near-ring NG satisfy S-Artinian condition?
7. Let $Z_{12}$ be the near-ring and $S_5$ be the group. Is the group near-ring $Z_{12}S_5$ S-semiprime?
8. Let $Z_5$ be the near-field and $S_7$ be the group. Will the near group ring $Z_5S_7$ satisfy the S-Artirian condition. Is it S-Notherian? Justify your answer.
9. Let $Z_{12}$ be a near-ring and S(5) the S-semigroup. $Z_{12}S(5)$ is a S-semigroup near-ring. Is $Z_{12}S(5)$ a S-near-ring? Find S-ideals if any in $Z_{12}S(5)$. Can $Z_{12}S(5)$ have zero divisors, S-zero divisors?
10. Let $Z_5S(7)$ be the S-semigroup near-ring. Find whether $Z_5S(7)$ is S-regular or regular? Or is $Z_5S(7)$ a S-near-ring?

## 5.4 On a special class of Smarandache seminear-rings and their generalizations

The study of seminear-rings is itself very meagre. With the well-known example of $Z_n$, n a composite number where '×' is defined as usual multiplication on $Z_n$ modulo n and '.' defined a . b = a or a . b = 0 for all a, b ∈ $Z_n$. Then {$Z_n$, '×', .} becomes a seminear-ring. To obtain a nontrivial class of seminear-rings we have ventured to study group seminear-rings and semigroup seminear rings. Now here we proceed on to study and define Smarandache seminear-rings (S-seminear-rings) and study them. The study or the definition of S-seminear-rings is given in [99].

**DEFINITION [99]:** *A non-empty set N is said to be a Smarandache seminear-ring (S-seminear-ring) if (N, +, .) is a seminear-ring having a proper subset A, (A ⊂ N) such that A under the same operations of N is a near-ring, that is (A, +, .) is a near-ring.*

*Example 5.4.1:* Let $Z_{10} = \{0, 1, 2, …, 9\}$. $(Z_{10}, \times, .)$ is a seminear-ring. $(Z_{10}, \times, .)$ is a S-seminear-ring as A = $\{2, 4, 6, 8\}$ is such that $(A, \times, .)$ is a near-ring.

*Example 5.4.2*: Let $Z_8 = \{0, 1, 2, …, 7\}$. $(Z_8, '\times', '.')$ is a seminear-ring. Clearly $Z_8$ is also a S-seminear-ring as $(A, \times, .)$ with A = (1, 3) is a near-ring.

**THEOREM 5.4.1:** *Let $Z_n = \{0, 1, 2, …, n − 1\}$, n a composite number; $(Z_n, '\times', .)$ is a S-seminear-ring.*

*Proof*: For take A = $\{1, n − 1\} \subset Z_n$ then $(A, \times, .)$ is a near-ring. We assume n is a composite number for otherwise if n is a prime $(Z_p, \times, .)$ itself is a near-ring as $Z_p \setminus \{0\}$ under '$\times$' is a group.

We have now seen only examples of S-seminear-rings N in which the semigroup (N, $\times$) is a commutative semigroup. To obtain non-commutative S-seminear-ring we proceed on as follows:

*Example 5.4.3*: Let $(Z_6, \times, .)$ be a seminear-ring and G any group. We see the group seminear-ring $Z_6G$ is a non-commutative seminear-ring provided G is taken to be non-commutative group.

*Example 5.4.4*: Let $(Z_4, \times, .)$ be a seminear-ring and G = $S_3$. The group seminear-ring $Z_4S_3$ is a S-seminear-ring which is non-commutative.

**DEFINITION 5.4.1**: *Let G be any group, N a seminear-ring, the group seminear-ring, NG is a Smarandache group seminear-ring (S-group seminear-ring) if and only if N is a S-seminear-ring. Thus we define in this way. Suppose N is not a S-seminear-ring and G any group even if the group seminear-ring is a S-seminear-ring we still call it only as a S-seminear-ring and not as a S-group seminear-ring.*

We now proceed on to define yet another new class of seminear-rings viz. Smarandache semigroup seminear-rings. The concept of semigroup seminear-rings was introduced in Chapter 3.

**DEFINITION 5.4.2**: *Let N be a seminear-ring, and S a S-semigroup. We call the semigroup seminear-ring to be a Smarandache semigroup seminear-ring (S-semigroup seminear-ring) thus we say NS the semigroup seminear-ring is a S-semigroup seminear-ring if and only if the semigroup S is a S-semigroup.*

Here also we may or may not have the S-semigroup seminear-ring to be a S-seminear-ring. If S is not a S-semigroup and N a seminear-ring even if the semigroup seminear-ring NS, is a S-seminear-ring still we do not call NS a S-semigroup seminear-ring.

**THEOREM 5.4.2**: *Let S be a semigroup with 1 and N a S-seminear-ring. Then the semigroup seminear-ring NS is a S-seminear-ring.*

*Proof*: Straightforward. Hence left for the reader to prove.

**THEOREM 5.4.3:** *Let S be a S-semigroup and N a S-seminear-ring. The S-semigroup seminear-ring NS has subseminear-rings which are group seminear-rings and group near-rings.*

*Proof*: Follows from the very definition of these concepts.

**THEOREM 5.4.4:** *Let S be a S-semigroup, and N a S-seminear-ring. The S-semigroup seminear-ring. NS has nontrivial S-group seminear-rings.*

*Proof*: Obvious by the very definitions.

**THEOREM 5.4.5:** *Let N be a seminear-ring and S be S-semigroup having S-subsemigroup. Then the semigroup seminear-ring NS has S-subseminear-ring.*

*Proof*: Straightforward by the definitions.

**DEFINITION 5.4.3:** *Let N be a seminear-ring. If (N, +) has a S-normal subgroup then the seminear-ring N is said to have Smarandache normal subseminear-ring (S-normal subseminear-ring).*

***Example 5.4.5***: Can $\{Z_{10}, '\times', '.'\}$ have S-normal subseminear-ring?

**DEFINITION 5.4.4:** *Let (N, +, .) be a seminear-ring. If in the S-semigroup (N, +) every proper subset (A, +) which is a group is commutative then we say the S-seminear-ring N is Smarandache commutative (S-commutative). Thus if (N, +, .) is commutative and if (N, +) is S-semigroup then trivially (N, +, .) is S-commutative. Secondly if (N, +, .) is S-commutative seminear-ring then N need not in general be commutative.*

**DEFINITION 5.4.5***: Let (N, +, .) be a seminear-ring. If (N, +) is a S-semigroup and if N has at least one proper subset which is a subgroup that is commutative then we say the seminear-ring (N, +, .) is a Smarandache weakly commutative seminear-ring (S-weakly commutative seminear-ring).*

**THEOREM 5.4.6.**: *Let (N, +, .) be a seminear-ring which is S-commutative then N is S-weakly commutative.*

*Proof*: Proof is direct and the reader is expected to prove.

**DEFINITION 5.4.6**: *Let (N, +, . ) be a seminear-ring. (N, +) be a S-semigroup such that every proper subset A of N which is a group is a cyclic subgroup then we say the seminear-ring (N, +, .) is a Smarandache cyclic seminear-ring (S-cyclic seminear-ring). In particular (N, +, .) has at least one proper subset which is a cyclic group, we call (N, +, .) a Smarandache weakly cyclic seminear-ring (S-weakly cyclic seminear-ring).*

The following theorem is an easy consequence of the definitions, hence left as an exercise for the reader to prove.

**THEOREM 5.4.7:** *Let (N, +, .) be a S-cyclic seminear-ring then N is a S-weakly cyclic seminear-ring.*

The reader is requested to construct an example of a S-weakly cyclic seminear-ring which is not a S-cyclic seminear-ring.

**DEFINITION 5.4.7**: *Let N be any seminear-ring. If (N, +) has two subgroups say (A, +) and ($A_1$, +) such that A and $A_1$ are conjugate then we say the seminear-ring N is a Smarandache conjugate seminear-ring (S-conjugate seminear-ring).*

**DEFINITION 5.4.8:** *Let N and $N_1$ be any two S-seminear-rings. We say a map $\phi$ from N to $N_1$ is a Smarandache seminear-ring homomorphism (S-seminear-ring homomorphism) from A to $A_1$ where $A \subset N$ is a near-ring and $A_1 \subset N_1$ is a near-ring and $\phi(x + y) = \phi(x) + \phi(y)$, $\phi(xy) = \phi(x)\phi(y)$ where $x, y \in A$ and $\phi(x), \phi(y) \in A_1$ and it is true for all $x, y \in A$. We need not even have the map $\phi$ to be well-defined or even defined on the whole of N. The concept of Smarandache isomorphism (S-isomorphism) etc. are defined in a similar way.*

**Example 5.4.6:** Let N = ($Z_{12}$ , × , .) be a S-seminear-ring and $N_1$ = ($Z_4$ , × , .) be a S-seminear-ring $\phi$ : N to $N_1$ is a S-seminear-ring isomorphism as $\phi$: A $\rightarrow$ $A_1$ is a near-ring isomorphism, here A = {1, 11} and $A_1$ = {1, 3}. Thus |N| = 12 and |$N_1$| = 4, yet N and $N_1$ are S-isomorphic under this $\phi$.

We define level II Smarandache seminear-rings as follows:

**DEFINITION 5.4.9:** *Let (N, +, .) be a seminear-ring. We say N is a Smarandache seminear-ring of level II (S-seminear-ring II) if N has a proper subset P where P is a semiring.*

**Example 5.4.7:** Let N = $Z_6 \times Z^o$ be the S-mixed direct product of the seminear-ring ($Z_6$ , × , .) and the semiring ($Z^o$, +, .) where $Z^o = Z^+ \cup \{0\}$. Clearly N is a seminear-ring which is a S-seminear-ring II.

**DEFINITION 5.4.10:** *Let N be a S-seminear-ring II, we say N is a Smarandache strict seminear-ring (S-strict seminear-ring) if A is a strict semiring ( $A \subset N$ and A is semiring).*

**DEFINITION 5.4.11**: *Let N be a S-seminear-ring II we say N is a Smarandache commutative seminear-ring II (S-commutative seminear-ring II) if every proper subset A of N which is a semiring is a commutative semiring. If the S-seminear-ring II has at least one proper subset which is a semiring and is commutative we call the N a Smarandache weakly commutative seminear-ring II (S-weakly commutative seminear-ring II).*

**DEFINITION 5.4.12:** *Let N and $N_1$ be two S-seminear-ring II, we say a map $\phi : N$ to $N_1$ is a Smarandache seminear-ring homomorphism II (S-seminear-ring homomorphism II) if $\phi : A \rightarrow A_1$ is a semiring homomorphism where $A \subset N$ and $A_1 \subset N_1$ and A and $A_1$ are semirings.*

***Example 5.4.8:*** Let $N = Z^o \times Z_{12} \times Z_3$ be the S-mixed direct product, where $Z^o$ is a semiring, $Z_{12}$ a seminear-ring and $Z_3$ a near-ring. Clearly N is a S-seminear-ring II and as well as S-seminear-ring I. Thus we see a seminear-ring can simultaneously be a S-seminear-ring I as well as S-seminear-ring II.

***Example 5.4.9:*** Let $N = L \times Z_{18}$, S-mixed direct product of a semiring and seminear-ring where L is a distributive lattice given by the following diagram.

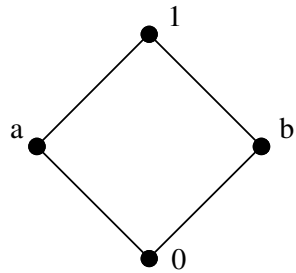

**Figure 5.4.1**

Clearly N is a seminear-ring which is a S-seminear-ring II. Thus we have also seen examples of S-seminear-rings II of finite order. Now we proceed on to define the concept of Smarandache pseudo seminear-ring.

**DEFINITION [99]:** *N is said to be a Smarandache pseudo seminear-ring (S-pseudo seminear-ring) if N is a near-ring and has a proper subset A of N which is a seminear-ring under the operations of N.*

We cannot obtain several examples of S-pseudo near-rings by using the S-mixed direct product.

***Example 5.4.10:*** Let (Z, +, .) be a near-ring. Clearly (Z, +, .) is a S-pseudo seminear-ring as $(Z^o, +, .)$ is a seminear-ring.

***Example 5.4.11:*** Let $M_{n \times n} = \{(a_{ij}) / a_{ij} \in Z, (Z, +, .)$ a near-ring$\}$. $M_{n \times n}$ is a S-pseudo seminear-ring as $M_{n \times n} = \{(a_{ij})/ a_{ij} \in Z^o = Z^+ \cup \{0\}\}$ is a seminear-ring.

**DEFINITION 5.4.13:** *Let N and $N_1$ be the S-pseudo seminear-ring. A mapping $\phi : N \rightarrow N_1$ is said to be a Smarandache pseudo seminear-ring homomorphism (S-pseudo seminear-ring homomorphism) if $\phi$ restricted from A to $A_1$ is a seminear-ring homomorphism where A and $A_1$ are proper subsets of N and $N_1$ which are seminear-rings. Thus $\phi$ need not even be defined on the whole of N.*

We define the concept of S-pseudo subseminear-ring and Smarandache pseudo ideals of a near-ring.

**DEFINITION 5.4.14**: *Let N be a near-ring if N has a proper subset A which is a subnear-ring A and if A itself a S-pseudo seminear-ring then we say A is a Smarandache pseudo subseminear-ring (S-pseudo subseminear-ring).*

**DEFINITION 5.4.15**: *Let N be a near-ring. A proper subset M of N is said to be a Smarandache pseudo ideal (S-pseudo ideal) if M is a S-ideal of the near-ring N.*

**THEOREM 5.4.8:** *Let N be a near-ring. If N has a S-pseudo subseminear-ring then N is a S-pseudo seminear-ring.*

*Proof*: Straightforward by the very definitions.

It is left for the reader to construct a S-pseudo seminear-rings which has no S-pseudo subseminear-rings.

**DEFINITION 5.4.16**: *Let N be a S-pseudo seminear-ring, if N has no proper S-pseudo subseminear-ring then we say N is a Smarandache pseudo simple seminear-ring (S-pseudo simple seminear-ring).*

The reader is assigned the work of finding an example of a S-pseudo simple seminear-ring.

<u>**PROBLEMS:**</u>

1. Is $Z_8S(3)$ a S-seminear-ring?
2. Prove or disprove the semigroup seminear-ring $Z_{12}S(4)$ where $(Z_{12}, \times, '.')$ is a seminear-ring is S-cyclic seminear-ring.
3. Is $Z_{12}S(4)$ S-weakly commutative seminear-ring?
4. Can $Z_{12}S(4)$ be at least a S-weakly cyclic seminear-ring?
5. Does $Z_{12}S(4)$ given in problem 2 have S-conjugate subgroup?
6. Is $Z_{12}S(4)$ given in problem 5 be S-conjugate seminear-ring?
7. Find an example of a S-simple pseudo seminear-ring.
8. Does $(Z, +, .)$ have S-pseudo subseminear-rings?
9. Is $(Z, +, .)$ a S-simple pseudo seminear-ring?
10. Can $(Z, +, .)$ have S-pseudo ideals? Justify your answer.

## 5.5 Some special properties in S-near -rings

In this section we introduce several interesting Smarandache properties in S-near-ring like Smarandache planar, Smarandache integral, Smarandache equiprime etc. we do not intend to give results or routine theorems for the reader but on the contrary expect the reader to work in these directions and develop effective theories and results on Smarandache notions.

**DEFINITION 5.5.1:** *Let N be a near-ring, A a non-empty subset of a S-N-subgroup V. The S-N-subgroup H of V, Smarandache centralizes (S-centralizes) a if $(h + u)a = ha + ua$ for all $h \in H$   $u \in A$ and $a \in N$.*

**THEOREM 5.5.1**: *If N is a near-ring and if a S-N-subgroup H of V, S-centralizes A then it centralizes A.*

*Proof*: Direct by the very definition. It is left as an exercise for the reader to give an example of a N-subgroup which centralizes but does not S-centralize.

**DEFINITION 5.5.2**: *A zero symmetric S-near-ring N is called Smarandache reliable (S-reliable) if whenever $\theta : I \rightarrow N$ is an epimorphism and I is a S-ideal of a S-near-ring A then for (x, y $\in$ I and a $\in$ A) x – y $\in$ Ker $\theta$ implies ax – ay $\in$ Ker $\theta$.*

Construct an example of a S-reliable $\theta$ symmetric S-near-ring N.

**DEFINITION 5.5.3**: *Let N be a S-near-ring we define a Smarandache equivalence relation (S-equivalence relation) on N by a $\equiv$ b if and only if na = nb for all n $\in$ P $\subset$ N where P is a semifield in N. The near-ring N is called Smarandache planar (S-planar) if there is more than three S-equivalence classes and every equation of the form xa = xb + c has unique solution x when a $\neq$ b.*

This S-planar near-ring will be utilized in Balanced Incomplete Block Designs (BIBD) of high efficiency and it has application in error correcting codes.

**DEFINITION 5.5.4**: *A Smarandache I-E S-semigroup (S-IE-S-semigroup) A is a S-semigroup which satisfies I (G) = E (G) (here G $\subset$ A and G is a group of the S-semigroup A) where I(G) is the S-near-ring generated by all inner automorphisms of G and E (G) is the S-near-ring generated by all the endomorphisms of G.*

**THEOREM 5.5.2**: *If A is a S-semigroup which is S-I-E semigroup then A has a proper subset which is a I-E group and conversely if A is a semigroup which has a subgroup of G which is I-E group then A is a S-I-E semigroup.*

*Proof*: Straightforward, hence left for the reader to prove. Several interesting properties of these S-I-E semigroups can be derived.

**DEFINITION 5.5.5**: *Let (G, +) be a S-semigroup written additively not necessarily abelian, with identity element denoted by 0. Let T $\subset$ A \ {0} where (A, +) is the subgroup of the S-semigroup (G, +). Define $\overline{\Gamma}$ by a $\overline{\Gamma}$ b = a if b $\in$ T; a $\overline{\Gamma}$ b = 0 if b $\notin$ T then (G, +, $\overline{\Gamma}$ ) is a Smarandache right zero symmetric near-ring (S-right zero symmetric near-ring).*

**DEFINITION 5.5.6**: *Let N be a S-near-ring and let $N_0$ be a subset of N such that $N_0$ generates B $\subset$ N where B is the near-field of N and $BN_0 \subset N_0$. Suppose $A_0$ is a subset of $N_0$ such that $BA_0 \subset A_0$ and $A_0 N_0 \subset A_0$, then $A_0$ is called the Smarandache $\sigma$ subset (S-$\sigma$-subset) of N and the S-subnear-ring generated by $A_0$ is called a Smarandache $\sigma$ subnear-ring (S-$\sigma$-subnear-ring) of N.*

*A S-$\sigma$- subnear-ring E of N is called essential in N if 0 $\neq$ A is an ideal of N implies E $\cap$ A $\neq$ 0.*

Several interesting results can be developed in this direction using these two Smarandache concepts. It is left for the researchers to do research in this direction and obtain some interesting results about them.

**DEFINITION 5.5.7:** *A subset S of a N-near-ring is called Smarandache integral (S-integral) if S has no non-zero divisors of zero. If H is an S-integral subset of N with the further property that $B^2 \subset H$ where $B \subset N$ is a near-field then N is said to be an Smarandache H integral near-ring (S-H-integral near-ring).*

**DEFINITION 5.5.8:** *Let N be a S-right near-ring and G a S-N-semigroup; G is called Smarandache equiprime (S-equiprime) if $PG \neq 0$, $P0 = 0$ (where $P \subset N$ and P is a near-field) and if for each $a \in P$ and $g$, $g^1 \in G$ with $Ga \neq 0$, $ang = ang^1$ for all $n \in P$ implies $g = g^1$.*

Thus the reader is to find some interesting relations between the S-equiprime ideals and S-equiprime N-semigroups.

**DEFINITION 5.5.9:** *Let N be a S-right near-ring and A an S-ideal or a S-left ideal of N. We define 3 Smarandache properties in the following :*

    i.    *A is S-equiprime if for any a, x, $y \in P \subset N$ such that $nx - any \in A$ for all $n \in P \subset N$ we have $a \in P$ or $x - y \in A$. ($P \subset N$ is a near-field).*

    ii.    *A is Smarandache strongly semiprime (S-strongly semiprime) if for each $a \in N \setminus A$ there exists a finite subset F of N such that if $x, y \in F \subset N$ and $afx - afy \in A$ for all $f \in F$ then $x - y \in A$.*

    iii.    *A is Smarandache completely equiprime (S-completely equiprime) if $a \in P \setminus A$ and $ax - ay \in A$ imply $x - y \in A$.*

**DEFINITION 5.5.10:** *A S-right near-ring N is a Smarandache left infra near-ring (S-left infra near-ring) if $(N, +, .)$ is a triple where N is a non-empty set and '+' and '.' are binary operations satisfying.*

    1.    *(N, +) is a S-semigroup.*
    2.    *(N, .) is a S-semigroup.*
    3.    *For all x, y, z in N $x . (y + z) = x . y - x . 0 + x . z$.*

**DEFINITION 5.5.11:** *A S-near-ring M which satisfies the property, A an S-ideal of B, B an S-ideal of N and $B / A \cong M$ implies, A is an S-ideal of N for all S-near-rings N, and S-subnear-rings A and B of N are called Smarandache F-near-ring (S-F-near-ring).*

**DEFINITION 5.5.12:** *A near-ring N is Smarandache left weakly regular (S-left weakly regular) if every $x \in N$ is of the form $x = ux$ for some $u_1$ in the S-principal ideal generated by x and is S-left weakly regular if every $a \in N$ is in the product of $\{Na\} * \{Na\} = \{$finite sums $\Sigma x_k y_k / x_k, y_k \in Na\}$.*

**DEFINITION 5.5.13:** *A Smarandache composition near-ring (S-composition near-ring) is a quadruple $\{C, +, o, .\}$ where $(C, +, o)$ and $(C, +, .)$ are S-near-rings such that $(a . b) o c = (a o b) . c$ for all $a, b, c \in C$.*

**DEFINITION 5.5.14:** *A non-zero S-ideal H of G is said to be Smarandache uniform (S-uniform) if for each pair of S-ideals $K_1$ and $K_2$ of G such that $K_1 \cap K_2 = (0)$, $K_1 \subset H$, $K_2 \subset H$ implies $K_1 = (0)$ or $K_2 = (0)$.*

**DEFINITION 5.5.15:** *An S-ideal H of G is said to have Smarandache finite Goldie dimension written as (S-FGD) if H does not contain an infinite number of non-zero S-ideals of G whose sum is direct.*

**DEFINITION 5.5.16:** *Let G be a S-N subgroup with S-FGD then every non-zero ideal of G contains a S-uniform ideal.*

**DEFINITION 5.5.17:** *Let N be a S-near-ring an element $a \in P$ is called Smarandache normal element (S-normal element) of N if aN = Na. If aN = Na for every $a \in P$ then N is called a Smarandache normal near-ring (S-normal near-ring).*

*Let Sn(N) denote the set of all S-normal elements of N. N is called a S-normal ring if and only if Sn (N) = P where $P \subset N$ and P a near-field.*

**DEFINITION 5.5.18:** *Let N be a S-near-ring. A S-subnear-ring S of N is called S-normal subnear-ring with respect to a non-empty subset T of N if ts = st for every $t \in T$. In particular we have nS = Sn for every $n \in P \subset N$ we say S is a Smarandache normal subnear-ring (S-normal subnear-ring) of N.*

**DEFINITION 5.5.19:** *Let N be a S-near-ring ; a S-subsemigroup S of N is called a Smarandache normal subsemigroup (S-normal subsemigroup) of S if nS = Sn for all n $\in P \subset N$.*

**DEFINITION 5.5.20:** *Let N be a S-near-ring. A S-subnear-ring M of N is Smarandache invariant (S-invariant) with respect to an element $a \in P \setminus M$ if aM $\subset$ M and Ma $\subset$ M. M is said to Smarandache strictly invariant (S-strictly invariant) with respect to a if aM = M and Ma = M.*

**DEFINITION 5.5.21:** *A commutative S-near-ring N with identity is called a Smarandache Marot near-ring (S-Marot near-ring) if each regular S-ideal of N is generated by regular elements where by regular S-ideal we mean an ideal with regular elements and the non-zero divisors of the near-ring.*

**DEFINITION 5.5.22:** *Let K and I be S-ideals of the S-N-subgroup G. K is said to be the Smarandache complement (S-complement) of I if the following two conditions hold good.*

1. *$K \cap I = (0)$ and*
2. *$K_1$ is an S-ideal of G such that $K \subset K_1$ then $K_1 \neq K$ imply $K_1 \cap I \neq (0)$.*

**DEFINITION 5.5.23:**

i.    *A subset S of the S-N-subsemigroup G is said to be Smarandache small (S-small) in G if S + K = G where K is an S-ideal of G imply K = G.*

ii.    *G is said to be Smarandache hollow (S-hollow) if every proper S-ideal of G is S-small in G and*

iii.   *G is said to have Smarandache finite spanning dimension (S-finite spanning dimension) if for any decreasing sequence of Smarandache NS-subsemigroups (S-NS-subsemigroups) $X_0 \supset X_1 \supset X_2 \supset \ldots$ of G such that $X_i$ is an S-ideal of $X_{i-1}$ there exists an integer such that $X_j$ is S-small in G for all $j \geq k$.*

**DEFINITION 5.5.24**: *Suppose H and K be any two S-N-subsemigroup of G. Then K is said to be Smarandache supplement (S-supplement) for H if H + K = G and H + $K_1 \neq$ G for any proper S-ideal $K_1$ of K.*

Now we have just defined several of the Smarandache concepts in near-ring so that an interested researcher will develop them.

**DEFINITION 5.5.25**: *Let N be a near-ring we say N is Smarandache ideally strong (S-ideally strong) if for every S-subnear-ring of R is a S-ideal of R.*

**DEFINITION 5.5.26**: *Let N be near-ring $\{I_k\}$ be the collection of all S-ideal of R. R is said to be a Smarandache I\* near-ring (S-$I^*$ near-ring) if for every pair of ideals. $I_1$, $I_2 \in \{I_k\}$ we have for every $x \in N \setminus \{I_1 \cup I_2\}$ the S-ideal generated by x and $I_1$ and x and $I_2$ are equal i.e. $\langle x \cup I_1 \rangle = \langle x \cup I_2 \rangle$ we can have several examples of these.*

## PROBLEMS :

1.    Let $Z_{28}$ be a S-near-ring. Can we have a S-equivalence relation on $Z_{28}$ ?
2.    Can $Z_{28}$ be S-planar?
3.    Is $Z_{28}$ planar?
4.    Give an example of a S-σ-near-ring.
5.    Give an example of a σ near-ring which is not a S-σ-near-ring.
6.    Give an example of a S-composition near-ring.
7.    Can Z be made into a Smarandache composition near-ring?
8.    Give an example of a S-normal near-ring.
9.    Does their exist an example of a normal near-ring which is not a S-normal near-ring?
10.   Give an example of a S-left weakly regular near-ring.
11.   Give an example of S-ideal H of G which has S-FGD.
12.   Does $Z_{29}$ the near-ring have S-normal element?
13.   Give an example of a S-H-integral near-ring.
14.   Can N = $Z_{15} \times Z$ be a S-H-integral near-ring?
15.   Give an example of a S-equiprime N-subsemigroup.
16.   Can N = $Z_{12} \times Z_{19}$ be a S-infra near-ring?
17.   Will N = $Z_{12} \times Z_{19}$ be a S-F near-ring?
18.   Give an example of a S-F near-ring.
19.   Find two S-subsemigroup K and H of which K is a S-supplement for H.

**Chapter Six**

# SMARANDACHE SEMINEAR-RINGS

This chapter specially serves the purpose of introducing Smarandache notions on seminear-rings. We define two levels of Smarandache seminear-ring and study them. This chapter has five sections, in section one we define Smarandache seminear-rings of type I and study some of its properties. For these S-seminear-rings we define and study homomorphism and its ideals in section two. Section three is devoted to the introduction and study of Smarandache seminear-rings of level II. In section four we define Smarandache pseudo seminear-rings and study several of its Smarandache properties.

In the final section we introduce two special classes of Smarandache seminear-rings using groups and semigroups. These definitions help us to find more and more new classes of S-seminear-rings. At the end of each section some problems are given for a reader to solve.

## 6.1 Definition and properties of Smarandache seminear-rings.

In this section we define a new concept called Smarandache seminear-rings (S-seminear-rings) and define some interesting properties in them like Smarandache subseminear-ring, Smarandache simple seminear-ring and Smarandache commutative seminear-rings. We prove if a seminear-ring has a Smarandache subseminear ring then the seminear-ring is a S-seminear-ring. But even if N is S-seminear-ring its subseminear-ring in general need not be a S-subseminear-ring.

**DEFINITION 6.1.1**: *A set S with two binary operations '+' and '.' is called a Smarandache seminear-ring (S-seminear-ring) if (S, +) and (S, .) are Smarandache semigroups and for all $s, s^I, s^{II} \in S, (s + s^I) s^{II} = ss^{II} + s^I s^{II}$.*

**Example 6.1.1**: Let $Z^+$ be the set of positive integers. $(Z^+, +)$ is a semigroup. Define '.' on $Z^+$ by $a.b = a$ for all $a, b \in Z^+$. $(Z^+, '+', '.')$ is a seminear-ring which is not a S-seminear-ring. Thus this example proves all seminear-rings are not in general S-seminear-ring.

Now we will see an example of a S-seminear-ring.

**Example 6.1.2**: Let $N = Z_2 \times (Z^+ \cup \{0\})$ clearly N is a S-seminear-ring. Here $Z^+ \cup \{0\}$ is a seminear ring with usual '+' and '.' defined by $a.b = a$ for all $a, b \in Z^+$. $P = Z_2 \times \{0\}$ is a subgroup of N under '+' so N is a S-semigroup; $\{1\} \times \{0\}$ is a subgroup of N under '.'. Hence N is a S-seminear-ring.

**DEFINITION 6.1.2**: *Let (N, +, .) be a seminear-ring. A proper subset P of N is said to be a Smarandache subseminear-ring (S-subseminear-ring) if (P, +, .) is a S-seminear-ring.*

**THEOREM 6.1.1**: *Let (N, +, .) be a seminear-ring. If N has a S-subseminear-ring then N is a S-seminear-ring.*

*Proof*: Follows from the fact that if N has P to be a S-subseminear-ring then (P, +) is a S-semigroup and (P, .) is a S-semigroup. Since P ⊂ N, (N, +) is a S-semigroup and (N, .) is a S-semigroup.

Thus we see in the definition of S-subseminear-ring we need not put the condition N is a S-seminear-ring. Still we may have examples in which a S-seminear-ring may not have proper S-subseminear-rings.

In view of this we define the following:

**DEFINITION 6.1.3**: *Let (N, +, .) be a S-seminear-ring. If N has no proper S-subseminear-rings we say N is a Smarandache simple seminear-ring (S-simple seminear-ring).*

**DEFINITION 6.1.4**: *Let (N, +, .) be a seminear-ring. We say (N, +, .) is a Smarandache commutative seminear-ring (S-commutative seminear-ring) if every S-subseminear-ring P, P ⊂ N is a commutative seminear-ring.*

**THEOREM 6.1.2**: *Let N be a commutative seminear-ring which has S-subseminear-rings then N is a S-commutative seminear-ring.*

*Proof*: Obvious.

It is left as an exercise to find S-commutative seminear-ring which is not a commutative seminear-ring.

Now we proceed on to define Smarandache ideals in seminear-rings.

**DEFINITION 6.1.5**: *Let (N, +, .) be a S-seminear-ring. If N has atleast one S-subseminear-ring which is commutative then we say (N, +, .) is a Smarandache weakly commutative seminear-ring (S-weakly commutative seminear ring).*

**DEFINITION 6.1.6**: *Let (N, +, .) be a seminear-ring. We say N has a Smarandache hyper subseminear-ring (S-hyper subseminear-ring) S,*

1. *If (S, +) is a S-semigroup.*
2. *If A is a proper subset S which is a subsemigroup of S and A contains the largest group of (S, +).*
3. *(A, .) is a S-subsemigroup or a group.*

**THEOREM 6.1.3**: *Let N be a seminear-ring. Every S-hyper subseminear-ring is a S-subseminear-ring, but every S-subseminear-ring in general is not a S-hyper subseminear-ring.*

*Proof*: One way is straightforward by the definition. To prove the other part the reader is requested to construct examples.

Now we define Smarandache dual hyper subseminear-ring in the following.

**DEFINITION 6.1.7**: *Let (N, +, .) be a seminear-ring. We say a proper subset S of N is a Smarandache dual hyper subseminear-ring (S-dual hyper subseminear-ring) if*

1. *(S, .) is a S-semigroup.*
2. *If A is a proper subset of S which is a subsemigroup of S and A contains the largest group of (S, .).*
3. *(A, +) is a S-subsemigroup or a group.*

The reader is advised to obtain examples and some more properties about these concepts.

Now we proceed on to define Smarandache normal subseminear-ring.

**DEFINITION 6.1.8**: *Let (N, +, .) be a seminear-ring. A proper subset A of N which is a group under the operations of both '+' and '.' ({A \ {0}} is a group under '.') is called a Smarandache normal subseminear-ring (S-normal subseminear-ring) if $x A \subseteq A$ and $Ax \subset A$ or $xA = \{0\}$ and $Ax = \{0\}$ for all $x \in S$.*

**DEFINITION 6.1.9**: *If a seminear-ring (N, +, .) has no S-normal subseminear-ring then we call the seminear ring Smarandache pseudo simple (S-pseudo simple).*

**DEFINITION 6.1.10**: *Let (N, +, .) be a seminear-ring. If A is a S-normal subseminear-ring of N then we define the Smarandache quotient seminear-ring (S-quotient seminear-ring) of the seminear-ring of N by N / A = {Ax / x ∈ N}.*

*It is not known under what condition N / A will be a near-ring or a S-near-ring or a seminear-ring or a S-seminear-ring.*

It is left as an open problem for the reader to solve.

**DEFINITION 6.1.11**: *Let (N, +, .) be a S-seminear-ring, we say P is a Smarandache maximal S-subseminear-ring (S-maximal S-subseminear-ring) if there is a S-subseminear-ring M such that P ⊂ M then P = M is the only possibility.*

**DEFINITION 6.1.12**: *Let (N, +, .) be a seminear-ring. If N has only one maximal subseminear-ring we call N a Smarandache maximal seminear-ring (S-maximal seminear-ring).*

**THEOREM 6.1.4**: *If (N, +, .) is a S-maximal seminear-ring then N is a S-seminear-ring.*

*Proof*: Follows from the very definition.

<u>**PROBLEMS**</u>:

1.   Can a S-seminear-ring of order 5 exist?
2.   What is the order of the smallest S-seminear-ring?
3.   Is N = $Z_2 \times Z_7$ where $Z_2$ and $Z_7$ are near-rings, a S-seminear-ring?

4. Will N = $Z_2 \times Z_2$ have S-subseminear-rings?
5. Find a S-subseminear of N = $Z_2 \times Z_7$ given in problem 3.
6. Give an example of a S-commutative seminear-ring which is not a commutative seminear-ring.
7. Give an example of a S-weakly commutative seminear-ring which is not a S-commutative seminear-ring.
8. Let N = $Z \times Z_2$ be a S-near-ring.
   a. Does N have S-hyper subseminear-ring?
   b. Can N have S-dual hyper subseminear-ring?
9. Can N = $Z \times Z_2 \times Z_3$ have S-normal subseminear-ring?
10. Give an example of a S-pseudo simple seminear-ring.
11. Is N = $Z_9 \times Z^+$ a S-seminear-ring? Justify your claim.
12. Give an example of a S-maximal seminear-ring.

## 6.2 Homomorphism and ideals of a S-seminear-ring

In this section we define the notion of S-ideals and S-seminear-rings and obtain some Smarandache notions by way of definitions and properties. These S-seminear rings are got when the semigroups under '+' and '.' are replaced by S-semigroups.

**DEFINITION 6.2.1**: *Let N and $N_1$ be two S-seminear-rings. A map $\theta$ from N to $N_1$ is said to be Smarandache seminear-ring homomorphism (S-seminear-ring homomorphism) if*

1. *$\theta(x + y) = \theta(x) + \theta(y)$ for all x, y $\in$ A, A an additive subgroup of (N, +) and $\theta(x)$, $\theta(y) \in A_1$, $A_1$ an additive subgroup of (N, +).*
2. *$\theta(xy) = \theta(x)\,\theta(y)$ for all x, y $\in$ B, B a multiplicative subgroup (N, .) and $\theta(x)$, $\theta(y) \in B_1$,   $B_1$ a subgroup of $(N_1, .)$.*

**THEOREM 6.2.1**: *Let (N, +, .) be a S-seminear-ring. A proper subset P of N is a near-ring if and only if*

1. *(P, +) is a subgroup of (N, +).*
2. *(P, .) is a subsemigroup of (N, .).*

*Proof*: Left for the reader to prove.

**Example 6.2.1**: Let $Z_6 = \{0, 1, 2, 3, 4, 5\}$. Define '×' as usual multiplication modulo 6 'o' as a o b = a for all a, b $\in$ $Z_6$. Clearly $(Z_6, \times)$ is a S-semigroup for A = $\{1, 5\}$ is a subgroup of $Z_6$. $(Z_6, .)$ is a S-semigroup for every singleton {x} is a group as x o x = x.

**DEFINITION 6.2.2**: *Let (N, +, .) be a S-seminear-ring. A proper subset P of N is said to be a Smarandache left ideal (S-left ideal) of N if*

1. *(P, +) is a S-subsemigroup of (N, +).*
2. *$n(n_1 + i) + n_r\, n^I \in P$ for each i $\in$ P and n, $n_1 \in A \subset N$ where (A, +) is a subgroup.*

**DEFINITION 6.2.3**: *Let P be a nonempty subset of a S-seminear-ring N. P is called an Smarandache ideal (S-ideal) of N if*

*1.  P is a S-left- ideal.*
*2.  PA ⊂ P (for A ⊂ N, A is a subgroup of (N, +)).*

**DEFINITION 6.2.4**: *Let N be a S-seminear-ring. N is said to be Smarandache left bipotent (S-left bipotent) if $Pa = Pa^2$ for all a ∈ N (P ⊂ N is subgroup of (N, +)).*

**DEFINITION 6.2.5**: *Let N be a S-seminear-ring N is said to be Smarandache dually left bipotent (S-dually left bipotent) if $B + a = B + a^2$ for all a ∈ N (B ⊂ N is a subgroup of (N, .)).*

**DEFINITION 6.2.6**: *Let N be a S-seminear-ring. We say N is Smarandache strongly bipotent (S-strongly bipotent) if*

*1.  $Pa = Pa^2$.*
*2.  $P + a = P + a^2$.*

*where (P, +) is subgroup of (N, +) and (P, .) is a subgroup of (N, .).*

From the above definitions we have the following theorem which is left for the reader to prove as an exercise.

**THEOREM 6.2.2**: *Let N be a S-seminear-ring which is S-strongly bipotent then N is S-left bipotent and S-dually left bipotent.*

**DEFINITION 6.2.7**: *A S-seminear-ring N is said to be a Smarandache s-seminear-ring (S-s-seminear-ring) if a ∈ Pa for each a ∈ N and P a subgroup of (N, +).*

**DEFINITION 6.2.8**: *Let N be a S-seminear-ring. An additive subgroup A of N is called the Smarandache N-subgroup (S-N-subgroup) of N related to P if PA ⊂ A and AP ⊂ A where (P, .) is a subgroup of (N, .).*

$$PA = \{pa \ / \ p \in P \ and \ a \in A\}.$$

*Thus by the very definition $S_1N$ subgroups may or may not exist for all S-seminear-ring.*

**Example 6.2.2**: Let $Z_{12}$ = {0, 1, 2, 3, …, 11} be a S-seminear-ring under the operation '×' and '.'. Take A = {1, 5, 7, 11} is a subgroup of ($Z_{12}$, ×). Take P = {1} a subgroup of ($Z_{12}$, .); A is a $S_1N$-subgroup of $Z_{12}$.

If we take $P_1$ = {9} a subgroup of ($Z_{12}$, .) A is not $S_1N$-subgroup. From this we see the term "related to P" is important for a subgroup may fail to be $S_1N$ subgroup if we change P. This is evident from the above example. So a natural question of interest will be does there exist an additive subgroup in S-seminear-ring N such that A is a $S_1N$ subgroup related to all subgroups of (N, .). If such seminear-rings exists characterize them. This is given as a suggested problem in chapter X.

On similar dual lines one can define $S_1N$ subgroups replacing the operation '+' by '.' and '.' and '+' and study the related questions. It is left for the reader to define Smarandache dual N subgroups, ($S_1DN$ subgroups) and illustrate them by examples.

<u>**PROBLEMS:**</u>

1. Give an example of a S-seminear-ring N in which (P, +) is a subgroup but (P, .) is not even a semigroup, where P is a proper subset of N.
2. Prove ($Z_{24}$, ×, .) is a S-seminear-ring. Show (P, +) and (P, .) are group and S-semigroup respectively. Show (P, .) a group exists in $Z_{24}$.
3. Find S-left ideals in $Z_{24}$ given in example 2.
4. Is (Z, ×, .) a S-seminear-ring? Justify. Will (Z, ×, . ) be a seminear-ring?
5. Prove (Q, ×, .) is a S-seminear-ring. Is (Q, ×, .), S-strongly bipotent or just S-left bipotent?

## 6.3 Smarandache seminear-rings of level II

In section 6.1 we saw S-seminear-rings N defined as (N, +) and (N, .) to be S-semigroups. We defined various Smarandache properties associated with them. Now we define in this section S-seminear-rings more in a conventional way, these S-seminear-rings will be distinguished from type I by putting II. By default of notions we do not denote S-seminear-ring I by just S-seminear-ring. We obtain several interesting results about them.

**DEFINITION [99]**: *A nonempty set N is said to be a Smarandache seminear-ring II (S-seminear-ring II) if (N, +, .) is a seminear-ring having a proper subset A (A ⊂ N) such that A under the same binary operations of N is a near-ring that is (A, +, .) is a near-ring.*

***Example 6.3.1***: Let $Z_{18}$ = {0, 1, 2 ,…, 17} be the set of integers modulo 18. $Z_{18}$ under usual multiplication '×' modulo 18 is a semigroup. Define '.' an operation on $Z_{18}$ as a.b = a for all a, b ∈ $Z_{18}$. Clearly ($Z_{18}$, ×, .) is a seminear-ring. ($Z_{18}$, ×, .) is a S-seminear-ring for take A = {1, 5, 7, 11, 13, 17}; (A, ×, .) is a near-ring hence the claim.

**THEOREM 6.3.1**: *All seminear-rings need not in general be S-seminear-rings II.*

*Proof*: By an example.

Let $Z^+$ = {set of positive integers}. $Z^+$ under '+' is a semigroup. Define 'o' a binary operation on $Z^+$ as a o b = a for all a, b ∈ $Z^+$. Clearly $Z^+$ is a seminear-ring which is not a S-seminear-ring.

**THEOREM 6.3.2**: *Let N be a S-seminear-ring II then N is a S-seminear-ring I provided (N, .) is a S-semigroup.*

*Proof*: Obvious by the very definitions as (N, +) is a S-semigroup and (N, .) is a S-semigroup.

**Example 6.3.2**: Let $M_{n \times n} = \{(a_{ij}) \,/\, a_{ij} \in Q\}$. Define matrix multiplicative as an operation on $M_{n \times n}$. $(M_{n \times n}, \times)$ is a semigroup. Define 'o' on $M_{n \times n}$ as A o B = A for all A , B $\in$ $M_{n \times n}$, $(M_{n \times n}, \times, .)$ is S-seminear-ring. For the set of all n × n matrix with |A| ≠ 0 is a group if we denote the collection by $A_{n \times n}$ then $A_{n \times n} \subset M_{n \times n}$. Clearly $(A_{n \times n}, \times, .)$ is a near-ring.

**Example 6.3.3**: Let $Z_{24} = \{0, 1, 2, …, 23\}$ be the set of integers modulo 24. Define usual multiplication '×' on $Z_{24}$. $(Z_{24}, \times)$ is a semigroup. Define 'o' on $Z_{24}$ as a o b = a for all a $\in$ $Z_{24}$. $Z_{24}$ is a S-seminear-ring as A = {1, 5, 7, 11, 13, 17, 19, 23} is a group under '×' as $(A_1, \times, .)$ is a near-ring. Thus $Z_{24}$ is a S-seminear-ring. In view of these examples we propose problems in chapter 10.

**Example 6.3.4:** Let $Z_4 = \{0, 1, 2, 3\}$, clearly $(Z_4, \text{'×'}, \text{'.'})$ is a S-seminear-ring as {A = {1, 3}, ×, . } is a near-ring.

**Example 6.3.5:** $Z_9 = \{0, 1, 2, …, 8\}$ is a S-seminear-ring as {A = {1, 8}, ×, .} is a near-ring.

From this example we see 9 is not a prime, still in case $n = p^{\alpha}$; where p is a prime, we see $Z_n$ is a S-seminear-ring.

**Example 6.3.6:** Let $Z_{25} = \{0, 1, …, 24\}$. $Z_{25}$ is a S-seminear-ring as (A = {1, 24}, ., ×) is a near-ring.

In view of this we have the following theorem.

**THEOREM 6.3.3**: *Let ( $Z_{p^2}$ , ×, .) be a seminear-ring. ( $Z_{p^2}$, ×, .) is a S-seminear-ring.*

*Proof*: To prove we have to show a near-ring as a substructure of ( $Z_{p^2}$ , ×, .). Take A = ({1, $p^2 - 1$}, ×, .) is a near-ring. Hence the claim.

**THEOREM 6.3.4**: *Let $Z_{p^n}$ = {0, 1, …, $p^n$-1} be a seminear-ring under '×' and '.'.*
*Clearly $Z_{p^n}$ is a S-seminear-ring for take {A = ({1, $p^n$-1}, ×, .)}; A is near-ring.*

*Proof:* Simple number theoretic arguments will yield the result.

Thus we have a natural class of S-seminear-rings. In fact this natural class of S-seminear-rings is also trivially S-seminear-rings of level I by taking singletons x in $Z_n$ where $x^2 = x$.

Now it may happen that many near-rings may contain subsets which are seminear-rings we choose to call these S-pseudo seminear-rings which will be dealt in the next section.

1.  Find a seminear-ring other than $Z^+$ which is not a S-seminear-ring.
2.  Prove $Z_{19} = \{0, 1, 2, \ldots, 18\}$ is S-seminear-ring.
3.  Does there exist an example of a S-seminear-ring I which is never a S-seminear-ring II? Justify your answer.
4.  Illustrate by an example that there exists a S-seminear-ring II which is not a S-seminear-ring I.
5.  Let $N = Z_3 \times Z_{14}$ be a S-seminear-ring. How many proper subsets in N are near-rings?
6.  Let $Z_{28} = \{0, 1, 2, \ldots, 27\}$ be a seminear-ring. Find the number of near-rings in $Z_{28}$.
7.  Is $Z_2 = \{0, 1\}$ a S-seminear-ring? Justify your answer.

## 6.4 Smarandache pseudo seminear-ring

In this section we just define the concept of S-pseudo seminear-rings and obtain some interesting results about them. This is an interesting case when a near-ring has the general structure as a proper subset. To find a natural class of such S-pseudo seminear-rings we define Smarandache mixed direct products in near-rings.

**DEFINITION 6.4.1**: *N is said to be Smarandache pseudo seminear-ring (S-pseudo seminear-ring) if N is a near-ring and has a proper subset A of N such that A is a seminear-ring under the operations of N.*

***Example 6.4.1***: Let Z be the set of integers under usual '+' and 'o' defined by a o b = a for all a, b $\in$ Z, (Z, '+', 'o') is a near-ring. Taking A = $Z^+$; ($Z^+$, +, .) is a seminear-ring So Z is a S-pseudo seminear-ring.

***Example 6.4.2***: Let Z[x] be the polynomial ring over the ring of integers. Define '+' on Z[x] as usual addition of polynomials. Define an operation 'o' on Z[x] as p(x) . q(x) = p(x) for all p(x), q(x) $\in$ Z[x]. Clearly (Z[x], +, .) is a S-pseudo seminear-ring for ($Z^+$ [x], +, .) is a seminear-ring.

**DEFINITION 6.4.2**: *Let N and $N_1$ be two S-pseudo seminear-rings. h : N $\rightarrow$ $N_1$ is a Smarandache pseudo seminear ring homomorphism (S-pseudo seminear-ring homomorphism) if h restricted from A to $A_1$ is a seminear-ring homomorphism.*

**DEFINITION 6.4.3**: *Let N be a S-seminear-ring II. An additive subgroup A of N is called a Smarandache N-subgroup II (S-N-subgroup II) if NA $\subset$ N (AN $\subset$ A) where NA = {na / n $\in$ N and a $\in$ A}.*

Thus we see in case of S-seminear-rings we have the concept of S-N-subgroup. We define S-left ideal for S-seminear-ring.

**DEFINITION 6.4.4:** *Let (N, +, .) be a S-seminear-ring. A proper subset I of N is called Smarandache left ideal (S-left ideal II) in N if*

    i.   *(I, +) is a normal subgroup of A $\subset$ N where (A, +) is a group.*

    ii.  *n ($n_1$ + i) + $n_r$ $n_1$ $\in$ I for each i $\in$ I, n, $n_1$ $\in$ A where $n_r$ denotes the unique inverse of n.*

**DEFINITION 6.4.5**: *A nonempty subset I of (N, +, .); N a S-seminear-ring is called a Smarandache ideal II (S-ideal II) in N if*

    1.  *I is a S-left ideal.*

    2.  *IA $\subset$ I (A $\subset$ N; A is a near-ring).*

**DEFINITION 6.4.6:** *Let N be a S-seminear-ring. N is said to be Smarandache left bipotent II (S-left bipotent II) if Na = $Na^2$ for every a $\in$ A, where A $\subset$ N is a near-ring.*

**DEFINITION 6.4.7:** *Let N be a S-seminear ring. N is said to be a Smarandache s-seminear-ring II (S-s-seminear-ring II) if a $\in$ Na for each a $\in$ A $\subset$ N; where A is a near ring.*

**DEFINITION 6.4.8**: *Let N be a S-seminear-ring, N is said to be a Smarandache regular II (S-regular II) if for each a in A $\subset$ N (A a near-ring) there exists x in A such that a = axa.*

The following result is left for the reader to prove.

**THEOREM 6.4.1:** *Let N be a S-seminear-ring, then N is S-regular II if and only if for each a ($\neq$ 0) in A, A $\subset$ N, there exists an idempotent e such that Aa = Ae (A a near-ring).*

**DEFINITION 6.4.9**: *A S-seminear-ring N is said to be Smarandache strictly duo (S-strictly duo) if every SN-subgroup (S-left ideal) is also a right SN-subgroup (S-right ideal).*

**DEFINITION 6.4.10**: *A S-seminear-ring N is called Smarandache irreducible II or Smarandache simple II (S-irreducible II or S-simple II) if it contains only trivial SN-subgroups (S-ideals) (0) and A; (A $\subset$ N; A the only near-ring which is irreducible).*

**DEFINITION 6.4.11**: *An S-ideal P ($\neq A_1$) is called a Smarandache strictly prime (S-strictly prime) if for any two SN-subgroups A and B of $A_1$ ($A_1 \subset$ N, $A_1$ a near-ring of the S-seminear-ring) AB $\subset$ P then, A $\subset$ P or B $\subset$ P.*

**DEFINITIONS 6.4.12**: *A S-left ideal B of a S-seminear-ring; N is called Smarandache strictly essential (S-strictly essential) if B $\cap$ K $\neq$ (0) for every non-zero SN-subgroup K of N.*

**DEFINITION 6.4.13**: *An element x in N (N a S-seminear-ring) is said to be Smarandache singular (S-singular) if there exists a non-zero S-strictly essential left ideal C in N such Cx = {0}.*

**DEFINITION 6.4.14**: *A nonzero S-seminear-ring N is said to be Smarandache subdirectly irreducible (S-subdirectly irreducible) if the intersection of all non-zero S-ideals II of N is non-zero.*

**Example 6.4.3**: Let $Z_{14} = \{0, 1, 2, \ldots, 13\}$ be a seminear-ring under the operation '×' and '.' where × is the usual multiplication modulo 14 and a o b = a for all a, b ∈ $Z_{14}$. A = {1, 3, 5, 9, 11, 13} has a SN-subgroup.

It is left to the reader to find interesting properties about S-seminear-rings.

**PROBLEMS;**

1. Prove $Z_{27}$ is a S-seminear-ring II.
2. Find S-ideals II in $Z_{28}$.
3. How many SN-subgroups II does the S-seminear-ring $Z_{30}$ have?
4. Find a S-seminear-ring homomorphism from $Z_{12}$ to $Z_{20}$.
5. Find all S-ideals II and SN-subgroups II of $Z_{210}$.

## 6.5 Miscellaneous properties of some new classes Smarandache seminear-ring.

In this section we recollect the various properties of seminear-ring and adopt them to S-seminear-rings. We study and introduce the two new classes of seminear-rings got using groups and semigroups. We study about units and zero divisors in these structure. We introduce the concept of semi-idempotents in case of S-seminear-rings. Just we recall the definition of group seminear-rings [84].

**DEFINITION 6.5.1**: *Let (N, +, .) be a seminear-ring such that (N, +) is a commutative semigroup and (N, .) is a semigroup, G be any group. The group G over the seminear-ring N; called the group seminear-ring, NG is a seminear-ring; NG consists of all formal sums of the form $\alpha = \Sigma \alpha(g) g$; $\alpha(g) \in N$ and $g \in G$ such that supp $\alpha = \{g / \alpha(g) \neq 0\}$ i.e. the support of $\alpha$ is finite satisfying the following operational rules.*

1. *$\Sigma \alpha(g) g = \Sigma \beta(g) g \Leftrightarrow \alpha(g) = \beta(g)$ for all $g \in G$.*
2. *$\Sigma \alpha(g) g + \Sigma \beta(g) g = \Sigma(\alpha(g) + \beta(g)) g$ for all $g \in G$.*
3. *$(\Sigma \alpha(g) g) (\Sigma \beta(g^1) g^1) = \Sigma \mu(g_1) g_1$ where $\mu(g) = \Sigma \alpha(g) \beta(g^1)$ with $gg^1 = g_1$ provided the coefficient are distributive from the left and right; otherwise we assume the product term is in NG.*
4. *$\Sigma \alpha(g) gn = \Sigma \alpha(g) ng \ \forall \alpha(g), n \in N$.*

*Since $1 \in G$ we see $N \subseteq NG$. G may or may not be contained in NG. For all $n \in N$ and $g \in G$ we have $ng = gn$. It can be verified NG is a seminear-ring.*

**THEOREM 6.5.1**: *Let $N = Z^+ \cup \{0\}$ be a seminear-ring and $G$ any group. The group seminear-ring NG has no nontrivial divisors of zero.*

*Proof*: Follows from the fact if $\alpha = \Sigma \, \alpha_i \, g_i$; $\alpha_i \in Z^+$ and $\beta = \Sigma \, \beta_j \, h_j$; $\beta_j \in Z^+$, $\alpha\beta \neq 0$ as in $Z^+ \cup \{0\}$; $0h = 0$ for all $h \in N$. Hence the claim.

Thus we see in case of group rings RG, every element of finite order in the group G will give way to a zero divisor. This is not in general true for group seminear-rings, evident from theorem 6.5.1.

Now we proceed on to define semigroup seminear-rings. If in the definition of group seminear-ring we replace the group G by a semigroup we get the semigroup seminear-ring.

**DEFINITION 6.5.2**: *Let $N$ be a seminear-ring and $S$ a semigroup with unit. The semigroup seminear-ring NS is defined analogous to group seminear-rings i.e. in the definition 6.5.1 replace the group G by a semigroup S with unit we get semigroup seminear-ring NS.*

It is easily proved NS is also a seminear-ring.

**THEOREM 6.5.2**: *Let $Z^+ \cup \{0\} = N$ be a seminear-ring and $S$ any semigroup with unit. $\beta\alpha = \alpha\beta = 1$ if and only if $\Sigma \, \alpha_i = 1$ and $\Sigma \, \beta_i = 1$ where $\alpha = \Sigma \, \alpha_i \, s_i$ and $\beta = \Sigma \, \beta_j \, s_j$ [only under product $\alpha\beta = \Sigma \, \alpha_i \, s_i \quad \Sigma \, \beta_j \, s_i = \Sigma \, \alpha_i$ for all $\alpha_i \, \beta_i \in NS$ (as in $Z^+$ we have $ab = a \; \forall a, b \in Z^+$)].*

*Proof*: Using the definitions the result is straightforward.

We now proceed on to define the concept of semi-idempotents in near-rings.

**DEFINITION 6.5.3**: *Let $N$ be a seminear- ring we say $\alpha \in N$ is a semi-idempotents of $N$ if $N(\alpha^2 - \alpha)$ does not contain $(\alpha^2 - \alpha)$.*

**DEFINITION 6.5.4:** *A seminear-ring $N$ is strongly subcommutative if every left ideal of it is a right quasi reflexive. (They define ideal in seminear-rings in a different way).*

**DEFINITION [2]:** *Let $R$ denote a seminear-ring. A non-empty subset $I$ of $R$ is an ideal if $x + y$, $rx$, $xr$ are in $I$ for all $x, y \in I$ and $r \in R$.*

In this book we call this ideal as a common ideal.

**DEFINITION [2]:** *A ideal $I$ is called strongly idempotent (SI-seminear-ring) if for every ideal $I$ is such that $\langle I^2 \rangle = I$.*

For more about these concepts refer [2].

We just introduce here the notion of chain seminear-ring.

**DEFINITION 6.5.5**: *A seminear-ring S is said to be a left chain seminear-ring if the set of all left ideals of S is totally ordered by inclusion.*

**DEFINITION 6.5.6**: *A seminear-ring N is said to be a chain seminear-ring if the set of ideals of N is totally ordered by inclusion.*

**DEFINITION 6.5.7**: *Let $P \subset S$ be an ideal of S, P is called prime if for all ideals I, J $\in$ S, IJ $\subset$ P implies I $\subset$ P or J $\subset$ P.*

**DEFINITION 6.5.8**: *An ideal $I \subseteq S$ is called semiprime if and only if for all ideals J $\subset$ S, $J^2 \subset$ I implies J $\subset$ I.*

**DEFINITION 6.5.9**: *Let S be a seminear-ring. The prime radical of S is the intersection of all prime ideals of S.*

As our main motivation is only the study of S-seminear-ring we have only just recalled the definitions and concepts about seminear-rings and by no means we deeply go into the subject as we assume the author to be well versed in the theory of near-rings. Hence we now proceed on to define Smarandache analogous of each and every concept mentioned here. We also by no means claim that we have exhausted all the properties about seminear-rings. We enlist here only those notions and definitions which was available to us.

Now first we study when do we call any group seminear-ring to be a Smarandache group seminear-ring.

**DEFINITION 6.5.10**: *A group seminear-ring NG of a group G over a near-ring NS is called a Smarandache group seminear-ring (S-group seminear-ring) if and only if N is a S-seminear-ring II.*

**THEOREM 6.5.3**: *Let NG be a S-group seminear-ring then N has a near-ring P such that P $\not\subset$ N but P $\subset$ NG.*

*Proof*: From the very construction and definition of S-group seminear-ring NG we have N to be S-seminear-ring so N contains a near-ring P $\subset$ N. Clearly PG the group near-ring which is a near-ring. Hence the claim.

*Example 6.5.1*: Let $Z_{24}$ be a seminear-ring with operation '$\times$' and 'o' where '$\times$' is the usual multiplication modulo 24 and a o b = a for all a, b $\in$ $Z_{24}$. A = {1, 5, 7, 11, 13, 17, 19, 23} is a near-ring. Let G = $S_3$. The group seminear-ring $Z_{24}$G contains AG to be a near-ring, hence the claim. Find zero divisors, idempotents, and units in $Z_{24}$G.

**DEFINITION 6.5.11**: *Let N be a seminear-ring and S a semigroup. The semigroup seminear-ring NS is said to be a Smarandache semigroup seminear-ring (S-semigroup seminear-ring) if and only if S is a S-semigroup.*

Thus it is an interesting problem to find when S-semigroup seminear-rings are S-seminear-rings.

***Example 6.5.2:*** Let N = Z$^+$ be a seminear-ring under '×' and 'o'. S = S(n) be the symmetric semigroup i.e. the semigroup of all mapping of a set of n elements to itself. NS is a Smarandache semigroup seminear-ring as S(n) is a S-semigroup.

**DEFINITION 6.5.12**: *Let N be a seminear-ring. We say N is Smarandache strongly subcommutative seminear-ring (S-strongly subcommutative seminear-ring) if every S-right ideal of it is a right quasi reflexive.*

**DEFINITION 6.5.13**: *Let R be a S-seminear-ring. We call a non-empty subset I of R to be a Smarandache common ideal (S-common ideal) related to A, (A ⊂ R and A is a near-ring) if*

i.      *x + y ∈ I, for all x, y ∈ I.*
ii.     *ax, anx,  xa ∈ I for all x ∈ I and a ∈ A.*

**THEOREM 6.5.4**: *If N is a S-seminear-ring. Every common ideal of N is a S-common ideal of N. But all S-common ideals of N need not be common ideals of N.*

*Proof*: Follows from the definitions. The reader is expected to construct an example of an S-common ideal which is not a common ideal.

***Example 6.5.3:*** Let $Z_{24}$ = {0, 1, 2, …, 23} be a S-seminear-ring under '×' and '.' Does this seminear-ring have S-common ideals and common ideals?

**DEFINITION 6.5.14**: *Let N be a seminear-ring. N is said to satisfy Smarandache chain seminear-ring condition (S-chain seminear-ring condition) if the set of S-ideals of N is totally ordered by inclusion.*

The reader is requested to formulate interesting results and properties about these S-chain conditions on seminear-rings.

**DEFINITIONS 6.5.15**: *Let N be a seminear-ring. An S-ideal I ⊆ N is called Smarandache semiprime ( S-semiprime ) if and only if for  all S-ideals J ⊂ S, $J^2$ ⊂ I implies J ⊂ I.*

**DEFINITION 6.5.16**: *Let P ⊂ N, N a seminear-ring. An S-ideal of N, P is called Smarandache prime (S-prime) if for all S-ideals I, J ∈ N,   IJ ⊂ P implies I ⊂ P or J ⊂ P.*

Motivated by these definitions the reader is expected to study the notions of S-zero divisors, S-idempotents, S-units, S-nilpotents and S-semi-idempotents which can be thought of as element wise study. Find also substructure like S-ideals, S-common ideals, S-chain condition for S-group seminear-rings, S-semigroup seminear-rings, S-seminear-ring and characterize them. It is upto the researchers to find innovative and ingenous results about S-seminear-rings. Several suggested problems proposed in chapter X may help the reader with more information.

We at this point wish to state several modular or non-modular lattices may contain subsets which are seminear-rings. Thus we can have a class of seminear-rings derived from modular or non-modular lattices as distributive lattices are trivially seminear-ring.

<u>**PROBLEMS:**</u>

1. Is the semigroup seminear-ring $Z_{12}$ S(5) a S-seminear-ring? Justify your answer.
2. Prove $Z_{27}S_7$ is a S-group seminear ring. Find S-ideals and S-zero divisors in them.
3. Find S-common ideals of the S-group seminear-ring $Z_{10}S_3$.
4. Does the S-seminear-ring $Z_{210}$ have
    i. Common ideals.
    ii. S-common ideals.
    iii. S-ideals which are common ideals.
    iv. Does it satisfy S-chain conditions?
5. Can the seminear-ring $Z_{24}$ satisfy
    a. Chain conditions.
    b. S-chain conditions.
    c. Does it have S-semiprime ideals.
    d. Does it have S-prime ideals.

**Chapter Seven**

# SOME APPLICATIONS OF SMARANDACHE NEAR-RINGS AND NEAR-RINGS

This chapter is a maiden effort to bring the two probable existing applications of near-rings and seminear-rings. We are not widely speaking about the applications of near-rings we only discuss here the two applications one in automatons and the other in the construction of error correcting codes. This chapter has three sections. In section one we recall the basic notions of an automaton, semi-automaton, Smarandache automaton and group automaton and give some interesting illustrations of them.

Section two completely gives the definition of Smarandache semigroup semi-automaton and illustrate them with examples and also give the associated syntactic near-ring. We introduce the state graph for these S-S-semigroup semi-automaton and extend this notion to the group automaton and Smarandache S-semigroup automaton (S-S-semigroup automaton) and define the concept of Smarandache syntactic near-ring. The final section is just only the introduction of Smarandache-planar near-rings used in the construction of  Balanced incomplete block design, (BIBD) we show the introduction of a S-planar near-rings leads to several BIBDs for any given S-planar near-ring.

## 7.1 Basics on automaton and on semi-automaton

In this section we just recall the definitions of automaton, semi-automaton, Smarandache automaton, Smarandache semi-automaton and group semi-automaton and finally the concept of syntactic near-ring given by [18]. We illustrate them with explicit examples so that it would become easy in the next section when we introduce the Smarandache equivalence of a syntactic near-ring and its probable application to group semi-automaton. As this section deals with application all basic relevant information are recalled.

**DEFINITION 7.1.1**: *A semi-automaton is a tripe $S = (Z, A, \delta)$ consisting of two nonempty sets Z and A and a function $\delta : Z \times A \to Z$, Z is called the set of states and A the input alphabet and $\delta$ the next state function of S.*

**DEFINITION 7.1.2**: *An automaton is a quintuple $A = (Z, A, B, \delta, \lambda)$ where $(Z, A, \delta)$ is a semi-automaton, B is a nonempty set called the output alphabet and $\lambda : Z \times A \to B$ is the output function.*

A (semi) automaton is finite if all sets Z, A and B are finite. Several types of automaton are studied.

We usually describe the (semi) automaton only by tables or by graphs.

Description of automaton:

Let A = {$a_1$, $a_2$, ..., $a_n$}, B = {$b_1$, $b_2$, ..., $b_n$} and Z = {$z_1$, $z_2$ ,..., $z_k$} where A is the input alphabet, B is the output alphabet and Z is the set of states. We describe $\delta$ the next state function from $Z \times A \rightarrow Z$ by the following transition table.

| $\delta$ | $a_1$ | **...** | $a_n$ |
|---|---|---|---|
| $z_1$ | $\delta(z_1, a_1)$ | $\cdots$ | $\delta(z_1, a_n)$ |
| $\vdots$ | $\vdots$ | | $\vdots$ |
| $z_k$ | $\delta(z_k, a_1)$ | $\cdots$ | $\delta(z_k, a_n)$ |

($\delta(z_i, a_j) \in$ Z).

The output table in case of automaton is given by $\lambda$, where $\lambda : Z \times A \rightarrow B$ where B is the output alphabet is given by the following output table.

| $\lambda$ | $a_1$ | **...** | $a_n$ |
|---|---|---|---|
| $z_1$ | $\lambda(z_1, a_1)$ | $\cdots$ | $\lambda(z_1, a_n)$ |
| $\vdots$ | $\vdots$ | | $\vdots$ |
| $z_k$ | $\lambda(z_k, a_1)$ | $\cdots$ | $\lambda(z_k, a_n)$ |

where $\lambda(z_i, a_i) \in$ B.

The graphical representation which is called a state graph are drawn by taking the set of states $z_1$, $z_2$, ..., $z_k$ as, 'discs'.

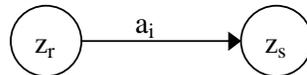

if $\delta(z_r, a_i) = z_s$. In case of an automaton, we have $\lambda(z_r, a_i)$ also as the output function so that the state graph in this case is

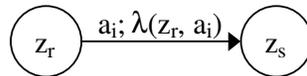

We illustrate by the following example the parity check automaton by both the tables and state graphs.

***Example (Parity-check Automaton)***: Let Z = Z = {$z_0$, $z_1$}, A = B = {0, 1} and

| $\delta$ | 0 | 1 |
|---|---|---|
| $z_0$ | $z_0$ | $z_1$ |
| $z_1$ | $z_1$ | $z_0$ |

| $\lambda$ | 0 | 1 |
|-----------|---|---|
| $z_0$ | 0 | 1 |
| $z_1$ | 0 | 1 |

The state graph of this automaton

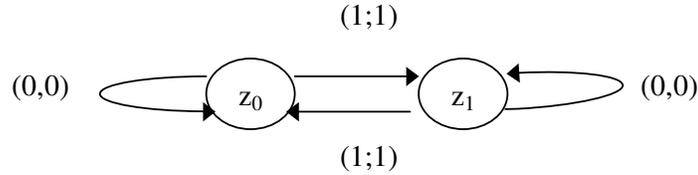

This automaton has liminations for a given input it can give an output. It cannot do any sequential operations so using free semigroups as automaton which can perform sequence of operations acting on sequence of elements. To achieve this we define the following:

$\overline{A}$ = free semigroup generated by the input alphabets A together with the empty sequence $\Lambda$.

$\overline{B}$ = free semigroup generated by the output alphabets B together with the empty sequence $\Lambda$.

Clearly the number of states Z in a machine cannot be adjusted as it is constructed or designed with a fixed number of states.

We extend the next state function $\delta$ and the output function $\lambda$ from $Z \times A$ to $Z \times \overline{A}$ by defining for $z \in Z$ and $a_1, a_2 \ldots, a_r \in A$ where $\delta : Z \times A \to Z$ and $\overline{\delta} . Z \times \overline{A} \to Z$ by

$$\overline{\delta}(z, \Lambda) = z$$
$$\overline{\delta}(z, a_1) = \delta(z, a_1)$$
$$\overline{\delta}(z, a_1 a_2) = \delta(\overline{\delta}(z, a_1), a_2)$$
$$\vdots$$
$$\overline{\delta}(z, a_1 a_2 \ldots a_n) = \delta(\overline{\delta}(z, a_1 a_2 \ldots a_{n-1}), a_n)$$

and

$\lambda : Z \times A \to B$ by $\overline{\lambda} : Z \times \overline{A} \to \overline{B}$ defined or extended by

$$\overline{\lambda}(z, \Lambda) = \Lambda$$
$$\overline{\lambda}(z, a_1) = \lambda(z, a_1)$$
$$\overline{\lambda}(z, a_1 a_2) = \lambda(z, a_1)\overline{\lambda}(\delta(z, a_1), a_2)$$
$$\vdots$$
$$\overline{\lambda}(z, a_1 a_2 \ldots a_n) = \lambda(z, a_1)\overline{\lambda}(\delta(z, a_1), a_2 \ldots a_r).$$

In this way we obtain functions

$$\overline{\delta} : Z \times \overline{A} \to Z \quad \text{and}$$

$$\overline{\lambda} : Z \times \overline{A} \to \overline{B}.$$

Thus the semi-automaton $S = (Z,\ A,\ \delta)$ is generalized to new semi-automaton $\overline{S} = (Z, \overline{A}, \overline{\delta})$. Similarly the automaton $A = (Z,\ A,\ B,\ \delta,\ \lambda)$ is generalized to the new automaton $\overline{A} = (Z, \overline{A}, \overline{B}, \overline{\delta}, \overline{\lambda})$; we can describe the operation by

$$
\begin{aligned}
z_1 &= z \\
z_2 &= \delta\,(z_1,\ a_1) \\
z_3 &= \overline{\delta}\,(z_1,\ a_1 a_2) \\
&= \overline{\delta}\,(\delta\,(z_1,\ a_1),\ a_2) \\
&= \delta\,(z_2,\ a_2) \\
&\quad \vdots
\end{aligned}
$$

The sectional graph of the automaton is:

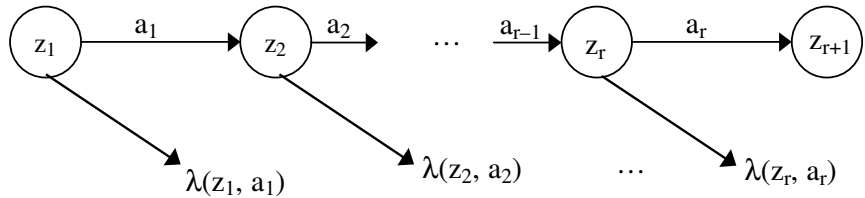

**Figure 7.1.1**

Now using the concept of semi-automaton and automaton we define the concept of Smarandache semi-automaton and Smarandache automaton as follows:

We would advise the reader to go through the definition of free groupoids and Smarandache free groupoids defined in chapter I.

**DEFINITION [98]:** $\gamma_s = (Z,\ \overline{A}_s,\ \overline{\delta}_s)$ *is said to be a Smarandache semi-automaton if* $\overline{A}_s = \langle A \rangle$ *is the free groupoid generated by A with $\Lambda$ the empty sequence adjoined with it and* $\overline{\delta}_s$ *is the function from $Z \times \overline{A}_s \to Z$. Thus the Smarandache semi-automaton (S semi-automaton) contains* $\gamma = (Z,\ \overline{A}_s, \overline{\delta}_s)$ *as a new semi-automaton which is a proper substructure of $\gamma_s$.*

Or equivalently we define a S-semi-automaton as one which has a new semi-automaton as a substructure.

Clearly the S-semi-automaton is a generalized structure than the new semi-automaton as all free semigroups are in $\overline{A}$ where $\overline{A}$ is the free groupoid generated by A.

**DEFINITION [98]:** $\overline{Y}'_s = (Z_1, \overline{A}_s, \overline{\delta}'_s)$ *is called the Smarandache subsemi-automaton (S-subsemi-automaton) of* $\overline{Y}_s = (Z_2, \overline{A}_s, \overline{\delta}'_s)$ *denoted by* $\overline{Y}_s \subseteq \overline{Y}_s$ *if $Z_1 \subset Z_2$ and* $\overline{\delta}'_s$ *is* *the restriction of $\overline{\delta}_s$ on $Z_1 \times \overline{A}_s$ and $\overline{Y}_s'$ has a proper subset $\overline{H} \subset \overline{Y}_s'$ such that $\overline{H}$ is a new semi-automaton.*

We just illustrate them by an example.

***Example 7.1.1***: Let $Z = Z_4 = (2, 2)$ and $A = Z_3$ (1, 2) we define $\delta$ (z, a) = z * a (mod 4) '*' as in Z. The table for the semi-automaton is given by

| $\delta$ | 0 | 1 | 2 |
|---|---|---|---|
| 0 | 0 | 2 | 0 |
| 1 | 2 | 0 | 2 |
| 2 | 0 | 2 | 0 |
| 3 | 2 | 0 | 2 |

The graph for it is

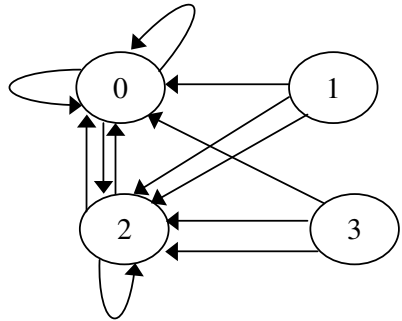

**Figure 7.1.2**

Thus this has Smarandache subsemi-automaton $Z_1$ given by $Z_1 = \{0, 2\}$ states.

**<u>Remark</u>**: A machine equipped with this S-semi-automaton can use any new automaton as per need. Now proceed on to recall the definition of group semi-automaton.

**DEFINITION [98]:** $\overline{K}_s = (Z, \overline{A}_s, \overline{B}_s, \overline{\delta}_s, \overline{\lambda}_s)$ *is defined to be a Smarandache automaton (S-automaton) if* $\overline{K} = (Z, \overline{A}_s, \overline{B}_s, \overline{\delta}_s, \overline{\lambda}_s)$ *is the new automaton and* $\overline{A}_s$ *and* $\overline{B}_s$, *the Smarandache free groupoids so that* $\overline{K} = (Z, \overline{A}_s, \overline{B}_s, \overline{\delta}_s, \overline{\lambda}_s)$ *is the new automaton got from K, and* $\overline{K}$ *is strictly contained in* $\overline{K}_s$.

Thus S-automaton enables us to adjoin some more elements which is present in A and freely generated by A, as a free groupoid; that will be the case when the compositions may not be associative.

Secondly, by using S-automaton we can couple several automaton as

$$Z \quad = \quad Z_1 \cup Z_2 \cup ... \cup Z_n$$
$$A \quad = \quad A_1 \cup A_2 \cup ... \cup A_n$$
$$B \quad = \quad B_1 \cup B_2 \cup ... \cup B_n$$

$$\lambda = \lambda_1 \cup \lambda_2 \cup ... \cup \lambda_n$$
$$\delta = \delta_1 \cup \delta_2 \cup ... \cup \delta_n.$$

where the union of $\lambda_i \cup \lambda_j$ and $\delta_i \cup \delta_j$ denote only extension maps as '$\cup$' has no meaning in the composition of maps, where $K_i = (Z_i, A_i, B_i, \delta_i, \lambda_i)$ for i = 1, 2, 3, ..., n and $\overline{K} = \overline{K}_1 \cup \overline{K}_2 \cup ... \cup \overline{K}_n$. Now $\overline{K}_s = (\overline{Z}_s, \overline{A}_s, \overline{B}_s, \overline{\lambda}_s, \overline{\delta}_s)$ is the S-automaton.

A machine equipped with this S-automaton can use any new automaton as per need.

***Example 7.1.2:*** Let $Z = Z_4$ (3, 2), $A = B = Z_5$ (2, 3). $K = (Z, A, B, \delta, \lambda)$ is a S-automaton defined by the following tables where $\delta$ (z, a) = z * a (mod 4) and $\lambda$ (z, a) = z * a (mod 5).

| $\lambda$ | 0 | 1 | 2 | 3 | 4 |
|---|---|---|---|---|---|
| 0 | 0 | 3 | 1 | 4 | 2 |
| 1 | 2 | 0 | 3 | 1 | 4 |
| 2 | 4 | 2 | 0 | 3 | 1 |
| 3 | 1 | 4 | 2 | 0 | 3 |

| $\delta$ | 0 | 1 | 2 | 3 | 4 |
|---|---|---|---|---|---|
| 0 | 0 | 2 | 0 | 2 | 0 |
| 1 | 3 | 1 | 3 | 1 | 3 |
| 2 | 2 | 0 | 2 | 0 | 2 |
| 3 | 1 | 3 | 1 | 3 | 1 |

We obtain the following graph:

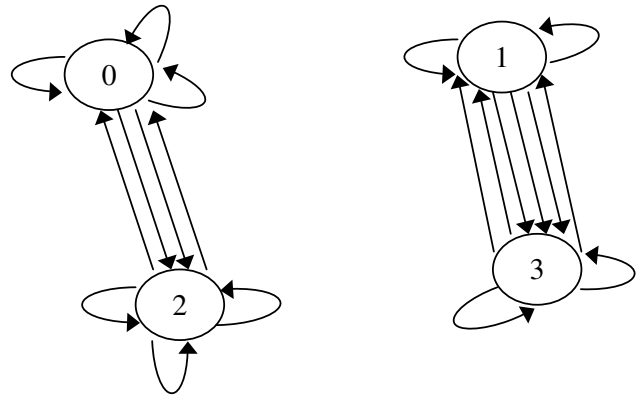

**Figure 7.1.3**

Thus we see this automaton has two Smarandache subautomatons given by the sates {0, 2} and {1, 3}.

**DEFINITION [98]:** $\overline{K}_s^{'} = (Z_1, \overline{A}_s, \overline{B}_s, \overline{\delta}_s^{'}, \overline{\lambda}_s^{'})$ *is called Smarandache sub-automaton (S-subautomaton) of* $\overline{K}_s = (Z_2, \overline{A}_s, \overline{B}_s, \overline{\delta}_s, \overline{\lambda}_s)$ *denoted by* $\overline{K}_s^{'} \subseteq \overline{K}_s$ *if* $Z_1 \subseteq Z_2$ *and* $\overline{\delta}_s^{'}$ *and* $\overline{\lambda}_s^{'}$ *are the restriction of* $\overline{\delta}_s$ *and* $\overline{\lambda}_s$ *respectively on* $Z_1 \times \overline{A}_s$ *and has a proper subset* $\overline{H} \subset \overline{K}_s^{'}$ *such that* $\overline{H}$ *is a new automaton.*

**DEFINITION [98]:** *Let* $\overline{K}_1$ *and* $\overline{K}_2$ *be any two S-automaton where* $\overline{K}_1 = (Z_1, \overline{A}_s, \overline{B}_s, \overline{\delta}_s, \overline{\lambda}_s)$ *and* $\overline{K}_2 = (Z_2, \overline{A}_s, \overline{B}_s, \overline{\delta}_s, \overline{\lambda}_s)$. *A map* $\phi : \overline{K}_1$ *to* $\overline{K}_2$ *is a Smarandache automaton homomorphism (S-automaton homomorphism*

*) if* $\phi$ *restricted from* $K_1 = (Z_1, A_1, B_1, \delta_1, \lambda_1)$ *and* $K_2 = (Z_2, A_2, B_2, \delta_2, \lambda_2)$ *denoted by* $\phi_r$, *is an automaton homomorphism from* $K_1$ *to* $K_2$. $\phi$ *is called a Smarandache monomorphism (epimorphism or isomorphism) if there is an isomorphism* $\phi_r$ *from* $K_1$ *to* $K_2$.

**DEFINITION [98]:** *Let* $\overline{K}_1$ *and* $\overline{K}_2$ *be two S-automatons, where* $\overline{K}_1 = (Z_1, \overline{A}_s, \overline{B}_s, \overline{\delta}_s, \overline{\lambda}_s)$ *and* $\overline{K}_2 = (Z_2, \overline{A}_s, \overline{B}_s, \overline{\delta}_s, \overline{\lambda}_s)$.

*The Smarandache automaton direct product (S-automaton direct product) of* $\overline{K}_1$ *and* $\overline{K}_2$ *denoted by* $\overline{K}_1 \times \overline{K}_2$ *is defined as the direct product of the automaton* $K_1 = (Z_1, A_1, B_1, \delta_1, \lambda_1)$ *and* $K_2 = (Z_2, A_2, B_2, \delta_2, \lambda_2)$ *where* $K_1 \times K_2 = (Z_1 \times Z_2, A_1 \times A_2, B_1 \times B_2, \delta, \lambda)$ *with* $\delta((z_1, z_2), (a_1, a_2)) = (\delta_1 (z_1, a_2), \delta_2 (z_2, a_2)), \lambda ((z_1, z_2), (a_1, a_2)) = (\lambda_1 (z_1, a_2), \lambda_2 (z_2, a_2))$ *for all* $(z_1, z_2) \in Z_1 \times Z_2$ *and* $(a_1, a_2) \in A_1 \times A_2$.

**_Remark:_** *Here in* $\overline{K}_1 \times \overline{K}_2$ *we do not take the free groupoid to be generated by* $A_1 \times A_2$ *but only free groupoid generated by* $\overline{A}_1 \times \overline{A}_2$.

Thus the S-automaton direct product exists wherever a automaton direct product exists.

We have made this in order to make the Smarandache parallel composition and Smarandache series composition of automaton extendable in a simple way.

**DEFINITION [98]:** *A Smarandache groupoid* $G_1$ *divides a Smarandache groupoid* $G_2$ *if the corresponding semigroups* $S_1$ *and* $S_2$ *of* $G_1$ *and* $G_2$ *respectively divides, that is, if* $S_1$ *is a homomorphic image of a sub-semigroup of* $S_2$.

*In symbols* $G_1 \mid G_2$. *The relation divides is denoted by '|'.*

**DEFINITION [98]:** *Let* $\overline{K}_1 = (Z_1, \overline{A}_s, \overline{B}_s, \overline{\delta}_s, \overline{\lambda}_s)$ *and* $\overline{K}_2 = (Z_2, \overline{A}_s, \overline{B}_s, \overline{\delta}_s, \overline{\lambda}_s)$ *be two S-automaton. We say the Smarandache automaton* $\overline{K}_1$ *divides the S-automaton* $\overline{K}_2$ *if in the automatons* $K_1 = (Z_1, A, B, \delta_1, \lambda_1)$ *and* $K_2 = (Z_2, A, B, \delta_2, \lambda_2)$, *if* $K_1$ *is the homomorphic image of a sub-automaton of* $K_2$.

*Notationally* $K_1 \mid K_2$.

**DEFINITION [98]:** *Two S-automaton* $\overline{K}_1$ *and* $\overline{K}_2$ *are said to be equivalent if they divide each other. In symbols* $\overline{K}_1 \sim \overline{K}_2$.

We proceed on to define direct product S-automaton and study about them. We can extend the direct product of semi-automaton to more than two S-automatons. Using the definition of direct product of two automaton $K_1$ and $K_2$ with an additional assumption we define Smarandache series composition of automaton.

**DEFINITION [98]:** *Let $K_1$ and $K_2$ be any two S- automatons where $\overline{K}_1$ = $(Z_1, \overline{A}_s, \overline{B}_s, \overline{\delta}_s, \overline{\lambda}_s$ ) and $\overline{K}_2$ = $(Z_2, \overline{A}_s, \overline{B}_s, \overline{\delta}_s, \overline{\lambda}_s$ ) with an additional assumption $A_2 = B_1$.*

*The Smarandache automaton composition series (S-automaton composition series) denoted by $\overline{K}_1 \Vdash \overline{K}_2$ of $\overline{K}_1$ and $\overline{K}_2$ is defined as the series composition of the automaton $K_1 = (Z_1, A_1, B_1, \delta_1, \lambda_1)$ and $K_2 = (Z_2, A_2, B_2, \delta_2, \lambda_2)$ with $\overline{K}_1 \Vdash \overline{K}_2$ = $(Z_1 \times Z_2, A_1, B_2, \delta, \lambda)$ where $\delta((z_1, z_2), a_1) = (\delta_1(z_1, a_1), \delta_2(z_2, \lambda_1(z_1, a_1))$ and $\lambda((z_1, z_2), a_1) = (\lambda_2(z_2, \lambda_1(z_1, a_1))$ $((z_1, z_2) \in Z_1 \times Z_2, a_1 \in A_1)$.*

*This automaton operates as follows: An input $a_1 \in A_1$ operates on $z_1$ and gives a state transition into $z_1' = \delta_1(z_1, a_1)$ and an output $b_1 = \lambda_1(z_1, a_2) \in B_1 = A_2$. This output $b_1$ operates on $Z_2$ transforms a $z_2 \in Z_2$ into $z_2' = \delta_2(a_2, b_1)$ and produces the output $\lambda_2(z_2, b_1)$.*

*Then $\overline{K}_1 \Vdash \overline{K}_2$ is in the next state $(z_1', z_2')$ which is clear from the following circuit:*

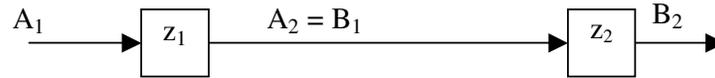

Now a natural question would be do we have a direct product, which corresponds to parallel composition, of the two S-automatons $\overline{K}_1$ and $\overline{K}_2$.

Clearly the S-direct product of automatons $\overline{K}_1 \times \overline{K}_2$ since $Z_1$, and $Z_2$ can be interpreted as two parallel blocks. $A_1$ operates on $Z_i$ with output $B_i$ ($i \in \{1,2\}$), $A_1 \times A_2$ operates on $Z_1 \times Z_2$, the outputs are in $B_1 \times B_2$. The circuit is given by the following diagram:

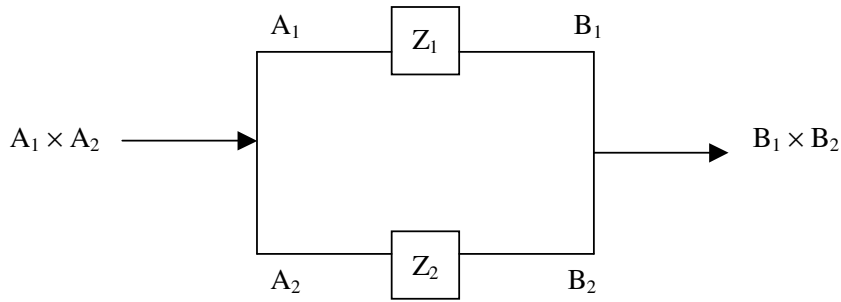

**Figure 7.1.4**

**DEFINITION [18]:** *A group semi-automaton is an ordered triple $S = (G, I, \delta)$ where $G$ is a group; $I$ is a set of inputs and $\delta : G \times I \to G$ is a transition function.*

For more about these refer [18]. We illustrate this by an example.

***Example 7.1.3:*** Let $G = S_3$, where $S = \{p_0, p_1, \ldots, p_5\}$ $I = \{0, 1, 2,\}$ be the set of inputs. The transition function $\delta : G \times I \to G$ defined by

$\delta (p_0, 0) = p_0$

$\delta (p_i, 0) = p_i$ for $i = 1, 2, \ldots, 5$

$\delta (p_i\ 1) = p_{i+1}$ and $\delta (p_5\ 1) = 0$, for $i = 1, 2, \ldots, 4$

$\delta (p_i\ 2) = p_{i+2}$, for $i = 1, 2, 3$ and $\delta (p_4, 2) = p_0\ \delta = (p_5, 2) = p_1$

Clearly, $(G, I, \delta)$ is a group semi-automaton.

Now using the notions of [27] we define the following.

**DEFINITION [27]**: *A semi-automaton $S = (Q, \times, \delta)$ with state set $Q$ input monoid $X$ and state transition function $\delta$; $\delta : Q \times X \to Q$ is called a group semi-automaton if $Q$ is an additive group. $\delta$ is called additive if there is some $x_0 \in X$ with $\delta(q, x) = \delta(q, x_0) + \delta(0, x)$ and $\delta(q - q^1, x_0) = \delta(q, x_0) - \delta(q^1, x_0)$ for all $q, q^1 \in Q$ and $x \in X$. Then there is some homomorphism*

$\psi : Q \to Q$ *and some map*

$\alpha : \times \to Q$ *with $(x_0) = 0$ and $\delta(q, x) = \psi(q) + \alpha(x)$.*

*Let $\delta_x$ for a fixed, $x \in X$ be the map $Q \to Q$; $q \to \delta(q, x)$ then $\{\delta_x \, / \, x \in X\}$ generates a subnear-ring $N(S)$ of near-ring $(M(Q), +, .)$ of all mappings on $Q$. This near-ring $N(S)$ is called the syntactic near-ring.*

Thus both [18] and [27] have defined group semi-automaton and syntactic near-ring.

**PROBLEMS:**

1. Describe the semi-automaton with $A = (0, 1, 2, 3, 4)$; $Z_3 = (0, 1, 2)$ and $\delta : Z_3 \times A \to Z_3$ by $\delta (z, a_i) = za_i \pmod 3$ by the table and graph.

2. Describe the following automaton $A = (Z, A, B, \delta, \lambda)$ given by $Z = Z_7$; $A = B = Z_5$, $\delta (z, a) = z + a \pmod 7$, $\lambda (z, a) = z.a \pmod 5$, by a state graph and by the tables.

3. Convert the automaton $A$ given in problem 2 into a new automaton.

4. Let $A = B = Z_5$ $(3, 2)$ be a S-groupoid of order 7 and $Z = Z$ $(2, 2)$. Describe the S-automaton where $\delta (z, a) = z * a \pmod 4$ and $\lambda (za) = z * a \pmod 5$ by its graph.

5. For the group $G = \langle g \, / \, g^7 = 1 \rangle$ and $I = \{0, 1, 2, 3\}$, set of input with $\delta : G \times I \to G$ defined by $\delta(g^i, 0) = g^{i+0}$ $\delta(g^i, 3) = g^{i+1}$, $\delta(g^i, 2) = g^{i+2}$ and $\delta(g^i, 3) = g^{i+3}$ for $g^i \in G$ describe the group semi-automaton what is the associated syntactic near-ring.

## 7.2 Smarandache S-semigroup semi-automaton and the associated Smarandache syntactic near-ring

In this section we introduce the concept of Smarandache semigroup semi-automaton analogous to the group semi-automaton. Further for each group semi-automaton an associated syntactic near-ring is also given. Here we define the Smarandache analogous of it.

**DEFINITION 7.2.1**: *Let G be S-semigroup. A Smarandache S-semigroup semiautomaton (S-S-semigroup semiautomaton) is an ordered triple S(S) = (G, I, δ) where G is a S-semigroup; I is a set of inputs and δ : G × I → G is a transition function.*

Clearly all group semi-automatons are trivially S-S-semigroup semiautomaton. Thus the class of S-S semigroup semiautomaton properly contains the class of group semi-automaton.

***Example 7.2.1***: Let G = {0, 1, 2, … , 5} = $Z_6$ and I = {0, 1, 2} be the set of states; define the transition function δ : G × I → G by δ(g, i) = g.i (mod 6) thus we get the S-semigroup semi-automaton. We give the following state graph representation for this S-S semigroup semi-automaton.

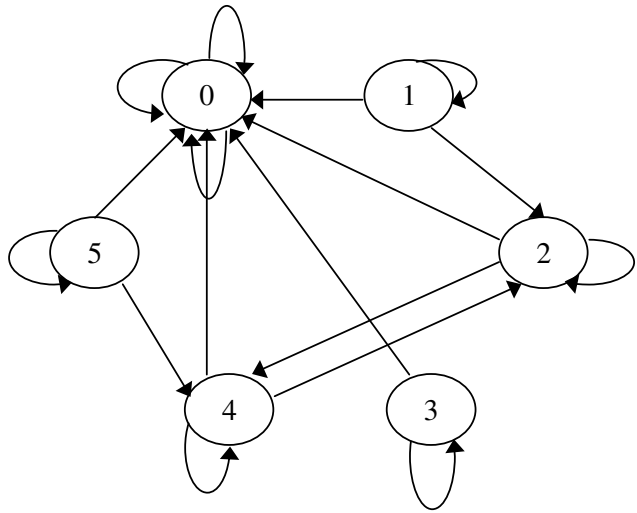

**Figure 7.2.1**

Now we see when we change the function $δ_1$ for the same S-semigroup, G and I = {0, 1, 2} we get a different S-S-semigroup semiautomaton which is given by the following figure; here the transition fuctions. δ : G × I → G is defined by $δ_1$ (g, i) = (g + i) mod 6. We get a nice symmetric S-S-semigroup semiautomaton.

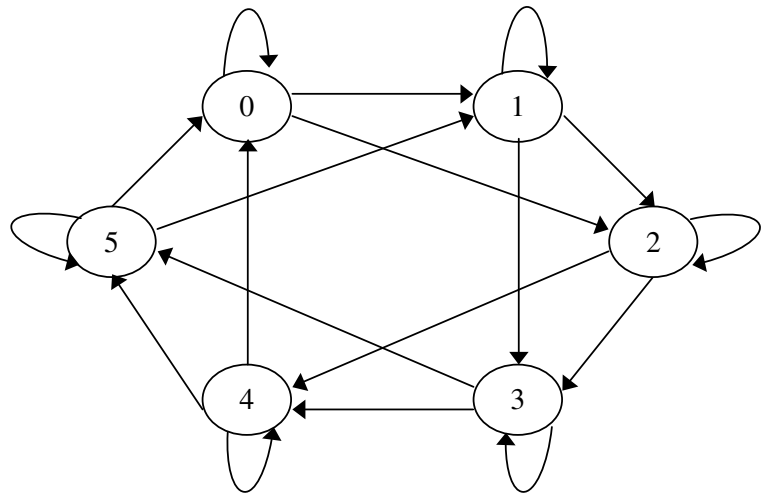

**Figure 7.2.2**
120

**DEFINITION 7.2.2**: *Let S (S) = (G, I, δ) where G is a S-semigroup I the set of input alphabets and δ the transition fiction as in case of S-S semigroup semi-automaton. We call a proper subset B ⊂ G with B a S-subsemigroup of G and $I_1 ⊂ I$ denoted (B, $I_1$, δ) = S(SB) is called the Smarandache subsemigroup semi-automaton (S-subsemigroup semiautomaton) if δ : B × $I_1$ → B. A S-S-semigroup semi-automaton may fail to have sometimes S-S-subsemigroup subsemi-automaton.*

**Example 7.2.2**: In example 7.2.1 with S = {$Z_6$, I, δ} when we take B = {0, 2, 4} we get a S-S-semigroup subsemi-automaton

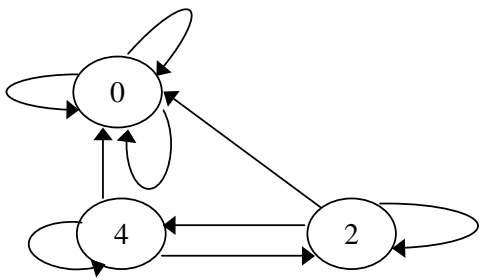

**Figure 7.2.3**

as B = {0, 2, 4} is a S-subsemigroup of $Z_6$.

Now to the best of our knowledge we have not seen the definition of group automaton we define it and then proceed on to get corresponding near-rings.

**DEFINITION 7.2.3**: *A group automaton is an ordered quintuple A = {G, A, B, δ, λ} where G is a group, A is a set of inputs, B is a set of outputs δ : G × A → G and λ : G × A → B is the transition and output functions respectively.*

The Smarandache analogous of this would be as follows.

**DEFINITION 7.2.4**: *A Smarandache S-semigroup automaton (S-S semigroup automaton) is an ordered quintuple S(A) = (G, A, B, δ, λ) where G is a S-semigroup and A set of input alphabets and B set of output alphabets, δ and λ are the transition and output functions respectively.*

*If δ restricted form G × A → G maps a subgroup T ⊂ G into itself δ : T × A → T then we can say the S-S-semigroup automaton has (T ⊂ G, A, B, δ, λ) to be a group automaton.*

We call this level I S-S-semigroup automaton. Now we proceed on to define S-S-semigroup semiautomaton and S-S-semigroup automaton of level II.

**DEFINITION 7.2.5**: *Let S = (G, I, δ) be an ordered triple where G is only a S-semigroup. If B ⊂ G and B is a group under the operations of G, and if $S^I$ = (B, I, δ) is a group semi-automaton with B a group, I is a set of inputs and δ : B × I → B then*

*S = (G, I, δ) is called the Smarandache S-semigroup semi-automaton of level II (S-S-semigroup semiautomaton II). Similarly we define Smarandache S-semigroup automaton of level II (S-S-semigroup automaton II).*

It is left for the reader to find any relation existing between the two levels of automaton. Now we proceed on to see how these concepts will produce for us the S-near-rings.

**DEFINITION 7.2.6**: *Let $S = (Q, X \delta)$ be a S-S-semigroup semi-automaton where Q is a S-semigroup under addition '+' i.e. $G \subset Q$ is such that G is a S-semigroup under '+' and G a proper subset of Q Let $\delta : G \times X \to G$ is well defined i.e. restriction of δ to G is well defined so $(G, X, \delta) = C$ becomes a semi-automaton C is called additive if there is some $x_0 \in X$ with*

$$\delta (q, x) = \delta (q, x_0) + \delta (0, x) \text{ and}$$
$$\delta (q - q^1, x_0) = \delta (q_1 x_0) - \delta (q^1, x_0) \text{ for all } q, q^1 \in G \text{ and } x \in X.$$

*Then there is some homomorphism*

$$\psi : G \to G \text{ and some map}$$
$$\alpha : X \to G \text{ with } \psi(x_0) = 0 \text{ and } \delta(q, x) = \psi(q) + \alpha(x).$$

*Let $\delta_x$ for a fixed $x \in X$ be the map from $G \to G$; $q \to \delta (q, x)$. Then $\{\delta_x / x \in X\}$ generates a subnear-ring N (P) of the near-ring (M (G), +, .) of all mappings on G; this near-ring N(P) is called as the Smarandache syntactic near-ring (S-syntactic near-ring).*

The varied structural properties are enjoyed by this S-syntactic near-ring for varying S-semigroups. The most interesting feature about these S-syntactic near-ring is that for one S-S-semigroup semi-automaton we have several group semiautomaton depending on the number of valid groups in the S-semigroup. This is the vital benefit in defining S-S-semigroup automaton and S-syntactic near-ring. So for a given S-S-semigroup semi-automaton we can have several S-syntactic near-ring.

Thus the Smarandache notions in this direction has evolved several group automaton for a given S-semigroup.

**PROBLEMS:**

1. Let $S = (G, I, \delta)$ where $G = S_5$ define : $G \times I \to G$ in any suitable way where I = {0, 1} and obtain state graph of S.
2. Let $S(S) = \{S(5), I = \{0, 1, 2\}, \delta\}$ study the S-S-semigroup semi-automaton.
3. Give an example of a syntactic near-ring.
4. Find an example of S-syntactic near-ring. Does it have S-ideals which have any impact on the S-S-semigroup semi-automaton?
5. Let $S = \{Z_{12}, I, \delta\}$, $I = \{0, 1, 2, 3, 4\}$, $\delta = Z_{12} \times I \to Z_{12}$, $\delta (x_i, k) = x_i.k$ (mod 12), and $Z_{12}$ is a S-semigroup under '$\times$' mod 12. Does this have S-S-semigroup subsemi-automaton?

## 7.3 Applications of near-rings to error correcting codes and their Smarandache analogue

The study of how experiments can be organized systematically so that statistical analysis can be applied in an interesting problem which is carried out by several researchers. In the planning of experiments it often occurs that results are influenced by phenomena outside the control of the experimenter. The introduction of balanced incomplete block design (BIBD) helps in avoiding undesirable influences in the experiment. In general, if we have to test the effect of r different conditions with m possibilities for each conditions this leads to a set of r orthogonal latin squares.

A planar near-ring can be used to construct balanced incomplete block designs (BIBD) of high efficiency; we just state how they are used in developing error correcting codes as codes which can correct errors is always desirable than the ones which can only detect errors. In view of this we give the following definition. Just we recall the definition of BIBD for the sake of completeness.

**DEFINITION 7.3.1**: *A balanced incomplete block design (BIBD) with parameters ($v$, $b$, $r$, $k$, $\lambda$) is a pair (P, B) with the following properties.*

>   i.  *P is a set with $v$ elements.*
>   ii.  *B = ($B_1$, …, $B_b$) is a subset of p (P) with b elements.*
>   iii.  *Each $B_i$ has exactly k elements where $k < v$ each unordered pair (p, q) with p, q $\in$ P, p $\neq$ q occurs in exactly $\lambda$ elements in B.*

*The set $B_1$, …, $B_b$ are called the blocks of BIBD. Each a $\in$ P occurs in exactly r sets of B. Such a BIBD is also called a ($v$, $b$, $r$, $k$, $\lambda$) configuration or 2 − ($v$, $k$, $\lambda$) tactical configuration or design. The term balance indicates that each pair of elements occurs in exactly the same number of block, the term incomplete means that each block contains less than $v$ - elements. A BIBD is symmetric if $v = b$.*

*The incidence matrix of a ($v$, $b$, $r$, $k$, $\lambda$) configuration is the v × b matrix A = ($a_{ij}$) where*

$$a_{ij} = \begin{cases} 1 \ if \ \ i \in B_j \\ 0 \ otherwise \end{cases},$$

*here i denotes the $i^{th}$ element of the configuration. The following conditions are necessary for the existence of a BIBD with parameters $v$, $b$, $r$, $k$, $\lambda$.*

>   1.  *bk = rv.*
>   2.  *r(k − 1) = $\lambda$ (v − 1).*
>   3.  *b $\geq$ v.*

*Recall a near-ring N is called planar (or Clay near-ring) if for all equation x o a = x o b + c. (a, b, c $\in$ N, a $\neq$ b) have exactly one solution x $\in$ N.*

***Example 7.3.1***: Let $(Z_5, '+', *)$ where '+' is the usual '+' and '*' is $n * 0 = 0$, $n * 1 = n * 2 = n$, $n * 3 = n * 4 = 4n$ for all $n \in N$. Then $1 \equiv 2$ and $3 \equiv 4$. N is planar near-ring; the equation $x * 2 = x * 3 + 1$ with $2 \neq 3$ has unique solution $x = 3$.

A planar near-ring can be used to construct BIBD of high efficiency where by high efficiency 'E' we mean $E = \lambda v / rk$; this E is a number between 0 and 1 and it estimates the quality of any statistical analysis if $E \geq 0.75$ the quality is good.

Now, how does one construct a BIBD from a planar ring? This is done in the following way.

Let N be a planar ring. Let $a \in N$. Define $g_a : N \to N$ by $n \to n \circ a$ and form $G = \{g_a / a \in N\}$. Call $a \in N$ "group forming" if $a \circ N$ is a subgroup of $(N, +)$. Let us call sets a . $N + b$ ($a \in N*$, $b \in N$) blocks. Then these blocks together with N as the set of "points" form a tactical configuration with parameters

$$(\upsilon, b, r, k, \lambda) = (\upsilon, \frac{\alpha_1 \upsilon}{|G|} + \alpha_2 \upsilon, \alpha_1 + \alpha_2 |G|, |G|, \lambda)$$

where $\upsilon = |N|$ and $\alpha_1$ ($\alpha_2$) denote the number of orbits of F under the group $G \setminus \{0\}$ which consists of entirely of group forming elements. The tactical configuration is a BIBD if and only if either all elements are group forming or just 0 is group forming.

Now how are they used in obtaining error correcting codes. Now using the planar near-ring one can construct a BIBD. By taking either the rows or the columns of the incidence matrix of such a BIBD one can obtain error correcting codes with several nice features. Now instead of using a planar near-ring we can use Smarandache planar near-ring. The main advantage in using the Smarandache planar near-ring is that for one S-planar near-ring we can define more than one BIBD. If there are m-near-field in the S-planar near-ring we can build m BIBD. Thus this may help in even comparison of one BIBD with the another and also give the analysis of the common features.

Thus the S-planar near-ring has more advantages than the usual planar near-rings and hence BIBD's constructed using S-planar near-rings will even prove to be an efficient error correcting codes.

<u>**PROBLEMS:**</u>

1.  Using $(Z_7, +, .)$ build a planar near-ring.
2.  Can $(Z_7, +, .)$ help in constructing S-planar near-ring? Justify your claim.
3.  For the planar near-ring $(Z_{11}, +, .)$ construct the associated BIBD.
4.  For the S-planar near-ring constructed using $Z_{24}$, how many associated BIBD's can be constructed?
5.  For the S-planar ring $N = (Z_{12}, +, .)$ using the BIBD associated with N construct error correcting codes.

**Chapter Eight**

# SMARANDACHE NON-ASSOCIATIVE NEAR-RING AND SEMINEAR-RING

This chapter is solely devoted to the study and introduction of Smarandache non-associative near-rings (S-NA-near-rings) and Smarandache non-associative seminear-rings (S-NA-seminear-ring). The study of non-associative near-rings started in the year 1978. They have called the loop near-ring which is a non-associative structure with respect to the operation '+' which is very unusual in rings.

This chapter is organized into five sections. In section one we just define S-NA-seminear-ring and define several of its Smarandache properties. In section two we construct several new types of S-NA-near-rings and seminear-rings using S-mixed direct product. Here we also define level II seminear-rings and study them. Section three is completely devoted to the introduction of Smarandache loop near-ring, Smarandache near loop ring and the identities satisfied by them. In section four we study the new S-NA-seminear-rings built using loops and groupoids. In the final section a few new notions which are studied for near-rings and seminear-rings are studied in case of S-seminear-rings and S-near-rings.

## 8.1 Smarandache non-associative seminear-ring and its properties

In this section we just define the notion of S-NA-seminear-ring and introduce several properties enjoyed by these structures to S-seminear-ring. This study is interesting and innovative as this gives a new class of non-associative seminear-rings. To the best our knowledge no one has introduced non-associative seminear-rings; only study of seminear-ring have been carried out by researchers.

**DEFINITION 8.1.1**: *Let (N, '+', '.') be a non-empty set endowed with two binary operation '+' and '.' satisfying the following:*

    a. *(N, +) is a semigroup.*
    b. *(N, .) is a groupoid.*
    c. *(a + b) c = a.c + b.c for all a, b, c ∈ N; (N, '+', '.') is called the right seminear-ring which is non-associative.*

*If we replace (c) by a. (b + c) = a.b + a.c for all a, b, c ∈ N Then (N, '+', '.') is a non-associative left seminear-ring.*

We in this text denote by (N, +, .) a non-associative right seminear-ring and by default of notation call N just a non-associative seminear-ring. The theory of right or left seminear-ring runs completely parallel in both cases. (N, .) will be only groupoids and never semigroups for we assume non-associativity of seminear-rings.

***Example 8.1.1***: Let $Z_2 = \{0, 1\}$ be a near-ring and G a groupoid. The groupoid near-ring $Z_2G$ is a non-associative near-ring. In view of this we have the following.

***Example 8.1.2***: All non-associative near-rings are non-associative seminear-rings.

***Example 8.1.3***: Let $Z_{12} = \{0, 1, 2, \dots, 11\}$ be a seminear-ring with operation '+' and '.'. (G, .) be a groupoid given by the following table:

| . | 0 | 1 | 2 | 3 | 4 | 5 |
|---|---|---|---|---|---|---|
| 0 | 0 | 3 | 0 | 3 | 0 | 3 |
| 1 | 1 | 4 | 1 | 4 | 1 | 4 |
| 2 | 2 | 5 | 2 | 5 | 2 | 5 |
| 3 | 3 | 0 | 3 | 0 | 3 | 0 |
| 4 | 4 | 1 | 4 | 1 | 4 | 1 |
| 5 | 5 | 2 | 5 | 2 | 5 | 2 |

$Z_{12}G$ is a non-associative seminear-ring. For $Z_{12}G$ under multiplication is only a groupoid. We define substructures in non-associative seminear-ring.

**DEFINITION 8.1.2**: *Let (N, +, .) be a seminear-ring which is not associative. A subset P of N is said to be a subseminear-ring if (P, +, .) is a seminear-ring.*

**DEFINITION 8.1.3**: *Let N be a non-associative seminear-ring. An additive subsemigroup A of N is called the N-subsemigroup (right N-subsemigroup) if $NA \subseteq A$ ($AN \subset A$) where $NA = \{na \, / \, n \in N, \, a \in A\}$.*

**DEFINITION 8.1.4**: *A non-empty subset I of N is called left ideal in N if*

    *i.       (I, +) is a normal subsemigroup of (N, +).*
    *ii.      $n (n_1 + i) + nn_1 \in I$ for each $i \in I$, $n, n_1 \in N$.*

**DEFINITION 8.1.5**: *Let N be a non-associative seminear-ring. A nonempty subset I of N is called an ideal in N if*

    *i.       I is a left ideal.*
    *ii.      $IN \subset I$.*

**DEFINITION 8.1.6**: *A non-associative seminear-ring N is called left bipotent if $Na = Na^2$ for a in N.*

**DEFINITION 8.1.7**: *A non-associative seminear-ring N is said to be a s-seminear-ring if $a \in Na$ for each a in N.*

The following definitions about strictly prime ideals would be of interest when we develop Smarandache notions.

**DEFINITION 8.1.8**: *An ideal P ($\neq N$) is called strictly prime if for any two N-subsemigroups A and B of N such that $AB \subset P$ then $A \subset P$ or $B \subset P$.*

**DEFINITION 8.1.9**: *An ideal B of a non-associative (NA for short) seminear-ring N is called strictly essential if B ∩ K ≠ {0} for every non-zero N-subsemigroup K of N.*

**DEFINITION 8.1.10**: *An element x in N is said to be singular if there exists a nonzero strictly essential left ideal A in N such that Ax = {0}.*

Several other analogous results existing in near-rings and seminear-rings can also be defined for non-associative seminear-rings. Now we proceed on to define Smarandache notions.

**DEFINITION 8.1.11**: *Let (N, +, .) be a non associative seminear-ring. N is said be a Smarandache non associative seminear-ring of level I (S-NA seminear-ring I) if*

       *i.  (N, +) is a S-semigroup.*
      *ii.  (N, .) is a S-groupoid.*
    *iii.  (a + b) c = a.c + b.c for all a, b, c ∈ N.*

**DEFINITION 8.1.12**: *Let (N, +, .) be a non-associative seminear-ring. N is said to have a Smarandache subseminear-ring (S-subseminear-ring) P ⊂ N if P is itself a S-seminear-ring.*

**DEFINITION 8.1.13**: *Let (N, +, .) be a non-associative seminear-ring. An additive S-subsemigroup A of N is called the Smarandache N-subsemigroup (S-N-subsemigroup) (right N-subsemigroup) if NA ⊂ A (AN ⊂ A) where Na = {na / n ∈ N, a ∈ A}.*

**DEFINITION 8.1.14**: *A non-empty subset I of N is called Smarandache left ideal (S-left ideal) in N if*

    *1.  (I, +) is a normal S-subsemigroup of (N, +).*
    *2.  n (n₁ + i) + nn₁ ∈ I for each i ∈ I, n, n₁ ∈ N.*

**DEFINITION 8.1.15**: *Let N be a NA-seminear-ring. A nonempty subset I of N is called an Smarandache ideal (S-ideal) in N if*

     *i.    I is a S-left ideal.*
    *ii.   IN ⊂ I.*

**DEFINITION 8.1.16**: *A S-NA-seminear-ring N is called Smarandache left bipotent (S-left bipotent) if Na = Na² for every a in N.*

**THEOREM 8.1.1**: *Every S-NA seminear-ring N has a nontrivial seminear-ring P if and only if (P, +) is a semigroup and (P, .) is a semigroup.*

*Proof*: If N has P to be a seminear-ring where P ⊂ N then (P, +) and (P, .) are semigroups. (N, .) is a S-groupoid which is always possible. Converse is straightforward.

**THEOREM 8.1.2**: *Every SNA seminear-ring N has a subset P which is a near-ring if and only if (P, +) is a group and (P, .) is a semigroup.*

*Proof*: This is also possible when N is a SNA-seminear-ring as (N, +) is a S-semigroup, so has a subgroup with respect + and (N, .) is a S-groupoid so has a subset which is a semigroup. If the same subset N happens to be group under '+' and semigroup '.' then we have the above theorem to be true.

It is pertinent to mention that only this definition of SNA seminear-ring has paved way for such possibilities so we can define S-seminear-ring I of type A as follows.

**DEFINITION 8.1.17**: *Let (N, +, .) be a NA seminear-ring. N is said to be a Smarandache seminear-ring I of type A (S-seminear-ring I of type A) if N has a proper subset P such that (P, +, .) is an associative seminear-ring.*

It is very important to note by no means we claim that the two definition SNA-seminear-ring I and SNA seminear-ring I of type A are equivalent. But of course it is an interesting problem to find under what conditions the two definitions are equivalent which is enlisted as suggested problems in chapter ten.

**DEFINITION 8.1.18**: *Let N be a NA-seminear-ring. N is said to be a Smarandache NA seminear-ring I of type B (S-NA seminear-ring I of type B) if N has a proper subset P where P is a near-ring.*

We say in the definition of S-NA seminear-ring I we assumed (N, +) is a S-semigroup so we have (P, +) for a suitable proper subset P to be a group and since (N, .) is a S-groupoid we have (Q, .) is a semigroup, if P = Q we see N has P to be a near-ring. This has forced us to define S-NA seminear-ring I of type B. It is pertinent to mention here that by no means we assume or state that the definitions S-NA seminear-ring I, S-NA seminear-ring I type B are equivalent. Thus this definitions give some what conditionally equivalent conditions or the Smarandache structures gives way for such concepts i.e. a non-associative algebraic structure containing completely an associative structure or more richer algebraic structure.

***Example 8.1.4***: Let (N, +, .) be a seminear-ring and G be a groupoid. The groupoid seminear-ring NG is a non-associative seminear-ring. NG is S-NA-seminear-ring if and only if N is a S-seminear-ring.

We ask the reader to obtain interesting results about S-non-associative seminear-rings.

**DEFINITION 8.1.19**: *Let N and $N_1$ be two S-seminear-rings; we say a map $\phi$ from N to $N_1$ is a Smarandache non-associative seminear ring homomorphism (S-NA seminear-ring homomorphism) if*

$$\phi(x + y) = \phi(x) + \phi(y)$$
$$\phi(xy) = \phi(x)\,\phi(y) \text{ for all } x, y \in N.$$

**PROBLEMS:**

1. What is the smallest order of a non-associative seminear-ring?
2. Find the order of the smallest S-non-associative seminear-ring.

3. Let N = $Z_9$ define '×' and 'o' on $Z_9$ as usual multiplication modulo p and a.b = a for all a, b ∈ N. Let G be the groupoid given by the following table:

| * | a | b | c | d |
|---|---|---|---|---|
| a | a | b | a | b |
| b | c | d | c | d |
| c | b | a | b | a |
| d | d | c | d | c |

    a. Is $Z_9$G the group seminear-ring a S-NA seminear-ring.

    b. Find S-subseminear-ring of $Z_9$G.

    c. Does $Z_9$G have S-ideals?

4. Let $Z^+$ be the seminear-ring and G a groupoid given in problem 3. Is the groupoid seminear-ring $Z^+$G have a S-seminear-ring? Does $Z^+$G have strictly prime ideals?

## 8.2 Some special Smarandache NA seminear-rings of type II

In this section we introduce Smarandache NA seminear-ring of type II and define some interesting results about them. We see that S-NA seminear-rings of type I happen to have a very meager property for we in the definition of S-NA seminear-rings of type I replace the groupoid by S-groupoids this does not make a big change except one of the obvious properties that these non-associative seminear-rings contain associative seminear-rings as substructures. Except for this we were not able to achieve many more properties we proceed on to define Smarandache NA seminear-rings.

**DEFINITION 8.2.1**: *Let (N, +, .) be a NA seminear-ring, we say N is a Smarandache NA seminear-ring II (S-NA-seminear-ring II) if N has a proper subset P which is a associative seminear-ring.*

**THEOREM 8.2.1**: *Let N be a S-NA seminear-ring II then N is a S-NA-seminear-ring I of type A.*

*Proof*: Obvious by the very definition. Now using the S-mixed direct product definition of seminear-rings we can extend it to the case of non-associative seminear-rings. This method will help to build a class of S-NA seminear-rings of type II.

**DEFINITION 8.2.2**: *Let $N_1$, ..., $N_k$ be k seminear-rings where at least one of the seminear-rings is non-associative and at least one is associative. The direct product of these seminear-rings $N = N_1 \times ... \times N_k$ is called the S-mixed direct product of seminear-rings (S-mixed direct product of seminear-rings). The operation in N is carried out component wise. Clearly N is a S-seminear-ring II and S-seminear-ring I of type A.*

Now we define a still a strong mixed product.

**DEFINITION 8.2.3**: *Let $N_1$, $N_2$, ..., $N_k$ be a collection of near-rings, non-associative seminear-rings and seminear-rings. The product $N = N_1 \times N_2 \times ... \times N_k$ is called the Smarandache strong mixed direct product (S-strong mixed direct product) of NA seminear-rings. Under the component wise operation on N we get N to be a NA seminear-ring. Clearly N is a S-NA seminear-ring I of type A and B and also S-NA seminear-ring II.*

Since these S-NA seminear-rings are non-associative structures we would be justified in studying the cases when these S.NA seminear-rings satisfy the special identities like Bol, Moufang, alternative etc.

As we do not have a very large class of S.NA-seminear-rings or for that matter NA seminear-rings we are forced to build S.NA seminear-rings and NA seminear-rings only by using loops and groupoids and these newly constructed structure will be called as loop seminear-rings and groupoid seminear-rings. These will be non-associative seminear-rings even if we take the seminear-rings to be associative. We study when these new class of NA seminear-rings satisfies these identities.

To the best of the authors knowledge no researcher has studied such types of identities in these structures. Further in view of the author the very study of non-associative near-rings and seminear-rings is very meager. It is all the more important to state especially these non-associative seminear-rings do not find any place in any text books on near-rings. Even the very concept of seminear-ring is dealt very sparingly in some of the text books. Further as contrary to the case of ring theory where there are several books, in case of near-rings there are only less than a dozen text books.

Thus we have tried to give the definitions as well, we have mentioned them in the bibliography. Now we proceed on to define the Smarandache substructures and Smarandache special elements in these NA seminear-rings. Throughout this section we only assume S-seminear-ring of level II or S-seminear-ring of level I type B.

**DEFINITION 8.2.4**: *Let N be a S-NA seminear-ring II. An additive S-semigroup A of N is said to be a S-left N subsemigroup (right-N subsemigroup) if $PA \subset A$ and ($AP \subset A$). where P is a proper subset of N and P is a seminear-ring which is associative. $PA = \{pa \mid p \in P$ and $a \in A\}$.*

**DEFINITION 8.2.5**: *Let N and $N_1$ be any two S-NA-seminear-rings A mapping $\phi$ from N to $N_1$ is called a Smarandache-NA seminear-ring homomorphism (S-NA seminear-ring homomorphism) if $\phi$ maps every $p \in P \subset N$ (p a associative seminear-ring-associative) into a unique element $\phi(p) \in P_1 \subset N_1$ where $P_1$ is an associative seminear-ring such that $\phi(p + p^I) = \phi(p) + \phi(p^I)$ and $\phi(p_1 \, p_2) = \phi(p_1) \, \phi(p_2)$ for every $p_1, p_2 \in P \subseteq N$.*

It is important to note that $\phi$ need not be defined on the whole of N it is sufficient if $\phi$ is defined on a subset P of N where P is an associative seminear-ring. Thus if $P \subset N$ and $P_1 \subset N_1$ and $\phi : N \to N_1$ is such that $\phi$ is one to one and on to from P to $P_1$ the two S NA seminear-rings would become isomorphic even if they are not having same number of elements in them.

**DEFINITION 8.2.6**: *Let N be a S-NA-seminear-ring. An additive S-semigroup A of N is said to be a Smarandache left ideal (S-left ideal) of N if it is an ideal of the semigroup (N, +) with the conditions.*

$$p_1 (p_2 + a) - p_1 p_2 \in A \text{ for each } a \in A; p_1, p_2 \in P \subset N; P \text{ a seminear-ring.}$$

**DEFINITION 8.2.7**: *A subset I of N is called a S-ideal if it is a S-left ideal and IP ⊂ I where P ⊂ N is the associative seminear-ring.*

**DEFINITION 8.2.8**: *A S-NA seminear-ring N has Smarandache IFP (S-IFP) (insertion of factor property) if for a, b ∈ P, ab = 0 implies a.p.B = 0 for all p ∈ P ⊂ N since N is non-associative we have to restrict our selves only to the associative substructure to define IFP property as a (nb) ≠ (an) b in general for all a, n, b ∈ N.*

*Now for S-NA seminear-ring of level I type B we define all the above concepts with the only modification we replace the seminear-ring P by the near-ring in all the definitions.*

Certainly we cannot say S-ideals in S-NA seminear-ring of level II will be the same as S-ideals of SNA seminear-ring I of type B even if N is both a S-NA seminear-ring II and S-NA-seminear-ring I type B. So all results will be different. The reader can verify them. We give only example.

***Example 8.2.1***: Let $N = N_1 \times N_2 \times N_3$ where $N_1$ is any non-associative seminear-ring, $N_2 = Z_{24}$ where '×' is the usual multiplication modulo 24 and '.' is a . b = a for all a, b ∈ $Z_{24}$ is an associative seminear-ring and $N_3 = Z_{12}$ where '+' is the operation i.e. addition modulo 12 and a.b = a for all a, b ∈ $N_3$; clearly $N_3$ is an associative near-ring. N by the very definition is a non-associative seminear-ring. Now N is a S-NA seminear-ring II for it has $N_2$ as a proper subset which is an associative seminear-ring. N is also a S-NA seminear-ring I type B as $N_3$ is a proper subset which is a near ring of N. Thus we see simultaneously N is S-NA seminear-ring II and S-NA seminear-ring I type B. It is left for the reader to define S-ideals in them and see how they are different sets behaving in a distinct way, in case of S-mixed direct product of seminear-rings II or S-NA seminear-ring I of type B or S-NA seminear-ring I of type A. This is clearly evident from the example.

***Example 8.2.2***: Let $N = Z \times Z_{15} \times N_1$ where Z is a near-ring in which (Z, +) is a group and for all x, y ∈ Z; x.y = x, ($Z_{15}$, '×', '.') is a seminear ring under usual multiplication modulo 15 and '.' as a.b = a for all a, b ∈ $Z_{15}$ and $N_1$ some non-associative seminear-ring. N is a S-NA seminear-ring II and S-NA seminear-ring I of type B.

From these examples we are not still in a position to find natural examples of NA seminear-ring, we have a class of associative seminear-ring given by $Z_n$ where ($Z_n$, x) is a semigroup under multiplication modulo n and a.b = a for all a, b ∈ $Z_n$. Thus the section 4 is completely devoted to finding a class of NA seminear-rings.

**PROBLEMS:**

1. What is the order of the smallest S-NA seminear-ring?

2. Prove the groupoid near-ring $Z_2G$, G given by the following table is a S-NA seminear-ring.

| . | 0 | 1 | 2 | 3 | 4 | 5 |
|---|---|---|---|---|---|---|
| 0 | 0 | 2 | 4 | 0 | 2 | 4 |
| 1 | 2 | 4 | 0 | 2 | 4 | 0 |
| 2 | 4 | 0 | 2 | 4 | 0 | 2 |
| 3 | 0 | 2 | 4 | 0 | 2 | 4 |
| 4 | 2 | 4 | 0 | 2 | 4 | 0 |
| 5 | 4 | 0 | 2 | 4 | 0 | 2 |

3. Let $Z_2 = \{0, 1\}$ be the near ring. Let L be loop given by the following table:

| . | e | $a_1$ | $a_2$ | $a_3$ | $a_4$ | $a_5$ |
|---|---|---|---|---|---|---|
| e | e | $a_1$ | $a_2$ | $a_3$ | $a_4$ | $a_5$ |
| $a_1$ | $a_1$ | e | $a_3$ | $a_5$ | $a_2$ | $a_4$ |
| $a_2$ | $a_2$ | $a_5$ | e | $a_4$ | $a_1$ | $a_3$ |
| $a_3$ | $a_3$ | $a_4$ | $a_1$ | e | $a_5$ | $a_2$ |
| $a_4$ | $a_4$ | $a_3$ | $a_5$ | $a_2$ | e | $a_1$ |
| $a_5$ | $a_5$ | $a_2$ | $a_4$ | $a_1$ | $a_3$ | e |

Is the loop near-ring $Z_2L$ a S-NA seminear-ring I? Is $Z_2L$ at least a non-associative seminear-ring? Justify your answer.

4. Find S-ideals & S-N-subsemigroups for seminear-ring $Z_2G$ given in problem 3.

5. Is $Z_2L$ a S-NA seminear-ring I type B? Justify your claim.

## 8.3 Smarandache Non-Associative near-rings

In this section we introduce the concept of Smarandache non-associative ring and also study about the near loop ring introduced in chapter 3 to be S-near-ring. Several interesting Smarandache concepts are introduced. Here also we once again mention loop near-rings and near loop rings are different for the former has '+' to be non-associative where as in the later '+' is associative but '.' happens to be non-associative and the near loop rings are built using a loop and a near-ring. Finally we also introduce the identities newly to the near-ring which are non-associative. We introduce the concept of Smarandache right loop-half groupoid near-ring which is the most generalized concept of loop near-ring.

**DEFINITION 8.3.1**: *The system N = (N, '+', '.' 0) is called a Smarandache right loop half groupoid near-ring (S-right loop half groupoid near-ring) provided.*

*1. (N, +, 0) is a Smarandache loop.*

*2. (N, '.') is a half groupoid.*

*3. $(n_1.n_2).n_3 = n_1.(n_2.n_3)$ for all $n_1, n_2\ n_3 \in N$ for which $n_1.n_2, n_2.n_3, (n_1.n_2).n_3$ and $n_1.(n_2.n_3) \in N$.*

4. $(n_1 + n_2) \, n_3 = n_1.n_3 + n_2.n_3$ *for all* $n_1, n_2, n_3 \in N$ *for which* $(n_1 + n_2).n_3,$ $n_1.n_3 \, , \, n_2.n_3 \in N$. *If instead of (4) in N the identity* $n_1.(n_2 + n_3) = n_1.n_2 + n_1.n_3$ *is satisfied then we say N is a Smarandache left half groupoid near-ring (S-left half groupoid near-ring).*

We just say (L, +) is a S-loop if L has a proper subset P such that (P, +) is an additive group.

**DEFINITION 8.3.2**: *A Smarandache right loop near-ring (S-right loop near-ring) N is a system (N, +, .) of double composition '+' and '.' such that*

1. *(N, +) is a S-loop.*
2. *(N, .) is a S-semigroup.*
3. *The multiplication '.' is right distributive over addition i.e. for all* $n_1, n_2, n_3 \in N$ $(n_1 + n_2).n_3 = n_1.n_3 + n_2.n_3$.

**THEOREM 8.3.1**: *Let N be a S-right loop near-ring then N is a S-right loop half groupoid near-ring.*

*Proof*: Obvious by the very definitions.

**THEOREM 8.3.2**: *Let N be a S-right loop half groupoid near-ring, then N is not in general a S-right loop near-ring.*

*Proof*: Obvious. (N, .) is only a half groupoid so it can never be a S-semigroup.

***Example 8.3.1***: Every S-near-ring is a S-loop near-ring.

**THEOREM 8.3.3**: *Let N be a S-loop near-ring. Then N contains a proper subset P such that (P, +) is a group and (P, .) is a semigroup so that P is a near-ring.*

*Proof*: By the very definition of S-loop near-ring we have (N, +) is a S-loop so N has a proper subset P which is a group under '+' Now (P, +, .) is near-ring.

**DEFINITION 8.3.3**: *Let (G, +, 0) be a S-loop and* $\Delta$ *be a S-groupoid of G. A set of all endomorphism of G is called a Smarandache* $\Delta$ *centralizer (S-* $\Delta$ *centralizer) of G provided.*

1. *The zero endomorphism* $\hat{o} \in S$.
2. *S /* $\hat{o}$ *(complement of* $\hat{o}$ *in S) is a group of automorphisms of G.*
3. $\phi(\Delta) \subset \Delta$ *for all* $\phi \in S$.
4. $\phi, \psi \in S$ *and* $\phi(\omega) = \psi(\omega)$ *for some* $\hat{o} \neq \omega \in \Delta$ *imply* $\phi = \psi$.

**DEFINITION 8.3.4**: *Let (G, +, 0) be a S-loop. Let* $\Delta$ *be a subset of G which is a S-groupoid of G and S a S-* $\Delta$ *centralizer of G. A mapping T of G into itself is called a Smarandache* $\Delta$ *- transformation (S-* $\Delta$ *- transformation) of G over S provided T [* $\phi$ *($\omega$)] = $\phi$[T ($\omega$)] for all* $\omega \in \Delta$ *and* $\phi \in S$. *It is interesting to note that if* $\overline{0} \in \Delta$ *and T is a S -* $\Delta$ *transformation of G over S then T fixes 0 i.e. T ($\overline{0}$) = ($\overline{0}$). We shall denote the* *set of all S - Δ transformations of G over S by S(N (S, Δ)). Further we see for any endomorphism φ of a S-loop G, [φ (g)]$_r$ = φ (g$_r$) for all g ∈ G.*

The reader is requested to prove the following theorem:

**THEOREM 8.3.4**: *Let (G, +, 0) be a S-loop. Δ a subgroup of G which is a groupoid containing $\overline{0}$ and S a S-Δ-centralizer of G. Then S(N (S, Δ)) is not a S-loop near-ring or even a loop near-ring (The operation on S (N (S, Δ)) being usual addition and composition of mappings).*

**DEFINITION 8.3.5**: *A nonempty subset M of a S-loop (N, +, '.', 0) is said to be a Smarandache subloop near-ring (S-subloop near-ring) of N if and only if (M, '+', '.' 0) is a S-loop near-ring.*

**DEFINITION 8.3.6**: *A S-loop near-ring N is said to be Smarandache zero symmetric (S-zero symmetric) if and only if n0 = 0 for every n ∈ P ⊂ N where (P, +) is a group, here '0' is the additive identity. S-zero symmetric loop near-ring will be denoted by S(N$_o$).*

**DEFINITION 8.3.7**: *A S-loop near-ring N is said to be a S-near-ring if N contains a proper subset P such that (P, +, .) is a near-field.*

**DEFINITION 8.3.8**: *An element a of a S-loop near-ring N is said to be a Smarandache left (right) zero divisor (S-left (right) zero divisor) in N if there is an element b ≠ 0 in N with a.b = 0   (b . a = 0) and there exist x, y ∈ S \ {a, b, 0}, x ≠ y with*

> 1.  a . x = 0 (xa = 0).
> 2.  by = 0  (yb = 0).
> 3.  xy ≠ 0 (yx ≠ 0).

*An element which is both a S-right and S-left zero divisor in N is called a Smarandache two sided zero divisor or simply a Smarandache zero divisor (S-zero divisor) in N.*

**DEFINITION 8.3.9**: *A map θ form a S-loop near-ring N into a S-loop near-ring N$_1$ is called a Smarandache homomorphism (S-homomorphism) if*

> θ (n$_1$ + n$_2$) = θ (n$_1$) + θ (n$_2$)  and
> θ (n$_1$ n$_2$) = θ (n$_1$) θ (n$_2$)

*for every n$_1$, n$_2$ ∈ P ⊂ N where (P, +, .) is a near-ring and φ (n$_1$), φ (n$_2$) ∈ P$_1$ where (P$_1$, +, .) is a near-ring of N$_1$.*

**DEFINITION 8.3.10**: *Let N be a S-loop near-ring. An additive S-subloop A of N is called a N-subloop (right N-subloop) if NA ⊂ A (AN ⊂ A) where NA = {na / n ∈ N and a ∈ A}.*

**DEFINITION 8.3.11**: *A nonempty subset I of a S-near-ring N is called a Smarandache left ideal (S-left ideal) in N if*

1. *(I, +) is a S-normal subloop of (N, +).*
2. *n (n₁ + i) + n_r n₁ ∈ I for each i ∈ I and n, n₁ ∈ N where n_r denotes the unique right inverse of n.*

**DEFINITION 8.3.12**: *A nonempty subset I of N is called a Smarandache ideal (S-ideal) of N if*

1. *I is a S-left ideal.*
2. *IN ⊂ I.*

**DEFINITION 8.3.13**: *A S-loop near-ring N is said to be Smarandache left bipotent (S-left bipotent) if $Na = Na^2$ for every a in N.*

***Example 8.3.2***: Let N = {e, a, b, c, d} be given by the following table for addition '+'.

| + | e | a | b | c | d |
|---|---|---|---|---|---|
| e | e | a | b | c | d |
| a | a | b | e | d | c |
| b | b | c | d | a | e |
| c | c | d | a | e | b |
| d | d | e | c | b | a |

Let '.' be defined on N by x.y = x for every x, y ∈ N. (N, +, .) is a S-loop near-ring. For P = {e, c} is a S-loop. Clearly $Na^2$ = Na for every a ∈ N. Thus N is a S-left bipotent loop near-ring.

**DEFINITION 8.3.14**: *A S-loop near-ring N is said to be a Smarandache-s-loop near-ring (S-s-loop near-ring) if a ∈ Na for every a in N.*

***Example 8.3.3***: Let the S-loop near-ring (N, +, .) be given by the following table for '+' and '.'.

| + | 0 | a | b | c |
|---|---|---|---|---|
| 0 | 0 | a | b | c |
| a | a | 0 | c | b |
| b | b | c | 0 | a |
| c | c | b | a | 0 |

and

| . | 0 | a | b | c |
|---|---|---|---|---|
| 0 | 0 | 0 | 0 | 0 |
| a | 0 | a | b | c |
| b | 0 | b | 0 | 0 |
| c | 0 | c | b | c |

Every pair (0, a) is a group so (N, +) is a S-loop. (N, .) is a S-semigroup as {a} is a group. Thus (N, +, .) is a S-s-loop near-ring.

**DEFINITION 8.3.15**: *A S-loop near-ring N is said to be Smarandache regular (S-regular) if for each a in N there exists x in P (P ⊂ N), P a S-loop such that a = axa.*

**DEFINITION 8.3.16**: *Let N be a near-ring. An element e ∈ N is said to be Smarandache idempotent (S-idempotent) in N if*

1. *$e^2 = e$.*
2. *There exists $b ∈ N \setminus \{e\}$ such that $b^2 = e$ and eb = b (be = b) or (be = e or eb = e).*

**THEOREM 8.3.5**: *All S-idempotents in S-loop near-rings are also idempotents and non-conversely.*

*Proof*: Left for the reader to prove.

**THEOREM 8.3.6**: *Let N be a S-s-loop near-ring, then N is S-regular if and only if for each a (≠0) in N there exists a S-idempotent e such that Na = Ne.*

*Proof*: Left for the reader to prove.

**DEFINITION 8.3.17**: *A S-loop near-ring N is said to be Smarandache strictly duo (S-strictly duo) if every S-N-subloop (S-left ideal) is also a S-right N-subloop (S-right ideal).*

**DEFINITION 8.3.18**: *For any subset A of a S-loop near-ring N we define, $\sqrt{A} = \{x ∈ N / x^n ∈ A \text{ for some n}\}$.*

**DEFINITION 8.3.19**: *A S-loop near-ring N is called Smarandache irreducible (S-irreducible) (S-simple) if it contains only the trivial S-N-subloops (S-ideals) (0) and N itself.*

**DEFINITION 8.3.20**: *A loop near-ring N is called a S-loop near-field.*

1. *if N is a S-loop near-ring.*
2. *if P ⊊ N contains an identity and each non-zero element in P ⊂ N, {(P, +) a group} has a multiplicative inverse.*

**DEFINITION 8.3.21**: *An S-ideal P (≠N) is called Smarandache strictly prime (S-strictly prime) if for any two S-N-subloops (S-ideals) A and B of N such that AB ⊂ P then A ⊂ P or B ⊂ P.*

**DEFINITION 8.3.22**: *As left ideal B of a S-loop near-ring N is called Smarandache strictly essential (S-strictly essentail) if B ∩ K ≠ {0} for any non zero S-N-subloop K of N.*

*An element x in N is said to be Smarandache non-singular (S-non-singular) if there exists a S-strictly essential left ideal A in N such that Ax = {0}.*

**DEFINITION 8.3.23**: *An element x of a S-loop near-ring N is said to be Smarandache central (S-central) if xy = yx for all y ∈ P ⊂ N, (P, +, .) is a near-ring.*

**DEFINITION 8.3.24**: *A nonzero S-loop near-ring N is said to be Smarandache subdirectly irreducible (S-subdirectly irreducible) if the intersection of all the nonzero S-ideals of N is non-zero.*

Now we study when are loops over near-rings i.e. near loop rings interminology Smarandache near loop rings, as loop near-rings are not non associative near-ring analogous to non associative ring; we proceed to define here the concept of Smarandache non associative near-rings first in the following:

**DEFINITION 8.3.25**: *Let (N, +, .) be a nonempty set. N is said to be a Smarandache non-associative near-ring (S-NA-near-rings) if the following conditions are satisfied.*

    *1.  (N, +) is a S-semigroup.*
    *2.  (N, .) is a S-groupoid.*
    *3.  (a + b) c = ac + bc for all a, b, c ∈ N.*

If (N, +, .) has an element 1 such that a.1 = 1.a = a for all a ∈ N, we call N a S-non associative near-ring with unit. If (N, .) is a S-loop; barring 0 we call N a Smarandache non-associative division near-ring (S-NA-division near-ring).

If N \ {0} is commutive S-loop under the operation '.' we call N a Smarandache non-associative near-field (S-NA near-field).

**DEFINITION 8.3.26**: *Let N be a near-ring; L a loop the near loop ring NL is a Smarandache near loop ring (S-near loop ring) if*

    *1.  N is a S-near-ring.*
    *2.  L is a S-loop.*

**THEOREM 8.3.7**: *Let N be a S-near ring and L a S-Moufang loop. Then the near loop ring NL is a S-moufang near loop ring.*

*Proof*: Left for the reader as an exercise.

***Example 8.3.4***: Let $Z_2 = \{0, 1\}$ be a near-ring and L be a S-loop.

The near loop ring $Z_2L$ is a S-near-ring where the S-loop is given by the following table:

| .     | e     | $a_1$ | $a_2$ | $a_3$ | $a_4$ | $a_5$ |
|-------|-------|-------|-------|-------|-------|-------|
| e     | e     | $a_1$ | $a_2$ | $a_3$ | $a_4$ | $a_5$ |
| $a_1$ | $a_1$ | e     | $a_4$ | $a_2$ | $a_5$ | $a_3$ |
| $a_2$ | $a_2$ | $a_4$ | e     | $a_5$ | $a_3$ | $a_1$ |
| $a_3$ | $a_3$ | $a_2$ | $a_5$ | e     | $a_1$ | $a_4$ |
| $a_4$ | $a_4$ | $a_5$ | $a_3$ | $a_1$ | e     | $a_2$ |
| $a_5$ | $a_5$ | $a_3$ | $a_1$ | $a_4$ | $a_2$ | e     |

Now $Z_2L = \{1, a_1, a_2, a_3, a_4, a_5, \text{sums taken 2 at a time } \ldots\}$ where $1.e = e.1 = 1$ is the assumption made so that 1 acts as the identity. When we say $Z_2L$ is a near loop ring we assume $L \subset Z_2L$. Clearly $Z_2L$ is a non-associative near-ring and it is a S-near-ring.

**THEOREM 8.3.8**: *Let L be a loop and N be a S-near-ring. NL is a S-near-ring which is non-associative.*

*Proof*: Follows from the fact $N \subset NL$ and $L \subset NL$ and as L is non-associative and N is a S-near-ring NL is a S-near-ring.

**THEOREM 8.3.9**: *Let L be any loop and N any near-ring. The near loop ring NL is a S-non-associative near-ring.*

*Proof*: Follows from the fact (NL, +) is a group; hence trivially a S-semigroup and (NL, .) is a S-groupoid as $N \subset NL$ is a semigroup. Hence the claim.

Now we will study when the class of S-non-associative near-rings are Moufang, Bol, etc.

**THEOREM 8.3.10**: *Let N be a near-ring and L a Moufang loop then the near loop ring is a S-Moufang near-ring.*

*Proof*: Follows from the fact (N, .) is a associative semigroup and L is a Moufang loop.

<u>*Corollary*</u>: Let N be near-ring and L not a Moufang loop; N the near loop ring is not a S-Moufang near-ring.

Thus we can say if we take any near-ring N with 1 and L to be a special type of loop say Bol, alternative … then NL will be correspondingly the special type of near-ring Bol, alternative, …

But even if L is a Bruck loop the near loop ring NL will not be a Bruck loop. NL will be not a Bruck near-ring even if N is a near-field and L a Bruck loop.

Now we have introduced in chapter one a new class of loops $L_n$ which will be a class of non-associative S-loops which are alternative, power associative and their near loop rings are S-near loop rings if the near-ring is a S-near-ring.

**THEOREM 8.3.11**: *Let $Z_2L$ be a near loop near of the S-loop L over $Z_2$. Let L be a power associative loop. Then the near loop ring $Z_2L$ in general is not a power associative near-ring.*

*Proof*: Every element in $Z_2L$ need not generate a subgroup.

**DEFINITION 8.3.27**: *Let N be a non-associative loop near-ring. We say N is a additively power associative loop near-ring if every element under '+' of N generates* *a group. The near-ring N is said to be a multiplicatively power associative loop near-ring if every element different from 0 generates a group under product. A loop near-ring N is said to be power associative near-ring if*

> 1. *Every element $n \in N$ under '+' generates a group.*
> 2. *Every element $n \in N \setminus \{0\}$ under '×' or '.' generates a group.*

**THEOREM 8.3.12**: *A loop near-ring N is a power associative near-ring then*

> 1. *N is an additively power associative near-ring.*
> 2. *N is multiplicatively a power associative near-ring.*

*Proof*: Straightforward by the very definitions.

Now we proceed on to define Smarandache power associative near-ring Smarandache additively power associative near-ring and Smarandache multiplicatively power associative near-ring. Before we go for Smarandache notions it is important to observe that if N is not a loop near-ring but a non-associative near-ring then the concept "additively power associative" has no meaning. So while defining the concept of power associative for non-associative near-ring we define only just power associative; likewise even in case of near loop rings we define only power associative. Only for loop near-ring when we talk of Smarandache concepts we may talk of the notion of additively power associative.

**DEFINITION 8.3.28**: *Let N be a non-associative near-ring. We say N is power associative if every element of $N \setminus \{0\}$ generates a group with respect to.*

**DEFINITION 8.3.29**: *Let N be a non-associative seminear-ring, we say N is a power associative seminear-ring if every element in N generates a semigroup under multiplication.*

**DEFINITION 8.3.30**: *Let N be a non-associative seminear-ring we say N is a Smarandache power associative (S-power asoociative) if (N, .) has a proper subset P where P is a S-subgroupoid of N and P is power associative i.e. every element $p \in P$ generates a semigroup.*

It is important to note that even if (N, '+', '.') is a S-seminear-ring which is still non-associative it may be the case (N, .) has no proper S-subgroupoids so the non-associative seminear-ring may fail to become S-power associative.

**DEFINITION 8.3.31**: *Let (N, +, .) be a S-non-associative near-ring we say N is Smarandache power associative (S-power associative) if (N, .) has a proper subset which has a S-subgroupoid P such that every element $p \in P$ is power associative i.e., every element in p generates a semigroup.*

Now the four definitions of power associative for non-associative near-rings, non-associative seminear-rings S-NA-near-rings and S-NA-seminear-rings can be made on similar lines for "diassociative".

Thus we will say a non-associative near-ring (seminear-ring) is diassociative if every pair of elements generates a group (semigroup). This notion is modified to S-NA-near-rings and S-NA seminear-rings in an analogous way.

We stop at this stage and propose problems for the reader to solve and obtain class of such near-rings as this is the first time such type of study has been carried out for S-near-ring and S-seminear-ring which are non-associative. We give examples and a class of Smarandache right alternative near-ring (S-right alternative near-rings) (left alternative near-rings) and Smarandache WIP near-rings (S-WIP near-rings).

**DEFINITION 8.3.32**: *Let (N, +, .) be a non-associative near-ring (or a seminear-ring) N is said to be Smarandache NA WIP-seminear-ring (S-NA-WIP seminear-ring); if N has S-subgroup P (S-subgroupoid) such that every set of x, y, z ∈ P satisfies (xy) z = e imply x (yz) = e.*

Similar Smarandache notions can be defined for non-associative near-rings and non-associative seminear-rings like right (left) alternative, Bol, moufang etc.

**PROBLEMS:**

1. What is the order of the smallest loop near-ring?
2. Can a near loop ring be of order less than $2^5$?
3. Give examples of S-non-associative seminear-ring which is S-Moufang seminear-ring.
4. Let $Z_2 = \{0, 1\}$ be the near-ring. L be the loop given by the following table:

| . | e | $a_1$ | $a_2$ | $a_3$ | $a_4$ | $a_5$ |
|---|---|---|---|---|---|---|
| e | e | $a_1$ | $a_2$ | $a_3$ | $a_4$ | $a_5$ |
| $a_1$ | $a_1$ | e | $a_3$ | $a_5$ | $a_2$ | $a_4$ |
| $a_2$ | $a_2$ | $a_5$ | e | $a_4$ | $a_1$ | $a_3$ |
| $a_3$ | $a_3$ | $a_4$ | $a_1$ | e | $a_5$ | $a_2$ |
| $a_4$ | $a_4$ | $a_3$ | $a_5$ | $a_2$ | e | $a_1$ |
| $a_5$ | $a_5$ | $a_2$ | $a_4$ | $a_1$ | $a_3$ | e |

What are the Smarandache properties enjoyed by this near loop ring $Z_2L$?

5. Study the Smarandache properties satisfied by the groupoid seminear-ring $Z_2G$ where G is the groupoid given by the following table:

| . | e | $a_0$ | $a_1$ | $a_2$ | $a_3$ | $a_4$ | $a_5$ |
|---|---|---|---|---|---|---|---|
| e | e | $a_0$ | $a_1$ | $a_2$ | $a_3$ | $a_4$ | $a_5$ |
| $a_0$ | $a_0$ | e | $a_5$ | $a_4$ | $a_3$ | $a_2$ | $a_1$ |
| $a_1$ | $a_1$ | $a_4$ | e | $a_2$ | $a_1$ | $a_0$ | $a_5$ |
| $a_2$ | $a_2$ | $a_2$ | $a_1$ | e | $a_5$ | $a_4$ | $a_3$ |
| $a_3$ | $a_3$ | $a_0$ | $a_5$ | $a_3$ | e | $a_2$ | $a_1$ |
| $a_4$ | $a_4$ | $a_4$ | $a_3$ | $a_2$ | $a_1$ | e | $a_5$ |
| $a_5$ | $a_5$ | $a_2$ | $a_1$ | $a_0$ | $a_5$ | $a_4$ | e |

### 8.4 Smarandache loop seminear-rings and Smarandache groupoid seminear-rings.

In this section we introduce two new classes of non-associative seminear-rings which is define as Smarandache loop seminear-rings and Smarandache groupoid seminear-rings. We know the class of loop seminear-rings is a generalization of group seminear-rings and loop near-rings. In fact groupoid seminear-rings happens to be most generalized class of non-associative seminear-rings which contains all these classes of seminear-rings. We are more interested in the study of their Smarandache nature. The concept of classical identities plays a vital role when the structure under study happens to be non-associative. With the advent of a new class of loops of even order and new class of groupoids of all order built using integers we are able to give examples which are concrete.

One of the major drawbacks in the study of near-rings or any algebraic structure is the lack of concrete examples. Unless we have concrete examples it becomes difficult to make the subject or that algebraic structure profoundly attractive to many researchers. For instance, the research in ring theory or group theory or semigroup theory are carried out by more researchers when compared to the theory of loop, theory of groupoid. So in this section we introduce the Smarandache analogue of these concepts and try our level best to illustrate them whenever possible with examples.

**DEFINITION 8.4.1**: *Let (N, +, .) be a non-empty set, N endowed with two binary operation '+' and '.' satisfying the following conditions.*

1. *(N, +) is a semigroup.*
2. *(N, .) is a groupoid.*
3. *(a + b).c = a.c + b.c for a, b, c ∈ N.*

*(N, +, .) is defined as a non-associative seminear-ring (NA seminear-ring).*

But one of the major draw backs is that we do not have several examples of these seminear-rings; so it has become difficult to give examples using $Z_p$ alone or Z alone or any of the known sets. Thus to obtain examples of non-associative seminear-rings we define loop seminear ring and groupoid seminear rings which has been defined in chapter 3.

Now we proceed on to define Smarandache NA seminear-ring.

**DEFINITION 8.4.2**: *Let N be a non-empty set. Define two binary operations '+' and '.' satisfying the following conditions*

1. *(N, +) is a S-semigroup.*
2. *(N, .) is a S-groupoid.*
3. *(a + b).c = a.c + b.c for all a, b, c ∈ N then we call N a Smarandache NA seminear-ring of level I (S-NA seminear-ring of level I).*

Now we define S-NA seminear-ring of level II.

**DEFINITION 8.4.3**: *Let (N, +, .) be a non-associative seminear-ring we say N is a Smarandache NA seminear-ring of level II (S-NA-seminear-ring of level II). If N has a proper subset P ⊂ N where (P, +, .) is a non-associative near-ring.*

**DEFINITION 8.4.4**: *Let (N, +, .) be a NA near-ring. N is said to be a Smarandache pseudo NA seminear-ring (S-pseudo NA-seminear-ring), if N has a proper subset P such that P is a non-associative seminear-ring under the operations of N.*

Now still we have a level III and level IV definition of Smarandache NA seminear-rings.

**DEFINITION 8.4.5**: *Let N be a NA-seminear-ring. We say N is a Smarandache NA seminear-ring of level III (S-NA-seminear-ring III); if N has a proper subset P where P is an associative seminear-ring.*

Finally we define Smarandache seminear-ring of level IV and obtain the possible and probable relation between the four levels of S-NA seminear-rings.

**DEFINITION 8.4.6**: *Let (N, +, .) be a NA seminear-ring we say N is a Smarandache NA seminear-ring of level IV (S-NA-seminear-ring of level IV) if N has a proper subset P such that (P, +, .) is an associative near-ring.*

**THEOREM 8.4.1**: *Let (N, +, .) be a S-NA seminear-ring IV then (N, +, .) is a S-NA seminear-ring III.*

*Proof*: Follows directly by the very definitions of level III and level IV, S-NA seminear-rings.

It is left as an exercise for the reader to construct a S.NA seminear-ring of level III which is not a S-NA-semiring of level IV. Thus we see the class of S-NA seminear-rings of level IV is strictly contained in the class of S-NA seminear-rings of level III.

**THEOREM 8.4.2**: *Let (N, +, .) be a S.NA seminear-ring IV then (N, +, .) is a S-NA seminear-ring II.*

*Proof*: Obvious from the fact the class of all associative near-rings is contained in the class of non-associative near-rings (as for example we have every group is a loop and not vice versa).

Thus we see the class of SA seminear-ring IV is strictly contained in the class of S-NA seminear-ring II. It is easily seen that a S-NA seminear-ring II can in general never be a S-NA seminear-ring IV.

Now we proceed on to define concepts like S-subseminear-ring S-ideal etc.

**DEFINITION 8.4.7**: *Let (N, +, .) be a NA seminear-ring we say a proper subset S of N is called a Smarandache subseminear-ring I (II or III or IV) (S-subseminear-ring I or II or III or IV), if (S, +, .) is a S-NA seminear-ring I (II or III or IV respectively).*

From this we have the following.

**THEOREM 8.4.3**: *Let (N, +, .) be a NA-seminear-ring, if N has a S-NA subseminear-ring I (II or III or IV) then N is a S-NA seminear-ring I (II or III or IV).*

*Proof*: Follows from the very definition and the proof is straightforward.

**DEFINITION 8.4.8**: *Let (N, +, .) be a NA seminear-ring. A non-empty subset I of N is called a Smarandache left ideal I (S-left ideal I) in N if*

1. *(I, +) is a S-subsemigroup of (N, +).*
2. *$n(n_1 + i) + nn_1 \in I$ for each $i \in I$ and $n, n_1 \in N$.*

*If $IN \subset I$ then I is called the Smarandache ideal I (S-ideal I) of N.*

**DEFINITION 8.4.9**: *Let (N, +, .) be a NA seminear-ring. We call a nonempty subset I of N to be a Smarandache left ideal II (S-left ideal II) of N if*

1. *(I, +) is a S-subsemigroup.*
2. *$n(n_1 + I) + nn_1 \in I$ for all $n_1, n \in P \subset N$ (where P is a non-associative near-ring) and $i \in I$.*

*I is called a Smarandache ideal II (S-ideal II) if $IP \subset I$.*

Now we define Smarandache left ideal III and IV by replacing P in the definition by an associative seminear-ring or by an associative near-ring respectively. The reader is advised to obtain relations between S-ideals I, II, III and IV whenever they exist and obtain some interesting results in this direction.

Now we give examples so that the reader is made to understand the concept in a better way.

***Example 8.4.1***: Let $Z_6 = \{0, 1, 2, \ldots, 5\}$. Define '×' and '.' on $Z_6$ by; '×' is the usual multiplication modulo 6 and '.' is defined by a.b = a for all a, b $\in Z_6$. Let G be a groupoid with 1. The groupoid near-ring $Z_6G$ is a S-NA near-ring III as $Z_6 \subset Z_6G$.

**THEOREM 8.4.4**: *Let N be a seminear-ring and G any groupoid with 1 the groupoid seminear-ring NG is a S-NA seminear-ring III.*

*Proof*: Obvious by the very construction; we see $N \subset Na$ and N is a seminear-ring which is associative so NG is a S-NA seminear-ring III.

***Example 8.4.2***: Let N be any near-ring and G a groupoid with 1. The groupoid near-ring NG is a S-NA seminear-ring IV clearly from the fact $N \subseteq NG$.

**THEOREM 8.4.5**: *Let N be a near-ring and G any groupoid with 1. The groupoid near-ring is a S-NA seminear-ring IV.*

*Proof*: Straightforward by the very definition.

Now we proceed on to define Smarandache bipotent seminear-ring.

**DEFINITION 8.4.10**: *Let N be a seminear-ring. We say N is Smarandache left bipotent II (S-left bipotent II) if $Pa = Pa^2$ where P is a non-associative near-ring in N i.e. $P \subset N$ and $(P, +, .)$ is a near-ring. Smarandache left bipotent III (S-left bipotent III) if $Pa = Pa^2$ where P is an associative seminear-ring. Smarandache left bipotent IV (S-left bipotent IV) if $Pa = Pa^2$ where P is an associative near-ring. Now we have N to be a Smarandache left bipotent I (S-left bipotent I) if N is a S-NA seminear-ring with $Na = Na^2$.*

*Similarly we define Smarandache s-seminear-ring in the following.*

**DEFINITION 8.4.11**: *Let $(N, +, .)$ be a NA seminear-ring we say N is a Smarandache s-seminear-ring I (S-s-seminear-ring I) if $a \in Na$ where N is a S-NA near-ring I.*

*N is said to be Smarandache s-seminear-ring II (S-s-seminear-ring II), if $a \in Pa$ for each $a \in N$ where P is a non-associative near-ring contained in P.*

*N is said to be a Smarandache s-seminear-ring III (IV) (S-s-seminear-ring III (IV)) if $a \in Pa$, P is associative seminear-ring (associative near-ring). We see that if N is a S-s-seminear-ring II III or IV then it is a S-seminear-ring.*

The proof of the next theorem is straightforward and left for the reader to prove.

**THEOREM 8.4.6**: *Let N be a NA seminear-ring. If N is a Smarandache s-seminear-ring II, III or IV then N is a S-NA s-seminear-ring.*

Now for the seminear ring N to be S-s-seminear-ring I we need N to be S-seminear-ring. Hence above theorem has no relevance in case of S-NA seminear-ring I.

Obtain interesting interrelations between S-bipotent and S-s-seminear-rings.

**DEFINITION 8.4.12**: *A seminear-ring N is said to be Smarandache regular I (S-regular I) if for each $a \in N$ there exists a $x \in S$ (where S is a proper subset of N and S is a S-subseminear-ring of N) such that $a = axa$.*

*We define Smarandache regular seminear-ring II (III or IV) (S-regular seminear-ring II (III or IV)) by stating for each $a \in N$ their exist $x \in S$; S a non-associative near-ring (S an associative seminear-ring or S an associative near-ring) such that $a = axa$.*

All results related with this concept can be derived for S-seminear-rings. The following important conclusion is formulated as the theorem.

**THEOREM 8.4.7**: *Let N be a NA-seminear-ring, if N is S-regular II III or IV then clearly N is a S-seminear-ring.*

*Proof*: It is straightforward and hence left for the reader to prove.

**DEFINITION 8.4.13**: *Let N be a NA seminear-ring. N is said to have Smarandache insertion of factors property II, III or IV (S-IFP-II,III or IV for short) if for a, b ∈ N we have a.b = 0 implies anb = 0, for all n ∈ P; where P ⊂ N and P is a non-associative seminear-ring. S-IFP-III if P is an associative seminear-ring and S-IFP-IV if P is an associative near-ring.*

The following theorem is left as an exercise for the reader to prove.

**THEOREM 8.4.8**: *Let N be a NA seminear-ring if N satisfies a S-IFP II (III or IV) then N is a S-NA IFP seminear-ring II (III or IV).*

**DEFINITION 8.4.14**: *Let N and $N^I$ be any two S seminear-rings; we say a map $\theta$ from N to $N^I$ is a S-seminear-ring homomorphism I: If $\theta(x + y) = \theta(x) + \theta(y)$ and $\theta(xy) = \theta(x). \theta(y)$ for all x, y ∈ N. $\theta$ from N to $N^I$ is a S-seminear-ring homomorphism II if $\theta(x + y) = \theta(x) + \theta(y)$, $\theta(xy) = \theta(x)\theta(y)$, where $\theta$ is a near-ring homomorphism from P to $P^I$ where P ⊂ N is a non-associative near-ring and $P^I \subset N^I$ is a non-associative near-ring contained in $N^I$; $\theta$ need not be even defined on whole of N or $N^I$. $\theta$ from N to $N^I$ is said to be S-seminear-ring homomorphism III if $\theta : P \to P^I$ is a seminear-ring homomorphism where P ⊂ N and $P^I \subset N^I$ are associative seminear-rings properly contained in N and $N^I$ respectively. $\theta$ need not be defined entirely on N. $\theta$ from N to $N^I$ is a S-seminear-ring homomorphism IV if $\theta : P \to P^I$ is a near-ring homomorphism from P to $P^I$ where P and $P^I$ are proper subsets of N and $N^I$ respectively and they are associative near-rings.*

*Now S-isomorphism I, II, III and IV can be defined as in the case of other algebraic structures. The uniqueness about S-homomorphism is the map need not be defined on the entire seminear-ring, if it is sectionally well defined on a substructure it is enough. Further we see that two NA seminear-rings can be S-isomorphic even without their cardinality being equal. Yet the Smarandache property always allows us for such identifications or to get isomorphism without being really identical.*

*Thus it is enough if a substructure happens to be identical, Smarandache notions declare them to be Smarandache identical.*

Now we proceed on to define Smarandache N-subsemigroup.

**DEFINITION 8.4.15**: *Let N be a NA seminear-ring we say A ⊂ N, A a S-subsemigroup to be a Smarandache N-(right) subsemigroup I (II or III or IV) if*

1. *NA ⊂ A (AN ⊂ A) where N is a S-seminear-ring I.*
2. *A ⊂ N is said to be S-(right) N-subsemigroup II if PA ⊂ A (AP ⊂ A) (where P ⊂ N, P is a NA near-ring).*
3. *A ⊂ N is said to be a S-(right) N-subsemigroup III if PA ⊂ A (AP ⊂ A) (where P is an associative seminear-ring).*
4. *A ⊂ N is said to be a S-(right) N-subsemigroup IV if PA ⊂ A (AP ⊂ A) if P is an associative near-ring.*

The following theorem is direct.

**THEOREM 8.4.9**: *If N is a seminear-ring having a S-right N-subsemigroup (I, II, III or IV) then N is a S-seminear-ring (I II III or IV) respectively.*

*Proof*: Left as an exercise for the reader to prove.

At this juncture the author leaves it open for the reader to obtain interesting and innovative inter relative results about this concept.

**DEFINITION 8.4.16**: *A non-associative seminear-ring N is said to be Smarandache irreducible (Smarandache simple) (S-irreducible or S-simple) if it contains only the trivial S-N-subsemigroups (S-ideals).*

**DEFINITION 8.4.17**: *A S-left ideal B of a N.A seminear-ring N is called Smarandache essential (S-essential) if B ∩ K ≠ {0} for every non-zero S-N-subsemigroup K of N.*

**DEFINITION 8.4.18**: *An element x in N is said to be Smarandache singular (S-singular) if there exists a nonzero S-strictly essential left ideal A of N such that Ax = {0}.*

**DEFINITION 8.4.19**: *An element x of a NA seminear-ring of N is said to be Smarandache central (S-central) if xy = yx for every y ∈ P; P ⊂ N and P is a non-associative near-ring or P is an associative seminear-ring or P is an associative near-ring.*

In view of this we have the following theorem:

**THEOREM 8.4.10**: *Let N be a NA seminear-ring if x ∈ N is S-central then N is a S-seminear-ring of one of the types.*

*Proof:* Left for the reader to prove as the result is straightforward.

**DEFINITION 8.4.20**: *A nonzero NA seminear –ring N is said to be Smarandache subdirectly irreducible (S-subdirectly irreducible) if the intersection of all nonzero S-ideals of N is non-zero.*

Throughout this section by NL we denote the loop seminear-ring of a loop L over a seminear-ring N, N need not be even non-associative we take N to be only an associative seminear-ring but we see by all means the loop seminear-ring NL happens to be only a non-associative seminear-ring. We take only loops to be from the new class of loops in $L_n$, defined in chapter one.

We define the generalized class of non-associative seminear-ring namely groupoid seminear-rings which we denote by NG; here N is a seminear-ring and G is a groupoid. We will discuss only the groupoids from the new class of groupoids mentioned in the chapter one. Before we enumerate these properties we define S-idempotents, S-units and S-zero divisors in NA seminear-rings. We do not demand the seminear-ring to be S-seminear-ring. Even if we have defined those earlier we recall it so that it will be used by the reader in finding S-zero divisors, S-units, S-idempotents etc in case of the groupoid seminear-ring NG and the loop seminear-ring NL.

**DEFINITION 8.4.21**: *Let N be a seminear-ring not necessarily associative. We say a, b* ∈ *N is a Smarandache zero divisor (S-zero divisor) if a.b = 0 and there exists x, y ∈ N* \ *{a, b, 0}; x ≠ y with*

> *1. ax = 0 or    xa = 0,*
> *2. by = 0 or    yb = 0 and*
> *3. xy ≠ 0 or    yx ≠ 0.*

**THEOREM 8.4.11**: *If x, y ∈ N is a S-zero divisor then it is a zero divisor.*

*Proof*: Direct by the very definition.

The reader is requested to construct an example of a zero divisor in a non-associative near-ring N which is not a S-zero divisor.

***Example 8.4.3***: Let N = {0, 1} be a near-ring and G be a groupoid given by the following table:

| . | e | $a_1$ | $a_2$ | $a_3$ | $a_4$ | $a_5$ |
|---|---|---|---|---|---|---|
| e | e | $a_1$ | $a_2$ | $a_3$ | $a_4$ | $a_5$ |
| $a_1$ | $a_1$ | e | $a_5$ | $a_4$ | $a_3$ | $a_2$ |
| $a_2$ | $a_2$ | $a_3$ | e | $a_1$ | $a_5$ | $a_4$ |
| $a_3$ | $a_3$ | $a_5$ | $a_4$ | e | $a_2$ | $a_1$ |
| $a_4$ | $a_4$ | $a_2$ | $a_1$ | $a_5$ | e | $a_3$ |
| $a_5$ | $a_5$ | $a_4$ | $a_3$ | $a_2$ | $a_1$ | e |

Does the groupoid near-ring $Z_2G$ have zero divisors?

**DEFINITION 8.4.22**: *Let N be a non-associative near-ring. We say an element 1 ≠ x ∈ N is said to be Smarandache unit (S-unit) in N if there exists y ∈ N with xy = 1. There exists a, b ∈ N \ {x, y, 1} such that*

> *i.      xa = y (or ax = y) or*
> *ii.     yb = x (or by = x) and*
> *iii.    ab = 1.*

*Clearly it can be proved every S-unit is a unit and not conversely.*

Does the group seminear-ring $Z_2G$ have S-units?

Now we proceed on to define Smarandache idempotents (S-idempotents).

**DEFINITION 8.4.23**: *Let N be a non-associative seminear-ring. An element 0 ≠ x ∈ N is a S-idempotent of N if*

> *1.  $x^2 = x$.*
> *2.  There exists a y ∈ S \ {x} such that*

*i.*     $y^2 = x$ *and*

*ii.*    $xy = y\ (yx = y)$ *or* $xy = x\ (yx = x).$

*'or' in the mutually exclusive sense.*

The study of the existence of S-idempotents in seminear-rings happens to be an interesting problem as the notion of central idempotents is used to define decomposition of S-seminear-rings IV.

**THEOREM 8.4.12**: *Let N be a S-NA seminear-ring IV. If P ⊂ N be the near-ring and e ∈ P is an idempotent of P then we get pierce decomposition which we choose to call as Smarandache pierce decomposition (S-pierce decompostition) of N. i.e. every p ∈ P ⊂ N has a decomposition of the form p = (p − pe) + pe.*

*Proof*: Straightforward.

The S-pierce decomposition of a NA seminear-ring, N is defined only when N is a S-seminear-ring of level IV and e is a nontrivial idempotent in P, P ⊂ N; P a near-ring contained in N. Several results in this direction are left for the reader to develop. We end this section by introducing Smarandache identities like Bruck, Bol, Moufang, WIP and alternative (left) right for the seminear-ring. We give definition for one of the identity, the definition for all other identities can be analogously defined in a similar way. We define only for Moufang identity, it is stated in a similar way for Bol, WIP alternative etc.

**DEFINITION 8.4.24**: *Let N be a seminear-ring. We say N is a Smarandache Moufang seminear-ring I (S-Moufang seminear-ring I) if N is a S-seminear-ring of level I and N is itself a Moufang seminear-ring.*

**DEFINITION 8.4.25:** *Let N be a seminear-ring we say N is a Smarandache Moufang seminear-ring of level II (S-Moufang seminear-ring of level II) if for all x, y z ∈ P, (P ⊂ N where P is a non-associative near-ring) we have (xy) (zx) = (x (yz)) x .*

We see S-Moufang seminear-ring cannot be defined in case of S-seminear-ring of level of III and IV in an analogous way. Thus we do not define S-Moufang seminear-ring III or IV. Now in a very similar way level I and II Smarandache Bol, WIP, alternative seminear-rings are defined.

***Example 8.4.4***: Let $Z_2 = \{0, 1\}$ be the near-field and G be a groupoid given by the following table:

| . | e | $g_0$ | $g_1$ | $g_2$ | $g_3$ | $g_4$ | $g_5$ |
|---|---|---|---|---|---|---|---|
| e | e | $g_0$ | $g_1$ | $g_2$ | $g_3$ | $g_4$ | $g_5$ |
| $g_0$ | $g_0$ | e | $g_3$ | $g_0$ | $g_3$ | $g_0$ | $g_3$ |
| $g_1$ | $g_2$ | $g_5$ | e | $g_5$ | $g_2$ | $g_5$ | $g_2$ |
| $g_2$ | $g_2$ | $g_4$ | $g_1$ | e | $g_1$ | $g_4$ | $g_1$ |
| $g_3$ | $g_3$ | $g_3$ | $g_0$ | $g_3$ | e | $g_3$ | $g_0$ |
| $g_4$ | $g_4$ | $g_2$ | $g_5$ | $g_2$ | $g_5$ | e | $g_5$ |
| $g_5$ | $g_5$ | $g_2$ | $g_4$ | $g_1$ | $g_4$ | $g_1$ | e |

We declare $1.g = g.1 = g$ for all $g \in G$. Is $Z_2G$ S-right alternative or S-left alternative?

***Example 8.4.5:*** $Z_9 = \{0, 1, 2, \ldots, 8\}$ define '×' multiplication modulo 9. $a.b = a$ for all $a, b \in Z_9$ $\{Z_9, \text{'×'}, \text{'.'}\}$ is a seminear-ring and G be a groupoid given by the following table:

| .     | e     | $g_0$ | $g_1$ | $g_2$ | $g_3$ | $g_4$ | $g_5$ | $g_6$ | $g_7$ |
|-------|-------|-------|-------|-------|-------|-------|-------|-------|-------|
| e     | e     | $g_0$ | $g_1$ | $g_2$ | $g_3$ | $g_4$ | $g_5$ | $g_6$ | $g_7$ |
| $g_0$ | $g_0$ | e     | $g_6$ | $g_4$ | $g_2$ | $g_0$ | $g_6$ | $g_4$ | $g_2$ |
| $g_1$ | $g_1$ | $g_2$ | e     | $g_6$ | $g_4$ | $g_2$ | $g_0$ | $g_6$ | $g_4$ |
| $g_2$ | $g_2$ | $g_4$ | $g_2$ | e     | $g_6$ | $g_4$ | $g_2$ | $g_0$ | $g_6$ |
| $g_3$ | $g_3$ | $g_6$ | $g_4$ | $g_2$ | e     | $g_6$ | $g_4$ | $g_2$ | $g_0$ |
| $g_4$ | $g_4$ | $g_0$ | $g_6$ | $g_4$ | $g_2$ | e     | $g_6$ | $g_4$ | $g_2$ |
| $g_5$ | $g_5$ | $g_2$ | $g_0$ | $g_6$ | $g_4$ | $g_2$ | e     | $g_6$ | $g_4$ |
| $g_6$ | $g_6$ | $g_4$ | $g_2$ | $g_0$ | $g_6$ | $g_4$ | $g_2$ | e     | $g_6$ |
| $g_7$ | $g_7$ | $g_6$ | $g_4$ | $g_2$ | $g_0$ | $g_6$ | $g_4$ | $g_2$ | e     |

$Z_9G$ is a non-associative seminear-ring. Clearly $e.g_i = g_i.e = g_i$ serves as the identity.

Now $Z_9$ is a S-NA seminear-ring III. Obtain S-units, S-zero divisors; S-idempotents and compare them with units, zero divisors and units.

***Example 8.4.6:*** Let $Z_8 = \{0, 1, 2, \ldots, 7\}$ be a seminear-ring under '×' and '.'. $L_5(3)$ be a loop of order 6. Let $Z_8L_5(3)$ be the loop seminear-ring. $Z_8L_5(3)$ is a NA seminear-ring. It is a S-seminear-ring III. Obtain S-ideals and S-subsemigroups in them.

**PROBLEMS:**

1. Let $Z_3L_5(4)$ be a NA seminear-ring. Prove $Z_3L_5(3)$ is a S-seminear-ring of level IV. Does $Z_3L_5(4)$ have S-zero divisors?
2. Prove $Z_3L_5(4)$ and $Z_3G$ where G is a groupoid given by the following table:

| .     | e     | $g_0$ | $g_1$ | $g_2$ |
|-------|-------|-------|-------|-------|
| e     | e     | $g_0$ | $g_1$ | $g_2$ |
| $g_0$ | $g_0$ | e     | $g_2$ | $g_1$ |
| $g_1$ | $g_1$ | $g_2$ | e     | $g_0$ |
| $g_2$ | $g_2$ | $g_1$ | $g_0$ | e     |

   are S- isomorphic IV?
3. Can $Z_3G$ given in problem 2 be S-regular? Substantiate your claim.
4. Let $\{Z_{12}, \times, \text{'.'}\}$ be a seminear-ring. $L_7(5)$ be the loop of order 10. For the loop seminear-ring $Z_{12}L_9(5)$ Find S-ideals, S-units and S-N-subsemigroups.
5. Is $Z_2L_9$ (5) a NA seminear ring? Justify your claim.
6. Give an example of a S-seminear-ring III which is not a S-seminear-ring IV.
7. Does there exist an example of a S-seminear-ring II which is not a S-seminear-ring III?

## 8.5 New notions on S-NA-near-rings and S-NA seminear-rings

In this section we introduce several concepts in S-NA-near-rings and S-NA seminear-rings, like Smarandache quasi non-associative near-ring (seminear-ring) S-quasi N-subgroups, Smarandache quasi N-subsemigroups Smarandache quasi ideals etc and concepts like Smarandache quasi bipotent elements and several analogous properties done in case of near-rings. Finally study of Smarandache SNP rings of three levels are introduced and this section concludes with the concept of Smarandache quasi regular elements in non-associative near-ring. The concept of Smarandache Jacobson radical is also introduced and it is left for the reader to develop.

**DEFINITION 8.5.1**: *Let N be a non-associative near-ring we say N is a Smarandache quasi non-associative near-ring (S-quasi non-associative near-ring) if N has a proper subset which is a ring under the operations of N.*

**DEFINITION 8.5.2**: *Let N be a non-associative seminear-ring we say N is a Smarandache quasi non-associative seminear-ring (S-quasi non-associative seminear-ring) if N has a proper subset P such that P is a semiring.*

**Example 8.5.1**: Let N = $Z_8 \times Z_{12}$ be the S-mixed direct product of a ring and a seminear-ring where $Z_8$ is the ring of integers modulo 8 and $Z_{12}$ is a seminear-ring under the operation '×' and '.'. G be any groupoid with unit. NG is the groupoid seminear-ring which is a S-quasi non-associative near-ring.

**Example 8.5.2**: Let N = $Z^o \times Z_8$ (where $Z^o = Z^+ \cup \{0\}$ is a semiring and $Z_8$ is a seminear-ring) be the S-mixed direct product of a semiring and a seminear-ring. Clearly NG is a S-quasi seminear-ring where G is any groupoid with identity. Thus the concept of S-mixed direct product has helped us to construct non-trivial and non-abstract examples of such structures. Now for these S-quasi near-rings and S-quasi seminear-rings we can define in an analogous way the concept of S-quasi ideals, S-quasi N-subsemigroup and S-quasi N-subgroup.

**DEFINITION 8.5.3**: *Let (N, +, .) be a S-quasi seminear-ring. We call a non-empty subset I to be a Smarandache quasi left ideal (S-quasi left ideals) in N if*

      *a.  (I, +) is a S-semigroup.*
      *b.  $n (n^1 + i) + nn^1 \in I$ for each $i \in I$ and $n, n^1 \in P$; P a semiring in N.*
          *We say I is a S-quasi ideal if $IP \subset I$.*

**DEFINITION 8.5.4**: *Let (N, +, .) be a S-quasi near-ring. We say a non-empty subset I of N to be a Smarandache quasi left ideal (S-quasi left ideal) in N if*

*a.  (I, +) is a subgroup.*

*b.  $n (n^1 + i) + nn^1 \in I$ for each $i \in I$ and $n, n^1 \in R$, $R \subset N$ and R a ring. We say I is a S-quasi ideal if I is a S-quasi left ideal of N and $IR \subset I$.*

**DEFINITION 8.5.5**: *Let N be a S-quasi near-ring (S-quasi seminear-ring). We say N is Smarandache quasi bipotent (S-quasi bipotent) if $Pa = Pa^2$ where $P \subset N$ and P is a ring ($P \subset N$ and P is a semiring) for every a in N.*

**DEFINITION 8.5.6**: *Let N be a S-quasi near-ring (S-quasi seminear-ring) N is said to be a Smarandache quasi s-near-ring (S-quasi s-near-ring) (Smarandache quasi s-seminear-ring, in short, S-quasi s-seminear-ring) if $a \in Pa$ for each a in N where P is a proper subset N which is a ring ($P \subset N$ and P is a semiring).*

**DEFINITION 8.5.7**: *Let N be a S-quasi near-ring (S-quasi near-ring) N is said to be Smarandache quasi regular (S-quasi regular) if for each a in N there exists x in P; P $\subset N$, P a ring ($P \subset N$ and P a semiring) such that $a = a(xa) = (ax)a$.*

**DEFINITION 8.5.8**: *A S-quasi near-ring (Smarandache quasi seminear-ring) N is called Smarandache quasi irreducible (S-quasi irreducible) (Smarandache quasi simple) if it contains only the trivial S-quasi N-subgroups (S quasi N-subsemigroups).*

**DEFINITION 8.5.9**: *Let N be a S-quasi near-ring (S-quasi seminear-ring) an element x is said to be quasi central if $xy = yx$ for all $y \in R$; $R \subset N$ is a ring (or $R \subset N$ and R is a semiring).*

**DEFINITION 8.5.10**: *Let N be a S-quasi near-ring (or S-quasi seminear-ring) N is said to be Smarandache quasi subdirectly irreducible (S-quasi subdirectly irreducible) if the intersection of all nonzero S-quasi ideals of N is nonzero.*

**DEFINITION 8.5.11**: *Let N be a S-quasi near-ring (S-quasi seminear-ring) N is said to have Smarandache quasi insertion of factors property (S-quasi IFP) if a, b $\in$ N, ab = 0 implies arb = 0 where $r \in R$, $R \subset N$ and R is a ring (or $r \in R$, $R \subset N$, R is a semiring).*

Note: We can take a (nb) or (an) b in all cases it should vanish that is anb = 0.

Now we define the concept of Smarandache quasi divisibility and Smarandache divisibility.

**DEFINITION 8.5.12**: *Let N be a non-associative S-near-ring we say N is Smarandache weakly divisible (S-weakly divisble) if for all x, y $\in$ N there exists a z $\in$ P; P $\subset$ N where P is an associative ring or P is a near-field such that xz = y or zx = y.*

**DEFINITION 8.5.13**: *Let N be a non-associative S-seminear-ring we say N is Smarandache weakly divisible (S-weakly divisible) if for all x, y $\in$ N there exists z $\in$ P, P $\subset$ N where P is an associative seminear-ring such that xz = y or zx = y.*

**DEFINITION 8.5.14**: *Let N be a S-quasi near-ring (S-quasi seminear-ring). We say N is Smarandache quasi weakly divisible (S-quasi weakly divisible) if for all x, y $\in$ N there exists z $\in$ R; R a ring $R \subset N$ (R a semiring $R \subset N$) such that xz = y or yz = x.*

**DEFINITION 8.5.15**: *Let N be a S-near-ring I (or S-seminear-rings). We say N is Smarandache strongly prime (S-strongly prime) I if for each a ∈ N \ {0} there exists a finite S-semigroup F such that a(Fx) ≠ 0 for all x ∈ N \ {0}.*

*In case N is S-seminear-ring II (III or IV) we say N is a Smarandache strong prime II (III or IV) (S-strong prime II (III or IV)) if for each a ∈ N \ {0} there exists a finite S-semigroup F such that a Fx ≠ 0 for all x ∈ P \ {0}, P ⊂ N; P is a associative seminear-ring / P is a associative near-ring.*

**DEFINITION 8.5.16**: *Let N be a non-associative right near-ring and A an S-ideal or a S-left ideal of N. We define three properties as follows.*

    *a.*     *A is Smarandache equiprime (S-equiprime) if for any a, x, y ∈ N such that a(nx) − a(ny) ∈ A for all n ∈ N or (an) x − (an) y ∈ Y we have a ∈ A or x − y ∈ A.*

    *b.*     *A is Smarandache strongly semiprime (S-strongly semiprime) if for each a finite subset F of N such that if x, y ∈ N and*

$$(af) x − (af) y ∈ A \text{ or}$$
$$a (fx) − a(fy) ∈ A \text{ or}$$
$$(af) x − a (fy) ∈ A \text{ or}$$
$$a (fx) − (af) y ∈ A$$

    *for all f ∈ F then x − y ∈ A.*

    *c.*     *A is Smarandache completely equiprime (S-completely equiprime) if a ∈ N \ A and ax − ay ∈ A imply x − y ∈ A.*

The study of properties in S-non associative near-rings has validity when the non associative identity plays a role in the definition.

**DEFINITION 8.5.17**: *Let S be a non-empty subset of a S-right near-ring N which is non-associative. Define left and right Smarandache polar subsets (S-polar subsets) of N by*

$$SL( S ) = \left\{ x \left| \begin{array}{l} x( Ns ) = 0 \text{ or} \\ ( xN )s = 0 \text{ for all } s ∈ S \end{array} \right. \right\}$$

*and*

$$SR( S ) = \left\{ y \left| \begin{array}{l} ( sN )y = 0 \text{ or} \\ s( Ny ) = 0 \text{ for all } s ∈ S \end{array} \right. \right\}.$$

*Suppose $SP_L (N)$ is the set of S-left polar subsets of N and $SP_R(N)$ is the set of S-right polar subsets of N one need to test whether $SP_L(N)$ and $SP_R(N)$ are complete bounded lattices.*

Several of the results which are discussed in chapter 5 can verbatim be proved in case of non-associative near-ring. Hence to avoid repetition we have not recalled them.

**DEFINITION 8.5.18:** *Let P be a seminear pseudo ring (SNP-ring) we say P is a Smarandache SNP-ring I (S-SNP-ring I) if P has a proper subset $T \subset P$ such that T is a seminear-ring. Smarandache SNP-ring II (S-SNP-ring II) if P has a proper subset R $\subset P$ such that R is a near-ring. Smarandache SNP-ring III (S-SNP-ring III) if P has a proper $W \subset P$ such (W, $\oplus$, $\odot$) is a semiring. Thus we have 3 levels of S-SNP rings. A Smarandache SNP subring (S-SNP subring) is defined as a proper subset U of P such that (U, $\oplus$, $\odot$) is a S-SNP-ring.*

**DEFINITION 8.5.19:** *Let (P, $\oplus$, $\odot$) be a SNP-ring. A proper subset I of P is called a Smarandache SNP- ideal (S-SNP-ideal) if*

1. *for all p, q $\in$ I, p $\oplus$ q $\in$ I.*
2. *0 $\in$ I.*
3. *For all p $\in$ I and r $\in$ P we have p $\odot$ r or r $\odot$ p $\in$ I.*
4. *I is a S-SNP-ring.*

**DEFINITION 8.5.20:** *Let (R, $\oplus$, $\odot$) be a quasi SNP-ring. R is said to be a Smarandache quasi SNP-rings (S-quasi SNP-ring) if and only if R is a S-SNP-ring.*

**DEFINITION 8.5.21:** *Let (R, $\oplus$, $\odot$) and (R$_1$, $\oplus$, $\odot$) be any two S-SNP-ring, we say a map $\phi$ is a Smarandache SNP-homomorphism I (II or III) (S-SNP homomorphism I, II or III) if $\phi$: S to S$_1$ where $S \subset R$ and $S_1 \subset R_1$ are seminear-ring (or near-ring or semiring) respectively and $\phi$ is a seminear-ring homomorphism from S to S$_1$ (or near-ring homomorphism from S to S$_1$ or a semiring homomorphism from S to S$_1$). $\phi$ need not be defined on the entire set R or R$^1$ it is sufficient if it is well defined on S to S$_1$.*

*Now we proceed on to define Smarandache right quasi regular element. We just recall that an element x $\in$ R, R a ring is said to be right quasi regular if there exist y $\in$ R such that x o y = x + y – xy = 0 and left quasi regular if there exist y$^1$ $\in$ R such that y$^1$ o x = 0 = y$^1$ + x – y$^1$x.*

*Quasi regular if it is right and left quasi regular simultaneously. We say an element x $\in$ R is Smarandache right quasi regular (S-right quasi regular) if there exist y and z $\in$ R such that x o y = x + y – xy = 0, x o z = x + z – xz = 0 but y o z = y + z – yz ≠ 0 and z o y = y + z – zx ≠ 0.*

*Similarly we define Smarandache left quasi regular (S-left quasi regular) and x will be Smarandache quasi regular (S-quasi regular) if it is simultaneously S-right quasi regular and S-left quasi regular.*

The study of the quasi regular concept happens to be an interesting study in case of near-rings and seminear-rings. We leave the reader to obtain results and define S-quasi regular elements in non-associative near-rings and seminear-rings; as the definition of S-quasi regularity does not involve any associative or non-associative elements. This study is a routine and the reader is expected to obtain some interesting results about them.

If we define S-non-associative ring to have a substructure which is a non-associative ring then we have several interesting results.

**THEOREM 8.5.1**: *Let $N$ be a S-near-ring having a proper subset $P$ of $N$ to be a commutative ring with unit and of characteristic 0. $L$ any loop of finite order. Then the near loop ring $NL$ has a right quasi regular element $x = \Sigma \alpha_i m_i$ ($m_i \in L$) $\alpha_i \in P \subset N$ is right quasi regular then $\Sigma \alpha_i \neq 1$.*

*Proof*: Let $y = \Sigma \beta_i h_j$, $\beta_i \in P$ and $h_j \in L$ be the right quasi inverse of x then $x + y - xy = 0$ that is $\Sigma \alpha_i m_i + \Sigma \beta_j h_j - (xy) = 0$. Equating the coefficients of the like terms and adding these coefficients we get $\Sigma \alpha_i + \Sigma \beta_j - \Sigma \alpha_i \Sigma \beta_j = 0$ or $\Sigma \alpha_i = \Sigma \beta_j (\Sigma \alpha_i - 1)$. Now if $\Sigma \alpha_i = 1$ then $\Sigma \alpha_i = 0$ a contradiction. Hence $\Sigma \alpha_i \neq 0$.

*Example 8.5.3*: Let $N = R \times Z_m$ be a near-ring, where $(Z_m, +, .)$ is a right near-ring and $R$ the field of reals. Clearly $N$ is a near-ring which is a S-right near-ring III. L be a loop given by the following table:

| . | e | a | b | c | d |
|---|---|---|---|---|---|
| e | e | a | b | c | d |
| a | a | e | c | d | b |
| b | b | d | a | e | c |
| c | c | b | d | a | e |
| d | d | c | e | b | a |

Clearly NL is a near loop ring ZL $\subset$ NL is a loop ring. Take e + b $\in$ ZL; (e + b) o (e + c) = e + b + e + c – (e + b) (e + c) = 0 so e + c is the right quasi inverse of e + b. Clearly e + d is a left quasi inverse of e + b. Note that e + d $\neq$ e + c as d and c are distinct elements in L. Thus we see in near loops rings the left quasi inverse and right quasi inverse of an element need not in general be the same.

**THEOREM 8.5.2**: *Let $N = Z_2 \times Z_{15}$ where $Z_2$ is the prime field of characterize two and $Z_{15}$ is a near-ring. Let $L$ be any loop. $NL$ be the near loop ring. If $x \in Z_2 L \times \{0\} \subset (Z_2 \times Z_{15})$; $L$ is right quasi regular then $|supp\ x|$ is an even number.*

*Proof*: The proof is easily obtained by simple calculations.

Now we proceed on to define Smarandache Jacobson radical of a near loop ring.

**DEFINITION 8.5.22**: *Let $N = N_1 \times N_2$ where $N_1$ is a field of characterize zero and $N_2$ is any near-ring. $NL$ be the near loop ring of the loop $L$ over the near-ring $N$.*
*We define SJ (P) to be the Smarandache Jacobson radical (S-Jacobson radical) of NL if $P \subset NL$ is a non-associative ring and J(P) denotes the usual Jacobson radical of the non-associative ring P.*

*Example 8.5.4*: Let $N = Z \times Z_{18}$ be the mixed direct product of the ring Z and the near-ring $Z_{18}$ L any finite loop, NL the near loop ring of the loop L over the near-ring N. Clearly ZL $\subset$ NL and ZL is a non-associative ring. If $x = \Sigma \alpha_i h_i \in ZL$ such that $\Sigma$ $\alpha_i \neq 0$ then $x \notin SJ(ZL)$. It is left for the reader the verify, as the conclusion derived is straightforward.

**Example 8.5.5**: Let L be any finite loop. $N = Z_9 \times Z_7$ be S-mixed direct product of the near-ring $Z_9$ and the prime field of characterize 7, $Z_7$, N is S-near-ring. NL is the near loop ring of the loop L over the near-ring N. If $x \in S(J(Z_7L))$ such that $x = \Sigma \alpha_i m_i$ and $\Sigma \alpha_i \neq 0 \pmod 7$ then $x \notin S(J(Z_7 L))$.

Using these two examples several interesting results can be generalized in case of loop rings. NL where N has subsets which are fields of characteristic 0 or characteristic p or Z.

**PROBLEMS :**

1.   Prove or disprove in general $SP_L(N)$ and $SP_R(N)$ are bounded lattices.
2.   Prove $N = Z_{12} \times Z_7$ is a S-seminear-ring. Let NL be the near loop ring and L be a loop given by the following table:

|       | e     | $g_1$ | $g_2$ | $g_3$ | $g_4$ | $g_5$ |
|-------|-------|-------|-------|-------|-------|-------|
| e     | e     | $g_1$ | $g_2$ | $g_3$ | $g_4$ | $g_5$ |
| $g_1$ | $g_1$ | e     | $g_3$ | $g_5$ | $g_2$ | $g_4$ |
| $g_2$ | $g_2$ | $g_5$ | e     | $g_4$ | $g_1$ | $g_3$ |
| $g_3$ | $g_3$ | $g_4$ | $g_1$ | e     | $g_5$ | $g_2$ |
| $g_4$ | $g_4$ | $g_3$ | $g_5$ | $g_2$ | e     | $g_1$ |
| $g_5$ | $g_5$ | $g_2$ | $g_4$ | $g_1$ | $g_3$ | e     |

   a.   Is NL a S-near loop ring?
   b.   Is NL a S-seminear-ring II or III or IV? Substantiate your claim.
   c.   Find S-ideals, S-N- subsemigroup and S-N subgroup.
   d.   Does NL have S-zero divisors, S-units or S-idempotents? Justify your answer.

3.   Let $N = Z_3 \times Z_{12}$ where N is the mixed direct product of the field $Z_3$ and the near-ring $Z_{12}$. L a loop given by the following table:

| .   | e   | a   | b   | c   | d   |
|-----|-----|-----|-----|-----|-----|
| e   | e   | a   | b   | c   | d   |
| a   | a   | e   | c   | d   | b   |
| b   | b   | d   | a   | e   | c   |
| c   | c   | b   | d   | a   | e   |
| d   | d   | c   | e   | b   | a   |

   NL the near loop ring of the loop L over the near-ring N. Find SJ $(Z_3L)$. Is J $(Z_3L)$ = SJ $(Z_3L)$? Find quasi regular elements of NL. Does NL have S-quasi regular elements?

4.   Let $N = Z \times Z_{15}$ be the S-mixed direct product of the ring Z and $Z_{15}$ the near-ring. L be a loop given by the following table:

| . | e | $g_1$ | $g_2$ | $g_3$ | $g_4$ | $g_5$ | $g_6$ | $g_7$ |
|---|---|---|---|---|---|---|---|---|
| e | e | $g_1$ | $g_2$ | $g_3$ | $g_4$ | $g_5$ | $g_6$ | $g_7$ |
| $g_1$ | $g_1$ | e | $g_5$ | $g_2$ | $g_6$ | $g_3$ | $g_7$ | $g_4$ |
| $g_2$ | $g_2$ | $g_5$ | e | $g_6$ | $g_3$ | $g_7$ | $g_4$ | $g_1$ |
| $g_3$ | $g_3$ | $g_2$ | $g_6$ | e | $g_7$ | $g_4$ | $g_1$ | $g_5$ |
| $g_4$ | $g_4$ | $g_6$ | $g_3$ | $g_7$ | e | $g_1$ | $g_5$ | $g_2$ |
| $g_5$ | $g_5$ | $g_3$ | $g_7$ | $g_4$ | $g_1$ | e | $g_2$ | $g_6$ |
| $g_6$ | $g_6$ | $g_7$ | $g_4$ | $g_1$ | $g_5$ | $g_2$ | e | $g_2$ |
| $g_7$ | $g_7$ | $g_4$ | $g_1$ | $g_5$ | $g_2$ | $g_6$ | $g_3$ | e |

Find for the near loop ring NL

    a. S (J (ZL)).
    b. J (ZL).
    c. Find S-quasi regular elements of NL.

5. Give an example of a S-weakly divisible near-ring.

6. Is $Z_{18}G$ S-quasi bipotent where $Z_{18}$ is a seminear-ring and G a groupoid given by the following table:

| . | e | $a_1$ | $a_2$ | $a_3$ | $a_4$ | $a_5$ |
|---|---|---|---|---|---|---|
| e | e | $a_1$ | $a_2$ | $a_3$ | $a_4$ | $a_5$ |
| $a_1$ | $a_1$ | e | $a_3$ | $a_5$ | $a_2$ | $a_4$ |
| $a_2$ | $a_2$ | $a_5$ | e | $a_4$ | $a_1$ | $a_3$ |
| $a_3$ | $a_3$ | $a_4$ | $a_1$ | e | $a_5$ | $a_2$ |
| $a_4$ | $a_4$ | $a_3$ | $a_5$ | $a_2$ | e | $a_1$ |
| $a_5$ | $a_5$ | $a_2$ | $a_4$ | $a_1$ | $a_3$ | e |

7. Give an example of a S-quasi subdirectly irreducible near-ring.

**Chapter Nine**

# FUZZY NEAR-RINGS AND SMARANDACHE FUZZY NEAR-RINGS

This chapter has three sections. In the first section we just recall all the definitions about fuzzy ideals and other known fuzzy concepts in near-rings. In the second section we give five new types of fuzzy near-rings which are formed in an unconventional way. In the final section we introduce the concept of Smarandache fuzzy near-rings, ideals and roots.

The concept of fuzzy subset was introduced by Zadeh L.A. [105]. Later several authors has developed it to fuzzy near-rings like [1, 25, 40, 46, 89, 94, 95, 98]. Till date there are only about a dozen papers in fuzzy near-rings. Here we introduce and study these with reference to Smarandache fuzzy near-rings we give the basic definitions needed to define this concept.

## 9.1 Basic notions on fuzzy near-rings

In this section we just recall the basic known and available fuzzy notions about near-rings. To the best of our knowledge we have only a very few papers on fuzzy near-ring that too mainly dealing with fuzzy ideals.

**DEFINITION 9.1.1**: *Let X be a nonempty set. A mapping $\mu: X \rightarrow [0, 1]$ is called a fuzzy subset of X.*

**DEFINITION 9.1.2:** *Let $\sigma$ and $\theta$ be two fuzzy subsets of N. We define the product of these fuzzy subsets (denoted by $\sigma \, o \, \theta$) as follows.*

$(\sigma \, o \, \theta)(x) = \sup_{x=yz} \{ min(\sigma(y), \theta(z)) \}.$ *(If x can be expressible as a product yz.)*

$(\sigma \, o \, \theta)(x) = 0,$ *(otherwise)* .

**DEFINITION 9.1.3**: *Let X and Y be any nonempty sets and f, a function of X into Y. Let $\mu$ and $\sigma$ be fuzzy subsets of X and Y respectively. Then f($\mu$), the image of $\mu$, under f, is a fuzzy subset of Y defined by*

$$(f(\mu)y) = \begin{cases} \sup \mu(x) \; if \; f^{-1}(y) \neq \phi \\ 0 \; if \; f^{-1} = 0 \end{cases}$$

*and $f^{-1}(\sigma)$ the pre image of $\sigma$ under f is a fuzzy subset of X defined by $(f^{-1}(\sigma))(x) = \sigma(f(x))$ for all $x \in X$.*

*Notation*: Let μ and σ be two fuzzy subsets of N then we write $\mu \subseteq \sigma$ if $\mu(x) \subseteq \sigma(x)$ for all $x \in N$.

**DEFINITION 9.1.4**: *Let $\mu$ be a non-empty fuzzy subset of a near-ring N (that is $\mu(x) \neq 0$ for some $x \in N$) then $\mu$ is said to be a fuzzy ideal of N if it satisfies the following conditions:*

    a.  $\mu(x + y) \geq min\ \{\mu(x),\ \mu(y)\}$.
    b.  $\mu(-x) = \mu(x)$.
    c.  $\mu(x) = \mu(y + x - y)$.
    d.  $\mu(xy) \geq \mu(x)$ and
    e.  $\mu\ \{x\ (y + i) - xy\ \} \geq \mu(i)$ for all x, y, i ∈ N.

1.  *If $\mu$ is a fuzzy ideal of N then $\mu(x + y) = \mu(y + x)$*
2.  *If $\mu$ is a fuzzy ideal of N then $\mu(0) \geq \mu(x)$ for all $x \in N$. The above two statements can be easily verified for if we put $z = x + y$, then $\mu(x + y) = \mu(z) = \mu(-x + z + x)$ since $\mu$ is a fuzzy ideal $\mu(-x + x + y + x) = \mu(y + x)$ (since $z = x + y$).*

*Likewise for the second statement $\mu(0) = \mu(x - x) \geq min\ \{\mu(x),\ \mu(-x)\} = \mu(x)$  since $\mu$ is a fuzzy ideal (since $\mu(-x) = \mu(x)$ by the very definition of fuzzy ideal).*

**DEFINITION 9.1.5**: *Let I be an ideal of N. we define $\lambda_I : N \rightarrow [0, 1]$ as*

$$\lambda_I(x) = \begin{cases} 1 & \text{if } x = 1 \\ 0 & \text{otherwise} \end{cases}$$

*$\lambda_I(x)$ is called the characterize function on 1.*

**THEOREM [25]**: *Let N be a near-ring and $\lambda_I$, the characteristic function on a subset I of N. Then $\lambda_I$ is a fuzzy ideal of N if and only if I is an ideal of N.*

**DEFINITION [25]**: *Let $\mu$ be a fuzzy subset of X. Then the set $\mu_k$ of all $t \in [0, 1]$ is defined by $\mu_t = \{x \in N / \mu(x) \geq t\}$ is called the level subset of t for the near-ring N.*

**DEFINITION [25]**: *Let N be a near-ring and $\mu$ be a fuzzy ideal of N. Then the level subset $\mu_t$ of N for all $t \in [0, t]$, $t \leq \mu(0)$ is an ideal of N if and only if $\mu$ is a fuzzy ideal of N.*

**DEFINITION [25]**: *A fuzzy ideal $\mu$ of N is called fuzzy prime if for any two fuzzy ideals $\sigma$ and $\theta$ of N $\sigma \circ \theta \subseteq \mu$ implies $\sigma \subseteq \mu$ or $\theta \subseteq \mu$.*

Now we use the concept of fuzziness in Γ-near-rings as given by [1, 72].

**THEOREM [72]**: *If $\mu$ is a fuzzy ideal of a near-ring N and $a \in N$ then $\mu(x) \geq \mu(a)$ for all $x \in \langle a \rangle$.*

*Proof*: Please refer [1, 72].

**DEFINITION [72]:** *Let μ and σ be two fuzzy subsets of M. Then the product of fuzzy subset ( σ o τ)(x) = sup {min (σ (y), τ (z))} if x is expressible as a product x = yz where y, z ∈ M and (σ o τ) (x) = 0 otherwise.*

**DEFINITION [72]:** *A fuzzy ideal μ of N is said to have fuzzy IFP if μ (a n b) ≥μ(ab) for all a, b, n ∈ N.*

**DEFINITION [72]:** *Let μ be a fuzzy ideal of N, μ has fuzzy IFP if and only if $μ_k$ is a IFP ideal of N for all $0 ≤ k ≤ 1$.*

**DEFINITION [72]:** *N has strong IFP if and only if every fuzzy ideal of N has fuzzy IFP.*

**THEOREM [72]:** *If μ is a fuzzy IFP-ideal of N then $N_μ$ = {x ∈ N / μ (x) = μ(0)} is an IFP ideal of N.*

*Proof*: We have μ (0) ≥ μ (x) for all x ∈ N; write t = μ (0). Now $N_μ = μ_t$ and by definitions $μ_t$ has IFP. Therefore $N_μ$ is an IFP ideal of N.

*Notation*: Let μ be a fuzzy ideal of N. For any s ∈ [0, 1] define $β_s$ : N → [0, 1] by

$$_μβ_s(x) = \begin{cases} s \; if \; μ(x) ≥ s \\ μ(x) \; if \; μ(x) < s. \end{cases}$$

Since $β_s$ depends on μ, we also denote $β_s$ by $_μβ_s$.

**_Results:_**

1. $β_s$ (x) ≤ s for all x ∈ N.
2. $β_s$ is a fuzzy ideal of N.
3. If μ (0) = t, then s ≥ t if and only if $μ_t = (β_s)_t$.

The proof of the above 3 statements are left as an exercise for the reader to prove.

We can still equivalently define as a definition or prove it as a theorem.

**DEFINITION [75]:** *μ is a fuzzy IFP ideal of N if and only if $β_s$ is a fuzzy IFP ideal for all s ∈ [0, 1].*

**_Result_**: A fuzzy ideal μ has IFP if and only if $N_{(βs)}$ has IFP for all s ∈ [0, 1].

The above result is assigned as an exercise for the reader to prove.

**DEFINITION [75]:** *Let μ : M→ [0, 1]. μ is said to be a fuzzy ideal of M if it satisfies the following conditions:*

1. *μ (x + y) ≥ min {μ (x), μ (y)}.*
2. *μ (-x) = μ (x).*
3. *μ (x) = μ (y + x − y).*

4.  $\mu(x \, \alpha \, y) \geq \mu(x)$ and
5.  $\mu\{(x \, \alpha \, (y + z) - x \, \alpha \, y\} \geq \mu(z)$ for all $x, y, z \in M$ and $\alpha \in \Gamma$.

*The following result is left as an exercise for the reader to prove.*

**THEOREM 9.1.1**: *Let $\mu$ be a fuzzy subset of M. Then the level subsets $\mu_t = \{x \in M / \mu(x) \geq t\}$, $t \in im\mu$, are ideals of M if and only if $\mu$ is a fuzzy ideal of M.*

All results true in case of fuzzy ideals of N are true in case of fuzzy ideals of M with some minor modifications. The following result can be proved by routine application of definitions.

**THEOREM 9.1.2**: *Let $M$ and $M^I$ be two $\Gamma$-near-rings, $h : M \rightarrow M^I$ be an $\Gamma$-epimorphism and $\mu$, $\sigma$ be fuzzy ideals of M and $M^I$ respectively then*

1.  $h(h^{-1}(\sigma)) = \sigma$.
2.  $h^{-1}(h(\mu)) \supseteq \mu$ and
3.  $h^{-1}(h(\mu)) = \mu$ if $\mu$ is constant on ker h.

**DEFINITION [72]:** *A fuzzy ideal $\mu$ of M is said to be a fuzzy prime ideal of M if $\mu$ is not a constant function; and for any two fuzzy ideals $\sigma$ and $\Gamma$ of M, $\sigma \circ \Gamma \subseteq \mu$, implies either $\sigma \subset \mu$ or $\Gamma \subset \mu$.*

Using these definitions it can be proved.

**THEOREM [72]**: *If $\mu$ is a fuzzy prime ideal of M then $M_\mu = \{x \in M / \mu(x) = \mu(0)\}$ is a prime ideal of M.*

**PROPOSITION [72]:** *Let I be an ideal of M and $s \in [0, 1)$. Let $\mu$ be a fuzzy subset of M, defined by*

$$\mu(x) = \begin{cases} 1 \text{ if } x \in I \\ s \text{ otherwise} \end{cases}$$

*Then $\mu$ is a fuzzy prime ideal of M if I is a prime ideal of M.*

*Proof*: Using the fact $\mu$ is a non-constant fuzzy ideal of M we can prove the result as a matter of routine using the basic definitions.

**DEFINITION [72]:** *Let I be an ideal of M. Then $\lambda_I$ is a fuzzy prime ideal of M if and only if I is a prime ideal of M.*

**_Result 1_**: If $\mu$ is a fuzzy prime ideal of M then $\mu(0) = 1$.

*Proof*: It is left for the reader to prove.

**_Result 2_**: If $\mu$ is a fuzzy prime ideal of M then |Im $\mu$| = 2.

1. Prove the following 3 conditions are equivalent

   a. N has strong IFP.
   b. Every fuzzy ideal of N has fuzzy IFP and
   c. The ideal $N_{(\mu\beta s)}$ has IFP for all ideals $\mu$ of N and for all $s \in [0, 1]$.

2. Prove if $\mu$ is a fuzzy prime ideal of M then |Im $\mu$| = 2.
   What is the analogous result in case of a fuzzy prime ideal of a near-ring.

3. Let $\mu$ be a fuzzy ideal of N for any $s \in [0, 1]$ define $\beta_s : N \rightarrow [0, 1]$.
   Prove $\beta_s$ is a fuzzy ideal of N.

4. Prove if N has strong IFP then every fuzzy ideal of N has fuzzy IFP.

5. If $\mu$ is a fuzzy IFP ideal of N then is it true that $N_\mu = \{x \in N \,/\, \mu(x) = \mu\,(0)\}$ is an IFP ideal of N?

6. If $\mu$ is a fuzzy prime ideal of M then $M_\mu = \{x \in M \,/\, \mu(x) = \mu(0)\}$ is a prime ideal of M. Prove.

7. Let I be an ideal of M. Prove $\lambda_I$ is a fuzzy prime ideal of M if and only if I is a prime ideal of M.

8. Prove if $\mu$ is a fuzzy prime ideal of M then $\mu(0) = 1$.

9. Let N be a near-ring and $\mu$ be a fuzzy subset of N. Then the level subset $\mu_t$ of N for all $t \in [0, 1]$, $t \leq \mu(0)$ is an ideal of N if and only if $\mu$ is a fuzzy ideal of N.

10. Prove if $\mu$ is a fuzzy IFP ideal of N then $\beta_s$ is a IFP fuzzy ideal for all $s \in [0,1]$.

11. Prove if $\mu$ is a fuzzy ideal of a near-ring N and $a \in N$ then $\mu\,(x) \geq \mu\,(a)$ for all $x \in \langle a \rangle$.

12. Prove if $\mu$ is a fuzzy ideal of M then $\mu\,(x + y) = \mu\,(y + x)$ for all x, y $\in$ M.

## 9.2 Some special classes of fuzzy near-rings

In this section the author solely introduces the concept of fuzzy complex near-rings, fuzzy near matrix ring, fuzzy polynomial near-ring, special fuzzy near-ring and fuzzy non-associative complex near-rings and studies them.

All these concepts are far from the conventional way of defining them. Hence we in this section define these five types of fuzzy near-rings and study some of its interesting properties.

**Definition 9.2.1**: *Let $P_{n \times n}$ denote the set of all $n \times n$ matrices with entries form [0, 1] i.e. $P_{n \times n} = \{(a_{ij}) \,/\, a_{ij} \in [0, 1]\}$ for any two matrices A, B $\in P_{n \times n}$ define $\oplus$ as follows:*

$$A = \begin{pmatrix} a_{11} & a_{12} & \ldots & a_{1n} \\ a_{21} & a_{22} & \ldots & a_{2n} \\ \vdots & \vdots & & \vdots \\ a_{n1} & a_{n2} & \ldots & a_{nn} \end{pmatrix}$$

*and*

$$B = \begin{pmatrix} b_{11} & b_{12} & \dots & b_{1n} \\ b_{21} & b_{22} & \dots & b_{2n} \\ \vdots & \vdots & & \vdots \\ b_{n1} & b_{n2} & \dots & b_{nn} \end{pmatrix}$$

$$A \oplus B = \begin{pmatrix} a_{11}+b_{11} & \dots & a_{1n}+b_{1n} \\ a_{21}+b_{21} & \dots & a_{21}+b_{2n} \\ \vdots & & \vdots \\ a_{n1}+b_{n1} & \dots & a_{nn}+b_{nn} \end{pmatrix}$$

*where*

$$a_{ij}+b_{ij} = \begin{cases} a_{ij}+b_{ij} \ \text{if} \ a_{ij}+b_{ij} < 1 \\ 0 \ \text{if} \ a_{ij}+b_{ij} = 1 \\ a_{ij}+b_{ij}-1 \ \text{if} \ a_{ij}+b_{ij} > 1 \end{cases}.$$

*Clearly* $(P_{n \times n}, \oplus)$ *is an abelian group and*

$$0 = \begin{pmatrix} 0 & 0 & \dots & 0 \\ 0 & 0 & \dots & 0 \\ \vdots & \vdots & & \vdots \\ 0 & 0 & \dots & 0 \end{pmatrix}$$

*is zero matrix which acts as the additive identity with respect to* $\oplus$.

*Define* $\odot$ *on* $P_{n \times n}$ *as follows. For A, B* $\in P_{n \times n}$

$$A \odot B = \begin{pmatrix} a_{11} & \dots & a_{1n} \\ a_{21} & \dots & a_{2n} \\ \vdots & & \vdots \\ a_{n1} & \dots & a_{nn} \end{pmatrix} \odot \begin{pmatrix} b_{11} & \dots & b_{1n} \\ b_{21} & \dots & b_{2n} \\ \vdots & & \vdots \\ b_{2n} & \dots & b_{nn} \end{pmatrix} = \begin{pmatrix} a_{11}+\dots+a_{1n} & \dots & a_{11}+\dots+a_{1n} \\ a_{21}+\dots+a_{2n} & \dots & a_{21}+\dots+a_{2n} \\ \vdots & & \vdots \\ a_{n1}+\dots+a_{nn} & \dots & a_{n1}+\dots+a_{nn} \end{pmatrix}$$

*where* $a_{ij}.b_{ij} = a_{ij}$ *for all* $a_{ij} \in A$ *and* $b_{ij} \in B$. *Clearly* $(P_{n \times n}, \odot)$ *is a semigroup. Thus (A $\oplus$ B) $\odot$ C = A $\odot$ C $\oplus$ B $\odot$ C. Hence* $(P_{n \times n}, \oplus, \odot)$ *is a near-ring, which we call as the fuzzy near matrix ring.*

**THEOREM 9.2.1**: *The fuzzy matrix near-ring is a commutative near-ring.*

*Proof*: Straightforward.

**Theorem 9.2.2**: *The fuzzy matrix near-ring is not an abelian near-ring.*

*Proof*: For A, B ∈ $P_{n \times n}$ we have A $\odot$ B ≠ B $\odot$ A in general.

**Theorem 9.2.3**: *In {$P_{n \times n}$, $\oplus$, $\odot$} we have $I_{n \times n} \neq A$ where $I_{n \times n}$ is the matrix with diagonal elements 1 and rest 0.*

*Proof*: Left for the reader to prove.

**Definition 9.2.2**: *Let {$P_{n \times n}$, $\oplus$, $\odot$} be a fuzzy near matrix ring we say a subset I of $P_{nxn}$ is a fuzzy left ideal of $P_{nxn}$ if*

1. *(I, +) is a normal subgroup of $P_{nxn}$.*
2. *$n(n^I + i) + n_r n^I \in I$ for each $i \in I$ and $n_r, n, n^I \in N$ where $n_r$ denotes the unique right inverse of n.*

*All properties enjoyed by near-rings can be defined and will be true with appropriate modifications.*

Next we proceed on to define the concept of fuzzy complex near-rings.

**Definition 9.2.3**: *Let V = {a + ib / a, b ∈ [0,1]} define on V the operation called addition denoted by $\oplus$ as follows*

*For a + ib, $a_1 + ib_1 \in V$, $a + ib \oplus a_1 + ib_1 = a + a_1 + i(b + b_1)$ where $a \oplus a_1 = a + a_1$ if $a + a_1 < 1$ and $a + a_1 = a + a_1 - 1$ if $a + a_1 \geq 1$ where '+' is the usual addition of numbers. Clearly (V, $\oplus$) is a group. Define $\odot$ on V by $(a + ib) \odot (a_1 + ib_1) = a + ib$ for all $a + ib, a_1 + ib_1 \in V$. (V, $\odot$) is a semigroup. It is easily verified. (V, $\oplus$, $\odot$) is a near-ring, which we call as the fuzzy complex near-ring.*

*Further P = {a / a ∈ [0,1]} and C = {ib / b ∈ [0 1]} are fuzzy complex subnear rings of (V, $\oplus$, $\odot$).*

**Proposition 9.2.1**: V has non-trivial idempotent.

*Proof*: Left for the reader to prove.

**Theorem 9.2.4**: *Let {V, $\oplus$, $\odot$} be a fuzzy complex near-ring. Every nontrivial fuzzy subgroup of N is a fuzzy right ideal of V.*

*Proof*: Obvious by the fact that if N is a fuzzy subgroup of V then NV $\subseteq$ N.

It is an open question. Does V have nontrivial fuzzy left ideals and ideals. The reader is requested to develop new and analogous notions and definitions about these concepts.

Now a natural question would be can we have the concept of fuzzy non-associative complex near-ring; to this end we define a fuzzy non-associative complex near-ring.

**DEFINITION 9.2.4**: *Let W = {a + ib / a, b ∈ [0, 1] called the set of fuzzy complex numbers}. Define on W two binary operations ⊕ and ⊙ as follows:*

*(W, ⊕) is a commutative loop where for a + ib, c + id ∈ W define a + ib ⊕ c + id = a ~c + i ( b ~ d) where '~' is the difference between a and b. Clearly (W, ⊕) is a commutative loop.*

*Define ⊙ on W by (a + ib) ⊙ (c + id) = a + ib for all a + ib, c + id ∈ W. (W, ⊕, ⊙) is called the fuzzy complex non-associative near-ring. ([0, 1], ⊕, ⊙) ⊆ (W, ⊕, ⊙) is a fuzzy non-associative subnear-ring.*

Obtain interesting properties about these non-associative fuzzy complex near-rings.

Now we proceed on to define fuzzy polynomial near-rings.

**DEFINITION 9.2.5**: *Let R be the set of reals. The fuzzy polynomial near-ring $R[x^{[0, 1]}]$ consists of elements of the form $p_0 + p_1 x^{\gamma_1} + p_2 x^{\gamma_2} + ... + p_n x^{\gamma_n}$ where $p_0, p_1, ... p_m \in R$ and $\gamma_1, \gamma_2, ..., \gamma_n \in [0, 1]$ with $\gamma_1 < \gamma_2 < ... < \gamma_n$. Two elements $p(x) = q(x) \Leftrightarrow p_i = q_i$. and $\gamma_i = s_i$ where $p(x) = p_0 + p_1 x^{\gamma_1} + ... + p_n x^{\gamma_n}$ and $q(x) = q_0 + q_1 x^{s_1} + ... + q_n x^{s_n}$ Addition is performed as in the case of usual polynomials.*

*Define ⊙ on R $[x^{[0\ 1]}]$ by p(x) ⊙ q(x) = p(x) for p(x), q(x) ∈ R $[x^{[0, 1]}]$. Clearly {R[$x^{[0, 1]}$], +, ⊙} is called the fuzzy right polynomial near-ring. $x^0 = 1$ by definition.*

**DEFINITION 9.2.6**: *Let {R $[x^{[0, 1]}]$, +, ⊙} be a fuzzy polynomial near-ring. For any polynomial p(x) ∈ R $[x^{[0, 1]}]$ define the derivative of p(x) as follows.*

*If $p(x) = p_0 + p_1 x^{s_1} + ... + p_n x^{s_n}$*

$$\frac{dp(x)}{dx} = 0 + p_1 x^{s_1 \sim 1} + ... + s_n p_n x^{s_n \sim 1} = (s_1 p_1) x^{s_1 \sim 1} + ... + (s_n p_n) x^{s_n \sim 1}$$

*where '~' denotes the difference between $s_i$ and 1. Clearly if $p(x) \in R[x^{[0,1]}]$ then $\frac{dp(x)}{dx} \in R[x^{[0,1]}]$ Likewise successive derivatives are also defined i.e. product of $s_i p_i \in R$ as $s_i \in [0, 1]$ and $p_i \in R$ i.e. the usual multiplication of the reals.*

**Example 9.2.1**: Let R be reals R [$x^{[0, 1]}$] be a polynomial near-ring.

$$p(x) = 5 - 6x^{1/5} + 2x^{3/8} - 15x^{7/9}$$

$$\frac{dp(x)}{dx} = 0 - \frac{1}{5}6x^{4/5} + \frac{2 \times 3}{8}x^{5/8} - \frac{15 \times 7}{9}x^{2/9} = -\frac{6}{5}x^{4/5} + \frac{3}{4}x^{5/8} - \frac{35}{3}x^{2/9}.$$

The observation to be made is that no polynomial other than the polynomial x vanishes after differentiation.

**DEFINITION 9.2.7**: *Let $p(x) \in \{Rx^{[0,\ 1]}\}$ the fuzzy degree of $p(x)$ is $s_n$ where $p(x) = p_0 + p_1 x^{s_1} + \ldots + p_n x^{s_n}$; $s_1 < s_2 \ldots < s_n$ $(p_n \neq 0)$ deg $p(x) = s_n$. The maximal degree of any polynomial $p(x)$ can take is 1. Now it is important to note that as in the case polynomial rings we cannot say deg $[p(x).q(x)] = $ deg $p(x) + $ deg $q(x)$.*

*But we have always in fuzzy polynomial near-ring.*

$$deg\ (p(x)\ q(x)) = deg\ p(x)\ for$$
$$p(x),\ q(x) \in R\ [x^{[0,\ 1]}]$$

*as this degree for fuzzy polynomial near-rings is a fuzzy degree we shall denote them by f(deg $(p(x))$.*

**DEFINITION 9.2.8**: *Let $p(x) \in [R\ [x^{[0,\ 1]}]$ - $p(x)$ is said to have a root $\alpha$ if $p(\alpha) = 0$.*

**Example 9.2.2**: Let $p(x) = \sqrt{2} - x^{1/2}$ be a fuzzy polynomial in R $\{x^{[0,\ 1]}\}$. The root of $p(x)$ is 2 for $p(2) = \sqrt{2} - 2^{1/2} = 0$.

But as in case of root of polynomial of degree n has n and only n roots which is the fundamental theorem on algebra; we in case of fuzzy polynomial near-rings cannot say the number of roots in a nice mathematical terminology that is itself fuzzy.

A study of these fuzzy polynomial near-rings is left open for any interested researchers. We proceed on to define fuzzy polynomial near-rings when the number of variables is more than one x and y.

**DEFINITION 9.2.9**: *Let R be the reals x, y be two variables we first assume xy = yx. Define the fuzzy polynomial near-ring.*

$$R\ [x^{[0,\ 1]},\ y^{[0,\ 1]}]\ by = \left\{ \sum r_i x^{p_i} y^{q_i}\ /\ r_i \in R;\ p_i \in [0,1]\ \ q_i \in [0,1] \right\}$$

*Define '+' as in the case of polynomial and '.' by $p(xy) \cdot q(x,\ y) = p(x,\ y)$. Clearly $R[x^{[0,\ 1]},\ y^{[0,\ 1]}]$ is called as a fuzzy polynomial right near-ring.*

**DEFINITION 9.2.10**: *Let $\{R\ [x^{[0,\ 1]},\ y^{[0,\ 1]}]\ \oplus$ '.'$\}$ be a fuzzy polynomial near-ring in the variable x and y.*

*A fuzzy polynomial $p(x,\ y)$ is said to be homogenous of fuzzy degree t, $t \in [0,\ 1]$ if $p(x,\ y) = a_n x^{s_1} y^{t_1} + \ldots + b_n x^{s_p t_p}$ then $t_i \neq 0$, $s_i \neq 0$ for all $i = 1, 2, \ldots, p$ and $s_i + t_i = t$ for $i = 1, 2, \ldots, p$.*

**DEFINITION 9.2.11**: *Let R $\{x^{[0,\ 1]},\ y^{[0,\ 1]},\ \oplus,\ $'.'$\}$ be the fuzzy polynomial near-ring in the variables x and y. A symmetric fuzzy polynomial is a homogenous polynomial of fuzzy degree t, $t \in [0,\ 1]$ such that $p(x,\ y) = p_1 x^{t_1} y^{s_1} + \ldots + p_n x^{t_n} y^{s_n}$ where $t_1 < t_2 < \ldots < t_n$, $s_n < s_{n-1} < \ldots < s_1$ with $t_1 = s_n$, $t_2 = s_{n-1}$, $\ldots$, $t_n = s_1$ further $p_1 = p_n$, $p_2 = p_{n-1}$, $\ldots$*

For example p (x, y) = $3x^{1/2} y^{2/3} + 3x^{2/3} y^{1/2}$. p (x, y) = $x^r + y^r$ r ∈ [0, 1].

p(x, y) = $x^r + y^r + x^s y^t + y^t s^s$ where s + t = 1. s, t, r ∈ [0, 1]. We have like other polynomials we can extend the fuzzy polynomials to any number of variables say $X_1$, $X_2$, …, $X_n$. under the assumption $X_i X_j = X_j X_i$ and denote it by R [$X_1^{[0, 1]}$, $X_2^{[0, 1]}$, … $X_n^{[0, 1]}$] called the fuzzy polynomial near-ring in n variables. The reader is advised to develop new results on these fuzzy polynomial near-rings.

We have introduced the concept of complex near-ring and the non-associative complex near-ring now we just define yet another new notion called fuzzy non-associative near-ring.

**DEFINITION 9.2.12**: *Let {W, ⊕, ⊙} be the fuzzy non-associative complex near-ring. Let x be an indeterminate. We define the fuzzy non-associative polynomial near-ring as follows*

*W [x] = {Σ $p_i$ $x^i$ / $p_i$ ∈ W}; we say p(x), q(x) ∈ W[x] are equal if and only if every coefficient of same power of x is equal i.e. if p(x) = $p_0$ + $p_1$x + … + $p_n x^n$ and q(x) = $q_0$ + $q_1$x + … + $q_n x^n$. p(x) = q(x) if and only if $p_i$ = $q_i$. for i = 1, 2, …, n. Addition is performed as follows p(x) ⊕ q(x) = $p_0$ ⊕ $q_0$ + … + ($p_n$ ⊕ $q_n$) $x^n$ where ⊕ is the operation on W. For p(x), q(x) in W[x] define p(x) ⊙ q(x) = p(x). Clearly {W (x), ⊕, ⊙} is a fuzzy non-associative complex polynomial near-ring.*

*Now take $Z^o = Z^+ \cup \{0\}$. Let p: $Z^o$ → W be defined by p(0) = 0, p(x) = $\dfrac{1}{x}$ for 0 ≠ x ∈ Z.*

*Clearly p(z) is a fuzzy non-associative subnear-ring of W. Thus p is a fuzzy non-associative subnear-ring. Let G = $Z^o \times Z^o$. Define a map $p_0$ : G → W[x] by*

$$p(0, 0) = 0.$$
$$p(x, y) = \frac{1}{x} + \frac{1}{y}, x \neq 0, y \neq 0.$$
$$p(x, 0) = \frac{1}{x};$$
$$p(0, y) = \frac{1}{y}.$$

*Then the map p is a fuzzy non-associative complex subnear-ring of G.*

Several interesting research in this direction is thrown open for the reader.

Now we proceed on to define a special class of fuzzy near-ring.

**DEFINITION 9.2.13**: *Let P = [0, 1] the interval from 0 to 1. Define ⊕ and ⊙ on P as follows. For a, b ∈ P define a ⊕ b = a + b if a + b < 1, a ⊕ b = 0 if a + b = 1 and a ⊕ b = a + b −1 if a + b > 1. Thus ⊕ acts as modulo 1. Define ⊙ on a, b ∈ P = [0, 1] by a ⊙ b = a; clearly (a ⊕ b) ⊙ c = a ⊙ c + b ⊙ c = a ⊕ b. Clearly (P, ⊕) is a group* *and (P, $\odot$) is a semigroup. Hence (P, $\oplus$, $\odot$) is a right near-ring. We call {P, $\oplus$, $\odot$} the special fuzzy right near-ring.*

**DEFINITION 9.2.14**: *Let (P, $\oplus$, $\odot$) be a fuzzy near-ring. $P_0$ = {p $\in$ P / p.0 = 0} is called the fuzzy zero symmetric part and $P_c$ = {n $\in$ P / n.0 = n} is called the fuzzy constant part.*

**THEOREM 9.2.5**: *The special fuzzy right near-ring {P, $\oplus$, $\odot$} has no fuzzy invertible elements.*

*Proof*: Left for the reader to prove.

Let S = {r/p, 0 / 1 < r < p} is a fuzzy subnear-ring or to be more specific if S = {0, ¼, ½, ¾}; S is a fuzzy subnear-ring.

**DEFINITION 9.2.15**: *A fuzzy subnear-ring N of P is called fuzzy invariant if NP $\subseteq$ N and PN $\subseteq$ N we call a fuzzy subnear-ring N of P to be a fuzzy right invariant if NP $\subset$ N.*

The following theorem is left as an exercise to the reader.

**THEOREM 9.2.6**: *Every fuzzy subnear-ring N of P is fuzzy right invariant.*

**DEFINITION 9.2.16**: *The set P = {0, 1} with two binary operations $\oplus$ and $\odot$ is called fuzzy right seminear-ring if {P, $\oplus$} and {P, $\odot$} are right seminear-ring.*

All results can be easily extended in case of fuzzy seminear-ring.

***Example 9.2.3***: Let {P, $\oplus$, $\odot$} be the fuzzy semi-near-ring. Define $\oplus$ as p $\oplus$ q if p + q < 1 and p $\oplus$ q = 0 if p + q $\geq$ 1. Then (P, $\oplus$) is a semigroup. Define $\odot$ as p $\odot$ q = p for all p, q $\in$ P. Clearly {P, $\oplus$, $\odot$} is a special fuzzy semi-near-ring.

## PROBLEMS :

1.  Give an example of a fuzzy matrix near-ring $M_{2\times 2}$. Study all properties analogous to near-ring in them.

2.  Construct a finite fuzzy matrix near-ring $M_{2\times 2}$ by choosing only a finite number of elements from [0, 1]; study its properties.

3.  Does [$x^{[0\ 1]}$] have a polynomials which has no roots? Justify your answer at least with examples.

4.  Is it possible to construct a finite fuzzy complex near-ring? Justify or substantiate your answer.

5.  Construct a finite fuzzy seminear-ring using [0, 1] and study all analogous properties enjoyed by them like fuzzy ideals, normal semigroups, invariant elements.

## 9.3 Smarandache fuzzy near-rings

In this section we just introduce the notion of Smarandache fuzzy near-rings in a very unconventional way. We further introduce the Smarandache concepts to the fuzzy near-rings introduced in section 9.2 Finally we define Smarandache roots for Smarandache polynomial near-ring and leave a vast stretch of research for the reader.

**DEFINITION 9.3.1**: *Let N be a near-ring, let $\mu$ be a non-empty fuzzy subset of a near-ring. If $\mu$ is a fuzzy ideal of the near-ring N then we call N a Smarandache fuzzy near-ring (S-fuzzy near-ring) related to $\mu$ and denote it by $_\mu N$ .*

*Every fuzzy subset of a near-ring need not in general give a S-fuzzy near-ring. Also we can define S-fuzzy near-ring in any way one likes. We call this definition as the S-fuzzy near-ring of level I.*

**DEFINITION 9.3.2**: *The characteristic function $\lambda$ on a near-ring is a S-fuzzy near-ring related to $\lambda$. This S-fuzzy near-ring is unique and we call this as the Smarandache fuzzy characteristic near-ring (S-fuzzy characteristic near-ring).*

**THEOREM 9.3.1**: *Every near-ring N has a S-fuzzy characteristic (function) near-ring.*

*Proof*: Straightforward.

**THEOREM 9.3.2**: *Let N be a near-ring. N a S-fuzzy near-ring related to $\mu$ ($\mu$ a fuzzy ideal of N) then the level subsets $\mu_t$, $t \in [0, t]$, $t \leq \mu(0)$ are S-fuzzy near-rings.*

*Proof*: Left for the reader to prove.

**DEFINITION 9.3.3**: *A S-fuzzy near-ring associated with $\mu$ is said to have Smarandache fuzzy IFP (S-fuzzy IFP) if and only if $\mu(a n b) \geq \mu(ab)$ for all a, b, n $\in$ N.*

**DEFINITION 9.3.4**: *Let N be a S-fuzzy near-ring associated with $\mu$. $\mu N$ has fuzzy IFP if and only if $\mu_k$ is a S-fuzzy near-ring for all $1 \leq \mu \leq 1$.*

**DEFINITION 9.3.5**: *N is a Smarandache strong fuzzy IFP (S-strong fuzzy IFP) if and only if every S-fuzzy near-ring is a S-fuzzy IFP near-ring.*

**DEFINITION 9.3.6**: *$_\mu N$ is S-fuzzy IFP near-ring if and only if $\beta_s$ is a S-fuzzy IFP near-ring for all $s \in [0, 1]$ where*

$$\beta_s = \begin{cases} s \ if \ \mu(x) \geq s \\ \mu(x) \ if \ \mu(x) < s. \end{cases}$$

**DEFINITION 9.3.7**: *A S-fuzzy near-ring $_\mu N$ is said to be Smarandache fuzzy prime (S-fuzzy prime) if $\mu$ is not a constant function and for any $_\sigma N \,_\tau N$ S-fuzzy near-rings $\sigma \circ T \subseteq \mu$ implies either $\sigma \subset \mu$ or $T \subset \mu$.*

The author has only introduced just the concept of S-fuzzy near-ring. It is left for the reader to do research in this direction and develop the concept of S-fuzzy near-ring as such notion is very new and has a lot of scope for any innovative researcher.

Now we proceed on to define a new type of S-fuzzy near-rings using the concepts defined in section 9.2.

**DEFINITION 9.3.8**: *Let $\{P_{n \times n}, \oplus, \odot\}$ be a fuzzy near matrix ring we call $\{P_{n \times n}, \oplus, \odot\}$ a Smarandache fuzzy matrix near-ring (S-fuzzy matrix near-ring) if $P$ has a subset $S$, $S \subset P$ such that $\{S, \oplus, \odot\}$ is a near-field.*

**DEFINITION 9.3.9**: *The fuzzy complex near-ring $\{V, \oplus, \odot\}$ is called a Smarandache fuzzy complex near-ring (S-fuzzy complex near-ring) if $V$ has a proper subset $R$ such that $\{R, \oplus, \odot\}$ is a near-field.*

**DEFINITION 9.3.10**: *Let $\{W, \oplus, \odot\}$ be a non-associative fuzzy complex near-ring, we say $\{W, \oplus, \odot\}$ is a Smarandache non-associative fuzzy complex near-ring (S-non-associative fuzzy complex near-ring) if $\{W, \oplus, \odot\}$ has a non-associative fuzzy complex near division ring.*

**DEFINITION 9.3.11**: *Let $R[x^{[0,\,1]}]$ be the fuzzy polynomial near-ring. We say $R[x^{[0,\,1]}]$ is a Smarandache fuzzy polynomial near-ring (S-fuzzy polynomial near-ring) if $R[x^{[0\,1]}]$ has a proper subset $T$ which is a near-field.*

**DEFINITION 9.3.12**: *Let $R[x^{[0,\,1]}]$ be a S-fuzzy polynomial near-ring we say a polynomial $f(x)$ in $T$, $T$ a near-field, has a Smarandache root (S-root) if $f(\alpha) = 0$ for some $\alpha \in T$.*

Thus we see sometimes a polynomial which has roots may not in general have Smarandache roots.

It is left for the reader to construct examples and counter examples to prove the above statements.

**PROBLEMS :**

1.  Prove $\beta_s$ is a S-fuzzy near-ring.

2.  Prove $_\mu N$ where

$$\mu(x) = \begin{cases} 1 \text{ if } x \in \,_\mu N \\ s \text{ otherwise} \end{cases}$$

is a S-fuzzy prime near-ring.

3.  Is $P = \{a \,/\, a \in [0\ 1]\}$ a fuzzy complex near-field?

4.  Can $P_{n \times n}$ have a subset which is a near-field? Justify your answer.

5.  Can $\{W, \oplus, \odot\}$ have nontrivial division near-rings?

6.  Find subsets in $R[x^{[0, 1]}]$ which are near-fields.

7.  Find or give examples of polynomials in $R[x^{[0, 1]}]$ which has roots but no S-roots.

**Chapter Ten**

# SUGGESTED PROBLEMS

This chapter suggests 145 problems on near-rings, seminear-rings, S-near-rings and S-seminear-rings and their generalizations. The problems mainly are connected with group near-rings, semigroup near-rings, group seminear-rings and semigroup seminear-rings and their Smarandache analogues which can always be formulated once the concepts are mastered. The main motivation is that the study of this type is very meagre or absent. But this study is significant as group-rings and semigroup rings are very widely studied. The reader should try to study these concepts so that one can get several examples of both commutative and non-commutative near-rings and seminear-rings using varied types of groups and semigroups.

Thus this study will lead to a multifold research on the four concepts of semigroup near-rings, group near-rings, group seminear-rings and semigroup seminear-rings and mainly their Smarandache analogue.

The study of non-associative near-rings and seminear-rings leads to the concept of studying famous identities like Moufang, Bol, Bruck, WIP etc. Here also the existence of a class of non-associative near-ring or seminear-ring from the natural integers or in a natural way in unknown to the author. We can get classes of non-associative seminear-rings using loops and groupoids. By using the natural class of loops and groupoids described in Chapter I the researcher will certainly get classes of both non-associative near-rings and seminear-rings.

Finally the study of S-near-rings, S-seminear-rings both associative and non-associative is very recent (2002). Thus this book gives a list of definitions about S-seminear-rings. The book does not completely deal with all analogous properties in case of Smarandache notions. Several Smarandache problems are given for the researchers to solve. The set of suggested problems in this book is an added feature.

1.  Obtain necessary and sufficient condition for left (right) quotient near-ring N to be a S-quotient near-ring.

2.  Classify those near-rings N which fulfills left (right) ore condition also fulfills S-left (right) one condition.

3.  Classify those class of near-ring $F_x$ which are free and are also S-free.

4.  Obtain a necessary and sufficient condition on the near-ring N so that

    i.      All ideals are S-ideals or
    ii.     The internal direct sum is the same as S-internal direct sum.
    iii.    Prove 1 $\Leftrightarrow$ 2.

5.  Characterize those near-rings N which cannot contain ideals I such that

       i.        I is not a S-direct summand.

      ii.       I is not a S-strong direct summand.

     iii.       I is always a S-direct summand.

6. Let P(I) be the prime radical; find the necessary and sufficient condition on the near-rings so that P(I) = S (P (I)).

7. Find a necessary and sufficient condition on the ideal I of a near-ring N so that its prime radical is a S-prime radical.

8. Find conditions on the near-ring N so that every idempotent in N is an S-idempotent.

9. For what loops L will the near loop ring NL have proper subsets which are

       i.        Near-fields.

      ii.       Near domains.

10. Find a necessary and sufficient condition for a near-ring to be a Marot near-ring.

11. Classify all near loop rings which are Marot near loop rings. (By putting conditions on the near-ring N and on the loop L).

12. Define a nice class of Lields and define Smarandache Lields.

13. Let N be a loop near-ring and L a loop having n elements of order 2. Can the loop-loop near-ring have divisors of zero?

14. N a near-ring and L a loop. When will the near loop ring have nontrivial divisors of zero?

15. Give examples of loop-loop near-rings which are

       i.        S-loop near-ring.

      ii.       Bipotent-loop near-ring.

16. Characterize near loop rings which are S-near loop-rings.

17. Characterize those loop-loop near-ring NL which are bipotent loop near-ring.

18. Characterize those near-rings which are strongly subcommutative.

19. Obtain a new class of near domains.

20. Characterize those loops L so that the loop near domains, NL has atleast a subset $P \subset NL$ such that P is a near domain.

21. Let $Z_p$ be a near-ring and L a finite loop, $|L| = n$. Find conditions on order of L and on p so that the near loop ring $Z_pL$ has

    i.      Nontrivial divisors of zero.
    ii.     Nontrivial idempotents.
    iii.    Nontrivial nilpotents.

22.    Obtain a necessary and sufficient conditions on $Z_p$ and the loop L so that the S-near loop ring NL has the mod p envelope; $L^*$ to be a loop.

23.    Let $Z_2 = \{0, 1\}$ and L be a loop of even order. $Z_2L$ be the near loop ring. Can $L^*$ be a loop? Justify your claim.

24.    Find conditions on the near-ring N and on the loop L so that

    i.      NL is a Bol near-ring.
    ii.     NL is a Bruck near-ring.
    iii.    NL is a Mougang near-ring.
    iv.    NL is an Alternative near-ring.

25.    Let L be a any loop and N a near-ring. Find conditions on NL the near loop ring to have quasi regular elements.

26.    Does there exist a near loop ring for which every element other than 1 and 0 are quasi regular? Study this problem for S-near loop rings.

27.    Let $Z_m$ be the near-ring and $L_n(m) \in L_n$ be a loop. Can $Z_m L_n(m)$ have quasi regular elements? If so find all quasi regular elements of $Z_n L_n (m)$.

28.    Give a class of S-seminear-rings.

29.    Does their exist a seminear-ring in which the same subset P of N happens to be a S-hyper subseminear-ring as well as S-dual hyper subseminear-ring?

30.    Let $Z_n = \{0, 1, 2, \ldots, n\text{-}1\}$ be the set of integers modulo n; $n = p_1^{\alpha_1} p_2^{\alpha_2} \ldots p_t^{\alpha_t}$ where $p_1, p_2, \ldots, p_t$ are primes t > 1.

Define two binary operations '×' and 'o' on $Z_n$; '×' is the usual multiplication modulo n. 'o' is defined as a o b = a for all b, a $\in Z_n$. Let A = $\{1, q_1, q_2, \ldots, q_r\}$ where $q_1, q_2, q_3, \ldots, q_r$ are all odd primes different from $p_1, p_2, \ldots, p_t$ and $q_1, \ldots, q_n \in Z_n$.

Prove (A, ×) is always a group under '×'.

[Solution to this problem would result in a class of S-seminear-rings for every positive n, n a composite integers].

31.    Find finite S-pseudo seminear-ring.

32.    Does there exist infinite near-ring constructed using reals or integers, which are not S-pseudo seminear-rings?

33. Characterize those S-seminear-rings N so that the $S_1N$ subgroups happen to be $S_1N$ subgroup related to all subgroups (P, .) of (N, .).

34. Define $S_1DN$-subgroups; characterize those S-seminear-ring for which the $S_1DN$ subgroups happen to be a $S_1DN$ subgroups for all subgroups (P, +) of (N, +).

35. Can modular lattices have sublattices that are seminear-rings.

36. Find necessary and sufficient conditions on the seminear-ring N and on the group G so that the group seminear-ring has

   i.    S-zero divisors.
   ii.   S-semi-idempotents.
   iii.  S-idempotents.
   iv.   S-units.

37. Study the problem (36) in case of semigroup seminear-rings.

38. Find conditions for a S-semigroup seminear-ring to be a S-seminear-ring.

39. Find a necessary and sufficient condition on the NA seminear-ring N so that the two definitions.

   i.    S-NA seminear-ring I and
   ii.   S-NA seminear-ring I type A are equivalent.

40. Find a necessary and sufficient conditions on the S-NA-seminear-ring so that the two definitions.

   i.    S-NA seminear-ring I and
   ii.   S-NA seminear-ring I of type B are equivalent.

41. Find a new class of non-associative near-rings and seminear-rings other than using loops over near-rings (near loop rings) and loops over seminear-rings.

42. Find a necessary and sufficient condition for a non-associative near-ring to be

   i.    Power associative.
   ii.   Diassociative.

43. Study the same problem in case of non-associative seminear-rings.

44. Find a necessary and sufficient conditions for a S-NA seminear-ring (and S-NA near-ring) to be

   i.    S-power associative.
   ii.   S-diassociative.

45. Does there exists a NA near-ring N which is a Bruck near-ring i.e. it satisfies the two identities

      i.      $[x\,(yx)]\,z = x\,(y\,(xz))$.
      ii.     $(xy)^{-1} = x^{-1}\,y^{-1}$ for all x, y, z ∈ L.

[Hint: If we take $Z_2 = (0, 1)$ to be a near-ring and L a Bruck loop will we have $Z_2L$ to be a non-associative near-ring which is a S-Bruck near-ring but $Z_2L$ is need not be a Bruck near-ring. Apart from this we are not in a position to construct any non-associative near-rings which are Bruck near-rings].

46. Study of full ideals in the composition ring of polynomials over integers in case of polynomials near-rings. $Z[x]$ remains unknown obtain some interesting results about ideals in $Z[x]$.

47. Develop interesting properties on matrix near-rings.

48. Let $Z^o = Z^+ \cup \{0\}$.$\{Z^o, +, .\}$ is a seminear-ring where '+' is the usual '+' and '.' is defined by $a.b = a$ for all a, b ∈ $Z^o$. Take G any groupoid from the new class of groupoids. Study all properties like existence of S-zero divisors, S-idempotents the identities they satisfy etc., for the groupoid seminear-ring $Z^oG$.

49. Let G be a groupoid form the new class of groupoids given in chapter 1. (So order of G is finite say m). Let $Z_m$ be the seminear-ring m a composite number. Obtain conditions on m and n and make a study about the properties enjoyed by $Z_mG$ when

      i.     $(m, n) = 1$.
      ii.    $m = n$.
      iii.   $(m, n) = p$,    p a prime.
      iv.   $(m, n) = n$,    n a prime.
      v.    $(m, n) = 1$ and $m = p^t$ and $n = q^t$ where p and q are primes.

50. Give an example of a S-weak Bol seminear-ring which is never a S-strong Bol seminear-ring.

51. What groupoids we will have the groupoid seminear-ring to be

      i.     S-Bol seminear-ring.
      ii.    S-Moufang seminear-ring.
      iii.   S-P-seminear-ring.
      iv.   S-alternative seminear-ring;

      study in case of weak and strong structures.

52. Let $(Z_{12}, \times, .)$ be a seminear-ring and G be a groupoid given by the following table:

| · | $a_0$ | $a_1$ | $a_2$ | $a_3$ | $a_4$ | $a_5$ | $a_6$ | $a_7$ | $a_8$ | $a_9$ | $a_{10}$ | $a_{11}$ |
|---|---|---|---|---|---|---|---|---|---|---|---|---|
| $a_0$ | $a_0$ | $a_4$ | $a_8$ | $a_0$ | $a_4$ | $a_8$ | $a_0$ | $a_4$ | $a_8$ | $a_0$ | $a_4$ | $a_8$ |
| $a_1$ | $a_3$ | $a_7$ | $a_{11}$ | $a_3$ | $a_7$ | $a_{11}$ | $a_3$ | $a_7$ | $a_{11}$ | $a_3$ | $a_7$ | $a_{11}$ |
| $a_2$ | $a_6$ | $a_{10}$ | $a_2$ | $a_6$ | $a_{10}$ | $a_2$ | $a_6$ | $a_{10}$ | $a_2$ | $a_6$ | $a_{10}$ | $a_2$ |
| $a_3$ | $a_9$ | $a_1$ | $a_5$ | $a_9$ | $a_1$ | $a_5$ | $a_9$ | $a_1$ | $a_5$ | $a_9$ | $a_1$ | $a_5$ |
| $a_4$ | $a_0$ | $a_4$ | $a_8$ | $a_0$ | $a_4$ | $a_8$ | $a_0$ | $a_4$ | $a_8$ | $a_0$ | $a_4$ | $a_8$ |
| $a_5$ | $a_3$ | $a_7$ | $a_{11}$ | $a_3$ | $a_2$ | $a_7$ | $a_{11}$ | $a_3$ | $a_7$ | $a_{11}$ | $a_3$ | $a_{11}$ |
| $a_6$ | $a_6$ | $a_{10}$ | $a_2$ | $a_6$ | $a_{10}$ | $a_2$ | $a_6$ | $a_{10}$ | $a_2$ | $a_6$ | $a_{10}$ | $a_2$ |
| $a_7$ | $a_9$ | $a_1$ | $a_5$ | $a_9$ | $a_1$ | $a_5$ | $a_9$ | $a_1$ | $a_5$ | $a_9$ | $a_1$ | $a_5$ |
| $a_8$ | $a_0$ | $a_4$ | $a_8$ | $a_0$ | $a_4$ | $a_8$ | $a_0$ | $a_4$ | $a_8$ | $a_0$ | $a_4$ | $a_8$ |
| $a_9$ | $a_3$ | $a_7$ | $a_{11}$ | $a_3$ | $a_7$ | $a_{11}$ | $a_3$ | $a_7$ | $a_{11}$ | $a_3$ | $a_7$ | $a_{11}$ |
| $a_{10}$ | $a_6$ | $a_{10}$ | $a_2$ | $a_6$ | $a_{10}$ | $a_2$ | $a_6$ | $a_{10}$ | $a_2$ | $a_6$ | $a_{10}$ | $a_2$ |
| $a_{11}$ | $a_9$ | $a_1$ | $a_5$ | $a_9$ | $a_1$ | $a_5$ | $a_9$ | $a_1$ | $a_5$ | $a_9$ | $a_1$ | $a_5$ |

Study the groupoid near-ring $Z_{12}G$ in the context of problem 51.

53. Does {V, ⊕, ⊙} have fuzzy left ideal? (or two sided fuzzy ideal). (By notation ideals of V will be called only as fuzzy ideals).

54. Is V a fuzzy simple near-ring?

55. Does (V, ⊕, ⊙) have prime ideals?

56. Does the fuzzy non-associative near-ring have fuzzy ideals?

57. Does (W, ⊕, ⊙) the complex non-associative near-ring have invariant fuzzy subnear-rings?

58. Is (W, ⊕,⊙) a fuzzy principal ideal domain?

59. Obtain some interesting properties about S-fuzzy near-ring.

60. Given G is a finite group and N a near-ring when will the group near-ring be artirian?

61. Whether these conditions on G are also necessary for the group near-ring to be Noetherian. (This problem in case of group rings is a classical one and complete answer however is not known, whether the conditions on G are also necessary for the group ring to be Noetherian).

62. Let N be a near-ring built using Z or Q or R and G any group. Find conditions on G and N so that the group near-ring ZG

    i.     is the subdirect product of prime near-rings.
    ii.    Intersection of all prime ideals of ZG is 0.

63. When will the group near-ring $Z_pG$ have no non-zero nil ideals?

64.    Let N be a near-ring built using Z or Q or R and G any group, when will the group near-ring ZG (QG or RG) have non-zero nil ideals?

65.    Define a concept in near-rings analogous to semisimple rings.

66.    When are group near-rings semiprime?

67.    Obtain a necessary and sufficient condition for a group near-ring to have zero divisors (S-zero divisors).

68.    Obtain a necessary and sufficient condition for a group near-ring to be regular (S-regular).

69.    Let NG and $N_1$H be any two group near-rings when are they isomorphic (S-isomorphic). That is like isomorphism problem for group rings when are two group near-rings isomorphic?

70.    Let G be a group and N a near-ring. Obtain necessary and sufficient condition for the group near-ring NG to be a p-near-ring (S-p-near-ring).

71.    If S is a commutative semigroup with divisors of zero and elements of finite order and N a near-field. Can the semigroup near-ring NS be a Marot semigroup near-ring (S-Marot semigroup near-ring)?

72.    Characterize those semigroup near-rings NS which are p-near-rings.

73.    Let $Z_2 = (0, 1)$ be the near-field. S any semigroup obtain a necessary and sufficient condition on the semigroup S so that the semigroup near-ring $Z_2$S has always a mod p envelop S* to be a S-semigroup.

74.    Study the same problem i.e. problem (73) when $Z_2$ is replaced by $Z_p$, p an odd prime.

75.    Let NS and $N_1S_1$ be any two semigroup near-rings. Obtain necessary and sufficient condition so that these two near-rings are isomorphic and S-isomorphic.

76.    Characterize those semigroup near-rings NS which have

         i.     Units and S-units.
         ii.    Zero divisors and S-zero divisors.
         iii.   Idempotents  and S-idempotents.

77.    Obtain necessary and sufficient for the group near-ring $Z_p$G so that

         i.     G = G*.
         ii.    G* is a group.
         iii.   G* is only a semigroup.

78. Let G = ⟨g / $g^p$ = 1⟩, p a prime and $Z_m$ = {0, 1, 2, …, m-1} m a prime. Find G* for the group near-ring $Z_mG$. Find G* when m = p when will G* be a group?

79. Study the group near-ring $Z_sG$, when s is not a prime what is the structure of G* if |G| = s. will G* ever be a group.

80. Find necessary and sufficient condition for the group near-ring $Z_pG$ to have G* = G.

81. Obtain a necessary and sufficient condition for the group near-ring NG to be a S-near-ring satisfying some polynomial identities.

82. Characterize those group near-rings NG (semigroup near-rings NS) to be a chain near-ring.

83. Find conditions on the group and on the loop near-ring N so that NG has

    i.      S-zero divisors.
    ii.     Non-trivial S-ideals.
    iii.    Satisfies chain conditions.

84. Let $Z_n$ be a seminear-ring. $L_n(m) \in L_n$ be the new class of loops. Study for the loop seminear-ring $Z_nL_n(m)$ for varying m (n non-prime) to have

    i.      S-units.
    ii.     S-idempotents.
    iii.    S-zero divisors and units and idempotents and zero divisors.

85. $Z_nL_n(m)$ given in problem 84 is a S-seminear-ring III. Find S-ideals III, S-N-subsemigroups. Is $Z_nL_n(m)$ S-regular? Obtain any other interesting results on these classes of S-seminear-rings III?

86. Obtain a necessary and sufficient condition for the S-seminear-ring III , $Z_nL_n(m)$, defined in problem 84 to be a S-seminear-ring IV.

87. Let G be a groupoid from the new class of groupoids built using $Z_n$. {$Z_n$, '×', '.'} be a near-ring. For the groupoid near-ring $Z_nG$ study the properties mentioned in problems 84 and 85.

88. Let L be any finite loop and N a S-near-ring. Obtain a necessary and sufficient condition for the near loop ring NL to have SJ(NL) = J(NL).

89. Does there exist a near-ring N in which every right quasi regular element is S-right quasi regular?

90. Characterize those near-rings which has quasi regular elements but no S-quasi regular elements.

91. Let $N = Z_{11} \times Z_p$ (where $Z_n$ is the near-ring and $Z_p$ is a field of characteristic p, (n, p) = 1, n not a prime) be the S-near-ring. Which of the near loop rings $NL_n(m)$ has $SJ(NL_n(m)) = J(NL_n(m))$; where $L_n(m) \in L_n$?

92. Characterize those group near-rings which are only S-near rings and not S-group near-rings.

93. Characterize those group near-rings NG which are such that NG has a S-subnear-ring if and only if NG is a S-group near-ring.

94. Let N be a near-field and G a torsion free abelian group. Does the group near-ring NG have zero divisors? (This problem is solved in case of group rings).

95. Let G be a torsion free non-abelian group and N a near-field. Can the group near-ring NG have zero divisors? (This problem even for group rings is still unsettled).

96. Let H be a torsion free abelian subgroup of a group G and N a near-field. NH is a group near-ring contained in the group near-ring NG. If $\alpha \in NH \subset NG$; $\alpha \neq 0$. Can $\alpha$ be a zero divisor in NG? (Hint: This question has been settled in the case of group rings; when $\alpha$ is not a zero divisor).

97. Let N be a near-ring and G a finite group.

    i. Obtain a necessary and sufficient condition for the near-ring NG to have S-artinian condition?
    ii. Obtain a necessary and sufficient condition for the S-semigroup near-ring NG to satisfy S-artinian condition.

98. Let N be a near-field and G any group. When will the group near-ring NG be S-semiprime?

99. Let N be a S-near-ring so that the group near-ring NG is a S-group near-ring. Obtain a necessary and sufficient condition for the group near-ring NG to be S-semiprime.

100. Let N be a near-field of the form $Z_p$, p a prime. G be any group with no elements of order p. Can the group near-ring NG have non-zero nil ideals?

101. If we replace the near-ring $N = Z_p$ by S-p-near-rings and G any group with no elements of order p. Can the group near-ring have non-zero nil ideals or non-zero S-nil ideals?

102. Let N be a near-field, $Z_p$ or a S-p-near-ring and G a finite group. Will the following conditions be equivalent?

    i. NG is S-semiprime and
    ii. $\Delta(G)$ has no elements of order p.

103. Let N be a near-ring and G any group. When does the group near-ring NG satisfy standard polynomial identity?

104. Define the concept of S-near-ring satisfying standard polynomial identities.

105. Obtain a necessary and sufficient condition for a S-group near-ring to have zero divisors and S-zero divisors.

   (We know in case of group ring KG if K is a field of characteristic zero and G a finite group or a group having elements of finite order then the group ring KG has zero divisors.)

106. Let N be a near-ring and G a group. Let NG be the S-group near-ring. Obtain a necessary and sufficient condition for NG to have non-trivial idempotents and S-idempotents.

107. Let $N = Z_p$ be the near-field. $Z_pG$ be the group near-ring of a locally finite group G over $Z_p$. Will $Z_pG$ be regular? Can $Z_pG$ be S-regular?

108. Obtain a necessary and sufficient condition on the group near-ring NG so that NG is S-regular. NG is S-group near-ring and regular.

109. Obtain conditions on the near-ring N and on the S-semigroup S so that the semigroup near-ring NS has

   i.    non-trivial zero divisors and S-zero divisors.
   ii.   idempotents and S-idempotents.

110. Let NS be a S-semigroup near-ring. Obtain a necessary and sufficient condition for NS to be S-regular and NS to be regular.

111. Let N be a near-field. S a S-semigroup. Can the semigroup near-ring NS have non-zero nil ideals?

112. Let NS be a S-semigroup near-ring. Obtain necessary and sufficient condition for NS to be semiprime, S-semiprime.

113. Characterize those S-near-rings N which has N-uniform ideal.

114. Obtain a necessary and sufficient condition for a near-ring to be a S-Marot near-ring.

115. Obtain a necessary condition for a planar near-ring to be S-planar.

116. Find sufficient conditions for a near-ring to be a S-σ near-ring.

117. Obtain a necessary and sufficient condition so that every σ-near ring is a S-σ-near-ring.

118. Find conditions on a near-ring so that it has a S-uniform ideal of G.

119.    Find a necessary and sufficient condition on a S-near-ring N to be S-normal.

120.    Find a condition on the near-ring N so that all uniform ideals of G are S-uniform ideals of G.

121.    Characterize those S-composition near-ring.

122.    Obtain a necessary and sufficient condition so that every normal near-ring is S-normal and vice versa.

123.    Obtain a necessary and sufficient condition for a F-near-ring to be a S-F-near-ring and conversely.

124.    Let G be a finite S-semigroup of order n and Z be the near-ring. Does ZG have

        i.       S-essential ideals.
        ii.      S-quasi ideals.
        iii.    S-central elements.
        iv.     S-regular elements.
        v.      S-principal ideals.

125.    For the S-S-semigroup automaton level I and II derive some inter-relations or a set of S-S-semigroup automaton level I and II with no inter-relation entirely acting in a discrete way.

126.    Using S-mixed direct product find an example of

        i.       S-left distributive near-ring.
        ii.      S-strong IFP near-ring.
        iii.    S-affine near-ring.
        iv.     S-medial near-ring.

127.    Obtain condition for a near-ring N to be a S-n-ideal near-ring.

128.    Characterize those near-rings which are S-LSD near-rings.

129.    Characterize those group near-rings which are only S-near-rings and not S-group near-rings.

130.    Characterize those S-group near-rings which satisfy S-ACCL and S-ACCN.

131.    When are S-group near-rings (seminear-rings) S-artinian?

132.    Characterize those group near-rings (seminear-rings) which are

        i.       S-commutative.
        ii.      S-weakly commutative.
        iii.    S-cyclic.
        iv.     S-weakly cyclic.

133. Obtain necessary and sufficient condition so that a seminear-ring is a S-strict seminear-ring.

134. Characterize those seminear-rings which are S-pseudo seminear-rings.

135. Characterize those semigroup (group) near-rings which are S-simple.

136. Find a class of S-reliable near-rings.

137. Characterize

      i.     S-planar near-rings.
      ii.    S-left infra near-rings.
      iii.   S-strongly semiprime near-rings.
      iv.   S-completely equiprime near-rings.
      v.    S-integral near-rings.

138. Give a class of S-NA-near-rings which are S-planar.

139. Find groupoid near-rings $Z_nG$ which has S-quasi subnear-rings.

140. Let $N = Z_p L_n(m)$ be the near loop ring of the new class of loops $L_n(m) \in L_n$ over the near-field $Z_p$.

      i.     Is $Z_pL_n(m)$ a S-s-near loop ring?
      ii.    Is N a S-Moufang loop near-ring?
      iii.   Is N S-strictly duo?

141. Find a class of near-rings which are S-quasi near-rings.

142. Find the S-Jacobson radical for $Z_p L_n(m)$; $L_n(m) \in L_n$.

143. When is SJ $(Z_pL_n(m)) = J(Z_pL_n(m))$?

144. Find a class of S-bipotent seminear-ring.

145. Find those S-group near-rings and S-semigroup seminear-rings which satisfy standard polynomial identities.

# INDEX

## About the Author

Dr. W. B. Vasantha is an Associate Professor in the Department of Mathematics, Indian Institute of Technology Madras, Chennai, where she lives with her husband Dr. K. Kandasamy and daughters Meena and Kama. Her current interests include Smarandache algebraic structures, fuzzy theory, coding/ communication theory. In the past decade she has guided seven Ph.D. scholars in the different fields of non-associative algebras, algebraic coding theory, transportation theory, fuzzy groups, and applications of fuzzy theory to the problems faced in chemical industries and cement industries. Currently, six Ph.D. scholars are working under her guidance. She has to her credit 241 research papers of which 200 are individually authored. Apart from this she and her students have presented around 262 papers in national and international conferences. She teaches both undergraduate and post-graduate students at IIT and has guided 41 M.Sc. and M.Tech projects. She has worked in collaboration projects with the Indian Space Research Organization and with the Tamil Nadu State AIDS Control Society.

She can be contacted at vasantha@iitm.ac.in
You can visit her on the web at http://mat.iitm.ac.in/~wbv

**Definition:**

**Generally, in any human field, a *Smarandache Structure* on a set A means a weak structure W on A such that there exists a proper subset B ⊂ A which is embedded with a stronger structure S.**

**These types of structures occur in our everyday's life, that's why we study them in this book.**

**Thus, as a particular case:**

**A *Near-ring* is a non-empty set N together with two binary operations '+' and '.' such that (N, +) is a group (not necessarily abelian), (N, .) is a semigroup. For all a, b, c ∈ N we have (a + b) . c = a . c + b . c.**

**A *Near-Field* is a non-empty set P together with two binary operations '+' and '.' such that (P, +) is a group (not necessarily abelian), (P \ {0}, .) is a group. For all a, b, c ∈ P we have (a + b) . c = a . c + b . c.**

**A *Smarandache Near-ring* is a near-ring N which has a proper subset P ⊂ N, where P is a near-field (with respect to the same binary operations on N).**